\documentclass{article}

\usepackage[letterpaper, hmarginratio=1:1]{geometry}

\usepackage{amsmath, amsthm, amsfonts, thmtools, amssymb,tikz,hyperref,bm,stmaryrd,bbm}
\usetikzlibrary{arrows}

\usepackage[mathscr]{euscript}
\usepackage{mathtools}

\numberwithin{equation}{section}

\usepackage{xcolor}

\newcommand{\ve}{v^{\mathsf{eff}}}
\newcommand{\Lip}{\mathrm{Lip}}
\newcommand{\1}{\mathbf 1}

\newcommand{\Z}{\mathbb{Z}}

\DeclareMathOperator{\eig}{eig}
\DeclareMathOperator{\Tr}{Tr}

\DeclareMathOperator{\dr}{dr}
\DeclareMathOperator{\eff}{eff}

\DeclareMathOperator{\Mat}{Mat}
\DeclareMathOperator{\dist}{dist}
\DeclareMathOperator{\supp}{supp}

\DeclareMathOperator{\dg}{dg}
\DeclareMathOperator{\sd}{sd}

\DeclareMathOperator{\T}{\mathbf{T}}
\DeclareMathOperator{\Hop}{H}

\DeclareMathOperator{\Var}{Var}

\DeclareMathOperator{\E}{\mathbb{E}}

\newtheorem{thm}{Theorem}[section]
\newtheorem{prop}[thm]{Proposition}
\newtheorem{lem}[thm]{Lemma}
\newtheorem{cor}[thm]{Corollary}

\theoremstyle{remark}
\newtheorem{rem}[thm]{Remark}
\theoremstyle{definition}
\newtheorem{definition}[thm]{Definition}
\newtheorem*{definition*}{Definition}

\newtheorem{example}[thm]{Example}
\newtheorem{assumption}[thm]{Assumption}
\theoremstyle{plain}
\newtheorem{innerthm}{Theorem}
\newenvironment{theorema}[1]
{\renewcommand\theinnerthm{#1}\innerthm}
{\endinnerthm}

\renewcommand{\L}{\mathbf{L}}
\renewcommand{\H}{\mathbf{H}}

\newcommand{\R}{\mathbb{R}}
\newcommand{\A}{\mathfrak A}

\newcommand{\llb}{\llbracket}
\newcommand{\rrb}{\rrbracket}

\title{Hydrodynamic limit of the Toda lattice}
\author{Amol Aggarwal \and Matthew Nicoletti}
\date{}

\begin{document}

\maketitle

\begin{abstract}
In this paper we consider the Toda lattice under random,
spatially inhomogeneous initial data, meaning that its initial momenta and
interparticle spacings are independent and their laws may slowly vary in space. When the density is strictly positive (but sufficiently small), we identify the model's hydrodynamic limit under Euler scaling. More precisely, we prove that spatial averages of each of the model’s local conserved quantities converge to moments of a density $m_{\tau} (\lambda, q)$ solving the generalized hydrodynamic equation predicted by Doyon \cite{GHCS} and Spohn \cite{Spo20GGE}. To do so, we justify the notion that the Toda lattice can be thought of as a dense collection of many “quasiparticles” that interact through scattering. As in \cite{Agg25a}, we define these quasiparticles through spectral properties of the Toda Lax matrix, but in our inhomogeneous setting we must now show that the eigenvectors of this Lax matrix remain localized throughout the dynamics. The core of our work is to show that, with high probability, the quasiparticle trajectories converge to the characteristic curves of the generalized hydrodynamic equation.

\end{abstract}

\tableofcontents
	
		\section{Introduction}
		\label{subsec:over}

		\subsection{Generalized hydrodynamics and the main result}

		Describing the large time behavior of many-body 
		Hamiltonian dynamics is a fundamental question 
		in mathematical physics. One aspect of this is
		proving hydrodynamic limits of such systems under random
		(or generic) initial data. Since the work of Morrey \cite{Mor55}, 
		it is believed that microscopic degrees of freedom locally equilibrate, 
		and thus at large
		scales the system is described by fewer effective parameters: the densities of (locally) conserved quantities.

        Chaotic systems have only a few conserved quantities, 
		typically the number of particles, total momentum, and total energy. 
		The hydrodynamic limit of the system is thus expected
		to consist of three coupled partial differential equations (known as Euler equations)
		describing the densities of particles, momenta, and energy;
		for detailed predictions, see the book of Spohn \cite{Spo91}.
		Much in this direction has been established at vanishingly low densities 
		(where the Boltzmann equation takes the place of the Euler equations) 
		in works of Lanford \cite{Lan75},
		Bodineau--Gallagher--Saint-Raymond--Simonella \cite{SDHG}, and
		Deng--Hani--Ma \cite{DHM24},
		or in the presence of a small amount of added noise 
		driving the dynamics by Olla--Varadhan--Yau \cite{OVY93}. However, for genuinely deterministic chaotic dynamics at strictly positive density, the proof of the Euler scale hydrodynamic equations remains open.

		 In contrast, integrable systems have infinitely many local conserved quantities. 
		This suggests that their long-time limits should be described by equations 
		describing the joint evolution of infinitely many densities. 
		These are known as \emph{generalized hydrodynamic equations}, 
		and their form was predicted much more recently in the physics literature by
		Bertini--Collura--De Nardis--Fagotti \cite{BCDF16} and
		Castro-Alvaredo--Doyon--Yoshimura \cite{CDY16}.
		Mathematically, proofs of convergence to
		these generalized hydrodynamic equations 
		had been previously obtained for only two integrable systems, 
		both of a combinatorial nature. The first is the hard rods (one-dimensional hard spheres) 
		model in works of Boldrighini--Dobrushin--Sukhov \cite{OHRH}, Ferrari--Franceschini--Grevino--Spohn \cite{HRHF} and Ferrari--Olla \cite{FO25,MFHRSL};
		the second is a cellular automaton known as the box-ball system in works of 
		Ferrari--Nguyen--Rolla--Wang \cite{SDBS}, Croydon--Sasada \cite{CS21}, and Olla--Sasada--Suda \cite{OSS24}. 
		
		In this work, we prove convergence to the
		generalized hydrodynamic theory for the Toda lattice \cite{Tod67}, a Hamiltonian system \(t \mapsto (\bm p(t);\bm q(t))\), with
		momenta \(\bm p(t)=(p_i(t))\) and
		positions \(\bm q(t)=(q_i(t))\) indexed by $i \in \mathbb{Z}$, evolving under the equations
		\begin{equation*}
		\partial_tq_i(t)=p_i(t),
		\qquad
		\partial_tp_i(t)=e^{q_{i-1}(t)-q_i(t)}-e^{q_i(t)-q_{i+1}(t)}.
		\end{equation*}
        
        \noindent As these dynamics are translation-invariant, we fix $q_0 (0) = 0$.  The integrability of the Toda lattice was established by Flaschka \cite{Fla74} and
		Manakov \cite{Man74}, and its scattering shift was computed by Moser \cite{Mos75}. 
		Since then, the Toda lattice has been viewed as a canonical example of a completely integrable system. 

        In cases when its initial data approximates a periodic background at $\pm \infty$ (such as when it decays or approximates a piecewise smooth profile), the long-time behavior of the Toda lattice has been studied, often by inverse scattering or Riemann-Hilbert methods. Some earlier works in this direction include those by Venakides--Deift--Oba \cite{DOV91}, Bloch--Kodama \cite{DRL}, Deift--Kamvissis--Kriecherbauer--Zhou \cite{DKKZ96}, Deift--McLaughlin \cite{CLL}, Bloch--Golse--Paul--Uribe \cite{BGPU03}, and Kr\"{u}ger--Teschl \cite{LADI}.
		More recently, Mazzuca--Piorkowski \cite{MP26} used Riemann--Hilbert methods to study the long-time behavior of specific dense $N$-soliton solutions 
		with random scattering data, and Brollo \cite{TCIST} used the inverse scattering method to predict large scale interactions between solitons and radiation. In contrast, in the present work we directly sample the physical variables $(\bm p(0);\bm q(0))$. Our initial data is outside the class of $N$-soliton solutions and is rough everywhere on the line (including at $\pm \infty$). Thus, the above analyses do not directly apply, and our approach to asymptotic analysis is different.


		We now state an informal version of our main hydrodynamic limit result. Let $T>0$ be the large parameter representing the time
		scale. Sample the initial variables
		$(q_{i+1}(0)-q_{i}(0),p_i(0))_{i\in\mathbb Z}$ independently, with laws that may vary
		in a piecewise Lipschitz way in the variable $i/T$;  in addition,
		assume the mean particle spacing $\mathbb{E}[q_{i+1}(0)-q_{i}(0)]$
		is of order $\mathfrak A$, where the fixed constant $\mathfrak A$ is
		sufficiently large. See Assumption \ref{ass:entry-law-profiles} for more details on the initial data.\footnote{We remark that in Theorem \ref{thm:A-intro}, the
		initial data need not be
		sampled from slowly varying invariant measures; see 
		Example \ref{ex:equally-spaced-particles}.}

        Next, we let $\mathfrak{k}_i^{[n]}(t)$ be the contribution of the $i$-th particle to the $n$-th locally conserved quantity of the Toda lattice. To define it, let $\mathbf{L}(t)$ denote the Lax matrix associated with the Toda lattice (see Definition \ref{def:flandL}), and set  $\mathfrak{k}_i^{[n]}(t):=[\L(t)^n]_{ii}$. For example, if $n=0$ then $\mathfrak{k}_i^{[n]} (t) = 1$ counts the $i$-th particle, and if $n=1$ then $\mathfrak{k}_i^{[n]} (t) = p_i(t)$ is the $i$-th particle's momentum. 
        
        Let $\alpha_0(q)$ denote the
		average particle spacing near the macroscopic position $q \in \mathbb{R}$
		and $\varrho_0(\lambda,q)$ denote the
		local limiting eigenvalue density of $\L(0)$ around $q $ (see  Definition \ref{rho0}). Then set $m_0(\lambda,q)=\alpha_0(q)^{-1}\varrho_0(\lambda,q)$. Subject to this initial data, let $m_\tau$ denote the solution of the \emph{generalized
		hydrodynamic (GHD) equation}   
        	\begin{equation}\label{eqn:rhoGHDintro_pref}
		\partial_\tau m_\tau(\lambda,q) +\partial_q\bigl(\ve_\tau(\lambda;q)m_\tau(\lambda,q)\bigr)=0,
		\qquad \tau\geq0,\quad q,\lambda\in\mathbb R;
		\end{equation}

        \noindent where $\ve_\tau$ satisfies the relation 
		\begin{equation}\label{eqn:cra_pref}
		\ve_\tau(\lambda;q)
		=
		\lambda
		+2\int_{-\infty}^{\infty}
		\log|\lambda-\lambda'|\,m_\tau(\lambda',q)
		\bigl(\ve_\tau(\lambda';q)-\ve_\tau(\lambda;q)\bigr)\,d\lambda'.
		\end{equation}

        \noindent The existence and uniqueness of the pair 
		$(m_\tau,\ve_\tau)$ solving \eqref{eqn:rhoGHDintro_pref} and \eqref{eqn:cra_pref} is provided by Proposition \ref{prop:GHDint} below. By \eqref{eqn:cra_pref} $\ve_\tau (\lambda,q)$, and thus $\ve_\tau (\lambda,q) m_\tau (\lambda,q)$, can be viewed as a functional of $\lambda \mapsto m_\tau(\lambda,q)$. Thus, \eqref{eqn:rhoGHDintro_pref} resembles an Euler equation of the form $\partial_{\tau} m_{\tau} (\lambda,q) = \partial_q J(m_{\tau} (\lambda,q))$ often describing hydrodynamic scaling limits for interacting particle systems \cite{Spo91}. This can be viewed as an infinite system of coupled equations indexed by $\lambda \in \mathbb{R}$.

       Under this notation, we prove the following result that describes the hydrodynamic evolution of this Toda lattice's conserved densities; it was predicted by Doyon \cite{GHCS} and Spohn \cite{Spo20,Spo20GGE}.

		\begin{theorema}{A}[Theorem \ref{thm:local-charge-testfn-intro}]
		\label{thm:A-intro}

		For every fixed $\tau \geq 0$ and integer $n\geq0$, and every smooth,
		compactly supported function $g:\mathbb R\to\mathbb R$, we have
		\[
		\frac1T\sum_{i\in\mathbb Z}
		g\bigl(T^{-1}q_i(T\tau)\bigr)\mathfrak{k}_i^{[n]}(T\tau)
		\longrightarrow
		\int_{\mathbb R}\int_{\mathbb R}
		g(q)\lambda^n m_\tau(\lambda,q)\,d\lambda\,dq
		\]
		in probability as $T\to\infty$.
		\end{theorema}

		Important to a proper understanding of these hydrodynamics is not only
		Theorem \ref{thm:A-intro} itself, but also the physical meaning of the
		function $m_\tau$. This can be interpreted as the density of
		``quasiparticles'' in the Toda lattice, whose analysis will play
		a central role both in proving Theorem \ref{thm:A-intro} 
		and intuiting equations \eqref{eqn:rhoGHDintro_pref} and \eqref{eqn:cra_pref}.

		\subsection{Quasiparticles}
		Our analysis centers around justifying the notion that the Toda lattice
		can be viewed as a dense collection of objects called
		\emph{quasiparticles}, a picture which is believed to hold more generally
		for both classical and quantum integrable systems (see the review of Doyon
		\cite{Doy20}). 
        Under this picture, each quasiparticle possesses a time-independent spectral
		parameter $\lambda_k$ and a time-dependent location $Q_k(t)$; this system of quasiparticles should satisfy two properties which together can be used to study the original integrable dynamics. 
        
        First, local averages of charge densities should be reconstructible from the quasiparticles (for example, for the $n=0$ charge, the number of Toda particles in a spatial interval should approximate the number of quasiparticles in it). 
         In this way, proving Theorem \ref{thm:A-intro} would reduce to understanding the limit of local densities of these quasiparticles. 
        In our setup, after properly defining the quasiparticles, we will prove that their limiting density is the solution $m_\tau(\lambda,q)$ of \eqref{eqn:rhoGHDintro_pref}, in the sense that $m_\tau(\lambda,q)\,d\lambda\,dq \cdot T$ approximates the number of quasiparticles at time $T \tau$ whose locations and spectral parameters lie in $[Tq, T(q + dq)]$ and  $[\lambda, \lambda + d \lambda]$, respectively. 
        

        Second, the evolution of quasiparticles should be approximately explicit, which for the Toda lattice takes the form
        \begin{equation}\label{eqn:asymptotic-scattering-overview}
		Q_k(t)
		\approx Q_k(0)+\lambda_k t\\
		+2\sum_{j\ne k}\log|\lambda_k-\lambda_j|
		\left(
		\mathbbm{1}_{Q_j(0)<Q_k(0)}
		-\mathbbm{1}_{Q_j(t)<Q_k(t)}
		\right).
		\end{equation}
        The dynamics \eqref{eqn:asymptotic-scattering-overview} have the following interpretation, known as the 
		\emph{flea-gas algorithm}: Quasiparticle $Q_k$ moves with velocity $\lambda_k$ until it meets another
		quasiparticle $Q_j$. At that moment, $Q_k$ instantaneously moves forward
		or backward by $2\log|\lambda_k-\lambda_j|$, depending on whether it meets
		$Q_j$ from the right or from the left, respectively. This displacement can
		make it pass another quasiparticle, producing a cascade of such
		interactions. The quasiparticle then proceeds with velocity $\lambda_k$
		until its next encounter, when the same procedure repeats. Here,
		$2\log|\lambda_k-\lambda_j|$ is the Toda lattice's two-body scattering
		shift, describing the phase displacement produced when two quasiparticles
		pass through one another. For this reason, we refer to \eqref{eqn:asymptotic-scattering-overview} as the \emph{asymptotic scattering relation}. 
        
        Now, suppose we can interpret $\ve_{\tau}(\lambda; q)$ as the limiting instantaneous velocity at time $T\tau$ of a quasiparticle at spatial location $Tq$ with spectral parameter $\lambda$. Then, the evolution equation \eqref{eqn:rhoGHDintro_pref} for $m_{\tau}$ would hold by definition, and the relation \eqref{eqn:cra_pref} for $\ve_{\tau}(\lambda; q)$ would be heuristically suggested by \eqref{eqn:asymptotic-scattering-overview} (see \cite[Section 5.3]{Spo23} for details in the related context of the hard rods model).

        In the physics literature, the concept of viewing many-body integrable systems as dense collections of quasiparticles was posited (but without providing a microscopic definition of their locations) by Bertini--Collura--De Nardis--Fagotti \cite{BCDF16} and
		Castro-Alvaredo--Doyon--Yoshimura \cite{CDY16}, and the asymptotic scattering relation was predicted by Doyon--H\"ubner--Yoshimura \cite{DHY24,DH23,DHY23}. The flea-gas algorithm above was described and named by Doyon--Yoshimura--Caux \cite{DYC18}. Earlier works of Zakharov \cite{Zak71} and El--Kamchatnov--Pavlov--Zykov \cite{El03,EK05,EKPZ11} also came across such dynamics, and observed that they parallel the collision structure between solitons in classical integrable systems, such as the Korteweg--de Vries (KdV) equation.\\

        To set up this framework mathematically, one must define quasiparticles and their spectral parameters. Since the spectral parameters are time-independent, they should be given by conserved quantities of the dynamics. For the Toda lattice, these are the eigenvalues of its \emph{Lax matrix}, which is the symmetric tridiagonal matrix
		$\L(t)=[L_{ij}(t)]$ whose nonzero entries are $L_{ii}(t)=p_i(t)$ and $L_{i,i+1}(t)=L_{i+1,i}(t)=e^{(q_i(t)-q_{i+1}(t))/2}$. The empirical spectral density of $\mathbf{L}(t)$ therefore describes the distribution of quasiparticle spectral parameters in the Toda lattice. Spohn \cite{Spo20GGE} predicted formulas for this density (under a wide class of invariant measures), which were subsequently proven in works by  Grava, Guionnet, Kozlowski, Little, Mazzuca, and Memin \cite{Maz22,GM22,MM24,LDPC}, with \cite{GM22,LDPC} also providing large deviation principles for these densities.
        
        A precise definition of the quasiparticle locations is more recent. For the Toda lattice under a natural invariant measure (known as thermal equilibrium), \cite{Agg25a} showed that their positions indeed likely satisfy \eqref{eqn:asymptotic-scattering-overview}. Subsequently, Doyon \cite{Doy26} defined quasiparticles for the KdV equation under a qualitatively different initial condition given by a soliton gas, and verified the counterpart of \eqref{eqn:asymptotic-scattering-overview}. 
        
        Let us recall the proposal from \cite{Agg25a}. Each quasiparticle corresponds to a fixed
		eigenvalue $\lambda_k$ of $\mathbf{L}(t)$; denote the associated unit eigenvector of $\mathbf{L}(t)$ by $\mathbf u_k(t)=(u_k(i;t))_i$. Under the random initial data we consider for the Toda lattice, the entries of $\mathbf{L}(t)$ are sampled independently at $t=0$. Thus, results of 
		Kunz--Souillard \cite{KS80} and Schenker \cite{Sch09} imply that $\mathbf{L}(0)$ is \emph{exponentially localized}; this means that each eigenvector $\mathbf u_k(0)$ admits a ``center'' index $\varphi_0(k)$ such that a
		bound of the form
		$|u_k(i;0)|\leq C e^{-c|i-\varphi_0(k)|}$ holds for every $i$. Using this, \cite{Agg25a} defined the initial location of the quasiparticle associated with $\lambda_k$ to be $Q_k(0)=q_{\varphi_0(k)}(0)$; in this way, the center $\varphi_0 (k)$ is viewed as the index of the Toda particle associated with the $k$-th quasiparticle. 

        If the eigenvectors $\mathbf{u}_k (t)$ remain localized, admitting centers $\varphi_t (k)$, throughout the dynamics, then one could similarly define the location of the $k$-th quasiparticle to be $Q_k(t)=q_{\varphi_t(k)}(t)$ for any time $t \ge 0$. Under an invariant initial condition, as was considered in \cite{Agg25a}, this localization persists since it holds at $t=0$.  In this paper, the initial data is not stationary, in which setting the entries of $\mathbf{L}(t)$ become highly correlated. Thus, eigenvector localization is no longer transparent. An important part of the present work is to justify that 
		a form of exponential localization remains valid under the dynamics, enabling one to define quasiparticles (as above) for all times.

		Proving this point and controlling the motion of the
		quasiparticles introduced above constitutes the bulk 
		of the work in proving Theorem \ref{thm:A-intro}. 
		We now state our main result in this direction. Let $x_\tau(\lambda;q)$ denote the 
		solution of the ordinary differential equation
		\begin{equation}\label{eqn:cfintro}
		x_0(\lambda;q)=q,
		\qquad
		\frac{d}{d\tau}x_\tau(\lambda;q)
		=
		\ve_\tau\bigl(\lambda;x_\tau(\lambda;q)\bigr).
		\end{equation}
		The following theorem shows that the
		 trajectories above are the limits as $T \rightarrow \infty$ 
		of quasiparticle trajectories, 
		as is broadly predicted by generalized hydrodynamics. 

		\begin{theorema}{B}[Theorem \ref{thm:ev}]
		\label{thm:B-intro}
		With high probability, for all $k$ satisfying
		$|Q_k(0)|\leq T \log T$, we have
		\[
		\sup_{t \in [0, T \log T]} \left|
		T^{-1}Q_k(t)
		-x_{t/T}\bigl(\lambda_k;T^{-1}Q_k(0)\bigr)
		\right|
		\leq T^{-1/4}\log T.
		\]
		\end{theorema}

		\begin{figure}[t]
		\centering
		\includegraphics[width=0.92\textwidth]{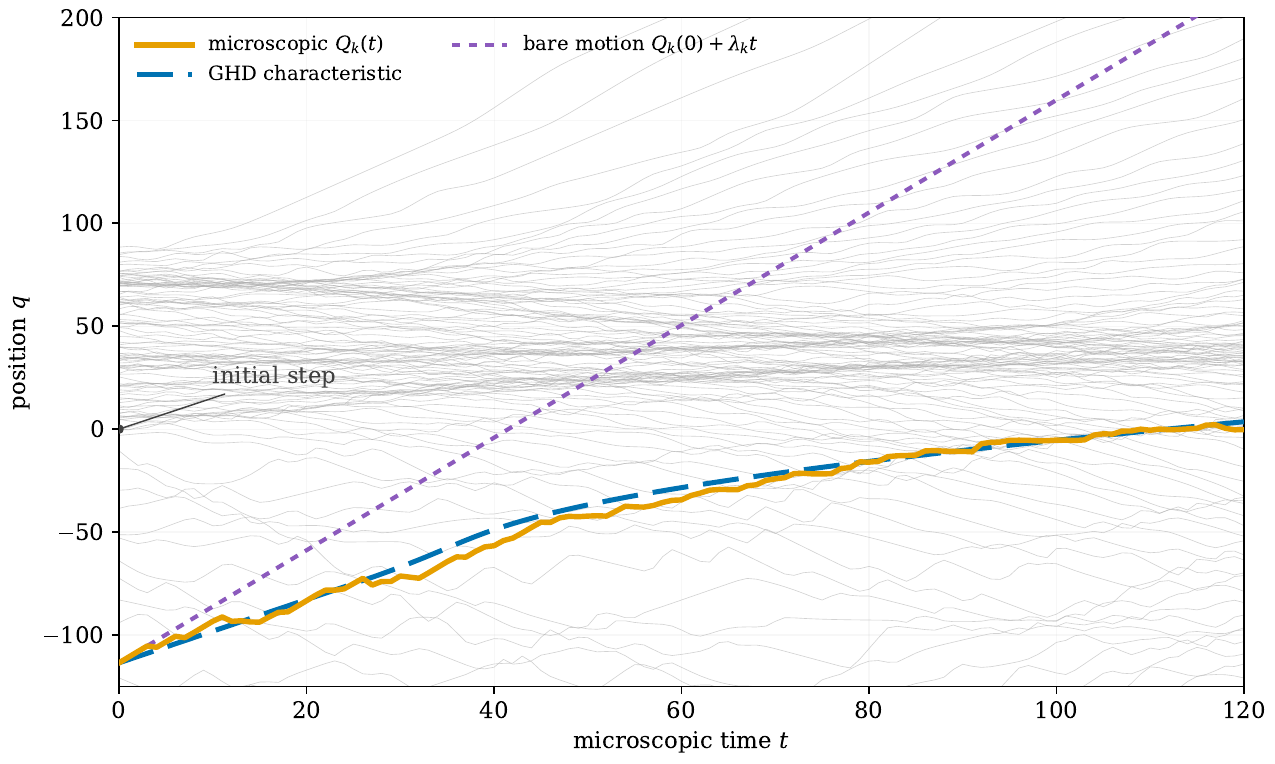}
		\caption{A numerical illustration of a quasiparticle in a Toda lattice with $600$ particles under Example
		\ref{ex:thermal-profiles}. The parameters $(\beta,\theta)$ there are
		$(0.2,0.2)$ to the left of the origin and $(2,2)$ to its right. The orange curve is the position $Q_k(t)$ of a fixed quasiparticle. The
		faint gray curves show the trajectories of the remaining Toda particles.
		The thick dashed blue curve is the GHD
		characteristic \eqref{eqn:cfintro} with the same initial position and spectral parameter
		$\lambda_k\approx 2.734$, while the purple curve is 
		the line of slope $\lambda_k$.}
		\label{fig:toda-step-ghd-characteristic}
		\end{figure}

		For the Toda lattice under thermal equilibrium,
		Theorem \ref{thm:B-intro} was
		established in \cite{Agg25}. Subsequent work of Aggarwal--Nicoletti
		\cite{AggNic26} addressed the finer question of analyzing the fluctuations of quasiparticles around their limiting trajectories. 
        
        The present work addresses the different issue of proving convergence to the generalized hydrodynamic
		evolution \eqref{eqn:rhoGHDintro_pref} for a broad class of nonstationary
		product initial profiles. In this setting, the local state and the effective
		velocity now vary along each trajectory; see Figure \ref{fig:toda-step-ghd-characteristic}. Central to the proof of Theorem \ref{thm:B-intro} is to prove the asymptotic scattering relation \eqref{eqn:asymptotic-scattering-overview} in our setup. 
		
        Before proceeding, let us mention that it would be of interest to understand the effect of integrability breaking for the Toda lattice. Consider perturbing the Toda Hamiltonian by
		$\lambda V$, where $0<\lambda\ll1$ and $V$ is a generic momentum-preserving
		interaction potential. Perturbative physics arguments 
		and numerical experiments predict that the GHD equation should remain the leading description on time scales $\tau \ll \lambda^{-2}$. On the longer 
		time scale $\tau \sim \lambda^{-2}$, a collision term should correct the GHD equation and
		drive a crossover toward conventional (chaotic) hydrodynamics, as discussed by
		Lopez-Piqueres--Ware--Gopalakrishnan--Vasseur \cite{LWGV21} and
		Doyon--Gopalakrishnan--M{\o}ller--Schmiedmayer--Vasseur \cite{DGMSV25}.
		Quantitatively pinpointing these phenomena for the Toda
		lattice remains open, both physically and mathematically; see
		\cite[Section IV.C]{DGMSV25}.

		\subsection{Proof strategy}
		


		We comment primarily on the proof of Theorem \ref{thm:B-intro},
		from which Theorem \ref{thm:A-intro} follows quickly; see Section
		\ref{subsec:mainpfs}. Before proving Theorem \ref{thm:B-intro},
		we must establish existence and uniqueness of the solution $(m_{\tau}, \ve_\tau)$ to 
		the GHD equations. H\"ubner--Doyon
		\cite{HubnerDoyon24} prove general well-posedness of the GHD
		equations \eqref{eqn:rhoGHDintro_pref}--\eqref{eqn:cra_pref} under
		certain assumptions on the scattering shift. The integrability
		assumption of \cite[Equation (3)]{HubnerDoyon24} does not hold
		for the Toda scattering shift (since $2 \log|\Lambda - \lambda| \rightarrow \infty$
		as $|\Lambda| \rightarrow \infty$),
		and thus we must modify their arguments to prove well-posedness in our setup,
		stated in Proposition \ref{prop:GHDint}; in particular, 
		in Lemma \ref{lem:H-contraction} we show the map \eqref{eqn:hop}, 
		whose fixed point is the desired solution,
		is a contraction on a certain space of functions with Gaussian decay. 
		Then, the core of the proof of Theorem \ref{thm:B-intro}
		is a bootstrap argument (outlined in items
		 \ref{item:proof-strategy-localization-to-scattering}--\ref{item:proof-strategy-trajectories-to-localization} below), 
		which can be viewed as a discrete
		analog of the fixed point approach taken in 
		\cite{HubnerDoyon24} and Proposition \ref{prop:GHDint}.

		As discussed in the previous subsection, a key step in the proof of
		Theorem \ref{thm:B-intro} is to establish eigenvector localization for $\L(t)$, 
		and the main difficulty in our setting is that the
		initial data are not stationary.
		The following three implications 
		are the central arguments of this work.
		\begin{enumerate}
			\item \label{item:proof-strategy-localization-to-scattering}
			Eigenvector localization implies the asymptotic
			scattering relation \eqref{eqn:asymptotic-scattering-overview}; see Proposition
			\ref{prop:discrete-continuum-scattering}, whose proof uses Proposition
			\ref{vn1n2out}.
			\item \label{item:proof-strategy-scattering-to-trajectories}
			The asymptotic scattering relation improves the approximation of the quasiparticle trajectories
			by the characteristic flow in Theorem \ref{thm:B-intro}; see Proposition
			\ref{prop:improv_M}.
			\item \label{item:proof-strategy-trajectories-to-localization}
			This improved quasiparticle approximation in turn improves the
			eigenvector localization estimate; see Proposition \ref{vni_out}.
		\end{enumerate}
		The apparent circularity of these implications is resolved by a 
		bootstrap argument that involves quantifying the improvements above.
		This allows us to simultaneously prove eigenvector localization, 
        		the asymptotic scattering relation, and Theorem \ref{thm:B-intro}.
		

%
		Before outlining proofs 
		of the three implications above, we provide a high level picture placing
		our arguments in the context of the inverse scattering method.
		To show Theorem \ref{thm:B-intro}, we replace the Toda lattice on $\mathbb{Z}$ 
		with one defined on
		the index interval $[N_1, N_2] \cap \mathbb{Z}$ with $-N_1,N_2 \sim T^2$
		(which will not affect the asymptotics on time scales of order $T \ll N_2-N_1$). 
		The classical inverse scattering method reconstructs the full
		Toda state $(\bm p(t);\bm q(t))$ from its eigenvalues and leading eigenvector
		entries, $\{(\lambda_k,  u_k(N_1;t))\}$, which evolve
		explicitly; this leads to an intricate high-dimensional
		Riemann--Hilbert problem, which we do not attempt to analyze.
		Instead, we encode the Toda state by its quasiparticle data $\{(\lambda_k, Q_k(t))\}$,
		and derive an approximately explicit mapping $\{(\lambda_k,  u_k(N_1;t))\} \mapsto \{(\lambda_k, Q_k(t))\}$ 
		(see \eqref{eqn:proof-strategy-transfer2})
		which has a natural continuum limit (see \eqref{eqn:Xtpref}); this is how 
		we establish items \ref{item:proof-strategy-localization-to-scattering} and 
		\ref{item:proof-strategy-scattering-to-trajectories} above. Establishing eigenvector 
		localization\footnote{The central role that eigenvector localization plays here may be 
		related to the well-definedness of the effective velocity $\ve$ 
		arising in the generalized hydrodynamic equation \eqref{eqn:rhoGHDintro_pref}.
		Indeed, in the notation used in Section \ref{subsec:GHD_2} below, 
		denoting the functions $\varsigma_0(\lambda)=1$ and $\varsigma_1(\lambda)=\lambda$, we have $\ve_\tau(\lambda;q)
		= \varsigma_1^{\dr,\tau,q}(\lambda) \varsigma_0^{\dr,\tau,q}(\lambda)^{-1}$,
		see \eqref{eqn:vedf}. For the Lax matrix obtained by
		extending the local environment at $(\tau,q)$,
		$(\varsigma_0^{\dr,\tau,q}(\lambda))^{-1}$ is proportional to the localization
		length for an eigenvector of eigenvalue $\lambda$; see
		\cite[Equation (A.6)]{Agg25}.
		Under this interpretation, delocalization would make $\ve_\tau(\lambda;q)$ meaningless, as its denominator would vanish.} is 
		a key step in 
		justifying this approximation; this is why we require step \ref{item:proof-strategy-trajectories-to-localization} above.

		
		Next, we will outline the proofs of items \ref{item:proof-strategy-localization-to-scattering}--\ref{item:proof-strategy-trajectories-to-localization}. Denote by 
		$\L^{[N_1,i-1]}(t)$ the principal block of $\L(t)$ corresponding 
		to indices in $[N_1, i-1]$.
		If $u_k(i;t)\ne0$, then since $\L(t)$ is a symmetric tridiagonal matrix (see Lemma \ref{lem:vkn1_formula}),
		\begin{equation}\label{eqn:proof-strategy-transfer}
		2\log|u_k(N_1;t)|
		=2\log|u_k(i;t)|+q_{N_1}(t)-q_i(t)
		-2\sum_{\mu\in\eig\L^{[N_1,i-1]}(t)}
		\log|\lambda_k-\mu|.
		\end{equation}
		A property of the Toda lattice due to Moser \cite{Mos75}
		is that the quantities $\mathcal{X}_k(t) \coloneqq -2 \log|u_k(N_1;t)| + q_{N_1}(t)$ evolve linearly, 
		$\mathcal{X}_k(t) = \mathcal{X}_k(0) + \lambda_k t$.
		Moreover, if we take $i = \varphi_t(k)$ to be the localization center of $\mathbf{u}_k(t)$, 
		since $|u_k(\varphi_t(k);t)|$ is not too small (see Definition \ref{def:loccent}) and
		since $Q_k(t) = q_{\varphi_t(k)}(t)$ by definition, the right hand side of the formula 
		\eqref{eqn:proof-strategy-transfer} states
		\begin{equation}\label{eqn:proof-strategy-transfer2}
		\mathcal{X}_k(t)  = Q_k(t) +  2\sum_{\mu\in\eig\L^{[N_1,i-1]}(t)}
		\log|\lambda_k-\mu| + \text{controlled error}.
		\end{equation}

		On the macroscopic scale,
		the continuum
		analog
		of \eqref{eqn:proof-strategy-transfer2} is
		heuristically (see \eqref{eqn:Xt}) defined as
		\begin{equation}\label{eqn:Xtpref}
		X_\tau(\lambda,q)
		=
		q+2 \int_{-\infty}^q \int_{\mathbb R}\log|\lambda-\mu|\,m_\tau(\mu,q')\,d\mu\,dq'.
		\end{equation}
		In this coordinate, the 
		characteristic (which, according to Theorem \ref{thm:B-intro}, should be the limiting quasiparticle trajectory)
		$x_\tau(\lambda;q)$ solving \eqref{eqn:cfintro} becomes a linear evolution,
		$X_\tau\bigl(\lambda,x_\tau(\lambda;q)\bigr)
		=
		X_0(\lambda,q)+\lambda \tau$ (see Lemma \ref{xt2}). Thus, 
		$X_\tau\bigl(\lambda_k, x_\tau(\lambda_k;q)\bigr)$
		evolves in the same way as $\mathcal{X}_k(t) $, and moreover 
		the right hand side of \eqref{eqn:proof-strategy-transfer2} is a discretization of the 
		right hand side of \eqref{eqn:Xtpref}.


		In what follows, we set $\mathfrak{M} = T^{3/4}$, and recall $\mathfrak{A}$ is a large enough
		but fixed constant
		representing the scale of average particle spacing. 
		First, we assume we already have a
		quantitative eigenvector localization estimate at time $t$, namely the bound 
		$|u_k(i;t)| \leq e^{\mathfrak{M} - \mathfrak{A}^{1/10} |i-\varphi_t(k)|}$.
		This localization estimate allows us (see Lemma \ref{lem:full-to-truncation}) to deduce that 
		all but $O(\mathfrak{M}/\mathfrak{A}^{1/10})$ eigenvalues of the block $\L^{[N_1,\varphi_t(k)-1]}(t)$ can be 
		closely approximated by distinct
		eigenvalues of the full matrix $\L(t)$ (as only eigenvalues with localization centers
		within $O(\mathfrak{M}/\mathfrak{A}^{1/10})$ of the 
		boundary of $[N_1, i-1]$ will not be accounted for). 
		This, together with the fact that $\{q_j(t)\}$ are approximately ordered
		(which we will prove during the bootstrap argument),
		turns \eqref{eqn:proof-strategy-transfer2} into
		\begin{equation}\label{eqn:Xktbetter}
		\mathcal{X}_k(t)
		=
		Q_k(t)
		+2\sum_{j:Q_j(t)<Q_k(t)}
		\log|\lambda_k-\lambda_j|
		+O(\mathfrak{M}/\mathfrak{A}^{1/10}).
		\end{equation}
		Subtracting this identity at time $0$ from the one at time $t$, together with the fact 
		that $\mathcal X_k(t) - \mathcal X_k(0) = \lambda_k t $,
		gives the asymptotic scattering relation
		\eqref{eqn:sc}, establishing implication
		\ref{item:proof-strategy-localization-to-scattering}. 
		The size of the error, $O(\mathfrak{M}/\mathfrak{A}^{1/10})$, will be
		important for obtaining an improvement in the next step.
        
		Next, if we assume the bound in Theorem \ref{thm:B-intro}, 
		namely $\left| T^{-1} Q_j(t) -x_{t/T}\bigl(\lambda_j;T^{-1}Q_j(0)\bigr) \right| \leq  \frac{\mathfrak{M}}{T} \log T$,
		then we may 
		turn the summation on the right hand side of \eqref{eqn:Xktbetter} into the 
		integral on the right hand side of \eqref{eqn:Xtpref} up to a controlled error (using Proposition \ref{prop:scat_conc} below).
		This in turn, together with the fact that $\mathcal{X}_k(t)$ and $X_\tau\bigl(\lambda_k, x_\tau(\lambda_k;q)\bigr)$ 
		obey the same linear evolution,
		allows us to bound $|x_\tau(\lambda;T^{-1} Q_k(0)) - T^{-1} Q_k(T\tau) | \leq \frac{\mathfrak{M}}{2 T} \log T$; 
		note the error is smaller by a factor of $2$ than the assumed weak version of Theorem
		\ref{thm:B-intro}, so this establishes implication
		\ref{item:proof-strategy-scattering-to-trajectories}. Similarly to the previous step, the improvement can be traced
		back to the fact that $\mathfrak{A}$ is large. Indeed, the factor
		$1/\mathfrak{A}^{1/10} $ appearing in the error in \eqref{eqn:Xktbetter} is a small constant, 
		and each newly accumulated error in step \ref{item:proof-strategy-scattering-to-trajectories} should heuristically be 
		bounded by $C \mathfrak{M}/\mathfrak{A}$; for example, 
		replacing $Q_j(t)$ by its limit $x_\tau(\lambda_j;T^{-1}Q_j(0))$ in the subscript of 
		the sum in \eqref{eqn:Xktbetter} 
		heuristically leads 
		to an error of size $O(\mathfrak{M}/\mathfrak{A})$,
		since we will see that only quasiparticles in a window of size $\mathfrak{M}$ around $Q_k (t)$ contribute to the error
		(and there are $O(\mathfrak{M}/\mathfrak{A})$ such particles, as their spacing is of order $\mathfrak{A}$).
		The argument proving implication
		\ref{item:proof-strategy-trajectories-to-localization} is similar; one again uses the transfer matrix formula
		\eqref{eqn:proof-strategy-transfer}, now together with the improved
		quasiparticle approximation estimate, to obtain the stronger localization estimate
		$|u_k(i;t)| \leq e^{\mathfrak{M}/2 - \mathfrak{A}^{1/10} |i-\varphi_t(k)|}$.
		

		These three self-improvement bounds together with 
		the short time continuity stated in Proposition \ref{prop:cont}
		show that the upper bound on $|x_{\tau} (\lambda; T^{-1} Q_k(0)) - T^{-1}Q_k(T\tau )|$ 
		cannot
		degenerate (see Theorem \ref{ttn} for a precise statement in terms of a certain stopping time),
		proving Theorem \ref{thm:B-intro}. A more detailed explanation of this argument is provided in 
		Section \ref{sec:outline} below.

	\subsection*{Acknowledgements}

	The authors thank Herbert Spohn for stimulating conversations and helpful comments on an initial draft of this paper. They also thank Alexei Borodin, Benjamin Doyon, Percy Deift, Tamara Grava, and Ken McLaughlin for valuable discussions. The work of Amol Aggarwal was partially supported by a Packard Fellowship for Science and Engineering. The work of Matthew Nicoletti was partially supported by an NSF postdoctoral fellowship, grant No. DMS 2402237.

	\subsection*{Statement of AI use}
	
	The authors conceived of the overarching structure of the proof presented in Section \ref{sec:outline}, and most of the statements of the ``lower lying'' results later in the paper needed to prove them. They used large language models to assist in writing the proofs of the latter ``lower lying'' results, in the following two ways. First, if the authors suspected that such a result could be shown directly from known results in a self-contained literature, the authors might ask GPT-5.6 Pro to provide a first draft of a proof, which the authors would polish (the proof of Lemma \ref{lem:local-energy-from-localization} was perhaps the most substantive example of this). Second, for more specialized results, the authors might write a proof outline to GPT-5.6 Pro, which would draft a more thorough argument (that the authors would then polish), by properly setting up notation and completing small gaps in the argument. Once the mathematical content was complete, the authors used GPT-5.6 Pro, and GPT-6 Astra when it was released, to proofread, edit, and check correctness of the work. GPT-6 Astra was also used to produce Figure \ref{fig:toda-step-ghd-characteristic}. As indicated above, all AI-assisted arguments and text were checked and revised by the authors, who take full responsibility for the contents of the paper.

	\section{Results}

	\subsection{Model definition}
		\label{subsec:assumptions}

		For any $-\infty\leq a\leq b\leq+\infty$, set
		$\llbracket a,b\rrbracket=[a,b]\cap\mathbb Z$.
		Let $-\infty\leq N_1<N_2\leq+\infty$ be integers or $\pm \infty$. The Toda lattice on
		$\llbracket N_1,N_2\rrbracket$ defines
	Hamiltonian dynamics
	\(t \mapsto (\bm p(t);\bm q(t))\), with
		momenta \(\bm p(t)=(p_{N_1}(t),p_{N_1+1}(t),\dots,p_{N_2}(t))\) and
		positions \(\bm q(t)=(q_{N_1}(t),q_{N_1+1}(t),\dots,q_{N_2}(t))\),
		given by the ordinary differential equations (ODEs)
		\begin{equation}\label{eqn:TD}
		\partial_tq_i(t)=p_i(t),
		\qquad
		\partial_tp_i(t)=e^{q_{i-1}(t)-q_i(t)}-e^{q_i(t)-q_{i+1}(t)},
		\end{equation}
		with the convention that $q_{N_1-1}\equiv -\infty$ and $q_{N_2+1} \equiv +\infty$
		if $N_1$ or $N_2$ is finite.
	When $N_1$ and $N_2$ are finite, the Picard Lindel\"{o}f theorem gives existence and uniqueness for all $t\geq0$ under arbitrary initial data; see, for example, the proof of \cite[Theorem 12.6]{OINL}.
	Next, we define the Flaschka variables.
	\begin{definition}\label{def:flandL}
	Denote the \emph{Flaschka variables} by $(\bm{a}(t); \bm{b}(t)) = (a_{N_1}(t),\dots, a_{N_2}(t); b_{N_1}(t),\dots, b_{N_2}(t))$, with 
	\begin{equation}\label{eqn:flv}
	a_i(t) = e^{\frac{1}{2} (q_i(t) - q_{i+1}(t))}; \qquad b_i(t) = p_i(t).
	\end{equation}
\noindent
	We denote the associated Lax matrix at any time $s \ge 0$ by $\L(s)= [L_{ij}(s)]_{i,j \in \llbracket N_1, N_2 \rrbracket}$, defined by
	\begin{equation}\label{eqn:L}
	L_{i, i+1}(s) = L_{ i+1, i}(s) = a_i(s); \qquad L_{i i}(s) = b_i(s).
	\end{equation}
All entries of $\L(s)$ not determined by \eqref{eqn:L} are equal to $0$. 
\end{definition}
We also may take $N_1 = -\infty$ and $N_2 = +\infty$, 
	which we refer to as the Toda lattice on $\mathbb{Z}$. 
For the Toda lattice on $\mathbb{Z}$, \cite[Proposition 4.7]{Agg25a} constructs
 the dynamics as limits of finite interval evolutions and proves their well posedness whenever
\begin{equation}\label{eqn:infinite-toda-growth}
\limsup_{j\in\mathbb Z}
\frac{|a_j(0)|+|b_j(0)|}{(|j|+1)^{\mathfrak p}}<\infty
\qquad\text{for some }\mathfrak p\in(0,1).
\end{equation}
Thus, under \eqref{eqn:infinite-toda-growth}, the Toda lattice is well defined on $\mathbb Z$.

		A fundamental feature of the Lax matrix is that its eigenvalues are preserved under the Toda dynamics. This was originally due to \cite{Fla74}; see also \cite[Section 2]{Mos75}.
	
	\begin{lem}[\cite{Fla74,Mos75}]

		\label{ltt}

		For any real numbers $t, t' \in \mathbb{R}_{\ge 0}$, we have $\eig \L (t) = \eig \L(t')$.

	\end{lem}

Thus, the eigenvalues $\eig \L(t) \coloneqq (\lambda_1 \ge \lambda_2 \ge . . . \ge \lambda_N)$ are preserved under the 
Toda dynamics, meaning that $\Tr(\L^k)$ is as well, for each $k \ge 0$. 
This motivates the following definition of local charge densities. The charges
are the diagonal entries of $\L(t)^m$; they are local quantities whose
total is preserved.
	
\begin{definition}[Local charges]\label{kti00}
	
 Fix an integer $m \geq 0$, index $i \in \llbracket N_1, N_2 \rrbracket$, and real number $t \geq 0$. Define the $m$-th \emph{local charge} $\mathfrak{k}_i^{[m]}(t)$ of $\L(t)$ at $i$ by setting
\begin{flalign} 
\label{mki} 
\mathfrak{k}_i^{[m]}(t) = [\L(t)^m]_{ii}.
\end{flalign}

\end{definition}

	We will study the Toda lattice under random initial data, in which the
	Flaschka variables are mutually independent and subject to the following
	conditions.
\begin{assumption}[Macroscopic entry law profiles]\label{ass:entry-law-profiles}
For each $y\in\mathbb R$, let $p_{\dg}(\cdot,y)$ be a probability density on $\mathbb R$. For 
each $y \in \mathbb{R}$, set
$\mu_{\dg}(dx,y):=p_{\dg}(x,y)\,dx$, and let $\mu_{\sd}(\cdot,y)$ be a probability measure on $\mathbb R_{>0}$.
Fix constants $\mathfrak K>1$, $\epsilon>0$, and $\mathfrak A>1$. Assume first that
$\int_0^\infty |\log x|\,\mu_{\sd}(dx,y)<\infty,
 y\in\mathbb R$,
and define the mean \emph{stretch}
\[
\widetilde\alpha(y):=-2\int_0^\infty\log x\,\mu_{\sd}(dx,y).
\]
Assume that
\begin{equation}\label{eqn:alph_bd}
\mathfrak K^{-1}\mathfrak A\le \inf_{y\in\mathbb R}\widetilde\alpha(y)
\le \sup_{y\in\mathbb R}\widetilde\alpha(y)
\le \mathfrak K\mathfrak A.
\end{equation}
Assume there is a finite set $\mathcal P = \{y_1,\dots, y_{\mathfrak{m}}\} \subset\mathbb R$, $\mathfrak{m} \in \mathbb{Z}_{\geq 0}$, such that, whenever
$y'<y$ and $[y',y]\cap\mathcal P =\varnothing$,
\begin{equation}\label{eqn:alpha-regularity}
|\widetilde\alpha(y)-\widetilde\alpha(y')|
\le  \mathfrak K \mathfrak A |y-y'|.
\end{equation}
Define the cumulative stretch map and the mean stretch in particle coordinates by
	\begin{equation}\label{eqn:Fstr}
	F_{\mathrm{str}}(u):=\int_0^u\widetilde\alpha(y)\,dy,
	\qquad
	\alpha(q):=\widetilde\alpha(F_{\mathrm{str}}^{-1}(q)).
	\end{equation}
	The inverse in \eqref{eqn:Fstr} is well defined by \eqref{eqn:alph_bd}.

For any $x\in\mathbb R$, we further assume the following uniform bounds on the entry laws:
	\begin{flalign}
		\label{munu} 
		\begin{aligned} 
			& \qquad \qquad \qquad  \displaystyle\sup_{y \in \mathbb{R}} p_{\dg}(x,y)  \le \mathfrak{K} e^{-x^2 / \mathfrak{K}}; \\
			&\displaystyle\sup_{y \in \mathbb{R}} \left( \int_{0}^{1} x^{-\epsilon/\widetilde\alpha(y)}\,\mu_{\sd}(dx,y)
			+\int_{0}^{\infty} e^{\epsilon x^2}\,\mu_{\sd}(dx,y)\right) \le \mathfrak K.
		\end{aligned} 
	\end{flalign}
	
	\noindent Assume also that
\begin{equation}\label{eqn:mudgc}
\sup_{X\in[-\mathfrak K,\mathfrak K]}
\inf_{y\in\mathbb R}
\inf_{|x-X|\le \mathfrak K^{-1}}
p_{\dg}(x,y)>\mathfrak K^{-1}.
\end{equation}

	\noindent Finally, for $y'<y$ satisfying $[y',y]\cap\mathcal P=\varnothing$, assume
	\begin{equation}\label{eqn:W2bd}
	\begin{aligned}
	W_2(\mu_{\dg}(\cdot;y), \mu_{\dg}(\cdot;y')) &\leq \mathfrak K |y-y'|; \\
	W_2(\mu_{\sd}(\cdot;y), \mu_{\sd}(\cdot;y')) &\leq \mathfrak K |y-y'|.
	\end{aligned}
	\end{equation}
		Here $W_2$ denotes the Wasserstein distance. 		
		\end{assumption}

\begin{assumption}[Infinite Toda lattice initial data]\label{ass:model_inf}
Adopt Assumption \ref{ass:entry-law-profiles}, and let $T\geq1$ be a real number (which will index a sequence of models). We consider the Toda lattice $(\bm p(t);\bm q(t))$ on $\mathbb Z$ and denote its associated Flaschka variables by $(\bm a(t);\bm b(t))$, as in Definition \ref{def:flandL}. Sample all entries mutually independently: for each $i\in\mathbb Z$, sample the diagonal entry $b_i(0)=L_{ii}(0)$ and the superdiagonal entry $a_i(0)=L_{i,i+1}(0)$ of $\L(0)$ according to $\mu_{\dg}(dx,i/T)$ and $\mu_{\sd}(dx,i/T)$, respectively. Set $p_i(0)=b_i(0)$ and $q_0(0)=0$; the identities $a_i(0)=e^{(q_i(0)-q_{i+1}(0))/2}$ then uniquely determine $\bm q(0)$.
\end{assumption}

We now comment briefly on Assumption \ref{ass:entry-law-profiles}. The law $\mu_{\sd}(\cdot,y)$ governs the positive superdiagonal entry $a_i(0)$ near lattice coordinate $y=i/T$, and thus governs the particle spacing $q_{i+1}(0)-q_i(0)$, while $\mu_{\dg}(\cdot,y)$ governs the diagonal entry $b_i(0)=p_i(0)$. The quantity $\widetilde\alpha(y)$ is the mean particle spacing $q_{i+1}(0)-q_i(0)$, with index $i\approx yT$; the quantity $\alpha(q)$ is the mean particle spacing $q_{i+1}(0)-q_i(0)$ around the macroscopic particle location $q\approx q_i(0)/T$ on $\mathbb R$. The latter statement is due to the fact, shown in Lemma \ref{lem:Estr}, that the cumulative stretch map $F_{\mathrm{str}}$ in
		\eqref{eqn:Fstr} is the continuum limit of the map
		$\mathbb R\ni i/T\mapsto q_i(0)/T\in\mathbb R$. Thus, \eqref{eqn:alph_bd} places the system in a uniformly positive density regime whose mean spacing is of order $\mathfrak A$. The two
moment bounds in the second line of \eqref{munu} provide tail decay on $ |\log X|$ and $X$, where $X \sim \mu_{\sd}(\cdot, y)$, uniformly in $y$. The estimate \eqref{munu} provides a uniform upper bound on the density $p_{\dg}$ of the diagonal law, and \eqref{eqn:mudgc} provides a uniform lower bound for this density on a small interval. The bounds \eqref{eqn:alpha-regularity} and \eqref{eqn:W2bd} state that the macroscopic profile varies in a Lipschitz manner away from finitely many possible discontinuity points. The parameter $\mathfrak A$ will be required to exceed a threshold depending only on $\mathfrak K$ and $\epsilon$ in the results below. Finally, we remark that the bounds in \eqref{munu} imply \eqref{eqn:infinite-toda-growth} almost surely under Assumption \ref{ass:model_inf}.

\begin{example}[Thermal profiles]\label{ex:thermal-profiles}
Fix $C>1$, let $\beta,\vartheta:\mathbb R\to\mathbb R_{>0}$ be
$C$-Lipschitz away from finitely many discontinuities and satisfy
\[
C^{-1}\leq\beta(y)\leq C,
\qquad
C^{-1}\leq \vartheta(y)\leq C,
\qquad y\in\mathbb R,
\]
and let $\theta(y) \coloneqq \mathfrak{A}^{-1}\vartheta(y)$. Consider the laws
\[
p_{\dg}(b,y)=C_{\beta(y)}e^{-\beta(y)b^2/2},
\qquad
\mu_{\sd}(da,y)=C_{\beta(y),\theta(y)}
a^{2\theta(y)-1}e^{-\beta(y)a^2}\,da,
\]
where the constants normalize the two measures. In this case,
\[
\widetilde\alpha(y)=\log\beta(y)-\psi(\theta(y)),
\]
where $\psi$ is the digamma function. Thus, for $\mathfrak A$ sufficiently
large, the mean stretch is comparable to $\mathfrak A$. As discussed in 
Section \ref{subsec:over}, the predictions motivating generalized hydrodynamics were originally
formulated for initial data drawn from a smoothly varying invariant measure;
this example is of that form, with the local invariant measure being thermal.
After choosing
$\mathfrak K$ sufficiently large depending only on $C$, then choosing
$\epsilon>0$ sufficiently small, and then choosing $\mathfrak{A} \geq \mathfrak{A}_0(\mathfrak{K}, \epsilon)$ large enough, 
these laws satisfy Assumption
\ref{ass:entry-law-profiles}.
\end{example}

\begin{example}[Equally spaced particles]\label{ex:equally-spaced-particles}
Fix $C>1$, and let $\beta:\mathbb R\to\mathbb R_{>0}$ be
$C$-Lipschitz away from finitely many discontinuities and satisfy
$C^{-1}\leq\beta(y)\leq C$ for every $y\in\mathbb R$. Set
\[
\mu_{\sd}(da,y)=\delta_{e^{-\mathfrak A/2}}(da),
\qquad
p_{\dg}(b,y)=\left(\frac{\beta(y)}{2\pi}\right)^{1/2}
e^{-\beta(y)b^2/2}.
\]
This gives initially equally spaced particles, $q_i(0)=\mathfrak A i$,
with independent Gaussian velocities whose inverse temperature varies
with the macroscopic location. Observe that this initial data is not a
smoothly varying invariant measure, since the initial particle positions are
fixed deterministically. After choosing $\mathfrak K$ sufficiently
large depending only on $C$ and $\epsilon>0$ sufficiently small, these laws
satisfy Assumption \ref{ass:entry-law-profiles}.
\end{example}

	\subsection{Hydrodynamic limit}
	
	\label{subsec:mainthm1}

Unless stated otherwise, constants may depend only on $\mathfrak K$ and
$\epsilon$.

We first define the density of states used in the limiting object. The proof of
the following lemma is given in Appendix \ref{app:frozen-dos}.
\begin{lem}\label{lem:frozen-dos-intro}
For $y\in\mathbb R$ and $n\geq1$, let $\mathbf{S}_y^{(n)}$ be the $n\times n$
symmetric tridiagonal matrix whose diagonal and superdiagonal entries are
sampled independently from $\mu_{\dg}(\cdot,y)$ and
$\mu_{\sd}(\cdot,y)$, respectively. There is a unique function
$\widetilde\varrho_0:\mathbb R^2 \rightarrow \mathbb R_{\geq0}$
such that $\lambda\mapsto\widetilde\varrho_0(\lambda,y)$ is a probability
density for every $y\in\mathbb R$ and, for every bounded continuous
$h:\mathbb R\to\mathbb R$,
\[
\frac1n\Tr h(\mathbf{S}_y^{(n)})
\longrightarrow
\int_{\mathbb R}h(\lambda)\widetilde\varrho_0(\lambda,y)\,d\lambda
\]
in probability as $n\to\infty$.
\end{lem}

We may now define the initial condition for the
macroscopic PDE which will describe the Toda lattice
when $\L(0)$ is sampled under Assumption
\ref{ass:model_inf}.
\begin{definition}
\label{rho0}
Let $\widetilde\varrho_0$ be the density supplied by Lemma
\ref{lem:frozen-dos-intro}, and set
\[
\varrho_0(\lambda,q):=\widetilde\varrho_0
(\lambda,F_{\mathrm{str}}^{-1}(q)),
\qquad
m_0(\lambda,q):=\alpha(q)^{-1}\varrho_0(\lambda,q).
\]
\end{definition}

\begin{rem}
The density $m_0(\lambda,q)$ is the unique one which satisfies
\begin{equation}\label{eqn:dos_def}
\frac{1}{T} \sum_{i \in \mathbb{Z}} f(q_i(0) T^{-1}) \mathfrak{k}_i^{[n]}(0)
\longrightarrow
\int_{-\infty}^{\infty}\int_{-\infty}^{\infty}
f(q)\lambda^n m_0(\lambda,q)\,d\lambda\,dq,
\end{equation}
for all smooth functions $f$ with compact support and all integers $n\geq0$,
where the convergence is in probability as $T\to\infty$. This follows from
Proposition \ref{prop:app-time-zero-dos} (combined with Lemmas \ref{rho0estimate}
and \ref{l0eigenvalues}) below.
\end{rem}

Under Assumption \ref{ass:entry-law-profiles}, the densities in Definition
\ref{rho0} satisfy the following. The proof is given in Appendix
\ref{app:frozen-dos}.
\begin{lem}
\label{rho0estimate}
There exist constants $C,c>0$ such that, for all
$(\lambda,q)\in\mathbb R^2$,
\begin{equation}\label{eqn:vrhob}
0\leq\varrho_0(\lambda,q)\leq C e^{-c\lambda^2},
\qquad
0\leq m_0(\lambda,q)\leq C\mathfrak A^{-1}e^{-c\lambda^2}.
\end{equation}
\end{lem}

To set up the hydrodynamic limit for the Toda lattice under Assumption
\ref{ass:model_inf}, the following proposition jointly (and uniquely)
constructs functions $m_t$ and $\ve_t$ satisfying
\eqref{eqn:rhoGHDintro_pref} and \eqref{eqn:cra_pref}. It encompasses a class of initial data
$m_0(\lambda,x)$ which includes those arising under Assumptions
\ref{ass:model_inf} and \ref{ass:model}.

If $\mathfrak c>0$ and $f:\mathbb R\to\mathbb R$, set
\begin{equation}\label{eqn:weighted-sup-norm}
\|f\|_{\mathfrak c}:=
\sup_{\lambda\in\mathbb R}e^{\mathfrak c\lambda^2}|f(\lambda)|.
\end{equation}
For a nonnegative function $\nu$, the space
$L^2(\mathbb R,\nu(\lambda)\,d\lambda)$ consists of the functions
$f:\mathbb R\to\mathbb R$ such that
$\int_{\mathbb R}|f(\lambda)|^2\nu(\lambda)\,d\lambda<\infty$.

The next proposition is proved in Section
\ref{subsec:ghd-uniqueness} below.
\begin{prop}[Well posedness of the GHD equations]\label{prop:GHDint}
Let $C_0,c_0>0$ and $\mathfrak c\in [c_0/10,c_0/2]$. There exist constants
$C=C(C_0,c_0)>1$ and
$\mathfrak A_0=\mathfrak A_0(C_0,c_0)>1$ such that the following holds
whenever $\mathfrak A\geq\mathfrak A_0$. Suppose $m_0(\lambda,x)$ satisfies
\begin{equation}\label{eqn:m0bd}
0\leq m_0(\lambda,x)\leq C_0\mathfrak A^{-1}e^{-c_0\lambda^2},
\qquad (x,\lambda)\in\mathbb R^2.
\end{equation}
There is a unique function $N_t(\lambda,x)$ with
$\|N_t(\cdot,x)\|_{\mathfrak c}<\infty$ for every $x\in\mathbb R$ and
$t\geq0$, having the following properties.
\begin{enumerate}
\item We have $N_0(\lambda,x)=\int_0^x m_0(\lambda,y)\,dy$.
\item We have
\begin{equation}\label{eqn:nlip}
e^{\mathfrak c\lambda^2}
\bigl|N_t(\lambda,x)-N_s(\lambda,y)\bigr|
\leq C\mathfrak A^{-1}\bigl(|x-y|+|t-s|\bigr)
\end{equation}
uniformly for $\lambda,x,y\in\mathbb R$ and $s,t\in\mathbb R_{\geq0}$.
Denoting the almost everywhere $x$-derivative
$m_t(\lambda,x):=\partial_xN_t(\lambda,x)$ for
any $(\lambda, t) \in \mathbb{R} \times \mathbb{R}_{\geq 0}$, we have $m_t(\lambda, x)\geq0$ 
for almost
every $x \in \mathbb{R}$. For every $t\geq0$ and almost every $x\in\mathbb R$, there is a
unique
$\ve_t(\cdot;x)\in L^2(\mathbb R,e^{-c_0\lambda^2/2}\,d\lambda)$ satisfying
\begin{equation}\label{eqn:cra}
\ve_t(\lambda;x)=\lambda+2\int_{-\infty}^{\infty}
\log|\lambda-\lambda'|\,m_t(\lambda',x)
\bigl(\ve_t(\lambda';x)-\ve_t(\lambda;x)\bigr)\,d\lambda'.
\end{equation}
\item For each $\lambda \in \mathbb{R}$,
\begin{equation}\label{eqn:GHDpot_unique}
\partial_tN_t(\lambda,x)=-\ve_t(\lambda;x)m_t(\lambda,x)
\end{equation}
holds almost everywhere in $(t,x)$. In particular, $(m_t,\ve_t)$ solves
\begin{equation}\label{eqn:rhoGHDintro}
\partial_t m_t(\lambda,x)
+\partial_x\bigl(\ve_t(\lambda;x)m_t(\lambda,x)\bigr)=0,
\qquad t\geq0,\ x,\lambda\in\mathbb R,
\end{equation}
in the weak sense.
\end{enumerate}
\end{prop}

			The main theorem of this work is the following. It says that averages of local charges are described as integrals against the solution of \eqref{eqn:rhoGHDintro}. Theorem \ref{thm:local-charge-testfn-intro} is proved in Section \ref{subsec:mainpfs} below.
	
		\begin{thm}[Hydrodynamic limit for local charges]\label{thm:local-charge-testfn-intro}
			Adopt Assumption \ref{ass:model_inf}, let $m_0$ be as in Definition
			\ref{rho0}, and let $N_t$ and $m_t$ be as in Proposition
			\ref{prop:GHDint}. Recall the local charges from Definition \ref{kti00}.
			There exists $\mathfrak A_0=\mathfrak A_0(\mathfrak K,\epsilon)>1$ such that the following holds whenever $\mathfrak A\ge\mathfrak A_0$. Let $g : \mathbb{R} \rightarrow \mathbb{R}$ denote a
		smooth function with compact support. Fix \(n\in\mathbb Z_{\geq 0}\) and
		\(\tau > 0\). Then we have the convergence in probability as $T\to\infty$,
		\[
		\frac{1}{T}
		\sum_{i \in \mathbb{Z}} g(q_i(T \tau) T^{-1}) \mathfrak{k}_i^{[n]}(T \tau) 
			\rightarrow \int_{-\infty}^{\infty} \int_{-\infty}^{\infty} g(q)  \lambda^n m_\tau(\lambda,q)\,d\lambda d q.
		\]
		\end{thm}

We have the following corollary (of the proof) of Theorem \ref{thm:local-charge-testfn-intro}. It would follow formally from Theorem 
\ref{thm:local-charge-testfn-intro} if one could set $g$ equal to the indicator of an interval of size of order $ T^{-1/6}$,
corresponding to a physical interval of size $T^{5/6}$. 
The exponent $5/6$ is not optimal, and we make no attempt to optimize it.
Corollary \ref{cor:local-charge-intro} is proved in Section
\ref{subsec:mainpfs} below.

		\begin{cor}\label{cor:local-charge-intro}
			Adopt the notation and assumptions of Theorem
			\ref{thm:local-charge-testfn-intro}. Fix \(n\in\mathbb Z_{\geq 0}\),
		\(\tau > 0\), \(5/6<a<a'<1\), $x \in \mathbb{R}$, and assume that
		$x$ is a continuity point of $y\longmapsto\int_{\mathbb R}\lambda^n m_\tau(\lambda,y)\,d\lambda$.
		Let \(J_T\subset
		\mathbb{R}\) be an interval such that
		\[
		T^a\le |J_T|\le T^{a'},
		\qquad
		\sup_{y\in J_T}\left|T^{-1}y-x\right|\longrightarrow 0.
		\]
		Then we have the convergence in probability as $T\to\infty$,
		\[
		\frac{1}{|J_T|}
		\sum_{i: q_i(T \tau) \in J_T} \mathfrak{k}_i^{[n]}(T \tau) 
		\rightarrow \int_{\mathbb R}\lambda^n m_\tau(\lambda,x)\,d\lambda.
		\]
		\end{cor}

	\subsection{Quasiparticles}
	\label{subsec:mainthm2}
	
		Next, we fix our running assumptions
		on the initial data for a finite
		version of the Toda lattice which we will study in order to deduce the
		results stated in the previous section. The following assumption is
		parallel to Assumptions \ref{ass:entry-law-profiles} and
		\ref{ass:model_inf}, except (in addition to being on a finite interval $\llbracket N_1, N_2 \rrbracket$, instead of the full line)
		it imposes
		the thermal profile of Example \ref{ex:thermal-profiles}, with $\beta$ and
		$\vartheta$ there taken to be constants here, sufficiently far from the origin.
		This serves mainly as a convenience; Theorem
		\ref{thm:ev} (and its proof) below shows
		that it does not asymptotically affect the behavior 
		 of the Toda lattice in intervals around the origin with length of order T.
	
	\begin{assumption}[Finite Toda lattice initial data]\label{ass:model}
	Let $\mu_{\dg}$ and $\mu_{\sd}$ be families of measures as in Assumption \ref{ass:entry-law-profiles}. 
	Fix a large real number \(T\ge 2\) and let \(N:=\lfloor 2 T^2 \rfloor \) be an integer.  Let
	\[
	N_1:=-\lfloor N/2\rfloor,\qquad N_2:=N_1+N-1.
	\]
	\noindent
	We consider the Toda lattice $(\bm{p}(t);\bm{q}(t))$ on the finite interval \(\llbracket N_1,N_2\rrbracket\), with initial data $(\bm{p}(0);\bm{q}(0))$. 
	Denote the associated Flaschka variables by $(\bm{a}(t); \bm{b}(t))$, defined as in Definition \ref{def:flandL}.
For each $y\in\mathbb R$, let $p_{\dg}^{(T)}(\cdot,y)$ be a probability density on $\mathbb R$, set $\mu_{\dg}^{(T)}(dx,y):=p_{\dg}^{(T)}(x,y)\,dx$, and let $\mu_{\sd}^{(T)}(\cdot,y)$ be a probability measure on $\mathbb R_{>0}$.
Assume that  
 $$\mu_{\dg}^{(T)}(dx,y) = \mu_{\dg}(dx,y) ,\qquad \mu_{\sd}^{(T)}(dx,y)= \mu_{\sd}(dx,y)  \qquad \text{for } |y| \leq (\log T)^5,$$
 and
$$\mu_{\sd}^{(T)}(dx,y) = C_{\mathfrak A}x^{2/\mathfrak A-1}e^{-x^2}\,dx,
\qquad
\mu_{\dg}^{(T)}(dx,y) = (2\pi)^{-1/2}e^{-x^2/2}\,dx
 \qquad \text{for } |y| > (\log T)^5.$$
Here $C_{\mathfrak A}>0$ is the normalization constant. Define
\[
\widetilde\alpha^{(T)}(y):=-2\int_0^\infty\log x\,\mu_{\sd}^{(T)}(dx,y).
\]

		We take the entries of $(\bm{a}(0); \bm{b}(0))$ to be mutually independent. We sample $a_i (0), i \in \llbracket N_1, N_2 -1 \rrbracket$ and $b_i (0),  i \in \llbracket N_1, N_2 \rrbracket$ under the probability laws $\mu_{\sd}^{(T)} (dx,i/T)$ and $\mu_{\dg}^{(T)} (dx, i/T)$, respectively. The initial data $(\bm p(0);\bm q(0))$ is then determined from $(\bm a(0);\bm b(0))$ by setting $q_0(0)=0$. Let $\L(s)$ be the associated Lax matrix, as in Definition \ref{def:flandL}, and set
		\[
		(\lambda_1,\lambda_2,\ldots,\lambda_N):=\eig\L(0),
		\qquad \lambda_1>\lambda_2>\cdots>\lambda_N.
		\]
		\end{assumption}

The choice $N\asymp T^2$
ensures that endpoints are sufficiently far from the bulk that 
boundary effects are not felt by the system on the hydrodynamic
space time scales under consideration. In principle, one could instead take
$N \sim T(\log T)^C$ for any sufficiently large constant $C>1$.
	

Next, we define \emph{quasiparticles}, which are the objects which allow us to prove the convergence stated in Theorem \ref{thm:local-charge-testfn-intro}. We proceed along the lines of
the definition given in \cite{Agg25a}.
\begin{definition}
\label{def:loccent}
    \label{center}
	Let $\mathbf{u} = (u(N_1),\dots, u(N_2)) \in \mathbb{R}^N$ be a unit vector. For any $\zeta \in \mathbb{R}_{\ge 0}$, we say that an index $\varphi \in \llbracket N_1, N_2 \rrbracket$ is a \emph{$\zeta$-localization center} of $\mathbf{u}$ if $|u(\varphi)| \ge \zeta$. 
	
	Next, let $\mathbf{M}$ be an $N \times N$ real symmetric matrix, with distinct eigenvalues and rows and columns indexed by $\llbracket N_1, N_2\rrbracket$. If $\lambda$ is an eigenvalue of $\mathbf{M}$, we call $\varphi$ a $\zeta$-localization center for $\lambda$ if $\varphi$ is a $\zeta$-localization center of the unit eigenvector $\mathbf{u}$ of $\mathbf{M}$ with eigenvalue $\lambda$. Further let $(\mathbf{u}_1, \mathbf{u}_2, \ldots , \mathbf{u}_N)$ denote an orthonormal eigenbasis for $\mathbf{M}$. We call a bijection $\varphi : \llbracket 1, N \rrbracket \rightarrow \llbracket N_1, N_2 \rrbracket$ a \emph{$\zeta$-localization center bijection} if $\varphi(j)$ is a $\zeta$-localization center for $\mathbf{u}_j$, for each $j= \llbracket 1, N \rrbracket$. By \cite[Lemma 2.7]{Agg25a}, any symmetric $N \times N$ matrix $\mathbf{M}$ admits a $(2N)^{-1}$-localization center bijection.
\end{definition}

While they exist, localization center bijections are not always unique. However, by Lemma \ref{bijectionl}, under Assumption \ref{ass:model}, at time $0$, they are unique up to a small error. 
One of the main tasks of this work is to show that eigenvectors remain sufficiently localized along the evolution of the dynamics 
even in the non stationary setup of the present work.

\begin{definition}[Quasiparticles] \label{def:quasi_intro}
Adopt Assumption \ref{ass:model}, and let $\zeta = e^{-100 (\log N)^{3/2}}$. For each $s\in\mathbb{R}_{\geq 0}$ and $k\in\llbracket 1,N\rrbracket$, let $\mathbf{u}_k(s)=(u_k(N_1;s),u_k(N_1+1;s),\ldots,u_k(N_2;s))\in\mathbb{R}^N$ denote a unit eigenvector of $\L(s)$ with eigenvalue $\lambda_k$.  
Under the orthonormal basis $(\mathbf{u}_1(s), \mathbf{u}_2(s),\ldots, \mathbf{u}_N(s))$ of $\L(s)$, let $\varphi_s:\llbracket 1,N\rrbracket \to \llbracket N_1,N_2\rrbracket$ be a $\zeta$ localization center bijection for $\L(s)$, and define the \emph{quasiparticle}
\begin{equation}\label{eqn:quasi_part_intro}
Q_j(s)=q_{\varphi_s(j)}(s),
\qquad
\text{for each }(j,s)\in\llbracket 1,N\rrbracket\times\mathbb{R}_{\ge 0}.
\end{equation}

\noindent For any $s \ge 0$ and $i \in \llbracket N_1, N_2 \rrbracket$, further define 
\begin{equation}
    \label{lambdait}
    \Lambda_i(s) = \lambda_{\varphi_s^{-1}(i)}, 
\end{equation}

\noindent which we view as the eigenvalue with localization center $i$ at time $s$. Also abbreviate
\[
\Lambda_i=\Lambda_i(0),
\qquad
\bar q_i=T^{-1}q_i(0).
\]

\end{definition}

The precise localization threshold $\zeta$ in Definition \ref{def:quasi_intro} is not essential. For the main results
one may fix any $\zeta=\zeta_N$ that decays at least polynomially fast, but more
slowly than stretched exponentially in $N$; the displayed choice
$e^{-100(\log N)^{3/2}}$ is a convenient representative.

The proof of Theorem \ref{thm:local-charge-testfn-intro} will follow from a
detailed understanding of the trajectories of the quasiparticle data
$(Q_i(t)/T,\lambda_i)$. The following theorem realizes
$m_\tau(\lambda,x)\,d\lambda\,dx$ as their limiting empirical measure; it is
proved in Section \ref{subsec:mainpfs} below.

		\begin{thm}\label{thm:main}
			Adopt Assumption \ref{ass:model}, let $m_0$ be as in Definition
			\ref{rho0}, and let $N_t$ and $m_t$ be as in Proposition
			\ref{prop:GHDint}. Recall the quasiparticle spectral parameter pairs
			$\{(Q_i(t),\lambda_i)\}_{i=1}^N$ from Definition
			\ref{def:quasi_intro}. There exists
			$\mathfrak A_0=\mathfrak A_0(\mathfrak K,\epsilon)>1$ such that the
			following holds whenever $\mathfrak A\ge\mathfrak A_0$. For any smooth,
			compactly supported function $f:\mathbb R^2\to\mathbb R$ and fixed
			$\tau>0$, we have the convergence in probability as $T\to\infty$,
		\[
		\frac1T\sum_{i=1}^N
		f\bigl(\lambda_i,Q_i(T\tau)T^{-1}\bigr)
		\longrightarrow
		\int_{-\infty}^{\infty}\int_{-\infty}^{\infty}
		f(\lambda,x)m_\tau(\lambda,x)\,dx\,d\lambda.
		\]
\end{thm}

	In order to prove Theorem \ref{thm:main}, we show that quasiparticles evolve according to an ordinary differential equation (ODE) determined by the
	effective velocity $\ve_{\tau}(\lambda; q)$
	characterized by Equation \eqref{eqn:cra}. 
		Indeed, we have the following definition of 
		a family of flows which represent the limiting, macroscopic trajectories of quasiparticles.
		\begin{definition}\label{def:xtd}
		Adopt the notation of Proposition \ref{prop:GHDint}, and recall the
		effective velocity $\ve_t$ from its second part.
		Fix \(\lambda\in\mathbb R\). A \emph{characteristic flow} is a family $(t,q)\mapsto x_t(\lambda;q), (t,q)\in\mathbb R_{\geq0}\times\mathbb R$, 
		such that \(x_0(\lambda;q)=q\), 
		 the map  $q\mapsto x_t(\lambda;q)$ is locally bi-Lipschitz uniformly for \(t\) in compact intervals, and
		  moreover,
		   \begin{equation} \label{eqn:xtode} 
		   \frac{d}{dt}x_t(\lambda;q) = v_t^{\eff}\bigl(\lambda;x_t(\lambda;q)\bigr) 
		   \end{equation} 
		  for almost every \((t,q)\in\mathbb R_{\geq0}\times\mathbb R\).
		\end{definition}

		The existence of a solution to \eqref{eqn:xtode} is provided by the
		following lemma, which is proved in Section
		\ref{subsec:ghd-uniqueness}.
		\begin{lem}
		\label{lem:xtode-wellposed-intro}
		Adopt the notation of Proposition \ref{prop:GHDint}, and recall the
		effective velocity $\ve_t$ from its second part. For every
		$\lambda \in\mathbb R$, there exists
		a unique characteristic flow $(t,q)\mapsto x_t(\lambda;q), (t,q)\in\mathbb R_{\geq0}\times\mathbb R$
		satisfying the properties of Definition \ref{def:xtd}.
		\end{lem}
		
		Using this definition, we can state our main theorem about the
		quasiparticles (from Definition \ref{def:quasi_intro}). Again, in what
		follows, the exponent $1/4$ is not optimal, and we make no attempt to
		optimize it. The theorem is proved in Section \ref{subsec:stop_outline}
		below.
		\begin{thm}\label{thm:ev}
		Adopt Assumption \ref{ass:model}. There exists $\mathfrak A_0=\mathfrak A_0(\mathfrak K,\epsilon)>1$ such that the following holds whenever $\mathfrak A\ge\mathfrak A_0$.
	There exists a constant $c = c(\mathfrak{A})> 0$ such that the following holds with probability at least $1-c^{-1} e^{-c(\log N)^2}$. For each $(k, t) \in \llbracket 1, N \rrbracket \times [0,T\log N]$ satisfying 
	$$\varphi_0(k) \in \llbracket - T (\log N)^{4}, T (\log N)^{4} \rrbracket,$$
	we have
	\begin{align}\label{eqn:ev}
		\left| Q_k (t) T^{-1}  - x_{t/T} (\lambda_k; Q_k(0) T^{-1}) \right|  \leq  T^{-1/4} \log T.  
		\end{align}
		\end{thm}

	\subsection{Notation}
	\label{notation}
	
	 Throughout, $c$ and $C$ denote small and large positive constants,
	 respectively, which may change between lines. Unless stated otherwise,
	 all constants may depend on $\mathfrak K$ and $\epsilon$ from Assumption
	 \ref{ass:entry-law-profiles}, even when this is not explicitly indicated.
	 For a bounded function $f : D \rightarrow \mathbb{R}$ defined on a set $D$, we write
	 $\|f\|_\infty:=\sup_{x\in D}|f(x)|$.
	 For a Lipschitz function $f:D\to\mathbb R$, where $D\subseteq\mathbb R$, we write
	 $\|f\|_{\Lip}:=\sup_{x,y\in D,\,x\ne y}|f(x)-f(y)|/|x-y|$.

	 Throughout the text, a vector $\bm{v} = (v_1, v_2, \ldots , v_n) \in \mathbb{C}^n$ is a 
	 unit vector if $\sum_{i=1}^n |v_i|^2 = 1$; we call it nonnegatively 
	 normalized if $v_j > 0$, where $j \in \llbracket 1, n \rrbracket$ is 
	 the minimal index such that $v_j \ne 0$. For any 
	 real symmetric $n \times n$ matrix $\mathbf{M}$,
	  let $\eig \mathbf{M} = (\lambda_1, \lambda_2, \ldots , \lambda_n)$ denote
	  the eigenvalues of $\mathbf{M}$, counted with multiplicity
	  and ordered so that $\lambda_1 \ge \lambda_2 \ge \cdots \ge \lambda_n$. 
	  
	  For a matrix $\mathbf{M}=[M_{ij}]$ whose rows and columns are indexed by
	$\llbracket N_1,N_2\rrbracket$ and an index interval
	$I=\llbracket a,b\rrbracket\subseteq\llbracket N_1,N_2\rrbracket$, write
	$\mathbf{M}^{I}=\mathbf{M}^{[a,b]}$ for the corresponding principal
	submatrix. When $\mathbf M=\L(t)$, we use the abbreviation
	$\mathbf M^I=\L^I(t)$. We use the empty matrix when $a>b$.

Throughout the text, we say an 
	 event has \emph{overwhelming probability} if its complement holds with probability at most $c^{-1}e^{-c(\log N)^2}$
		 for some $c = c(\mathfrak{A}, \mathfrak{K}, \epsilon)>0$.

	\section{Proof of quasiparticle limiting trajectories}
	\label{sec:outline}
	
As we will see in Section
\ref{subsec:mainpfs} below, the results in Sections
\ref{subsec:mainthm1} and \ref{subsec:mainthm2} follow from Theorem
\ref{thm:ev}, which prescribes the limiting trajectories of the
quasiparticles. Thus, in this section, we prove this theorem conditional on
various intermediate results that will be proved later in the paper.
Throughout the section, we adopt Assumption \ref{ass:model} and set
\begin{equation}\label{eqn:R}
R:=\mathfrak A^{1/10}.
\end{equation}
In addition, we use the notation defined in Definitions \ref{def:flandL}, \ref{def:loccent}, and \ref{def:quasi_intro}.
	
Sections \ref{subsec:lin} and \ref{subsec:ghd-uniqueness} below construct the
unique solution of the GHD equation, proving Proposition \ref{prop:GHDint};
Section \ref{subsec:tc} introduces
 the solution of the GHD equation associated to the initial data 
defined in Assumption \ref{ass:model};
Section \ref{subsec:stop_outline} introduces a certain 
	 stopping time and states Theorem \ref{ttn}, from which Theorem
	 \ref{thm:ev} directly follows. Sections
	 \ref{subsec:asymptotic_scattering} through \ref{subsec:ideas2} outline a
	 self improving argument to establish Theorem \ref{ttn}, reducing it to a
	 series of estimates whose proofs will constitute the bulk of this work.
		
	\subsection{Existence of solutions to the generalized hydrodynamic equation}
	\label{subsec:lin}
	In order to make sense of the trajectories $x_{\tau}$, defined in Definition \ref{def:xtd}, which 
	appear in the statement of Theorem \ref{thm:ev}, one must first 
	solve the GHD equation \eqref{eqn:rhoGHDintro}, so that the effective velocity 
	$\ve_{\tau}(\lambda ; x)$ is well defined for $\tau > 0$.
	 Existence and uniqueness of solutions 
	 of the GHD equation is the content of Proposition \ref{prop:GHDint}. In this subsection
	 we prove the existence part of that proposition, through Lemma
	 \ref{lem:int_ghd} below; uniqueness is proved in the next subsection. In
	 what follows, we use the following notation.
	\begin{definition}\label{def:Tdf}
	Let $\mathbf T$ denote the linear operator on functions of a single variable with integral kernel $2 \log|x - y|$. More precisely, if $f : \mathbb{R} \rightarrow \mathbb{R}$, then 
	\begin{equation}\label{eqn:Tdf}
	(\mathbf{T} f)(x) \coloneqq \int_{-\infty}^{\infty} 2 \log|x-y| f(y) dy
	\end{equation}
whenever the right hand side makes sense.
	\end{definition}
	
	We also adopt the following assumption throughout this subsection.
		\begin{assumption}\label{ass:sec2}
			Fix $\mathfrak K>1$.
			Let $\varrho_0 : \mathbb{R}^2 \rightarrow \mathbb{R}$ be a function such that for each $y \in \mathbb{R}$, \(\lambda \mapsto \varrho_0(\lambda,y)\) is a probability
			density. We assume for some constants $C_0, c_0>0$ it satisfies 
			\begin{equation}\label{eqn:varrho0bd}
			0\le \varrho_0(\lambda, y)\le C_0 e^{-c_0 |\lambda|^2}, \qquad (\lambda, y) \in \mathbb{R}^2.
			\end{equation}
			Let $\alpha_0 : \mathbb{R} \rightarrow \mathbb{R}$ and $\mathfrak{A} > 0$, and assume that for
			each \(x \in \mathbb{R}\),
			$\mathfrak K^{-1}\mathfrak A\leq\alpha_0(x)\leq\mathfrak K\mathfrak A$. We write
			\(m_0(\lambda,x)=\alpha_0(x)^{-1}\varrho_0(\lambda,x)\), for $\lambda,x \in \mathbb{R}$.
			We further assume
			$\mathfrak A\ge\mathfrak A_0(\mathfrak K, C_0,c_0)$, for some $\mathfrak A_0(\mathfrak K, C_0,c_0) > 1$.
			The threshold $\mathfrak A_0$ will be determined throughout this
			section. Unless stated otherwise, all constants in this section may
			depend on $\mathfrak K$, $C_0$, and $c_0$, even when this is not
			explicitly indicated.
		\end{assumption}
	Under Assumption \ref{ass:model_inf}, take $\alpha_0=\alpha$ from
	\eqref{eqn:Fstr} and let $\varrho_0$ be as in Definition \ref{rho0}.
	Then \eqref{eqn:alph_bd} and Lemma \ref{rho0estimate} show that
	$(\alpha_0,\varrho_0)$ satisfies Assumption \ref{ass:sec2}, so the results
	below apply to this choice.

	The following definition changes variables from the spatial coordinate $x$
	to a new coordinate $X$. In these coordinates, the function $N_t$ from
	Proposition \ref{prop:GHDint} admits the direct fixed point equation in Lemma
	\ref{lem:fixed_pointout}, and the characteristic flow is a linear
	evolution (see Lemma \ref{xt2}).
	\begin{definition}\label{def:initial-ghd-coordinates}
	Define
	\[
	X_0(\Lambda,x)
	=
	x+2\int_{0}^x\int_{-\infty}^{\infty}\log|\Lambda-\lambda|\,m_0(\lambda,y)\,d\lambda\,dy.
	\]
	Its $x$-derivative is
	\begin{equation}\label{eqn:J0}
	J_0(\Lambda,x)
	:=
	\partial_xX_0(\Lambda, x)
	=
	1+2\int_{-\infty}^{\infty}\log|\Lambda-\lambda|\,m_0(\lambda,x)\,d\lambda.
	\end{equation}
	Define $\varpi_0$ by pushforward:
	\begin{equation}\label{eqn:vpdef}
	\varpi_0(\Lambda;X_0(\Lambda,x))
	=
	\frac{m_0(\Lambda,x)}{J_0(\Lambda,x)}.
	\end{equation}
	Equivalently,
	$\varpi_0(\Lambda;X)\,d\Lambda\,dX
	=m_0(\Lambda,x)\,d\Lambda\,dx$
	under $X=X_0(\Lambda,x)$. Set
	\begin{equation}\label{eqn:N0}
	N_0(\lambda, x) =  \int_{0}^x m_0(\lambda,y)\,dy =  \int_{0}^{X_0(\lambda,x)} \varpi_0(\lambda,X)\,dX =:  N_0^X(\lambda, X_0(\lambda, x)).
	\end{equation}
	For $t\ge 0$, define the free transport of $\varpi_0$ by
	\begin{equation}\label{eqn:wtX}
	\varpi_t(\Lambda;X):=\varpi_0(\Lambda;X-t\Lambda).
	\end{equation}
	\end{definition}

	By Assumption \ref{ass:sec2},
	$(\T m_0(\cdot,x))(\Lambda)\ge -C\mathfrak A^{-1}$, since the
	negative part of the logarithmic kernel is uniformly integrable against
	the bound on $m_0$ obtained from \eqref{eqn:varrho0bd}. Thus, fixing
	$\mathfrak A\ge\mathfrak A_0$ sufficiently large, we have
	\[
	J_0(\Lambda,x)\ge 1-C\mathfrak A^{-1}\geq \frac{1}{2}> 0
	\]
	for all $(\Lambda,x)\in\mathbb R^2$. Hence $x\mapsto X_0(\Lambda,x)$ is
	strictly increasing for each $\Lambda$ and admits an inverse. Moreover,
	\[
	m_0(\Lambda,x)
	\le C_0\mathfrak K\mathfrak A^{-1}e^{-c_0|\Lambda|^2}.
	\]
	Consequently, at time zero,
	$0\le\varpi_0(\Lambda;X)\le
	2C_0\mathfrak K\mathfrak A^{-1}e^{-c_0|\Lambda|^2}$ for all
	$(\Lambda,X)\in\mathbb R^2$, while for every $t\ge0$,
	\begin{equation}\label{eqn:vpi0bd}
	0\le \varpi_t(\Lambda;X)
	\le 2C_0\mathfrak K\mathfrak A^{-1}e^{-c_0|\Lambda|^2},
	\qquad (\Lambda,X)\in\mathbb R^2.
	\end{equation}

		The function $N_t$ from Proposition \ref{prop:GHDint} will be defined by
		the solution to the fixed point equation \eqref{eqn:NtFP} below. Let
		$\mathfrak c\in [c_0/10,c_0/2]$, where $c_0$ is the constant in
	Assumption \ref{ass:sec2}, and set
	\begin{equation}\label{eqn:Bsp}
	\mathcal B
	:=
	\bigl\{f:\mathbb R\to\mathbb R:\ \|f\|_{\mathfrak c}<\infty\bigr\},
	\end{equation}
	equipped with the norm $\|\cdot \|_{\mathfrak c}$ defined in \eqref{eqn:weighted-sup-norm}.
	We have the following result, whose proof is given at the end of Section
	\ref{subsec:GHD_alt}.
	
	
	\begin{lem}\label{lem:fixed_pointout}
	Adopt Assumption \ref{ass:sec2}, and let $\mathcal B$ be the space above.
	For every $x\in \mathbb{R}$ and $t\ge 0$, there exists a unique function $N_t(\cdot,x)\in \mathcal B$ satisfying 
	\begin{equation}\label{eqn:NtFP}
	N_t(\lambda, x) = \int_{0}^{x+(\T N_t(\cdot, x))(\lambda) - \lambda t} \varpi_0(\lambda, X ) d X .
	\end{equation}
	Moreover, there is a constant $C>0$ such that $N_t(\lambda, x) $ satisfies
\begin{equation}\label{eqn:Nt-joint-lipschitz}
e^{\mathfrak c\lambda^2}
\bigl|N_t(\lambda,x)-N_s(\lambda,y)\bigr|
\leq C\mathfrak A^{-1}\bigl(|x-y|+|t-s|\bigr)
\end{equation}
uniformly for $\lambda,x,y\in\mathbb R$ and $s,t\in\mathbb R_{\geq0}$.
	\end{lem}

	The next lemma is proved in Section \ref{subsec:GHD_2}. It shows that the fixed point
	defined by Lemma \ref{lem:fixed_pointout} above
	indeed satisfies \eqref{eqn:GHDpot_unique}.
	\begin{lem}\label{lem:int_ghd}
	Adopt Assumption \ref{ass:sec2}, and let $N_t(\lambda,x)$ be the fixed point from Lemma \ref{lem:fixed_pointout}. Let
	$m_t(\lambda,x)=\partial_xN_t(\lambda,x)$ denote its weak $x$-derivative,
	which is defined for almost every $x \in \mathbb{R}$, for each $t \geq 0$ and $\lambda \in \mathbb{R}$. 
	For each 
	$(t, \lambda) \in \mathbb{R}_{\geq 0} \times \mathbb{R}$, $m_t(\lambda,x)\geq0$ 
	for almost every $x \in \mathbb{R}$.
	Moreover, for each $t \geq 0$, for almost every $x \in \mathbb{R}$ there is a
	unique function
	$\ve_t(\cdot;x)\in L^2(\mathbb{R},e^{-\frac{c_0}{2} \lambda^2} d\lambda)$ satisfying
	\eqref{eqn:cra}. With this choice of $\ve_t$, for each $\lambda \in \mathbb{R}$,
	\begin{equation}\label{eqn:int_ghd}
	\partial_t N_t(\lambda, x)
	=
	-\ve_t(\lambda;x)m_t(\lambda,x)
	\end{equation}
	for almost every $(t,x) \in \mathbb{R}_{\ge 0} \times \mathbb{R}$. 
	\end{lem}
	
	\subsection{Uniqueness of solutions to the generalized hydrodynamic equation}
	\label{subsec:ghd-uniqueness}

	In this subsection we prove Proposition \ref{prop:GHDint}. Its existence
	aspect follows from Lemma \ref{lem:int_ghd}, so its uniqueness aspect remains.
	To that end, we introduce coordinates associated with a candidate solution
	of \eqref{eqn:rhoGHDintro}.
	\begin{definition}\label{def:anyN}
	Adopt Assumption \ref{ass:sec2}, and let $\widetilde N_t$ denote any
	function satisfying the three properties listed in Proposition
	\ref{prop:GHDint}. Set
	$\widetilde m_t(\lambda,x):=\partial_x\widetilde N_t(\lambda,x)$ and define
	\begin{equation}\label{eqn:Xt}
	\widetilde X_t(\Lambda,x)
	:=x+(\T\widetilde N_t(\cdot,x))(\Lambda),
	\end{equation}
	\begin{equation}\label{eqn:s0dr}
	\widetilde J_t(\Lambda,x)
	:=1+(\T\widetilde m_t(\cdot,x))(\Lambda).
	\end{equation}
	If $\widetilde N_t=N_t$ is the fixed point from Lemma
	\ref{lem:fixed_pointout}, then we write $\widetilde X_t=X_t$ and
	$\widetilde J_t=J_t$.
	\end{definition}

	By Lemmas \ref{lem:Jlow} and \ref{xt2} below, the time $\tau$ flow $(\lambda, x) \mapsto x_{\tau}(\lambda, x)$ from Definition \ref{def:xtd} is well defined for any $\tau \geq 0 $.
	To state them, recall the definitions \eqref{eqn:Xt} and \eqref{eqn:s0dr}.
	The proof of the lemma below is in Section \ref{subsec:GHD_2}.
	\begin{lem}\label{lem:Jlow}
	Adopt the assumptions and notation of Definition \ref{def:anyN}.
	With $\widetilde J_t(\Lambda, x)$ as in \eqref{eqn:s0dr},
	\begin{equation}\label{eqn:avrbd2coarse}
	\widetilde J_t(\Lambda, x)
	\geq \frac12,
	\end{equation}
	and in particular, for every $t \geq 0$ and $\Lambda \in \mathbb{R}$, the map $x\mapsto \widetilde X_t(\Lambda,x)$ is bi-Lipschitz and invertible.
	 \end{lem}
	The following lemma is
	 proved in Section \ref{subsec:onsolghd}, and will be useful
	 in the proof of Theorem \ref{ttn}. Together with Lemma \ref{lem:Jlow},
	 it provides existence of the time $\tau$ flow $(\lambda, x) \mapsto x_{\tau}(\lambda, x)$ from Definition \ref{def:xtd},
	 and also provides a simple description of this flow, namely a linear evolution, in $X$ coordinates.
	\begin{lem}
		\label{xt2}
		 Adopt Assumption \ref{ass:sec2}, and let $\widetilde N_t$ denote any function
		 satisfying the three properties listed in Proposition
		 \ref{prop:GHDint}. Let $\widetilde{v}^{\eff}_t$ denote the corresponding effective
		 velocity.
		Define \(x_t(\lambda;q)\) by
		\begin{equation}\label{eqn:Xtlin}
		x_t(\lambda;q) = \widetilde X_t(\lambda, \cdot)^{-1}\left( X_0(\lambda,q)+\lambda t \right),
		\end{equation}
		for any $\lambda \in \mathbb{R}$ and $(t, q) \in \mathbb{R}_{\geq 0} \times \mathbb{R}$. Then,
		\((t,q) \mapsto x_t(\lambda;q)\) is a characteristic flow
		in the sense of Definition \ref{def:xtd}, with $\ve_t$ there replaced by
		$\widetilde v_t^{\eff}$.
		On the other hand, 
		if $(t,q) \mapsto x_t(\lambda;q)$ is a characteristic flow, then \eqref{eqn:Xtlin} holds.
		\end{lem}

	Using these ingredients, we are ready to provide the proof of Proposition \ref{prop:GHDint}. 
	\begin{proof}[Proof of Proposition \ref{prop:GHDint}]
	Existence follows from Lemmas \ref{lem:fixed_pointout} and \ref{lem:int_ghd}.
	For uniqueness, let $\widetilde N_t(\lambda,x)$ be any function satisfying
	the properties in the proposition, and let $x_s(\lambda;q)$ be defined as in
	Lemma \ref{xt2}. Define its pullback to $X$ coordinates by
	$\widetilde N_t^X(\lambda,X)
	:=\widetilde N_t\bigl(\lambda,
	\widetilde X_t(\lambda,\cdot)^{-1}(X)\bigr)$.
	The inverse exists by Lemma \ref{lem:Jlow}.
	
	 Equation
	\eqref{eqn:GHDpot_unique} implies that
	$s \mapsto \widetilde N_s(\lambda,x_s(\lambda; q)) $ is constant for almost every $q \in \mathbb{R}$,
	 and thus for every $q \in \mathbb{R}$ by continuity. 
	 Thus, using Lemma \ref{xt2} and \eqref{eqn:Xt}, we have
	\begin{align}
	\widetilde N_t(\lambda,x_t(\lambda; q)) &= N_0(\lambda,q) \notag \\
	&=   N_0^X(\lambda,X_0(\lambda, q)) \notag \\
	 &= N_0^X(\lambda, \widetilde X_t(\lambda, x_t(\lambda; q)) - \lambda t)\notag \\
	 &= N_0^X\bigl(\lambda, x_t(\lambda; q)+
	 (\mathbf{T}\widetilde N_t(\cdot,x_t(\lambda;q)))(\lambda)-t\lambda\bigr).
	 \label{eqn:m0push0}
	\end{align}
	By \eqref{eqn:N0}, upon replacing $x_t(\lambda;q)$ by any $x\in\mathbb R$
	(which can be done using \eqref{eqn:Xtlin} and the invertibility in Lemma
	\ref{lem:Jlow}), relation \eqref{eqn:m0push0} for $\widetilde N_t$
	coincides with the fixed point equation \eqref{eqn:NtFP} for $N_t$. Thus,
	by the uniqueness statement of Lemma \ref{lem:fixed_pointout}, we have
	$\widetilde N_t=N_t$, proving the proposition.
	\end{proof}
	
	\begin{rem}\label{rem:pushrem}
	Since, as shown in the proof of Proposition \ref{prop:GHDint},
	$N_s(\lambda,x_s(\lambda;x))$ is constant in $s\geq0$, we have
	\begin{equation}\label{eqn:m0push}
	m_{\tau} = (x_{\tau})_{\#} m_0.
	\end{equation}
	Here $(x_{\tau})_{\#}m_0$ denotes the density of the pushforward of the
	measure $m_0(\lambda,x)d\lambda dx$ under the map
	$(\lambda,x)\mapsto(\lambda,x_{\tau}(\lambda;x))$.
	\end{rem}
	\begin{proof}[Proof of Lemma \ref{lem:xtode-wellposed-intro}]
		Apply Lemma \ref{xt2} to the solution 
		furnished by Proposition \ref{prop:GHDint}. The first statement of the lemma
		provides existence, while the second statement guarantees uniqueness.
	\end{proof}

		\subsection{Limiting quantities under Assumption \ref{ass:model}}
	\label{subsec:tc} 
	Recall the notation $F_{\mathrm{str}}$, $\widetilde\varrho_0$,
	$\varrho_0$, $N_t$, $m_t$, and $x_t$ from \eqref{eqn:Fstr}, Lemma
	\ref{lem:frozen-dos-intro}, Definition \ref{rho0}, Lemma
	\ref{lem:fixed_pointout}, Proposition \ref{prop:GHDint}, and Lemma
	\ref{lem:xtode-wellposed-intro}, respectively,
	under Assumption \ref{ass:entry-law-profiles}. The $T$-dependent
	entry laws in Assumption \ref{ass:model} satisfy Assumption
	\ref{ass:entry-law-profiles}. We write $F_{\mathrm{str}}^{(T)}$,
	$\widetilde\varrho_0^{(T)}$, $\varrho_0^{(T)}$, $N_t^{(T)}$,
	$m_t^{(T)}$, and $x_t^{(T)}$ for the corresponding functions obtained by replacing
	$(p_{\dg},\mu_{\sd})$ with $(p_{\dg}^{(T)},\mu_{\sd}^{(T)})$.

	Heuristically, by \eqref{eqn:m0push} and since flows given in 
	Definition \ref{def:xtd} with spectral parameters $\lambda$ which are not too large cannot propagate too quickly,
	we expect $|m_{\tau}^{(T)} - m_{\tau}|$ to be small as long as
	 $\tau$ is not too large compared to $T$. We have the following lemma making this precise.
	 Its proof is in Section \ref{subsec:onsolghd}.

\begin{lem}
\label{lem:thermal-completion}
Under the notation above, fix $x_0\in\mathbb R$ and let
\(L>(\log T)^3\). Suppose that
\[
m_0^{(T)}(\lambda,x)=m_0(\lambda,x),
\qquad x\in[x_0-L,x_0+L].
\]
There are constants \(c,C>0\) such that, if
$0\le s\le (\log T)^2$
and
$|q-x_0|+\tau (1+|\lambda|)\le cL$,
then
\[
\left|
x_s^{(T)}(\lambda;q)-x_s(\lambda;q)
\right|
\le
Ce^{-c(\log T)^2}.
\]
\end{lem}

\subsection{Stopping time}
	\label{subsec:stop_outline}
	Recall that throughout this subsection we adopt and use the notation of Assumption
	\ref{ass:model} and Definition \ref{def:quasi_intro}, as well as the notation
	from Section \ref{subsec:tc} (such as $F_{\mathrm{str}}^{(T)}$ and
	$x_t^{(T)}$). We have the following definition of a stopping time
	$\mathfrak{T}$, which will be used in the proof of Theorem \ref{thm:ev}.
	Fix sufficiently large and small constants
	$C_1=C_1(\mathfrak K,\epsilon)>1$ and
	$c_1=c_1(\mathfrak K,\epsilon)\in(0,1)$, respectively. Their sizes, crucially, do not depend on the large parameter $\mathfrak{A}$ 
	(which will ultimately be taken large enough, in a way that depends on $c_1, C_1$). 
	For $i\in\llbracket N_1,N_2\rrbracket$, set
	\[
	\bar q_i:=T^{-1}q_i(0).
	\]

		\begin{definition} 
	\label{mt}
	Fix a universal constant $\delta\in(0,1/200)$, let
	$\mathfrak{M} \in [T^{3/4-\delta}, T^{3/4}]$, and let $R$ be as in \eqref{eqn:R}.
	For any $t=T\tau\geq0$, let $\mathfrak B(t)=\mathfrak B_{\mathfrak M}(t)$ denote the event on which
	the following three properties hold, with items \ref{item:T2} and
	\ref{item:T3} required for each $k \in \llbracket 1, N \rrbracket$ such that
	$\varphi_t (k) \in \llbracket N_1+T (\log N)^{10}, N_2 - T (\log N)^{10} \rrbracket$.
	
	\begin{enumerate}
	\item For every
$i>j,
i,j\in
\llbracket
N_1+T(\log N)^{10},\,
N_2-T(\log N)^{10}
\rrbracket$,
we have
\begin{equation}
\label{eqn:qbd}
q_i(t)-q_j(t)
\le
C_1\mathfrak A
\max\left\{
i-j,\frac{\mathfrak M}{R}
\right\}.
\end{equation}
Moreover, whenever $i-j>\frac{\mathfrak M}{R}$,
we have
\begin{equation}
\label{eqn:qbd-lower}
q_i(t)-q_j(t)
\ge
c_1\mathfrak A(i-j).
\end{equation}
\label{item:T1}

\item For all $j \in \llbracket N_1, N_2 \rrbracket$,
		\begin{equation}\label{eqn:uMR}
			|u_k (j;t)| \le e^{\mathfrak{M} - R |j-\varphi_t (k)|}.
		\end{equation} 
	\label{item:T2}

	\item We have 
			\begin{equation}\label{eqn:Qx}
				|Q_k (t) -T  x_\tau^{(T)} (\lambda_k; \bar q_{\varphi_0(k)})| \le \mathfrak{M} \log (|\lambda_k| + 2).
			\end{equation}
			\label{item:T3}
		
	\end{enumerate}
	Let $\mathfrak T_{\mathfrak M}=\mathfrak T\in[0,T\log N]$ denote the largest real number for which $\mathfrak B(s)$ holds for every $s\in[0,\mathfrak T]$.
	
	\end{definition}

Heuristically, $\mathfrak{T}$ is the first time that either the error of approximation of quasiparticles by their continuum limit becomes too large or an eigenvector of the Lax matrix becomes too delocalized. The main results of this work can be reduced
 to the following theorem, which states that $\mathfrak T$ is as large as possible. Its proof,
conditional on several intermediate results, is given in Section \ref{subsec:asymptotic_scattering} below.

	\begin{thm} 
		
		\label{ttn}
		
			Adopt Assumption \ref{ass:model}. There exists $\mathfrak A_0=\mathfrak A_0(\mathfrak K,\epsilon)>1$ such that the following holds whenever $\mathfrak A\ge\mathfrak A_0$. Let $\mathfrak{M} = T^{3/4}$ in Definition \ref{mt}. With overwhelming probability, $\mathfrak{T} = T \log N$. 
		
	\end{thm} 
	
	Theorem \ref{ttn} is the primary component in the proof of Theorem
	\ref{thm:ev}, but to establish the latter from the former we require a few
	additional facts. The first is the following proposition providing a finite
	speed of propagation for localization centers, which will be used to verify
	\eqref{eqn:uMR} directly when $|j-\varphi_t(k)|$ is sufficiently large. Its
	proof, which is similar to that of
	\cite[Lemma 5.2]{Agg25a}, is given in Section \ref{app:lemLocR}.

\begin{prop}\label{prop:outside-strip-confinement}
Adopt Assumption \ref{ass:model}, and let $R$ be as in \eqref{eqn:R}. There exists $\mathfrak A_0=\mathfrak A_0(\mathfrak K,\epsilon)>1$ such that whenever $\mathfrak A\ge\mathfrak A_0$,
the following holds with
overwhelming probability.
For any $k\in\llbracket1,N\rrbracket$, $s \in [0, T \log N]$,
$s'\in\{0,s\}$, and $m\in\llbracket N_1,N_2\rrbracket$ satisfying
$|m-\varphi_{s'}(k)| > T (\log N)^3$, we have
\begin{equation}\label{eqn:umphi}
        |u_k(m;s)|\le \exp\{-4R |m-\varphi_{s'}(k)|\}.
\end{equation}
\end{prop}

Before proving Theorem \ref{thm:ev}, we record the boundedness estimate
used in its proof.
\begin{definition}\label{adelta}
Fix a real number $A>0$, let $\mathscr I$ be an
index set, and let $\mathbf{M}=[M_{ij}]_{i,j\in\mathscr I}$ be a real
symmetric matrix. Define
\[
\mathsf{BND}_{\mathbf{M}}(A)
:=\left\{\max_{i,j\in\mathscr I}|M_{ij}|\le A\right\}
\cap\left\{\max_{\lambda\in\eig\mathbf{M}}|\lambda|\le A\right\}.
\]
\end{definition}
The proof of the following lemma is entirely analogous to that of
\cite[Lemma 3.15]{Agg25a} and is therefore omitted.
\begin{lem}[{\cite[Lemma 3.15]{Agg25a}}]
\label{l0eigenvalues}
For any real number $A\ge1$,
\[
\mathbb P\left[\bigcap_{t\in\mathbb R_{\ge0}}
\mathsf{BND}_{\L(t)}(A)\right]
\ge1-c^{-1}Ne^{-cA^2}.
\]
\end{lem}

In addition, we have the following lemma describing the spacing of Toda
particles at time $0$. Its proof, which is a quick consequence of the Bernstein
concentration inequality, is given in Section \ref{app:Estr}.
\begin{lem}\label{lem:Estr}
Adopt Assumption \ref{ass:model}. The event
\begin{equation}\label{eq:stretch-good}
\Bigl\{
\max_{i,j\in\llb N_1,N_2\rrb : i > j}
|\bar q_i - \bar q_j-(F_{\mathrm{str}}^{(T)}(i/T)-F_{\mathrm{str}}^{(T)}(j/T))|\le \mathfrak{A} \sqrt{i-j} T^{-1}(\log N)^2
\Bigr\}
\end{equation}
holds with overwhelming probability.
\end{lem}
	
	\begin{proof}[Proof of Theorem \ref{thm:ev}]
By \eqref{eqn:umphi} with $s'=0$,
with overwhelming probability the hypothesis
$\left|\varphi_0(k)\right|\le T(\log N)^4$ implies that
$\varphi_s(k)\in\llbracket N_1+T(\log N)^{10},N_2-T(\log N)^{10}\rrbracket$
for every $s\le T\log N$. Hence Theorem \ref{ttn} and \eqref{eqn:Qx} in the definition of $\mathfrak{T}$, with $\mathfrak M=T^{3/4}$, give with overwhelming
probability,
\begin{equation}\label{eqn:finite-T-quasiparticle-bound}
\left|T^{-1}Q_k(t)-x_{t/T}^{(T)}
\bigl(\lambda_k;T^{-1}Q_k(0)\bigr)\right|
\le T^{-1/4}\log(|\lambda_k|+2).
\end{equation}
Set $L:=\mathfrak K^{-1}\mathfrak A(\log T)^5$. By Assumption
\ref{ass:model}, Definition \ref{rho0}, and \eqref{eqn:alph_bd},
$m_0^{(T)}(\lambda,q)=m_0(\lambda,q)$ whenever $|q|\le L$.
On the event of Lemma \ref{lem:Estr}, the hypothesis on
$\varphi_0(k)$ gives
$|\bar q_{\varphi_0(k)}|
\le C\mathfrak A(\log N)^4=o(L)$.
Moreover, on the event of Lemma \ref{l0eigenvalues} with $A=\log N$,
$|\lambda_k|\le\log N$. For sufficiently large $T$, we have
$L>(\log T)^3$, $t/T\le\log N<(\log T)^2$, and
$(t/T)(1+|\lambda_k|)=o(L)$. Thus the hypotheses of Lemma
\ref{lem:thermal-completion} hold with $x_0=0$,
$q=\bar q_{\varphi_0(k)}=T^{-1}Q_k(0)$, and this value of $L$. That lemma
replaces $x^{(T)}$ by $x$ at an exponentially small cost; together with
\eqref{eqn:finite-T-quasiparticle-bound}, this yields \eqref{eqn:ev}.
		\end{proof}

	\subsection{Self improving estimates}
	\label{subsec:asymptotic_scattering}

	The next three propositions are at the heart of the proof of Theorem
	\ref{ttn}. The first, second, and third state that the estimates in the first,
	second, and third properties of Definition \ref{mt}, respectively, can be
	improved when $t\leq\mathfrak T$. They may be thought of as discrete
		counterparts of the contraction estimates used to construct $N_t$ (from
	Lemma \ref{lem:fixed_pointout}) in Section \ref{sec:GHD_sol} below. Propositions \ref{prop:qbdtfix}, \ref{vni_out}, and
	\ref{prop:improv_M} are proved in Sections \ref{subsec:qbd1N},
	\ref{subsec:ideas2}, and \ref{subsec:ideas_Qx}, respectively.
	
	\begin{prop}\label{prop:qbdtfix}
		Adopt Assumption \ref{ass:model} and recall Definition \ref{mt}. There exists $\mathfrak A_0=\mathfrak A_0(\mathfrak K,\epsilon)>1$ such that the following holds whenever $\mathfrak A\ge\mathfrak A_0$. Fix $t \in [0, T \log N]$. On an event of overwhelming probability,
	the following holds whenever $\mathfrak{B}(t)$ holds.
	For all $i > j$, $i, j \in \llbracket N_1 + T (\log N)^{4}, N_2-T (\log N)^{4} \rrbracket$,
	\begin{equation}\label{eqn:qbd1Nup}
	q_{i}(t) -q_j(t) \leq \frac{1}{2}C_1 \mathfrak{A} \max(|i-j|, \mathfrak{M}/(2 R)).
	\end{equation}
	Moreover, for all $i -j > \mathfrak{M}/(2 R)$, $i, j \in \llbracket N_1 + T (\log N)^{4}, N_2-T (\log N)^{4} \rrbracket$,
	\begin{equation}\label{eqn:qbd1Nlow}
		q_i(t)-q_j(t) \geq 2c_1\mathfrak A(i-j).
		\end{equation}
		\end{prop}

	\begin{prop}
		
		\label{vni_out} 
		
			Adopt Assumption \ref{ass:model} and recall Definition \ref{mt}. There exists $\mathfrak A_0=\mathfrak A_0(\mathfrak K,\epsilon)>1$ such that the following holds whenever $\mathfrak A\ge\mathfrak A_0$. Fix $t \in [0, T \log N]$.
		The following holds on an event of overwhelming probability. If $\mathfrak{B}(t)$ holds, $k \in \llbracket 1, N \rrbracket$
		satisfies $\varphi_t(k)\in\llbracket N_1+T(\log N)^{9},\,N_2-T(\log N)^{9}\rrbracket$,
		and $i \in \llbracket N_1, N_2 \rrbracket$ satisfies $|i - \varphi_t (k)| > \mathfrak{M}/(2R)$,
		then
		\[
		|u_k (i;t)| \le  \exp(-4R|i-\varphi_t(k)|).
		\]
	\end{prop}

	\begin{prop}\label{prop:improv_M}
		Adopt Assumption \ref{ass:model} and recall Definition \ref{mt}. There exists $\mathfrak A_0=\mathfrak A_0(\mathfrak K,\epsilon)>1$ such that the following holds whenever $\mathfrak A\ge\mathfrak A_0$. Fix $t \in [0, T \log N]$, and denote $\tau = t/T$.
	For a constant \(C>0\) independent of \(\mathfrak A\), on an overwhelmingly probable event the following
	holds. If $\mathfrak{B}(t)$ holds, and $k \in \llbracket 1, N \rrbracket$ satisfies $\varphi_t(k) \in \llbracket N_1 + \frac{1}{2}T (\log N)^{10}, N_2-\frac{1}{2}T (\log N)^{10} \rrbracket$, then
	\begin{flalign}
			\label{qkt} 
			\big| Q_k (t) - T x_\tau^{(T)} (\lambda_k; \bar q_{\varphi_0(k)}) \big| \le C R^{-1} \mathfrak{M} \log (|\lambda_k|+2).
	\end{flalign}
	\end{prop}
	
	Propositions \ref{vni_out}, \ref{prop:qbdtfix}, and \ref{prop:improv_M}
	suggest the following approach to prove Theorem \ref{ttn}. Assume
	$\mathfrak T<T\log N$ and set $t=\mathfrak T$. Then Proposition
	\ref{vni_out} implies that \eqref{eqn:uMR} holds with
	$e^{\mathfrak M-R|j-\varphi_t(k)|}$ replaced by
	$e^{\mathfrak M/2-2R|j-\varphi_t(k)|}$.
	Moreover, Proposition \ref{prop:qbdtfix} implies that \eqref{eqn:qbd} and \eqref{eqn:qbd-lower} also hold with $\mathfrak{M}$ replaced by $\mathfrak{M}/2$, 
	and Proposition \ref{prop:improv_M} implies a similar improvement on the
	bound \eqref{eqn:Qx}. Thus, if the Toda particles $q_i(t)$, quasiparticles
	$Q_k(t)$, and eigenvector entries $u_k(j;t)$ satisfy a weak continuity bound,
	the estimates \eqref{eqn:qbd}, \eqref{eqn:qbd-lower}, \eqref{eqn:uMR}, and
	\eqref{eqn:Qx} continue to hold slightly beyond $t=\mathfrak T$,
	contradicting the maximality of $\mathfrak{T}$. The proposition below provides exactly this
	continuity; its proof is given in Section \ref{subsec:proof-prop-cont}.
	
	\begin{prop}\label{prop:cont}
		Adopt Assumption \ref{ass:model} and recall Definition \ref{mt}. There exists $\mathfrak A_0=\mathfrak A_0(\mathfrak K,\epsilon)>1$ such that the following holds with overwhelming probability whenever $\mathfrak A\ge\mathfrak A_0$.
	\begin{enumerate}
	\item For every $i \in \llbracket N_1, N_2 \rrbracket$ and all $t,t' \geq 0$
satisfying $|t-t'| \leq 2 N^{-1}$,
\begin{equation}
\label{eqn:qbd_cont}
|q_i(t)-q_i(t')|
\le
4(\log N) /N .
\end{equation}

\item Suppose $t \in [0, T \log N] \cap N^{-1} \mathbb{Z}$, the event $\mathfrak{B}(t)$ holds, and $t' \in [0, T \log N]$ satisfies $|t'-t| \leq 2 N^{-1}$.
For all $k \in \llbracket 1, N \rrbracket$ such that
 $\varphi_{t'}(k) \in \llbracket N_1 + T(\log N)^{10}, N_2 - T (\log N)^{10} \rrbracket$, and each $j \in \llbracket N_1, N_2 \rrbracket$,
		\begin{equation}
			|u_k (j;t')| \le e^{\mathfrak{A} \mathfrak{M} - R |j-\varphi_{t'} (k)|}. 	\label{eqn:uMRct}
		\end{equation} 

\item Suppose $t = \tau T \in [0,T \log N] \cap N^{-1} \mathbb{Z}$, the event $\mathfrak{B}(t)$ holds, and $t' = \tau' T \in [0, T \log N]$ satisfies $|t'-t| \leq 2 N^{-1}$. For all $k \in \llbracket 1, N \rrbracket$ such that
	 $\varphi_{t'}(k) \in \llbracket N_1 + T(\log N)^{10}, N_2 - T (\log N)^{10} \rrbracket$,
			\begin{align} 
				|Q_k (t) -Q_k (t') | &\le \mathfrak{A}^2 \mathfrak{M} \log (|\lambda_k| + 2),  \label{eqn:Qxct}\\
				|T  x_\tau^{(T)} (\lambda_k; \bar q_{\varphi_0(k)}) -T  x_{\tau'}^{(T)} (\lambda_k; \bar q_{\varphi_0(k)})| &\le
				(\log N)^5/N. \label{eqn:Qxct2}
			\end{align}
	\end{enumerate}
	Moreover, in the setup of the second and third parts above,
	\begin{equation}\label{eqn:lc_weak}
	\varphi_t(k) \in \llbracket N_1 + T(\log N)^{10}/2, N_2 - T (\log N)^{10}/2 \rrbracket.
	\end{equation}

	\end{prop}

Next, we have the following lemma stating that the eigenvectors of $\L(0)$
are quantitatively localized; its proof is in Section \ref{app:lemLocR}.
This will ensure that \eqref{eqn:uMR} holds at time $0$.

	\begin{lem} 		
		\label{lem:time-0-locR}
			Adopt Assumption \ref{ass:model} and recall Definition \ref{mt}. There exists $\mathfrak A_0=\mathfrak A_0(\mathfrak K,\epsilon)>1$ such that the following holds with overwhelming probability whenever $\mathfrak A\ge\mathfrak A_0$. Fix any index $j \in \llbracket 1, N \rrbracket$, and let $\varphi_0(j) \in \llbracket N_1, N_2 \rrbracket$ be its $\zeta$-localization center from Definition \ref{def:quasi_intro}. Then, for
		 any $i \in \llbracket N_1, N_2 \rrbracket$, we have $| u_j (i;0) | \le e^{T^{2/3}-20 R |i-\varphi_0(j)|}$.
	\end{lem}

	With these preparations in place, we can give the proof of Theorem \ref{ttn}.
	
	\begin{proof}[Proof of Theorem \ref{ttn}]
			Using the constant $\delta$ from Definition \ref{mt},
			suppose $\mathfrak{M} \in [T^{3/4-\delta/2}, T^{3/4}]$.
			Denote by $\mathfrak{T} = \mathfrak{T}_{\mathfrak{M}}$
			 the stopping time from Definition \ref{mt},
			  and let
			 $\mathfrak{B}_{\mathfrak{M}}(t)$ denote the event from 
			 Definition \ref{mt} with its $\mathfrak{M}$-dependence made explicit. In addition,
			 suppose $\mathfrak{A}$ in Assumption \ref{ass:model} is chosen large enough. 
			 
			 First, note that $\mathfrak{B}_{\mathfrak{M}}(0)$ holds with overwhelming probability.
			 Item \ref{item:T3} holds since
			 \[
			 Q_k(0)=q_{\varphi_0(k)}(0)=T\bar q_{\varphi_0(k)}
			 =T x_0^{(T)}(\lambda_k;\bar q_{\varphi_0(k)}),
			 \]
			 while Items \ref{item:T1} and \ref{item:T2} hold by Lemmas
			 \ref{lem:Estr} and \ref{lem:time-0-locR}, respectively.
			 
			 Propositions \ref{prop:qbdtfix}, \ref{vni_out}, and \ref{prop:improv_M}, together with a union bound over $t\in[0,T\log N]\cap N^{-1}\mathbb Z$, give an overwhelmingly probable event on which, 
			 at every mesh time $t \in [0,\mathfrak T_{\mathfrak M}] \cap N^{-1} \mathbb{Z}$, the three bounds in Definition \ref{mt} hold with $\mathfrak M$ replaced by $\mathfrak M/2$; in other words,
			 $\mathfrak{B}_{\mathfrak{M}/2}(t)$ holds. 
			 Iterating these fixed time improvements
			 $\lceil100\log_2\mathfrak A\rceil$ times, we conclude that at every such mesh time the same bounds hold with $\mathfrak M$
			  replaced by $\mathfrak M\mathfrak A^{-100} \in [T^{3/4-\delta}, T^{3/4}]$; hence, the bounds in the conclusions of
			  Propositions \ref{prop:qbdtfix}, \ref{vni_out}, and \ref{prop:improv_M} hold for each 
			  $t \in [0,\mathfrak T_{\mathfrak M}] \cap N^{-1} \mathbb{Z}$.
			  Let us denote by $\mathsf{E}_0$ the event on which this is the case.

			Now let $\mathfrak{M} = T^{3/4}$. 
			Let $\mathsf{E}_1$ be the event provided by Proposition \ref{prop:cont} with $\mathfrak{M}$
			there replaced by $\mathfrak{M} 2^{-k}$ for $k = \lceil  100 \log_2 \mathfrak{A} \rceil$. In what follows we restrict to
			the overwhelmingly probable event $\mathsf{E}_{\rm final} = \mathsf{E}_0 \cap \mathsf{E}_1 $.
			
			 First, on these
		 	events,
			the combination of \eqref{eqn:qbd_cont} together with \eqref{eqn:qbd1Nup} (and with \eqref{eqn:qbd1Nlow}) implies that 
			\eqref{eqn:qbd} (and \eqref{eqn:qbd-lower}, respectively) in fact holds for all $t \in [0, \min(T \log N, \mathfrak{T} + \frac{1}{N})]$.

			Moreover, \eqref{eqn:uMRct}, applied with $\mathfrak{M}$ 
			replaced by $\mathfrak{M} \mathfrak{A}^{-100}$, 
			allows us to bound 
		\[
			| u_k(j; t') | \leq e^{\mathfrak{M} \mathfrak{A}^{-99}-R |j-\varphi_{t'}(k)|}
			\]
			 for all $|t-t'| \leq 2/N$, $ t\in [0, \mathfrak{T}] \cap \frac{1}{N} \mathbb{Z}$, and
			$\varphi_{t'}(k) \in \llbracket N_1 + T(\log N)^{10}, N_2 - T (\log N)^{10} \rrbracket$. This implies that
			\eqref{eqn:uMR} holds for times $t \in [0, \min(T \log N, \mathfrak{T} + \frac{1}{N})]$.

			Finally, \eqref{eqn:Qxct} and \eqref{eqn:Qxct2} applied with $\mathfrak{M}$ replaced by $\mathfrak{M} \mathfrak{A}^{-100}$,
			in combination with \eqref{qkt} (which holds by \eqref{eqn:lc_weak}), implies 
			that the bound \eqref{eqn:Qx} holds for $t \in [0, \min(T \log N, \mathfrak{T} + \frac{1}{N})]$. 
			
			This leads to a contradiction unless $\mathfrak{T} = T \log N$ on $\mathsf{E}_{\rm final}$. \qedhere
	\end{proof}

	\subsection{Proof of Proposition \ref{prop:improv_M}}
	\label{subsec:ideas_Qx}
	
	As explained in Section \ref{subsec:asymptotic_scattering}, the main results
	leading to Theorem \ref{ttn} are Propositions \ref{prop:qbdtfix},
	\ref{vni_out}, and \ref{prop:improv_M}. 
	Since \eqref{eqn:Qx} approximates the quasiparticles $(Q_k)$, 
	which constitute a permutation of the Toda particles $(q_i)$, 
	we will be able to deduce Proposition \ref{prop:qbdtfix} using \eqref{eqn:Qx} (and properties of the continuum trajectories $x_t^{(T)}$);
	 further discussion on this will be deferred to Section \ref{sec:qbd},
	 which contains the proof of Proposition \ref{prop:qbdtfix}. 
	 We will explain how Proposition \ref{prop:qbdtfix} is relevant for the 
	 proof of Proposition \ref{vni_out} in Section \ref{subsec:ideas2}.
	 In this
	subsection we state several intermediate propositions
	used to prove Proposition \ref{prop:improv_M},
	and then we provide the proof of the latter below. 
	The proofs of Proposition \ref{prop:improv_M} and
	Proposition \ref{prop:qbdtfix} are logically independent. 
	Throughout this subsection, recall $N_{\tau}^{(T)}$ from Section
	\ref{subsec:tc}.

	To show Proposition \ref{prop:improv_M}, we use the following result (due to
	\cite{Mos75}) indicating that the logarithms of the first entries of the Lax
	matrix eigenvectors are angle variables, in that they essentially evolve
	linearly under the Toda lattice.

	\begin{lem}[{\cite[Section 3]{Mos75}}]
	\label{uknt}
	For any $t \in \mathbb{R}_{\ge 0}$ and $k \in \llbracket 1, N \rrbracket$, we have
	\begin{equation}\label{eqn:ukN1lin}
		u_k (N_1; t)^2  = e^{-\lambda_k t}  u_k (N_1; 0)^2 \cdot \Bigg( \displaystyle\sum_{j=1}^N e^{-\lambda_j t}  u_j (N_1; 0)^2 \Bigg)^{-1}.
	\end{equation}
	\end{lem}

	In view of Lemma \ref{uknt} (and Lemma \ref{ltt}), to analyze the
	quasiparticles $(Q_k)$, we next seek to approximate them through the first
	entries of Lax matrix eigenvectors and the eigenvalues. This is done through
	the following proposition (which further expresses them in terms of the last
	entries, as well), to be proven in Section \ref{subsec:ideas1}. Below we
	recall $\mathfrak B(t)$ from Definition \ref{mt}.

	\begin{prop}
	\label{vn1n2out}
		Adopt Assumption \ref{ass:model}. There exists $\mathfrak A_0=\mathfrak A_0(\mathfrak K,\epsilon)>1$ such that the following holds whenever $\mathfrak A\ge\mathfrak A_0$. Fix $t \in [0, T \log N]$. There is a
		constant $C > 0$, independent of $\mathfrak{A}$, such that
		the following holds with overwhelming probability.
		If $\mathfrak{B}(t)$ holds and $k\in\llbracket 1,N\rrbracket$ satisfies
		$\varphi_t(k)\in
\llbracket N_1+2T(\log N)^{10}, N_2-2T(\log N)^{10}\rrbracket$,
		we have
		\begin{align}
			\left| 2 \log |u_k (N_1;t)| - q_{N_1} (t) + Q_k (t) + 2 \displaystyle\sum_{j: Q_j (t) < Q_{k} (t)} \log |\lambda_k - \lambda_j|\right| & \leq  C R^{-1} \mathfrak{M} \log (|\lambda_k|+2) ; \label{eqn:ukN1} \\
			\left| 2 \log |u_k (N_2;t)|  - Q_{k} (t) + q_{N_2} (t) + 2 \displaystyle\sum_{j: Q_j(t)> Q_k(t) } \log |\lambda_k - \lambda_j|   \right| &\leq  C R^{-1}  \mathfrak{M} \log (|\lambda_k|+2) . \label{eqn:ukN2}
		\end{align}
	\end{prop}

	Observe that, by identifying $Q_k(t)$ with $Tx$ and formally replacing the
	sum of $\log|\lambda_k-\lambda_j|$ with $T$ times an integral of
	$\log|\lambda_k-\lambda|$, the expression \eqref{eqn:ukN1} for
	$2\log|u_k(N_1;t)|$ can be viewed as a discrete counterpart (of a shift) of
	equation \eqref{eqn:Xt} defining $X_t(\lambda_k,x)$. Similarly, the linear
	evolution \eqref{eqn:ukN1lin} for $2\log|u_k(N_1;t)|$ can be viewed as a
	discrete counterpart of that of $X_t$ from equation \eqref{eqn:Xtlin}; as
	such, $X_t$ may be thought of as the continuum angle variables governing the
	hydrodynamic equation \eqref{eqn:rhoGHDintro}.

	The next proposition will be proven in Section \ref{subsec:quasi_thmpf}.
	Its first part \eqref{eqn:sc} provides a form of the asymptotic scattering
	relation mentioned in Section \ref{subsec:over}, and it is obtained by
	subtracting \eqref{eqn:ukN1} at $t=0$ from \eqref{eqn:ukN1} (and using a
	cancellation involving $q_{N_1}(0)$, $q_{N_1}(t)$, and the denominator on the right
	side of \eqref{eqn:ukN1lin}; see \eqref{eqn:Btbd} below). Its second part
	\eqref{eqn:scattering-continuum-cutoff} is obtained from replacing the
	$Q_j(t)$ by $T x_{t/T}^{(T)}(\lambda_j;Q_j(0)/T)$ as in \eqref{eqn:Qx}, and
	further replacing the sums in \eqref{eqn:sc} with integrals (which will be
	enabled by Proposition \ref{prop:scat_conc} below). Observe that the right
	side of \eqref{eqn:scattering-continuum-cutoff} is smaller than that of
	\eqref{eqn:Qx} by a factor of $R^{-1}$; this will eventually be due to the
	fact that the density of terms in the sums in \eqref{eqn:sc} is of order
	$\mathfrak A^{-1}<R^{-1}$.

	\begin{prop}[Discrete and continuum scattering relations]
	\label{prop:discrete-continuum-scattering}
	Adopt Assumption \ref{ass:model}. There exists
	$\mathfrak A_0=\mathfrak A_0(\mathfrak K,\epsilon)>1$ such that the following
	holds whenever $\mathfrak A\ge\mathfrak A_0$. Fix
	$t=T\tau\in[0,T\log N]$. There is a constant $C>0$, independent of
	$\mathfrak A$, such that, with overwhelming probability, whenever
	$\mathfrak B(t)$ holds and
	$|\varphi_0(k)|\le T(\log N)^{10}$, the following two estimates hold.
	\begin{enumerate}
	\item The quasiparticles satisfy
	\begin{align}
	\label{eqn:sc}
	\Bigg|Q_k(t)-Q_k(0)-\lambda_kt
	&-2\sum_{j:Q_j(0)<Q_k(0)}\log|\lambda_k-\lambda_j| \\
	&+2\sum_{j:Q_j(t)<Q_k(t)}\log|\lambda_k-\lambda_j|\Bigg|
	\le CR^{-1}\mathfrak M\log(|\lambda_k|+2).\notag
	\end{align}

	\item We have
	\begin{equation}
	\label{eqn:scattering-continuum-cutoff}
	\begin{aligned}
	&\Bigg|Q_k(t)-q_{\varphi_0(k)}(0)-\lambda_kt
	-2T\int_{\mathbb R}\log|\lambda_k-\lambda| \\
	&\hspace{28mm}\cdot\Bigl(N_0^{(T)}(\lambda,\bar q_{\varphi_0(k)})
	-N_\tau^{(T)}(\lambda,T^{-1}Q_k(t))\Bigr)\,d\lambda\Bigg| \\
	&\hspace{35mm}\le CR^{-1}\mathfrak M\log(|\lambda_k|+2).
	\end{aligned}
	\end{equation}
	\end{enumerate}
	\end{prop}
	
	Proposition \ref{prop:discrete-continuum-scattering} is stated for indices
	$k$ satisfying $|\varphi_0(k)|\le T(\log N)^{10}$. For larger indices,
	Assumption \ref{ass:model} indicates that the initial data
	$(\bm a(0);\bm b(0))$ is close to thermal equilibrium. The following lemma
	indicates that, in this region, the quasiparticles $(Q_k)$ continue to
	approximate their limiting trajectories. Its proof is given in Section
	\ref{app:thermal_edge} and proceeds by first reducing to the case of thermal
	equilibrium (using Lemma \ref{lem:edge-thermal-coupling} below), in which
	setting \eqref{eq:thermal-tail-quasiparticle} was essentially shown in
	\cite{Agg25}.

\begin{lem}[Quasiparticle approximation in the thermal region]
\label{lem:thermal-tail-quasiparticle}
Adopt Assumption \ref{ass:model}, fix \(t=T\tau\in[0,T\log N]\), and set
\(H_T:=T(\log N)^{10}\) and
\(\mathcal I:=
\llbracket N_1+\lfloor H_T/2\rfloor,-\lfloor H_T/2\rfloor\rrbracket
\cup
\llbracket \lfloor H_T/2\rfloor,N_2-\lfloor H_T/2\rfloor\rrbracket\).
There is a constant \(C>0\) such that, with overwhelming probability,
\begin{equation}
\label{eq:thermal-tail-quasiparticle}
\max_{k:\,\varphi_0(k)\in\mathcal I}
\left|
Q_k(t)
-
T x_\tau^{(T)}(\lambda_k; Q_k(0)/T)
\right|
\le
T^{1/2}(\log N)^C.
\end{equation}
\end{lem}

	\begin{proof}[Proof of Proposition \ref{prop:improv_M}]
	By Lemma \ref{lem:thermal-tail-quasiparticle} and Proposition \ref{prop:outside-strip-confinement}, 
	on an event of overwhelming probability the conclusion holds when
	$|\varphi_0(k)|>T(\log N)^{10}$. We may therefore assume that
	$|\varphi_0(k)|\le T(\log N)^{10}$. Set
	\[
	a:=T^{-1}Q_k(t),
	\qquad
	x_*:=x_\tau^{(T)}(\lambda_k;\bar q_{\varphi_0(k)}),
	\qquad
	L_k:=\log(|\lambda_k|+2).
	\]
	The second part of Proposition \ref{prop:discrete-continuum-scattering}
	gives
	\begin{align*}
	Q_k(t)=q_{\varphi_0(k)}(0)+\lambda_kt
	+2T\int_{\mathbb R}\log|\lambda_k-\lambda|
	\Bigl(N_0^{(T)}(\lambda,\bar q_{\varphi_0(k)})
	-N_\tau^{(T)}(\lambda,a)\Bigr)\,d\lambda
	+O(R^{-1}\mathfrak ML_k).
	\end{align*}
	On the other hand, Lemma \ref{xt2}, applied to $x_\tau^{(T)}$,
	gives
	\begin{align*}
	Tx_*=q_{\varphi_0(k)}(0)+\lambda_kt
	+2T\int_{\mathbb R}\log|\lambda_k-\lambda|
	\Bigl(N_0^{(T)}(\lambda,\bar q_{\varphi_0(k)})
	-N_\tau^{(T)}(\lambda,x_*)\Bigr)\,d\lambda.
	\end{align*}
	Subtracting the last two displays and using the definition of
	$X_\tau^{(T)}$, we obtain
	\[
	T\big|X_\tau^{(T)}(\lambda_k,a)
	-X_\tau^{(T)}(\lambda_k,x_*)\big|
	\le CR^{-1}\mathfrak ML_k.
	\]
	The lower bound $\partial_xX_\tau^{(T)}\ge \frac12$ from Lemma
	\ref{lem:Jlow} therefore implies
	\[
	T|a-x_*|
	=
	\big|Q_k(t)-Tx_\tau^{(T)}(\lambda_k;\bar q_{\varphi_0(k)})\big|
	\le CR^{-1}\mathfrak M\log(|\lambda_k|+2),
	\]
	which proves the proposition.
	\end{proof}
		
		\subsection{Proof of Proposition \ref{vn1n2out}}
		\label{subsec:ideas1}

	To show Proposition \ref{vn1n2out}, we begin with the following known lemma
	expressing eigenvectors of tridiagonal random matrices in terms of the
	eigenvalues of their truncated submatrices.
	\begin{lem}[{\cite[Lemmas 3.9 and 3.10]{Agg25a}}]\label{lem:vkn1_formula}
		Let $\L$ be a symmetric, tridiagonal matrix indexed by $\llbracket N_1, N_2 \rrbracket$ with positive off diagonal entries. Define the differences of the $q_i$ by
		$q_i-q_j:=2\sum_{\ell=i}^{j-1}\log L_{\ell,\ell+1}$, for $N_1\le i<j\le N_2$.
		Let $\eig\L=(\lambda_1,\lambda_2,\ldots,\lambda_N)$ and, for each
		$k\in\llbracket1,N\rrbracket$, let
		$\mathbf v_k=(v_k(i))_{i=N_1}^{N_2}$ denote the nonnegatively normalized
		unit eigenvector of $\L$ with eigenvalue $\lambda_k$.
		For every $k \in \llbracket 1,N\rrbracket$ and every $i \in \llbracket N_1, N_2 \rrbracket$ with $v_k(i) \neq 0$,
	\begin{align} \label{eqn:vkn1_formula}
	2 \log |v_k (N_1)| &= 2 \log |v_k (i)| + q_{N_1} - q_{i} - 2 \displaystyle\sum_{\mu \in \eig \L^{[ N_1, i-1 ]}} \log |\lambda_k - \mu|\\
	2 \log |v_k (N_2)| &= 2 \log |v_k (i)| +q_{i} - q_{N_2} - 2 \displaystyle\sum_{\mu \in \eig \L^{[ i+1 , N_2]}} \log |\lambda_k - \mu| . \label{eqn:vkn2_formula}
	\end{align} 
	\end{lem}

		To show Proposition \ref{vn1n2out}, let $\L=\L(t)$ and
		$\varphi=\varphi_k(t)$. Using Lemma \ref{lem:vkn1_formula}, we have the exact
		identity
		\begin{equation}\label{eqn:vn1exact}
		2\log |u_k(N_1;t)|
		=
		2\log |u_k(\varphi;t)|
		+
		q_{N_1}(t)-q_{\varphi}(t)
		-
		2\sum_{\mu\in \eig  \L^{[N_1, \varphi-1]}(t)}\log|\lambda_k-\mu|.
\end{equation}
	If eigenvectors of $\L$ are sufficiently localized, which is
	the case if $t \leq \mathfrak{T}$ by \eqref{eqn:uMR}, we can replace the sum
	over $\mu\in \eig  \L^{[N_1, \varphi-1]}$ with a sum over eigenvalues of $\L(t)$, as follows.
	The lemma below is proved in Section \ref{subsec:IStm2}.

\begin{lem}\label{lem:cutLbd}
Adopt Assumption \ref{ass:model}. There exists $\mathfrak A_0=\mathfrak A_0(\mathfrak K,\epsilon)>1$ such that the following holds whenever $\mathfrak A\ge\mathfrak A_0$. Fix $t \in [0, T \log N]$. There is a
	constant $C > 0$, independent of $\mathfrak{A}$, such that
	the following holds with overwhelming probability. 
	If $\mathfrak{B}(t) $ holds and $k \in \llbracket 1, N \rrbracket$
	 satisfies $ \varphi_t(k) \in \llbracket N_1+2T(\log N)^{10}, N_2 -2 T(\log N)^{10}\rrbracket$, denoting $\varphi = \varphi_t(k) $,
	we have
	\begin{equation}
		\left| \displaystyle\sum_{\mu \in \eig \L^{[ N_1, \varphi-1 ]}(t)} \log |\lambda_k - \mu| -  \displaystyle\sum_{j: Q_j (t) < Q_{k} (t)} \log |\lambda_k - \lambda_j|\right|  \leq  C R^{-1} \mathfrak{M} \log (|\lambda_k|+2) . \label{eqn:itoQ1} 	\end{equation}
\end{lem}

	\begin{proof}[Proof of Proposition \ref{vn1n2out}]
	We only show \eqref{eqn:ukN1}, as the proof of \eqref{eqn:ukN2} follows by
	symmetry. This follows from \eqref{eqn:vn1exact}, Lemma
	\ref{lem:cutLbd}, and the fact that
	$2|\log|u_k(\varphi;t)||\leq2|\log\zeta|\leq(\log N)^2$.
\end{proof}

\subsection{Proof of Proposition \ref{vni_out}}
	\label{subsec:ideas2}

To show Proposition \ref{vni_out}, we again use Lemma
\ref{lem:vkn1_formula}.

%

Let $v_k(i) \coloneqq u_k (i; t)$ and $\varphi = \varphi_t(k)$. Let us assume without loss of generality that $i < \varphi$, as the case $i > \varphi$ is similar.
By subtracting \eqref{eqn:vkn1_formula} by \eqref{eqn:vkn1_formula} with
$i$ replaced by $\varphi$, we obtain (setting $\L = \L(t)$)
	\begin{equation} \label{eqn:vki_exact}
	 2 \log |v_k (i)|  = 2 \log |v_k (\varphi)| + q_{i} - q_{\varphi} + 2 \displaystyle\sum_{\mu \in \eig \L^{[N_1, i-1 ]}} \log |\lambda_k - \mu| - 2 \displaystyle\sum_{\mu \in \eig \L^{[N_1, \varphi-1 ]}} \log |\lambda_k - \mu|.
	\end{equation} 
	First, note that we may discard $ 2 \log |v_k (\varphi)| $ since $\varphi_t$ is chosen as a $\zeta$-localization center bijection
	for $\L(t)$, 
	and therefore 
	\begin{equation}\label{eqn:zetbd}
	2 \left| \log |v_k (\varphi)| \right| \leq  (\log N)^2.
	\end{equation}
	The following lemma, proven in Section \ref{subsec:ideas4}, bounds the difference of the two sums on the right side of \eqref{eqn:vki_exact}.
	
	\begin{lem}\label{lem:log_sum_dif}
		Adopt Assumption \ref{ass:model}. There exists $\mathfrak A_0=\mathfrak A_0(\mathfrak K,\epsilon)>1$ such that the following holds whenever $\mathfrak A\ge\mathfrak A_0$. Fix $t \in [0, T \log N]$, and let $\mathfrak{M}$ be as in Definition \ref{mt}.
		There is a constant $C> 0$, independent of $\mathfrak{A}$, such that the following holds on an event of overwhelming probability. Assume $\mathfrak{B}(t) $ holds,
		$\varphi = \varphi_t(k)\in\llbracket N_1+T(\log N)^{9},\,N_2-T(\log N)^{9}\rrbracket $,
		and let $i \in \llbracket N_1 , N_2  \rrbracket$ be such that $- T (\log N)^3 < i - \varphi_t (k) \le -(2R)^{-1} \mathfrak{M}$.
		Then, we have
		\begin{equation}\label{enq:iphilsbd}
		2 \displaystyle\sum_{\mu \in \eig \L^{[N_1, i-1 ]}(t)} \log |\lambda_k - \mu| - 2 \displaystyle\sum_{\mu \in \eig \L^{[N_1, \varphi-1 ]}(t)} \log |\lambda_k - \mu|
		\leq C |i - \varphi| .
		\end{equation}
	\end{lem}
	
	\begin{proof}[Proof of Proposition \ref{vni_out}]
	Suppose $i < \varphi$; the other case follows by symmetry. If $|i - \varphi| > T (\log N)^3$, then Proposition \ref{prop:outside-strip-confinement}
		 immediately implies the result, so we may assume otherwise.
	Proposition \ref{prop:qbdtfix} shows that on an event of
	overwhelming probability, whenever $\mathfrak B(t)$ holds,
	$\varphi-i>\mathfrak M/(2R)$, and $i$ and $\varphi$ belong to
	$\llbracket N_1+T(\log N)^7,N_2-T(\log N)^7\rrbracket$,
	\begin{equation}\label{eqn:qdlb}
	q_{i} - q_{\varphi} \leq 2c_1 \mathfrak{A} (i-\varphi).
	\end{equation}
	The assumption on $\varphi$ and the bound $|i-\varphi|\le T(\log N)^3$
	ensure that the required conditions on $i$ and $\varphi$ are satisfied.
	 
	 Starting with \eqref{eqn:vki_exact}, Lemma \ref{lem:log_sum_dif} and
	 the observations \eqref{eqn:zetbd} and \eqref{eqn:qdlb} yield
	\[
	2 \log |v_k (i)| \leq 2c_1 \mathfrak{A} (i-\varphi)+ (\log N)^2 + O\big( |i-\varphi| \big) \leq 10 R (i-\varphi).
	\]
	where this last inequality is due to the fact that
	$\mathfrak{A}$ is chosen large enough. This implies the result.
	\end{proof}

	\section{Preliminaries}
	\label{sec:prelim}
	We collect various preliminary facts about random Lax matrices at time $0$,
	as well as several general
	facts about tridiagonal symmetric matrices.

	\subsection{Eigenvalue differences}
	\label{subsec:loclems}

	For a real symmetric matrix $\mathbf{M}$ and $\eta>0$, define
	\[
	\mathsf{SEP}_{\mathbf{M}}(\eta)
	:=\left\{\min_{\substack{\nu,\nu'\in\eig\mathbf{M}\\ \nu\ne\nu'}}
	|\nu-\nu'|\ge\eta\right\}.
	\]
	We omit the proof of the following lemma, which is entirely analogous to
	that of \cite[Lemma 3.18]{Agg25a}.
	
	\begin{lem}[{\cite[Lemma 3.18]{Agg25a}}]
		
		\label{eigenvalues0}
		
		Adopt Assumption \ref{ass:model}. For any real number $\delta>0$,
		\[
		\mathbb P\bigl[\mathsf{SEP}_{\L(0)}(\delta)\bigr]
		\ge 1-c^{-1}(\delta N^3+e^{-cN^2}).
		\]
		
	\end{lem}


We next provide a technical lemma which will be useful in the proof of Proposition \ref{vni_out};
its statement uses the notation $\mathbf{M}^{[a,b]}$ defined in Section \ref{notation}.
Its proof appears in Section \ref{app:suffpf}.

	\begin{lem}\label{lem:suff}
	Let $\L=[L_{ij}]_{i,j\in\llbracket N_1,N_2\rrbracket}$ be a real symmetric tridiagonal matrix with $L_{i,i+1}>0$ for every $i\in\llbracket N_1,N_2-1\rrbracket$. Let $\lambda\in\eig\L$, and let $\varphi$ be a $\zeta$-localization center for $\lambda$, where $\zeta\in(0,1]$. Set
	\[
	\H_-:=\L^{[N_1,\varphi-1]},
	\qquad
	\H_+:=\L^{[\varphi+1,N_2]}.
	\]
	Let $\mathcal A_-,\mathcal A_+\subseteq\eig\L\setminus\{\lambda\}$ be disjoint, and suppose there are injections
	$
	\pi_-:\mathcal A_-\longrightarrow\eig \H_-$
	and
	$
	\pi_+:\mathcal A_+\longrightarrow\eig \H_+$.
	Define
	\[
	\mathcal U_\pm:=\eig \H_\pm\setminus\pi_\pm(\mathcal A_\pm),
	\qquad
	\mathcal B:=(\eig\L\setminus\{\lambda\})\setminus(\mathcal A_-\cup\mathcal A_+).
	\]
	If
	\begin{equation}\label{eq:matched-log-error}
	\sum_{\theta\in\mathcal A_-}
	\left|\log|\lambda-\pi_-(\theta)|-\log|\lambda-\theta|\right|
	+
	\sum_{\theta\in\mathcal A_+}
	\left|\log|\lambda-\pi_+(\theta)|-\log|\lambda-\theta|\right|
	\leq\mathcal E,
	\end{equation}
	then
	\begin{equation}\label{eq:unmatched-singular-bound}
	-\sum_{\mu\in\mathcal U_-}\log|\lambda-\mu|
	-\sum_{\mu\in\mathcal U_+}\log|\lambda-\mu| 
	=
	-\sum_{\theta\in\mathcal B}\log|\lambda-\theta|
	 + O(\left|\log\zeta\right|)+O(\mathcal E).
	\end{equation}
	The implicit constants in \eqref{eq:unmatched-singular-bound} are absolute.
	\end{lem}

\subsection{Localization at time zero}
\label{app:tz-localization}

In this section we state localization estimates for the eigenvectors of
$\mathbf{L}(0)$, where throughout this section we adopt Assumption
\ref{ass:model} and the notation of Definition \ref{def:quasi_intro}.

The following result gives exponential decay for fractional moments of the off diagonal resolvent entries of the Lax matrix for the open Toda lattice under the above initial condition. The bound \eqref{gijs} follows from\footnote{While it is assumed in \cite[Theorem 4]{Sch09} that the entries of $\L(0)$ are identically distributed, it is quickly verified (and, indeed mentioned below \cite[Equation (2.1)]{Sch09}) that the proof applies when the entries of $\L(0)$ are uniformly \emph{fluctuation regular}, as holds in our setting by \eqref{munu} and \eqref{eqn:mudgc}.} \cite[Theorem 4]{Sch09} (and the remark following it) when $z \in \mathbb{R}$ is real; together with \cite[Equation (B.8)]{FFCL}, this implies that it continues to hold for any complex $z \in \mathbb{C}$.

\begin{lem}[{\cite{Sch09,FFCL}}]
\label{gijexponential}
For any real number $s \in (0, 1)$, there exists a constant $c = c(s) > 0$ such that the following holds. For any $z \in \mathbb{C}$, denote $\mathbf{G} (z) = [G_{ij} (z) ] = (\L(0) - z\mathbf{I})^{-1}$. We have
\begin{flalign}
\label{gijs}
\displaystyle\sup_{z \in \mathbb{R}} \mathbb{E} \big[ |G_{ij} (z) |^s \big] \le c^{-1} e^{-c|i-j|}.
\end{flalign}
\end{lem}

The below corollary, which is due to \cite[Theorem A.1]{FFCL} (using Lemma \ref{gijexponential} to verify its hypotheses), will quickly imply that eigenvectors of $\L(0)$ are localized.

\begin{cor}[{\cite[Theorem A.1]{FFCL}}]
\label{uijexponential}
There exists a constant  $c = c(\mathfrak{A}) > 0$ so that
\begin{flalign*}
\displaystyle\max_{k \in \llbracket 1, N \rrbracket} \displaystyle\max_{i, j \in \llbracket N_1, N_2 \rrbracket} \mathbb{E} \big[ | u_k (i;0) u_k (j;0) | \big] \le c^{-1} e^{-c|i-j|}.
\end{flalign*}
\end{cor}

We omit its proof, which is entirely analogous to that of
\cite[Lemma 5.1]{Agg25a}, using Corollary \ref{uijexponential}, a Markov
bound, and the definition of a $\zeta$-localization center.

\begin{lem}
\label{bijectionl}
There exists a constant $c = c(\mathfrak{A}) >0$ such that the following holds with overwhelming probability. Assume that $\zeta \ge e^{-200(\log N)^{3/2}}$. Fix any index $j \in \llbracket 1, N \rrbracket$, and let $\varphi_j \in \llbracket N_1, N_2 \rrbracket$ denote any $\zeta$-localization center for $\mathbf{u}_j(0)$.
\begin{enumerate}
\item For any $i \in \llbracket N_1, N_2 \rrbracket$ with $| i - \varphi_j | \ge (\log N)^2/2$, we have $| u_j (i; 0) | \le e^{-c|i-\varphi_j|}$.
\item If $i \in \llbracket N_1, N_2 \rrbracket$ satisfies $| i - \varphi_j | \ge (\log N)^2/2$, then $i$ is not a $\zeta$-localization center for $\mathbf{u}_j(0)$.
\end{enumerate}
\end{lem}

\subsection{Eigenvalue matching}
\label{subsec:matchlem}

Next we describe a situation in which certain eigenvalues
 of a localized tridiagonal matrix $\L$ can be matched
 to those of another Lax matrix $\widetilde \L$, whose entries
 are close to those of $\L$
 on a subset of indices.
\begin{assumption} 
		
		\label{lmatrixl} 
		
		Fix constants $C_0, c_0>0$. Let $\L=[L_{ij}]_{i,j\in\llbracket N_1,N_2\rrbracket}$ and $\widetilde{\L}=[\widetilde L_{ij}]_{i,j\in\llbracket N_1,N_2\rrbracket}$ be real symmetric tridiagonal matrices.
		\noindent	In addition, assume there is an index set $\mathcal{D} \subseteq \llbracket N_1, N_2 \rrbracket$ defined as $\mathcal D:=\llbracket N_1, N_2 \rrbracket \setminus\mathcal I$ for some $\mathcal I=\llbracket m, \ell \rrbracket$ with $m < \ell$,
	and there is a real number $\delta \in [0, 1)$ satisfying
		\begin{flalign}
			\label{ll} 		
			\displaystyle\max_{(i,j) \notin \mathcal I^2}(|L_{ij}|+  |\widetilde{L}_{ij}|) \le 2 \log N; \qquad \displaystyle\max_{i,j \in \mathcal I} |L_{ij} - \widetilde{L}_{ij}| \le \delta.
		\end{flalign} 

		\noindent Let $\lambda_j$ denote an eigenvalue of $\L$, let $\zeta > 0$,
		and let $\varphi_j$ denote a $\zeta$-localization center for $\lambda_j$
		satisfying $\varphi_j\in\mathcal I$ and
		$\dist(\varphi_j,\mathcal D)>100$. Assume that a unit eigenvector
		$\mathbf u_j=(u_j(i))_{i=N_1}^{N_2}$ of $\L$ with eigenvalue $\lambda_j$
		satisfies
				\begin{equation}\label{eqn:lmatrix-boundary-decay}
		|u_j (k)| \le C_0 e^{-c_0 |k-\varphi_j|}
		\end{equation}
		for all $k$ such that
		$\max\{\dist(k,\mathcal D),\dist(k,\mathcal I)\}\leq1$.
		\end{assumption}

	We then have the following lemma, indicating that the eigenvalues of $\L$ (or $\widetilde{\L}$) with localization centers sufficiently distant from $\mathcal{D}$ are also nearly eigenvalues of $\widetilde{\L}$ (or $\L$, respectively). Its proof (which is similar to its counterpart \cite[Corollary 5.6]{Agg25a} addressing the case of thermal equilibrium) is in Section \ref{app:ftt}.

		\begin{lem}[Full to truncation comparison]
			\label{lem:full-to-truncation}
			Adopt Assumption \ref{lmatrixl}. There exists a constant
			$c=c(C_0,c_0)>0$ and an eigenvalue
			$\widetilde\lambda\in \eig \widetilde{\L}$ such that
			\[
				|\lambda_j-\widetilde\lambda|
				\le
				c^{-1} \zeta^{-2} \log N \bigl(\delta + e^{-c\dist(\varphi_j,\mathcal D)}\bigr),
			\]
			and $\varphi_j$ is an $N^{-1}\zeta$-localization center for $\widetilde\lambda$ with respect to $\widetilde{\L}$.
		\end{lem}

\subsection{Coupling to thermal equilibrium near the edge}
\label{subsec:coupl}

Throughout this subsection, we adopt Assumption \ref{ass:model}. In the lemma
below we consider a coupling between $\L(0)$ and a Lax matrix
$\L^{\rm th}(0)$ sampled from \emph{thermal equilibrium} as follows.

Consider the Toda evolution $\L^{\rm th}(t)$ through its Flaschka variables
$(\bm a^{\rm th}(t);\bm b^{\rm th}(t))$, with
$\bm a^{\rm th}(t)=(a_i^{\rm th}(t))_{i=N_1}^{N_2-1}$ and
$\bm b^{\rm th}(t)=(b_i^{\rm th}(t))_{i=N_1}^{N_2}$, whose initial data
$\{a_i^{\rm th}(0)\}_{i=N_1}^{N_2-1}\cup
\{b_i^{\rm th}(0)\}_{i=N_1}^{N_2}$ are sampled independently and according
to the laws
\[
a_i^{\rm th}(0)
\sim
C_{\mathfrak A} a^{2\mathfrak A^{-1}-1}e^{-a^2} da,
\qquad
b_i^{\rm th}(0)
\sim
\frac{1}{\sqrt{2\pi}}e^{-b^2/2} db.
\] 
The two densities above are chosen to match those of 
Assumption \ref{ass:model} in the region $|y| > (\log T)^5$. 

We couple $\L^{\rm th}(0)$ to $\L(0)$ by requiring
\[
a_i^{\rm th}(0)=a_i(0),
\qquad
b_i^{\rm th}(0)=b_i(0),
\qquad
|i| > T(\log T)^5;
\]
and sampling $\{a_i^{\rm th}(0)\}\cup\{b_i^{\rm th}(0)\}$ and
$\{a_i(0)\}\cup\{b_i(0)\}$ mutually independently and according to their
respective laws for $|i|\leq T(\log T)^5$. We use the convention that
$a_{N_2}(t)\equiv a_{N_2}^{\rm th}(t)\equiv0$.

The following lemma indicates that $\L^{\rm th}(t)$ is close to $\L(t)$ on
entries sufficiently far from the origin. We omit its proof, which follows
from using (the proof of) \cite[Lemma 4.5]{Agg25a} to compare
$(\bm a(t);\bm b(t))$ to $(\bm a^{\rm th}(t);\bm b^{\rm th}(t))$, with the
parameters $(A,T,K)$ there equal to
$(\log N,T\log N,T(\log N)^3)$ here, using Lemma
\ref{l0eigenvalues}.

\begin{lem}[{\cite[Lemma 4.5]{Agg25a}}]
\label{lem:edge-thermal-coupling}
Let 
$\mathcal T_-=
\llbracket
N_1,\,
-2T(\log T)^5
\rrbracket$,
and
$\mathcal T_+=
\llbracket
2T(\log T)^5,\,
N_2
\rrbracket$.
Then, with overwhelming
probability,
\[
\sup_{0\le t\le T\log N}
\;
\sup_{i\in\mathcal T_-\cup\mathcal T_+}
\left(
|a_i(t)-a_i^{\rm th}(t)|
+
|b_i(t)-b_i^{\rm th}(t)|
\right)
\le
e^{-T(\log N)^2}.
\]

\end{lem}

	\subsection{Concentration lemmas}
\label{subsec:conc}

In this section we state concentration estimates for the Toda Lax matrix at
time $t=0$. Here, we adopt Assumption \ref{ass:model}; recall Definition
\ref{def:quasi_intro}; and recall the notation from Section \ref{subsec:tc}.
In particular, we use the notation
\[
\Lambda_j:=\lambda_{\varphi_0^{-1}(j)},
\qquad
\bar q_j:=T^{-1}q_j(0),\qquad j \in \llbracket N_1,N_2\rrbracket.
\]

The next assumption imposes bounds on the test functions to be used in our
concentration estimates, which are stated in the two propositions below it.

\begin{assumption}[Test functions]
\label{ass:admissible-test}
Let
$h:\mathbb R^2\longrightarrow\mathbb R$
be bounded and locally Lipschitz. We interpret
\(\partial_qh\) as the almost everywhere derivative
 in the second variable.
Let \(A,B>1\). Assume that there is an interval
$I \subseteq[N_1/T,N_2/T]$ with $|I|\leq2(\log N)^{20}$ such that
\[
\operatorname{supp}h\subseteq\mathbb R\times I.
\]
Assume
\[
\|h\|_{\Lip,\lambda} :=
\sup_{q\in\mathbb R}
\sup_{\lambda\ne\lambda'}
\frac{|h(\lambda,q)-h(\lambda',q)|}
{|\lambda-\lambda'|}\le B,
\]
and
\[
\sup_{\lambda\in\mathbb R}
\int_{\mathbb R}|\partial_qh(\lambda,q)|\,dq\le A,
\qquad
\|h\|_\infty\le A.
\]
\end{assumption}

The proof of the proposition below is in Section \ref{app:time-zero-dos-proof}.

\begin{prop}
\label{prop:app-time-zero-dos}
Adopt Assumption \ref{ass:model}, and let \( h : \mathbb{R}^2 \rightarrow \mathbb{R}\)
satisfy Assumption \ref{ass:admissible-test}.
There is a constant $C>0$ such that, with overwhelming probability,
\begin{equation}\label{eqn:hiTconv}
\left| \sum_{i = N_1}^{N_2} h( \Lambda_i, i/T)
-
T \int_{-\infty}^{\infty}\int_{\mathbb R}
h(\lambda, y)\widetilde\varrho_0^{(T)}(\lambda;y)\,d\lambda\,dy \right| 
\leq (\log N)^C \left(  A T^{3/5} + B T^{2/5} \right)   .
\end{equation}
Moreover, if in addition $I$ from Assumption \ref{ass:admissible-test} satisfies $I \subset [-(\log N)^{20}, (\log N)^{20}]$,
\begin{equation}\label{eqn:hiTspace}
\left| \sum_{i = N_1}^{N_2} h( \Lambda_i, q_i(0)/T)
-
T \int_{-\infty}^{\infty}\int_{\mathbb R}
h(\lambda, y) m_0^{(T)}(\lambda;y)\,d\lambda\,dy \right| 
\leq  (\log N)^C \left(  A T^{3/5} +  B T^{2/5} \right)  .
\end{equation}
\end{prop}

The proof of the proposition below is in Section \ref{app:scat_conc}.

\begin{prop}[Scattering sums to integrals]
\label{prop:scat_conc}
Adopt Assumption \ref{ass:model}. Fix $s \in [0, T \log N]$. There is a constant \(C>0\), independent of \(\mathfrak A\), such that the following holds with overwhelming
probability. For \(j\in\llb 1,N\rrb\) and $\tau \geq0$, define
\[
Y_j(\tau):=x_\tau^{(T)}\bigl(\lambda_j;Q_j(0)T^{-1}\bigr).
\]
If $a \in \mathbb{R}$ and \(k\in\llb 1,N\rrb\) satisfy

	\[
	a\in\left[-(\log N)^{15},(\log N)^{15}\right], \qquad \varphi_0(k)\in
\llb -T(\log N)^{10},T(\log N)^{10}\rrb, 
	\]
		then writing \(i_k:=\varphi_0(k)\) and $\tau = s T^{-1}$, one has
\begin{equation}
\label{eq:time0-halfspace-conc}
\begin{aligned}
\Bigg|
\Bigg[
\sum_{\substack{i:\,\bar q_i<\bar q_{i_k}\\ i\ne i_k}}\log|\lambda_k-\Lambda_i|
-
\sum_{\substack{j:\,Y_j(\tau)<a\\ j\ne k}}\log|\lambda_k-\lambda_j|
\Bigg]
-T\mathcal I_{\tau,k}^{(T)}(a)
\Bigg|
\le
CT^{3/4-1/100}\log(|\lambda_k|+2).
\end{aligned}
\end{equation}
Here
\[
\mathcal I_{\tau,k}^{(T)}(a)
:=
\int_{\mathbb R}\log|\lambda_k-\lambda|\,
\left(N_0^{(T)}(\lambda,\bar q_{i_k})-N_\tau^{(T)}(\lambda,a)\right)\,d\lambda,
\]
where \(N_\tau^{(T)}(\lambda,x)\) denotes the fixed point from Lemma
\ref{lem:fixed_pointout} and Section \ref{subsec:tc}.
\end{prop}

		\section{Density functions solving the GHD equation}
		\label{sec:GHD_sol}
The goal of this section is to provide the missing proofs of results in
Sections \ref{subsec:lin} and \ref{subsec:ghd-uniqueness}. Throughout this
section, we adopt Assumption \ref{ass:sec2}, recalling the notation there and
from Definitions \ref{def:Tdf} and \ref{def:initial-ghd-coordinates}.

	\subsection{Construction of the fixed point}
	\label{subsec:GHD_alt}
	
	Recall from the discussion around \eqref{eqn:NtFP}, if
	\begin{equation}\label{eqn:N0X}
	N_0^X(\Lambda, X) \coloneqq \int_{0}^X\varpi_0(\Lambda;Y)\,dY,
	\end{equation}
	then the fixed point equation \eqref{eqn:NtFP} is
	\begin{equation}\label{eqn:Hf-primitive}
	N_t(\Lambda, x)
	=
	N_0^X\bigl(\Lambda, x+(\T N_t(\cdot,x))(\Lambda)-t\Lambda \bigr).
	\end{equation}
	For $x\in\mathbb R$ and $t\ge0$, define
	\begin{equation}\label{eqn:hop}
	\Hop_{x,t}[f](\Lambda)
	:=
	N_0^X\bigl(\Lambda, x+(\T f)(\Lambda)-t\Lambda\bigr),
	\end{equation}
	and recall the space $\mathcal{B}$ defined in \eqref{eqn:Bsp}.

	\begin{lem}[Contraction estimate]\label{lem:H-contraction}
	For every $x\in\mathbb R$ and $t\ge0$, the map $\Hop_{x,t}$ sends
	$\mathcal B$ to itself. Moreover, there exists $a\in(0,1)$ such that
	\begin{equation}\label{eqn:H-contraction}
	\|\Hop_{x,t}[f_1]-\Hop_{x,t}[f_2]\|_{\mathfrak c}
	\le a\|f_1-f_2\|_{\mathfrak c}
	\end{equation}
	for all $f_1,f_2\in\mathcal B$, uniformly in $x$ and $t$.
	\end{lem}

	\begin{proof}
	For $f\in \mathcal B$ we have the bound
	\begin{equation}\label{eqn:lpbd}
	|(\T f)(\Lambda)|
	\le
	C\log(|\Lambda|+2)\|f\|_{\mathfrak c},
	\end{equation}
	where $C=C(\mathfrak c)$, because $\int_{-\infty}^{\infty}|\log|\Lambda-\lambda||\,e^{-\mathfrak c\lambda^2}\,d\lambda
	\le C\log(|\Lambda|+2)$.
	Hence, using \eqref{eqn:vpi0bd},
	\begin{align*}
	|\Hop_{x,t}[f_1](\Lambda)-\Hop_{x,t}[f_2](\Lambda)|
	&\le
	2C_0\mathfrak K\mathfrak A^{-1}e^{-c_0\Lambda^2}|(\T(f_1-f_2))(\Lambda)|\\
	&\le
	C\mathfrak A^{-1}e^{-c_0\Lambda^2}\log(|\Lambda|+2)\|f_1-f_2\|_{\mathfrak c}.
	\end{align*}
	Multiplying by $e^{\mathfrak c\Lambda^2}$ and using $\mathfrak c<c_0$ gives
	\begin{equation}\label{eqn:Hf1f2}
	\|\Hop_{x,t}[f_1]-\Hop_{x,t}[f_2]\|_{\mathfrak c}
	\le
	C\mathfrak A^{-1}\|f_1-f_2\|_{\mathfrak c}.
	\end{equation}
	After choosing $\mathfrak A$ large, the right hand side is at most $a\|f_1-f_2\|_{\mathfrak c}$ for some $a<1$.

	The same estimate shows that $\Hop_{x,t}$ maps $\mathcal B$ to itself:
	\[
	|\Hop_{x,t}[f](\Lambda)|
	\le
	2C_0\mathfrak K\mathfrak A^{-1}e^{-c_0\Lambda^2}
	\bigl(|x| + t |\Lambda|+C\log(|\Lambda|+2)\|f\|_{\mathfrak c}\bigr).
	\]
	This proves the lemma.
	\end{proof}

%
	
	By Lemma \ref{lem:H-contraction} and Banach's fixed point theorem,
	we obtain a unique fixed point of the functional $\Hop_{x,t}$ from \eqref{eqn:hop}
	for each $t \geq 0, x \in \mathbb{R}$. 
	
	\begin{lem}\label{lem:ghd_lip}
	Adopt Assumption \ref{ass:sec2} and let $N_t(\lambda,x)$ be the fixed point described above.
	 There exists a constant $C>0$ such that
	 \eqref{eqn:Nt-joint-lipschitz} holds
	uniformly in $(t,s,x,y,\lambda)\in
	\mathbb R_{\ge0}^2\times\mathbb R^2\times\mathbb R$.
	\end{lem}
	
	\begin{proof}
	Let $N_{x,t}(\Lambda):=N_t(\Lambda,x)$. Since $N_{x,t}$ and $N_{y,s}$ are fixed points of $\Hop_{x,t}$ and $\Hop_{y,s}$, respectively,
	\[
	N_{x,t}-N_{y,s}=\Hop_{x,t}[N_{x,t}]-\Hop_{y,s}[N_{y,s}].
	\]
	Add and subtract \(\Hop_{x,t}[N_{y,s}]\). By the triangle inequality,
	\[
	\begin{aligned}
	\|N_{x,t}-N_{y,s}\|_{\mathfrak c}
	&\le
	\|\Hop_{x,t}[N_{x,t}]-\Hop_{x,t}[N_{y,s}]\|_{\mathfrak c}
	+
	\|\Hop_{x,t}[N_{y,s}]-\Hop_{y,s}[N_{y,s}]\|_{\mathfrak c}.
	\end{aligned}
	\]
	Let $a\in(0,1)$ be the contraction constant from Lemma
	\ref{lem:H-contraction}.
	The first term is controlled by Lemma \ref{lem:H-contraction}:
	\[
	\|\Hop_{x,t}[N_{x,t}]-\Hop_{x,t}[N_{y,s}]\|_{\mathfrak c}
	\le
	a\|N_{x,t}-N_{y,s}\|_{\mathfrak c}.
	\]
	For the second term, recall that for each \(\Lambda\), 
	defining $U=x+(\T N_{y,s})(\Lambda)-t\Lambda$ and
	$U'=y+(\T N_{y,s})(\Lambda)-s\Lambda$, we have
	\[
		\Hop_{x,t}[N_{y,s}](\Lambda)
		=
		\int_{0}^{U}\varpi_0(\Lambda;X)\,dX, \quad \Hop_{y,s}[N_{y,s}](\Lambda)
		=
		\int_{0}^{U'}\varpi_0(\Lambda;X)\,dX.
	\]
	So the two upper integration limits differ by
	$U-U'=(x-y)-(t-s)\Lambda$.
	Therefore, using \eqref{eqn:vpi0bd},
	\[
	\begin{aligned}
	\bigl|\Hop_{x,t}[N_{y,s}](\Lambda)-\Hop_{y,s}[N_{y,s}](\Lambda)\bigr|
	&\le
	\sup_{X\in\mathbb R}\varpi_0(\Lambda;X)
	\bigl(|x-y|+|\Lambda|\,|t-s|\bigr) \\
	&\le
	2C_0\mathfrak K\mathfrak A^{-1}e^{-c_0\Lambda^2}
	\bigl(|x-y|+|\Lambda|\,|t-s|\bigr).
	\end{aligned}
	\]
Multiplying by \(e^{\mathfrak c\Lambda^2}\) and using \(\mathfrak c<c_0\), we obtain
	\[
	\|\Hop_{x,t}[N_{y,s}]-\Hop_{y,s}[N_{y,s}]\|_{\mathfrak c}
	\le
	C\mathfrak A^{-1}\bigl(|x-y|+|t-s|\bigr).
	\]
	Hence
	\[
	\|N_{x,t}-N_{y,s}\|_{\mathfrak c}
	\le
	a\|N_{x,t}-N_{y,s}\|_{\mathfrak c}
	+
	C\mathfrak A^{-1}\bigl(|x-y|+|t-s|\bigr),
	\]
	and therefore
	\[
	\|N_{x,t}-N_{y,s}\|_{\mathfrak c}
	\le
	\frac{C}{1-a}\mathfrak A^{-1}\bigl(|x-y|+|t-s|\bigr).
	\]
	After absorbing \((1-a)^{-1}\) into the constant, this proves
	\[
	\|N_{x,t}-N_{y,s}\|_{\mathfrak c}
	\le
	C\mathfrak A^{-1}\bigl(|x-y|+|t-s|\bigr),
	\]
	which is equivalent to the claimed pointwise bound.
	\end{proof}

	We may now provide a proof of Lemma \ref{lem:fixed_pointout}.
	\begin{proof}[Proof of Lemma \ref{lem:fixed_pointout}]
	Lemma \ref{lem:H-contraction} and Banach's fixed point theorem give a
	unique solution of \eqref{eqn:Hf-primitive} in $\mathcal B$. The joint
	Lipschitz assertion follows from Lemma \ref{lem:ghd_lip}.
	\end{proof}

		\subsection{Properties of fixed point solution} 
	\label{subsec:GHD_2}
	First, define $X_t$ and $J_t$ as in \eqref{eqn:Xt} and \eqref{eqn:s0dr},
	respectively, with the $\widetilde N_t$ there equal to the fixed point $N_t$
	provided by Lemma \ref{lem:fixed_pointout}.
	
	\begin{lem}\label{lem:fixed_point_derivative0}
	Let $N_t(\lambda,x)$ be the fixed point from Lemma
	\ref{lem:fixed_pointout}.
	Denote the almost everywhere $x$-derivative of $N_t(\lambda,x)$ by $m_t(\lambda,x)$.
	There exists a constant $C> 0$, such that for every $t\ge 0$ and almost every $x\in \mathbb{R}$,
	\begin{equation}\label{eqn:avrbd}
	|m_t(\lambda,x)|
	\le
	C\mathfrak A^{-1}e^{-\frac{c_0}{2}\lambda^2},
	\end{equation}
	and
	\begin{equation}\label{eqn:Jtbd}
	|J_t(\lambda,x)-1|
	\le
	C\mathfrak A^{-1}\log(|\lambda|+2).
	\end{equation}
	Moreover,
	\begin{equation}\label{eqn:mt}
	m_t(\Lambda,x)
	=
	\varpi_t(\Lambda;X_t(\Lambda,x))
	J_t(\Lambda,x).
	\end{equation}
	\end{lem}
	\begin{proof}
	The $t=s$ case of \eqref{eqn:Nt-joint-lipschitz} gives
	\eqref{eqn:avrbd}. In particular, $x\mapsto N_t(\Lambda,x)$ is
	Lipschitz, hence absolutely continuous, for each $\Lambda$. Thus
	$x\mapsto\mathbf T N_t(\Lambda,x)$ is also absolutely continuous, with
	almost everywhere derivative $\mathbf Tm_t(\Lambda,x)$. Moreover, by \eqref{eqn:avrbd},
	\[
	|(\T m_t(\cdot,x))(\Lambda)|
	\le
	C\mathfrak A^{-1}\log(|\Lambda|+2).
	\]
	The definition \eqref{eqn:s0dr} of $J_t$ now gives
	\eqref{eqn:Jtbd}. Finally, differentiating \eqref{eqn:NtFP} gives \eqref{eqn:mt}.
	\end{proof}

 The next lemma shows that the fixed point $ x \mapsto N_t(\Lambda, x)$ is nondecreasing. Thus, in the case that $N_t$ is absolutely continuous
	(this indeed holds in our setting by Lemma \ref{lem:fixed_pointout}), the lemma shows that its almost everywhere derivative $m_t(\Lambda, x)$ is nonnegative.

	\begin{lem}\label{lem:mt_nonnegative}
	Adopt Assumption \ref{ass:sec2}, and let $N_t$ denote the fixed point from Lemma \ref{lem:fixed_pointout}. For every \(t\ge 0\) and \(\Lambda\in\mathbb R\), the map \(x\mapsto N_t(\Lambda, x)\) is nondecreasing. 
	\end{lem}

	\begin{proof}
	Fix \(h>0\) and define the forward difference quotient
	\[
	D_h(\Lambda,x)
	:=
	\frac{N_t(\Lambda, x+h)-N_t(\Lambda, x)}{h}.
	\]
	Also write
	$Y_x(\Lambda)
	:=
	x+(\T N_t(\cdot, x))(\Lambda)-t\Lambda$, so 
	the fixed point equation \eqref{eqn:NtFP} says
	$N_t(\Lambda, x)
	=
	\int_{0}^{Y_x(\Lambda)}
	\varpi_0(\Lambda,X)\,dX$.
	Therefore,
	\[
	\begin{aligned}
	N_t(\Lambda,x+h)-N_t(\Lambda,x)
	&=
	\int_{Y_x(\Lambda)}^{Y_{x+h}(\Lambda)}
	\varpi_0(\Lambda,X)\,dX  \\
	&=
	\bigl(Y_{x+h}(\Lambda)-Y_x(\Lambda)\bigr)
	\int_0^1
	\varpi_0\bigl(\Lambda,
	Y_x(\Lambda)+s(Y_{x+h}(\Lambda)-Y_x(\Lambda))
	\bigr)\,ds .
	\end{aligned}
	\]
	Since
	\[
	Y_{x+h}(\Lambda)-Y_x(\Lambda)
	=
	h+\T\bigl(N_t(\cdot,x+h)-N_t(\cdot,x)\bigr)(\Lambda)
	=
	h\bigl(1+(\T D_h(\cdot,x))(\Lambda)\bigr),
	\]
	we obtain
	\begin{equation}\label{eqn:Dh-fixed-point}
	D_h(\Lambda,x)
	=
	\eta_h(\Lambda,x)
	\left(
	1+(\T D_h(\cdot,x))(\Lambda)
	\right),
	\end{equation}
	where $\eta_h(\Lambda,x)
	:=
	\int_0^1
	\varpi_0\bigl(\Lambda,
	Y_x(\Lambda)+s(Y_{x+h}(\Lambda)-Y_x(\Lambda))
	\bigr)\,ds $.
	Note that by \eqref{eqn:vpi0bd},
	\[
	0\leq \eta_h(\Lambda,x)
	\leq
	2C_0\mathfrak K\mathfrak A^{-1}e^{-c_0\Lambda^2}.
	\]

	We claim that \(D_h(\Lambda,x)\geq 0\) for every \(\Lambda\). To see this, consider the map
	\[
	\mathcal L_h[g](\Lambda)
	:=
	\eta_h(\Lambda,x)\left(1+(\T g)(\Lambda)\right).
	\]
	The same estimate as the one leading to \eqref{eqn:Hf1f2} in the proof of Lemma~\ref{lem:H-contraction} shows that, if $g_1, g_2 \in \mathcal{B}$, then
	\[
	\|\mathcal L_h[g_1]-\mathcal L_h[g_2]\|_{\mathfrak c}
	\leq
	C\mathfrak A^{-1}\|g_1-g_2\|_{\mathfrak c},
	\]
	so \(\mathcal L_h\) is a contraction on \(\mathcal B\). Moreover, \(D_h\) is a fixed point of \(\mathcal L_h\), by \eqref{eqn:Dh-fixed-point}.

	Let \(g_0=0\), and define \(g_{\ell+1}=\mathcal L_h[g_\ell]\). We show by induction that
	\[
	g_\ell(\Lambda)\geq 0
	\qquad\text{and}\qquad
	\|g_\ell\|_{\mathfrak c}\leq C_1\mathfrak A^{-1}
	\]
	for every \(\ell\), where \(C_1>0\) is chosen sufficiently large as a function of $\mathfrak{c}, c_0, C_0, \mathfrak K$. The claim is trivial for \(\ell=0\). Suppose it holds for \(g_\ell\). Then
	\[
	\|g_{\ell+1}\|_{\mathfrak c}
	\leq
	2C_0\mathfrak K \mathfrak A^{-1}
	\left(1+C \|g_\ell\|_{\mathfrak{c}}\right)
	\leq
	C_1 \mathfrak A^{-1},
	\]
	where we used the estimate \eqref{eqn:lpbd} from the proof of Lemma~\ref{lem:H-contraction} and took \(C_1\) and \(\mathfrak A\) large enough that \(2C_0\mathfrak K (1 + C C_1 \mathfrak{A}^{-1}) < C_1\).

	It remains to check positivity. Since \(g_\ell\geq 0\), the only negative contribution to \((\T g_\ell)(\Lambda)\) comes from \(|\Lambda-\lambda|<1\). Thus
	\[
	\begin{aligned}
	(\T g_\ell)(\Lambda)
	&=
	2\int_{-\infty}^{\infty}\log|\Lambda-\lambda|\,g_\ell(\lambda)\,d\lambda \\
	&\geq
	-2\int_{|\Lambda-\lambda|<1}
	|\log|\Lambda-\lambda||\,g_\ell(\lambda)\,d\lambda \\
	&\geq
	-C\|g_\ell\|_{\mathfrak c}
	\geq
	-C C_1 \mathfrak A^{-1}.
	\end{aligned}
	\]
	Taking \(\mathfrak A\) sufficiently large, this gives
	$1+(\T g_\ell)(\Lambda)\geq \frac12$.
	Since \(\eta_h\geq 0\), we conclude
	$g_{\ell+1}(\Lambda)
	=
	\eta_h(\Lambda,x)\left(1+(\T g_\ell)(\Lambda)\right)
	\geq 0$.
	This proves the induction step.

	Since \(\mathcal L_h\) is a contraction, \(g_\ell\to D_h(\cdot,x)\) in \(\mathcal B\) as $\ell\to\infty$. Hence $D_h(\Lambda,x)\geq 0$
	for every \(\Lambda\) and every \(h>0\). Therefore \(x\mapsto N_t(\Lambda,x)\) is nondecreasing.
		\end{proof}

	We now prove Lemma \ref{lem:Jlow}.
	\begin{proof}[Proof of Lemma \ref{lem:Jlow}]
	The bound
	$0 \leq \widetilde N_t(\lambda, y)-\widetilde N_t(\lambda, x)
	\leq C \mathfrak{A}^{-1} e^{-\mathfrak{c} \lambda^2}(y-x)$ for $y>x$
	gives
	$\widetilde X_t(\Lambda,y)-\widetilde X_t(\Lambda,x)
	\geq \left(1-C\mathfrak A^{-1}\right)(y-x)
	\geq\tfrac12 (y-x)$ for $\mathfrak A$ sufficiently large and $y>x$.

	The bounds above prove the second part of the statement. The first follows
	from differentiating in $x$ to obtain the almost everywhere $(x,t)$ bound
	$\partial_x\widetilde X_t(\Lambda,x)
	=1+(\T \widetilde m_t(\cdot,x))(\Lambda)
	=\widetilde J_t(\Lambda,x)\geq\frac12$.
	Observe that we also have
	\begin{equation}\label{eqn:J-fixed-point-lower}
	J_t(\Lambda,x)\geq\frac12,
	\end{equation}
	by the same reasoning, using Lemmas \ref{lem:ghd_lip} and
	\ref{lem:mt_nonnegative} to verify the bound
	$0\leq N_t(\lambda,y)-N_t(\lambda,x)
	\leq C\mathfrak A^{-1}e^{-c_0\lambda^2/2}(y-x)$ for $y>x$.
	\end{proof}

For every $t\geq0$ and almost every $x\in\mathbb R$, set
$m_{t,x}(\lambda)\coloneqq m_t(\lambda,x)$ and
$J_{t,x}(\Lambda)\coloneqq J_t(\Lambda,x)$. For every $t\geq0$ and
$x\in\mathbb R$, set
\[
\mathfrak n_{t,x}(\lambda)\coloneqq
\varpi_t(\lambda;X_t(\lambda,x)).
\]
By \eqref{eqn:mt}, for every $t\geq0$
and almost every $x\in\mathbb R$,
\begin{equation}\label{eqn:mtxntx}
m_{t,x}(\lambda)=\mathfrak n_{t,x}(\lambda)J_{t,x}(\lambda).
\end{equation}
We also define 
\begin{equation}\label{eqn:HtxL2}
\mathcal H
:=
L^2(\mathbb R, e^{-\frac{c_0}{2} \lambda^2} \,d\lambda),
\end{equation}
where $c_0 > 0$ is the constant from Assumption \ref{ass:sec2}.
In Lemma \ref{lem:dressing_invertibility} below, 
	we will show that the operator $1-\mathbf{T} \mathfrak n_{t,x}$
	is bounded and invertible on the Hilbert space $\mathcal H$.
	Its inverse is known as the \emph{dressing operator},
	and its application to a function $f$ at space time point $(x,t)$ is
	denoted
	\[
	f^{\dr,t,x}:=(1-\mathbf{T}\mathfrak n_{t,x})^{-1}f.
	\]
	Finally, define the single variable functions
	\[
	\varsigma_0(\lambda)\coloneqq1\in\mathcal H,
	\qquad
	\varsigma_1(\lambda)\coloneqq\lambda\in\mathcal H.
	\]

\begin{lem}[Invertibility of the local dressing operator]\label{lem:dressing_invertibility}
Adopt Assumption \ref{ass:sec2}. For every \(t\geq0\) and
\(x\in\mathbb R\), define the operator
$\mathcal K_{t,x}:=\T \mathfrak n_{t,x}$
on \(\mathcal H\) by
\[
(\mathcal K_{t,x}f)(\Lambda)
=
2\int_{-\infty}^{\infty}
\log|\Lambda-\lambda|\,
\mathfrak n_{t,x}(\lambda)f(\lambda)\,d\lambda.
\]
If \(\mathfrak A\) is sufficiently large, then the operator norm
$\|\mathcal K_{t,x}\|_{\mathcal H\to\mathcal H}
\leq C\mathfrak A^{-1}<1$.
Consequently,
$1-\T\mathfrak n_{t,x}
=
1-\mathcal K_{t,x}$
is invertible on \(\mathcal H\). In addition, for every $t\geq0$ and almost
every $x\in\mathbb R$,
\begin{equation}\label{eqn:J1dr}
J_{t,x}=\varsigma_0^{\dr, t, x}.
\end{equation}
\end{lem}

\begin{proof}
	By the definition of \(\mathfrak n_{t,x}\) and \eqref{eqn:vpi0bd}, we have
	$0\leq \mathfrak n_{t,x}(\lambda)
	\leq
	2C_0\mathfrak K\mathfrak A^{-1}e^{-c_0\lambda^2}$.
	Set $w(\lambda):=e^{-c_0\lambda^2/2}$, and let
	\(U:\mathcal H\to L^2(\mathbb R,d\lambda)\) be the unitary map
	\[
	(Uf)(\lambda):=w(\lambda)^{1/2}f(\lambda).
	\]
	Then \(U\mathcal K_{t,x}U^{-1}\) is the integral operator with kernel
	\[
	L_{t,x}(\Lambda,\lambda)
	:=
	2\log|\Lambda-\lambda|\,
	w(\Lambda)^{1/2}\mathfrak n_{t,x}(\lambda)w(\lambda)^{-1/2}.
	\]
	The preceding Gaussian bound on $\mathfrak n_{t,x}$ gives
	\[
	|L_{t,x}(\Lambda,\lambda)|^2
	\le
	C\mathfrak A^{-2}
	(\log|\Lambda-\lambda|)^2\,
	e^{-c_0\Lambda^2/2}e^{-3c_0\lambda^2/2}.
	\]
	Since the logarithmic singularity is locally square integrable, we get
	\[
	\iint_{\mathbb R^2}|L_{t,x}(\Lambda,\lambda)|^2\,d\Lambda\,d\lambda
	\le
	C\mathfrak A^{-2}.
	\]
	Therefore $\|U\mathcal K_{t,x}U^{-1}\|_{\mathrm{HS}}
	\le
	C\mathfrak A^{-1}$,
	and hence
	$\|\mathcal K_{t,x}\|_{\mathcal H\to\mathcal H}
	\le C\mathfrak A^{-1}$.

	Taking \(\mathfrak A\) sufficiently large gives
	$\|\mathcal K_{t,x}\|_{\mathcal H\to\mathcal H}<1$.
	Thus \(1-\mathcal K_{t,x}\) is invertible with inverse given by 
	 the convergent series $\sum_{n=0}^{\infty} \mathcal K_{t,x}^n$.

	Finally, fix $t\geq0$ and an $x$ for which \eqref{eqn:mtxntx} holds.
	The bound \eqref{eqn:Jtbd} shows that $J_{t,x}\in\mathcal H$. Moreover,
	using \eqref{eqn:mtxntx},
	\[
	J_{t,x}
	=
	1+\T m_{t,x}
	=
	1+\T(\mathfrak n_{t,x}J_{t,x}).
	\]
	Equivalently,
	$(1-\T\mathfrak n_{t,x})J_{t,x}=1$.
	So $J_{t,x}=\varsigma_0^{\dr, t, x}$,
	which is \eqref{eqn:J1dr}.
\end{proof}

\begin{lem}\label{lem:effective-velocity-dressing}
	For every $t\geq0$ and almost every $x\in\mathbb R$, the effective velocity
	from Lemma \ref{lem:int_ghd} is given by
	\begin{equation}\label{eqn:vedf}
	\ve_t(\Lambda;x)=
	\frac{\varsigma_{1}^{\dr,t,x}(\Lambda)}{\varsigma_{0}^{\dr,t,x}(\Lambda)}.
	\end{equation}
	\end{lem}

	We next prove Lemmas \ref{lem:int_ghd} and
	\ref{lem:effective-velocity-dressing} together.
	
	\begin{proof}[Proofs of Lemmas \ref{lem:int_ghd} and
	\ref{lem:effective-velocity-dressing}]
	The statement $m_t(\lambda,x)\geq0$ from Lemma \ref{lem:int_ghd} follows from Lemma
	\ref{lem:mt_nonnegative}. Fix $t\geq0$ and an $x$ for which
	\eqref{eqn:mtxntx} and \eqref{eqn:J1dr} hold.
	
	Since $J_{t,x}\geq\frac12$ by \eqref{eqn:J-fixed-point-lower} and
	$\varsigma_1^{\dr,t,x}\in\mathcal H$, the right side of
	\eqref{eqn:vedf} defines an element $v(\cdot,x)\in\mathcal H$. It satisfies
	\begin{equation}\label{eqn:effective-velocity-linear-identity}
	v(\Lambda,x)\bigl(1+(\T m_t(\cdot,x))(\Lambda)\bigr)
	=
	\Lambda+(\T(v(\cdot,x)m_t(\cdot,x)))(\Lambda).
	\end{equation}
	Indeed, multiplying the definition of $v$ by
	$\varsigma_0^{\dr,t,x}=J_{t,x}$ and using
	$(1-\T\mathfrak n_{t,x})\varsigma_1^{\dr,t,x}=\varsigma_1$ together
	with $m_{t,x}=\mathfrak n_{t,x}J_{t,x}$ gives
	\eqref{eqn:effective-velocity-linear-identity}, which is equivalent to
	the integral equation in \eqref{eqn:cra}.
	Conversely, suppose that $\widetilde{\ve_t}(\cdot;x)\in\mathcal H$ satisfies
	\eqref{eqn:cra}. The Gaussian bound \eqref{eqn:avrbd} implies that
	$\T(m_{t,x}\widetilde{\ve_t}(\cdot;x))\in\mathcal H$. Reversing the
	preceding calculation, using \eqref{eqn:mtxntx}, gives
	\[
	(1-\T\mathfrak n_{t,x})
	\bigl(J_{t,x}\widetilde{\ve_t}(\cdot;x)\bigr)
	=
	\varsigma_1.
	\]
	The invertibility just proved gives
	$J_{t,x}\widetilde{\ve_t}(\cdot;x)=\varsigma_1^{\dr,t,x}$, and hence
	$\widetilde{\ve_t}(\cdot;x)=v(\cdot,x)$. This proves the existence and
	uniqueness stated in Lemma \ref{lem:int_ghd}. Thus its effective velocity
	$\ve_t(\cdot;x)$ is $v(\cdot,x)$, which also proves Lemma
	\ref{lem:effective-velocity-dressing}.
	
	Differentiating \eqref{eqn:Hf-primitive} with respect to $t$ and using
	\eqref{eqn:N0} and \eqref{eqn:Xt}, we obtain
	for each $\Lambda \in \mathbb{R}$ and almost every $(x, t)$,
	\begin{align} \notag
	\partial_t N_t(\Lambda,x)
	&=
	\varpi_0(\Lambda; x+(\T N_t(\cdot,x))(\Lambda)-t\Lambda)\,\partial_t \left( x+(\T N_t(\cdot,x))(\Lambda)-t\Lambda \right) \\
	&=
	\varpi_t(\Lambda; X_t(\Lambda,x))\,\partial_t \left( x+(\T N_t(\cdot,x))(\Lambda)-t\Lambda \right) \notag \\
	&=
	\varpi_t(\Lambda; X_t(\Lambda,x))\, \left( (\T \partial_t N_t(\cdot,x))(\Lambda)-\Lambda \right) .  \label{eqn:chainNt0}
	\end{align}
	Set
	$r_{t,x}:=\T\partial_tN_t(\cdot,x)-\varsigma_1$.
	The Gaussian bound \eqref{eqn:Nt-joint-lipschitz} on $\partial_tN_t$ shows that
	$r_{t,x}\in\mathcal H$, and \eqref{eqn:chainNt0} gives
	$\partial_tN_t(\cdot,x)=\mathfrak n_{t,x}r_{t,x}$. Therefore
	$(1-\T\mathfrak n_{t,x})r_{t,x}=-\varsigma_1$.
	Applying Lemma \ref{lem:dressing_invertibility}, we obtain
	\begin{equation}\label{eqn:intghdfin}
	\begin{aligned}
	\partial_t N_t(\lambda,x) &= -  \mathfrak n_{t,x}(\lambda)  (1- \T \mathfrak n_{t,x})^{-1} (\varsigma_1)(\lambda) \\
	&= -  \mathfrak n_{t,x}(\lambda) \varsigma_1^{\dr, t,x}(\lambda) \\
	&= - \ve_t(\lambda; x) \mathfrak n_{t,x}(\lambda) J_t(\lambda, x) \\
	&= - \ve_t(\lambda;x) m_t(\lambda, x).
	\end{aligned}
		\end{equation}
In the third equality, we have used \eqref{eqn:J1dr} and the expression \eqref{eqn:vedf} for $\ve_t$, and in the fourth we have used \eqref{eqn:mtxntx}.	
	This completes the proof of Lemma \ref{lem:int_ghd}.
	\end{proof} 
	
	Next, we have a lemma which bounds $ \ve_{t}$ and its Lipschitz constant as a function 
	of $\lambda$. The proof, which is similar to the proof of
\cite[Corollary 3.8]{Agg25}, is given in Appendix \ref{app:vebd}.

\begin{lem}\label{lem:vebd}
Adopt Assumption \ref{ass:sec2}. Then,
$\lambda \mapsto \ve_{t}(\lambda; x)$ is Lipschitz. There is a constant $C>0$
such that for all
$(t, \lambda) \in \mathbb{R}_{\geq 0} \times \mathbb R$, and almost every $x \in \mathbb{R}$,
\begin{equation}\label{eqn:vebd}
|\ve_t(\lambda;x)|\le C(1+|\lambda|);
\qquad
|\partial_\lambda\ve_t(\lambda;x)|\le C(1+|\lambda|).
\end{equation}
\end{lem}

\subsection{On solutions of GHD equations and characteristic flows}
\label{subsec:onsolghd}

	In the proof below, we adopt Assumption \ref{ass:sec2} and the notation
	from Definition \ref{def:anyN}.

\begin{proof}[Proof of Lemma \ref{xt2}]
Fix $\Lambda\in\mathbb R$ and set
$K_t(\Lambda,x):=\widetilde X_t(\Lambda,x)-t\Lambda
=x+(\T \widetilde N_t(\cdot,x))(\Lambda)-t\Lambda$.
The weighted Lipschitz estimate in Proposition \ref{prop:GHDint}
gives
$|\partial_x\widetilde N_t(\lambda,x)|+|\partial_t\widetilde N_t(\lambda,x)|
\leq C\mathfrak A^{-1}e^{-\mathfrak c\lambda^2}$ almost everywhere.
Thus, with $\widetilde m_t(\lambda,x):=\partial_x\widetilde N_t(\lambda,x)$, we have
$0\leq \widetilde m_t(\lambda,x)\leq
C\mathfrak A^{-1}e^{-\mathfrak c\lambda^2}$ almost everywhere.
This bound permits commuting a derivative through $\T$, and hence
\begin{equation}\label{eqn:Kx-general-solution}
\partial_xK_t(\Lambda,x)
=1+(\T \widetilde m_t(\cdot,x))(\Lambda) = \widetilde J_t(\Lambda, x)
\end{equation}
and, by \eqref{eqn:GHDpot_unique},
\begin{equation}\label{eqn:Kt-general-solution}
\begin{aligned}
\partial_tK_t(\Lambda,x)
&=(\T(\partial_t\widetilde N_t(\cdot,x)))(\Lambda)-\Lambda\\
&=-(\T(\widetilde v_t^{\eff}(\cdot;x)\widetilde m_t(\cdot,x)))(\Lambda)-\Lambda
\end{aligned}
\end{equation}
for each $\Lambda$ and almost every $(x,t)$. Rewriting
\eqref{eqn:cra} using \eqref{eqn:Kx-general-solution} and
\eqref{eqn:Kt-general-solution} yields
\begin{equation}\label{eqn:Kttransport-general}
\partial_tK_t(\Lambda,x)
=-\widetilde v_t^{\eff}(\Lambda;x)\,\partial_xK_t(\Lambda,x)
\end{equation}
for each $\Lambda$ and almost every $(x,t)$.

Let $(t,q)\mapsto y_t(\Lambda;q)$ be a characteristic flow, as in
Definition \ref{def:xtd}, with $\ve_t$ there replaced by
$\widetilde v_t^{\eff}$.
The chain
rule and \eqref{eqn:Kttransport-general} give
\begin{equation}\label{eqn:X-flow-chain-rule}
\begin{aligned}
\frac{d}{d t}\!\left[\widetilde X_t\bigl(\Lambda,y_t(\Lambda;q)\bigr)\right]
&=\Lambda+\partial_xK_t\bigl(\Lambda,y_t(\Lambda;q)\bigr)\\
&\quad\times\left[
\frac{d}{dt} y_t(\Lambda;q)
-\widetilde v_t^{\eff}\bigl(\Lambda;y_t(\Lambda;q)\bigr)
\right]
\end{aligned}
\end{equation}
for almost every $(t,q)$. Thus \eqref{eqn:X-flow-chain-rule} implies
\eqref{eqn:Xtlin}.


Finally, define
$x_t(\Lambda;q)$ as in the statement.
Then $K_t(\Lambda;x_t(\Lambda;q))=K_0(\Lambda,q)$, or equivalently
$\widetilde X_t(\Lambda,x_t(\Lambda;q))=X_0(\Lambda,q)+t\Lambda$.
Equation \eqref{eqn:Kx-general-solution} and Lemma \ref{lem:Jlow} give a
lower bound on $\partial_x\widetilde X_t(\Lambda,x)$, while
\eqref{eqn:Kx-general-solution} and the bound $m_t(\lambda,x)\leq
C\mathfrak A^{-1}e^{-\mathfrak c\lambda^2}$ established at the beginning
of the proof give
$\partial_x\widetilde X_t(\Lambda,x)\leq C\log(|\Lambda|+2)$. These bounds show that
$q\mapsto x_t(\Lambda;q)$ is locally bi-Lipschitz
in $q$. The same lower bound, together with the fact that
$\widetilde N_t$ is Lipschitz in time by \eqref{eqn:nlip}, shows that
$t \mapsto x_t(\Lambda;q)$ is Lipschitz. The display
\eqref{eqn:X-flow-chain-rule} above therefore gives
$\frac{d}{d t} x_t(\Lambda;q)=
\widetilde v_t^{\eff}(\Lambda;x_t(\Lambda;q))$ for almost every $(t,q)$.
\end{proof}


%
%

\begin{proof}[Proof of Lemma \ref{lem:thermal-completion}]
We first assume that $x_0=0$ and outline the modifications for $x_0\ne0$ at
the end of the proof.
We append a superscript \((T)\) to the objects constructed from
\(m_0^{(T)}\). More precisely, let $X_t^{(T)}$ and $J_t^{(T)}$ denote the functions
$X_t$ and $J_t$ from Definition \ref{def:anyN}, respectively, defined from
$N_t^{(T)}$ instead of $N_t$; define $\varpi_0^{(T)}$ similarly.
Analogously, recalling $\Hop_{x,t}$ from \eqref{eqn:hop},
let
$\Hop_{x,t}^{(T)}$ be defined as in \eqref{eqn:hop}, but replacing
$\varpi_0$ in the definition \eqref{eqn:N0X} of $N_0^X$ with
$\varpi_0^{(T)}$. 
Since the initial densities
agree on $[x_0-L,x_0+L]$, $N_0^{(T)}(\lambda, x)=N_0(\lambda, x)$ and
$X_0^{(T)}(\lambda, x)=X_0(\lambda, x)$ for $x \in [x_0-L,x_0+L]$ and any $\lambda \in \mathbb{R}$. 

Fix \(0<s\le\tau := (\log T)^2\) and \(|x-x_0|\le L/2\). Let
\(\Hop_{x,s}\) and \(\Hop_{x,s}^{(T)}\) be the fixed point maps
for $N_t$ and $N_t^{(T)}$ (as in Lemma \ref{lem:H-contraction} for the former,
and as defined above, for the latter).
	 For each $\Lambda$, let $q_\Lambda$ be such that
$x_s(\Lambda;q_\Lambda)=x$ (recalling from \eqref{eqn:Xtlin} that
$x_s(\Lambda;\cdot)$ is increasing and thus invertible, since
$X_s(\Lambda,\cdot)$ is increasing by \eqref{eqn:avrbd2coarse}). By
Lemma \ref{lem:vebd},
\begin{equation}\label{eqn:xsx0b}
|q_\Lambda-x|\le Cs(1+|\Lambda|).
\end{equation}
Thus, if \(|\Lambda|\le cL/\tau\), with \(c>0\) sufficiently small,
then \(q_\Lambda\in[x_0-L,x_0+L]\), and hence
\[
\Hop_{x,s}^{(T)}[N_s(\cdot,x)](\Lambda)
=
\Hop_{x,s}[N_s(\cdot,x)](\Lambda).
\]

Using the Gaussian estimate
\eqref{eqn:vrhob} (and its analog 
with $m^{(T)}$ in place of $m$) for \(|\Lambda|>cL/\tau\), the above gives
\[
\left\|
\Hop_{x,s}^{(T)}[N_s(\cdot,x)]
-
\Hop_{x,s}[N_s(\cdot,x)]
\right\|_{\mathfrak c}
\le
C\exp\left(-c\frac{L^2}{\tau^2}\right).
\]
Since Lemma \ref{lem:H-contraction}, applied with $\varpi_0^{(T)}$ in place
of $\varpi_0$, implies that $\Hop_{x,s}^{(T)}$ satisfies the contraction
bound \eqref{eqn:H-contraction} with the same $a\in(0,1)$, adding and subtracting
\(\Hop_{x,s}^{(T)}[N_s(\cdot,x)]\) from $\Hop_{x,s}^{(T)}[N_s^{(T)}(\cdot,x)] - \Hop_{x,s}[N_s(\cdot,x)]$, we therefore obtain
\[
\sup_{|x-x_0|\le L/2}
\left\|
N_s^{(T)}(\cdot,x)-N_s(\cdot,x)
\right\|_{\mathfrak c}
\le
C\exp\left(-c\frac{L^2}{\tau^2}\right).
\]
By \eqref{eqn:lpbd} and the previous inequality, after decreasing \(c\),
\begin{equation}
\label{eq:thermal-completion-X-short}
\left|
X_s^{(T)}(\lambda,x)-X_s(\lambda,x)
\right|
\le
C\exp\left(-c\frac{L^2}{\tau^2}\right)
\end{equation}
whenever \(|x-x_0|\le L/2\) and \(|\lambda|\le cL/\tau\).

Next, decrease \(c>0\) if necessary. The hypothesis
\[
|q-x_0|+s(1+|\lambda|)\le cL
\]
and \eqref{eqn:xsx0b} ensure that both
\(x_r(\lambda;q)\) and \(x_r^{(T)}(\lambda;q)\) remain in
\([x_0-L/2,x_0+L/2]\) for \(0\le r\le s\). Moreover,
Lemma \ref{xt2} and the equality of the initial $X$-coordinates gives
\[
X_s^{(T)}
\left(\lambda,x_s^{(T)}(\lambda;q)\right)
=
X_s\left(\lambda,x_s(\lambda;q)\right).
\]
Since \eqref{eqn:avrbd2coarse} gives
\(\partial_xX_s^{(T)}(\lambda,x)\ge \frac12\), it follows (using the display above) that
\[
\begin{aligned}
\frac12\left|
x_s^{(T)}(\lambda;q)-x_s(\lambda;q)
\right|
&\le
\left|
X_s^{(T)}
\left(\lambda,x_s(\lambda;q)\right)
-
X_s\left(\lambda,x_s(\lambda;q)\right)
\right|\\
&\le
C\exp\left(-c\frac{L^2}{\tau^2}\right),
\end{aligned}
\]
where the last bound is by \eqref{eq:thermal-completion-X-short}. Since
\(L/\tau>\log T\), this gives the claimed estimate.

If $x_0\ne0$, we must shift $N_t^{(T)}$ and thus also slightly modify the
definition of $\Hop_{x,t}^{(T)}$. Set
\[
\begin{gathered}
d(\lambda):=N_0^{(T)}(\lambda,x_0)-N_0(\lambda,x_0),
\qquad e(\Lambda):=(\T d)(\Lambda),
\\
\widehat N_t^{(T)}:=N_t^{(T)}-d,
\qquad
\widehat X_t^{(T)}:=X_t^{(T)}-e.
\end{gathered}
\]
This leaves $m_t^{(T)}$, $J_t^{(T)}$, and $x_t^{(T)}$ unchanged. With the
notation
$\widehat\varpi_0^{(T)}(\Lambda,Y)
:=\varpi_0^{(T)}(\Lambda,Y+e(\Lambda))$, define
\[
\widehat\Hop_{x,t}^{(T)}[f](\Lambda)
:=\int_{-e(\Lambda)}^{x+(\T f)(\Lambda)-t\Lambda}
\widehat\varpi_0^{(T)}(\Lambda,Y)\,dY-d(\Lambda).
\]
The original fixed point equation and
$\widehat X_t^{(T)}(\Lambda,x)
=x+(\T\widehat N_t^{(T)}(\cdot,x))(\Lambda)$ give
\[
\widehat\Hop_{x,t}^{(T)}[\widehat N_t^{(T)}(\cdot,x)](\Lambda)
=\widehat N_t^{(T)}(\Lambda,x).
\]
The contraction bound \eqref{eqn:H-contraction} remains valid with this
definition of $\widehat\Hop_{x,t}^{(T)}$, and the agreement of the initial
densities gives
$\widehat N_0^{(T)}(\lambda,x)=N_0(\lambda,x)$ and
$\widehat X_0^{(T)}(\lambda,x)=X_0(\lambda,x)$ for
$x\in[x_0-L,x_0+L]$. The proof then follows the previous argument when
$x_0=0$, upon replacing $N^{(T)}$, $X^{(T)}$, and $\Hop^{(T)}$ by
$\widehat N^{(T)}$, $\widehat X^{(T)}$, and $\widehat\Hop^{(T)}$,
respectively; we omit further details.
\end{proof}

\subsection{Flowing back in time}
\label{subsec:backward_flow}

Let $m^{(T)}$ denote the solution of the GHD equation corresponding 
to Assumption \ref{ass:model} as defined in Section \ref{subsec:tc}.
The following uniform lower bound on $m^{(T)}$ is proved in Section
\ref{app:vac}.

\begin{lem}
\label{lem:uniform-non-vacuum}
Adopt Assumption \ref{ass:model}. There exist constants $c,C>0$ and
$\Lambda_\star\in[-C,C]$, independent of $\mathfrak{A}$, such that, setting
$J_\star=[\Lambda_\star-c,\Lambda_\star+c]$,
\[
\inf_{\tau \ge 0} \inf_{x \in \mathbb{R}}  \inf_{\lambda \in J_{\star}} m_\tau^{(T)}(\lambda,x)\ge c\mathfrak A^{-1}.
\]
\end{lem}

In what follows, let $m_t(\lambda,x)$ be defined as in Lemma
\ref{lem:int_ghd}, under Assumption \ref{ass:sec2}. Assume in addition that there is a constant
\(c_*>0\) such that
\begin{equation}\label{eqn:mintlow}
\inf_{s\ge0}\inf_{x\in\mathbb R}
\int_{\mathbb R}m_s(\lambda,x)\,d\lambda
\ge c_*\mathfrak A^{-1}.
\end{equation}
Lemma \ref{lem:uniform-non-vacuum} shows that this holds for
the solution $m_t^{(T)}$ under
Assumption \ref{ass:model}.

In the case that \eqref{eqn:mintlow} holds, we have the following stronger lower bound on
$J_t(\lambda,x)$. We postpone its proof until the end of this subsection.

\begin{lem}\label{lem:fixed_point_derivative}
Adopt Assumption \ref{ass:sec2}, let $N_t(\lambda,x)$ be the fixed point from
Lemma \ref{lem:fixed_pointout}. If, for some $t \geq 0$, there exists a
constant \(c_*>0\) such that \eqref{eqn:mintlow} holds,
then for some $c>0$,
\begin{equation}\label{eqn:avrbd2}
J_t(\lambda,x)
\ge c\left(1+\mathfrak A^{-1}\log(|\lambda|+2)\right),
\qquad \lambda,x\in\mathbb{R}.
\end{equation}
\end{lem}

Let \(x_t\) denote the characteristic flow
from Definition \ref{def:xtd}. By \eqref{eqn:Xtlin},
$x_t(\lambda;\cdot)$ is increasing and thus invertible, since
$X_t(\lambda,\cdot)$ is increasing by \eqref{eqn:avrbd2coarse}. Set
\[
\Phi_t(\lambda,x)
:=
\bigl(x_t(\lambda;\cdot)\bigr)^{-1}(x).
\]

\begin{lem}[Characteristic flow bounds]
\label{lem:backward-flow-rectangle}
Let $m_t$ be as in Lemma \ref{lem:int_ghd}, and suppose that
\eqref{eqn:mintlow} holds. Fix \(\tau\ge0\) and \(L\ge1\). There exists a
constant $C>0$, depending only on the constants from Assumption
\ref{ass:sec2} and on $c_*$ in \eqref{eqn:mintlow}, such that, for any
\(\lambda,\lambda'\in[-L,L]\) and \(q,q',x,x'\in\mathbb R\),
\begin{align}
|x_\tau(\lambda;q)-x_\tau(\lambda;q')|
&\le C |q-q'|,
\label{eq:flow-space}\\
|x_\tau(\lambda;q)-x_\tau(\lambda';q)|
&\le
C\tau(1+L)|\lambda-\lambda'|,
\label{eq:flow-spec}\\
|\Phi_\tau(\lambda,x)-\Phi_\tau(\lambda,x')|
&\le C|x-x'|,
\label{eq:backflow-space}\\
|\Phi_\tau(\lambda,x)-\Phi_\tau(\lambda',x)|
&\le
C\tau(1+L)|\lambda-\lambda'|.
\label{eq:backflow-spectral}
\end{align}

Consequently, for intervals \(I\subseteq\mathbb R\) and
\(J\subseteq[-L,L]\), there is an interval \(I_0\) such that
\[
\{(q,\lambda):\lambda\in J,\ x_\tau(\lambda;q)\in I\}
\subseteq I_0\times J
\]
and
\[
|I_0|
\le
C\left(
|I|+\tau(1+L)|J|
\right).
\]
\end{lem}

\begin{proof}
Recall $J_s(\mu,y)=\partial_yX_s(\mu,y)$ from Definition
\ref{def:anyN} and, throughout this proof, set
$h(\mu):=1+\mathfrak A^{-1}\log(|\mu|+2)$.
The Jacobian estimates \eqref{eqn:Jtbd} and \eqref{eqn:avrbd2} give
\begin{equation}\label{eq:characteristic-J-comparison}
c h(\mu)\le J_s(\mu,y)\le C h(\mu).
\end{equation}
Moreover, Lemma \ref{xt2} gives
\begin{equation}\label{eq:characteristic-coordinate-identity}
X_s(\mu,x_s(\mu;q))-s\mu=X_0(\mu,q).
\end{equation}
Differentiating the identity in \(q\), and then its inverse in \(x\), yields
\[
\partial_qx_s(\mu;q)
=
\frac{J_0(\mu,q)}{J_s(\mu,x_s(\mu;q))},
\qquad
\partial_x\Phi_s(\mu,x)
=
\frac{J_s(\mu,x)}{J_0(\mu,\Phi_s(\mu,x))}.
\]
Thus \eqref{eq:characteristic-J-comparison}, followed by integration, proves
\eqref{eq:flow-space} and \eqref{eq:backflow-space}.

We next prove \eqref{eq:backflow-spectral}.  Fix \(x\in
\mathbb R\) and \(\lambda,\lambda'\in[-L,L]\), and set
\begin{equation}\label{qqpdf}
q:=\Phi_\tau(\lambda,x),
\qquad
q':=\Phi_\tau(\lambda',x).
\end{equation}
Define
\[
F_s(\mu):=X_s(\mu,x_s(\lambda';q'))-s\mu,
\qquad 0\le s\le\tau.
\]
Writing $y_s:=x_s(\lambda';q')$, the chain rule gives
\[
\begin{aligned}
\frac{d}{ds}F_s(\mu)
&=\partial_sX_s(\mu,y_s)+J_s(\mu,y_s)\dot y_s-\mu\\
&=J_s(\mu;y_s)
\bigl(\ve_s(\lambda';y_s)-\ve_s(\mu;y_s)\bigr).
\end{aligned}
\]
Here the second equality follows from \eqref{eqn:Kx-general-solution},
\eqref{eqn:Kttransport-general}, and the fact $\dot y_s = \ve_s(\lambda';y_s)$.
Hence \eqref{eqn:Jtbd} and \eqref{eqn:vebd} give, for
\(\mu,\lambda'\in[-L,L]\),
\[
|F_\tau(\mu)-F_0(\mu)|
\le
C\tau h(\mu)(1+L) |\mu-\lambda'|.
\]
Taking \(\mu=\lambda\) in this estimate and using
\eqref{eq:characteristic-coordinate-identity}, we obtain
\[
\begin{aligned}
|X_0(\lambda,q)-X_0(\lambda,q')|
&=|X_\tau(\lambda,x_\tau(\lambda;q))-\tau\lambda-X_0(\lambda,q')|\\
&=|X_\tau(\lambda,x_\tau(\lambda';q'))-\tau\lambda-X_0(\lambda,q')|\\
&=|F_\tau(\lambda)-F_0(\lambda)|\\
&\le C\tau h(\lambda)(1+L)|\lambda-\lambda'|.
\end{aligned}
\]
In the second equality we used $x_\tau(\lambda;q) = x = x_\tau(\lambda';q')$ from \eqref{qqpdf}. The lower bound in \eqref{eq:characteristic-J-comparison} gives
$|q-q'| \leq \frac{C}{h(\lambda)} |X_0(\lambda,q)-X_0(\lambda,q')|$; this 
together with the display above implies \eqref{eq:backflow-spectral}.

To prove \eqref{eq:flow-spec}, fix \(q\in\mathbb R\) and put
\[
z:=x_\tau(\lambda;q),
\qquad
\widetilde q:=\Phi_\tau(\lambda',z).
\]
Then \(x_\tau(\lambda';\widetilde q)=z\), so \eqref{eq:flow-space} and
\eqref{eq:backflow-spectral} imply
\[
|x_\tau(\lambda;q)-x_\tau(\lambda';q)|
\le C|\widetilde q-q|
\le
C^2\tau(1+L)|\lambda-\lambda'|.
\]

It remains to prove the rectangle conclusion.  It suffices to consider
nonempty bounded intervals \(I\) and \(J\), since the other cases are
immediate.  Replacing \(I\) by its closure, write \(I=[a,b]\) and set
\[
I_0:=
\left[
\inf_{\lambda\in J}\Phi_\tau(\lambda,a),
\sup_{\lambda\in J}\Phi_\tau(\lambda,b)
\right].
\]
The monotonicity of \(\Phi_\tau(\lambda,\cdot)\) gives the required
containment.  For any \(\lambda_0\in J\), \eqref{eq:backflow-space} and
\eqref{eq:backflow-spectral} give
\[
|I_0|
\le
\Phi_\tau(\lambda_0,b)-\Phi_\tau(\lambda_0,a)
+
C\tau(1+L)|J|
\le
C\left(
|I|+\tau(1+L)|J|
\right).
\]
\end{proof}

\begin{proof}[Proof of Lemma \ref{lem:fixed_point_derivative}]

To show \eqref{eqn:avrbd2}, first assume that $|\lambda|\ge 4$ is
sufficiently large and set
\[
I_\lambda:=[-|\lambda|^{1/2},|\lambda|^{1/2}].
\]
By \eqref{eqn:avrbd},
\[
\int_{I_\lambda^{\complement}}|\log|\lambda-\mu||m_t(\mu,x)\,d\mu
\le
C\mathfrak A^{-1}e^{-c|\lambda|}.
\]
On the other hand, if \(\mu\in I_\lambda\), then
\[
|\lambda-\mu|
\ge
|\lambda|-|\lambda|^{1/2}
\ge
\frac{|\lambda|}{2},
\]
and hence
\[
\log|\lambda-\mu|
\ge
\log\frac{|\lambda|}{2}
\ge
c\log(|\lambda|+2).
\]
Therefore
\[
\begin{aligned}
(\T m_t(\cdot,x))(\lambda)
&=
2\int_{\mathbb R}\log|\lambda-\mu|\,m_t(\mu,x)\,d\mu \\
&\ge
2\int_{I_\lambda}\log|\lambda-\mu|\,m_t(\mu,x)\,d\mu
-
2\int_{I_\lambda^{\complement}}|\log|\lambda-\mu||\,m_t(\mu,x)\,d\mu \\
&\ge
c\log(|\lambda|+2)\int_{I_\lambda}m_t(\mu,x)\,d\mu
-
C\mathfrak A^{-1}.
\end{aligned}
\]
Define
$\alpha_t(x)^{-1}:=\int_{\mathbb R}m_t(\mu,x)\,d\mu$. By
\eqref{eqn:mintlow}, $\alpha_t(x)^{-1}\ge c_*\mathfrak A^{-1}$, and
so the tail bound \eqref{eqn:avrbd} on $m_t(\lambda,x)$ above gives
\[
\int_{I_\lambda}m_t(\mu,x)\,d\mu
=
\alpha_t(x)^{-1}+O(\mathfrak A^{-1}e^{-c|\lambda|}).
\]
Hence
\[
(\T m_t(\cdot,x))(\lambda)
\ge
c\,\alpha_t(x)^{-1}\log(|\lambda|+2)-C\mathfrak A^{-1}.
\]
Thus, assuming $\mathfrak{A}$ is large enough, after decreasing the constant \(c\), we obtain
\[
1+(\T m_t(\cdot,x))(\lambda)
\ge 1/2+c\mathfrak A^{-1}\log(|\lambda|+2)
\]
for all sufficiently large \(|\lambda|\), thereby verifying
\eqref{eqn:avrbd2} if $|\lambda|$ is sufficiently large. For other
$\lambda$, \eqref{eqn:avrbd2} follows from \eqref{eqn:avrbd2coarse},
which proves the lemma in general.
\end{proof}

\section{Quasiparticle approximations}
\label{sec:qa}
		
	Throughout this section, we adopt Assumption \ref{ass:model}. We first prove
	Proposition \ref{prop:discrete-continuum-scattering}, thereby completing the
	proof of Proposition \ref{prop:improv_M} from Section
	\ref{subsec:ideas_Qx}. We then prove Theorems
	\ref{thm:local-charge-testfn-intro} and \ref{thm:main} using 
	Theorem \ref{ttn} as input.

	Throughout, set $\bar q_i=T^{-1}q_i(0)$, as in Section
	\ref{subsec:stop_outline}, and recall the notation from Definition
	\ref{def:quasi_intro}, the functions $N_\tau^{(T)}$ and $m_\tau^{(T)}$ from
	Section \ref{subsec:tc}, and the characteristic flows
	$x_\tau^{(T)}$ and $x_\tau$ from Definition \ref{def:xtd} under
	$m^{(T)}$ and $m$, respectively. We also recall
	$\mathfrak M$ from Definition \ref{mt} and
	$R=\mathfrak A^{1/10}$ from \eqref{eqn:R}.

	\subsection{Proof of Proposition \ref{prop:discrete-continuum-scattering}}
	\label{subsec:quasi_thmpf}

		First, we have the following technical bound. Its proof is in Section \ref{subsec:min}.
	
\begin{lem}\label{lem:sumlambdak01_half}
Let \(\mathfrak M\) be as in Definition \ref{mt}. Fix \(\tau\in[0,\log N]\) and set \(t=T\tau\).
There is a constant \(C>0\) such that, with overwhelming probability, the
following holds whenever
\(\mathfrak B(t)\) holds, simultaneously for every \(k\in\llb1,N\rrb\) and every
$a\in[-(\log N)^{100},(\log N)^{100}]$:
\begin{flalign}
\label{eqn:sumlambdak01_half}
\begin{aligned}
 \left| \displaystyle\sum_{\substack{i\ne k:\,
 x_{\tau}^{(T)}(\lambda_i; Q_i(0)T^{-1} ) < a}}  
\log |\lambda_k - \lambda_i| 
-
\displaystyle\sum_{\substack{i\ne k:\,
 T^{-1}Q_i(t)<a}}  
\log |\lambda_k - \lambda_i|   \right| 
\le
C R^{-1}\mathfrak M \log(|\lambda_k|+2).
\end{aligned}
\end{flalign}
\end{lem}

	\begin{proof}[Proof of Proposition \ref{prop:discrete-continuum-scattering}]
	 Fix $k$ such that $|\varphi_0(k)| \le T (\log N)^{10}$. We restrict
	 to the event on which the conclusion of
	 Proposition \ref{vn1n2out} holds (for $t$ there given by both $0$ and $t$).
	
	 To that end, set 
	\begin{flalign*}
		\mathfrak{C}(t) = \log \displaystyle\sum_{j=1}^N e^{-\lambda_j t} u_j (N_1; 0)^2.
	\end{flalign*}
	
	\noindent Then, Lemma \ref{uknt} gives 
	\begin{flalign*}
	2 \log |u_k (N_1;t)| = 2 \log |u_k (N_1;0)| - \lambda_k t - \mathfrak{C}(t). 
	\end{flalign*} 
	
	\noindent As $Q_j (t) = q_{\varphi_t (j)} (t)$, this with Proposition \ref{vn1n2out} 
	(using Proposition \ref{prop:outside-strip-confinement} to verify the hypothesis
	on $\varphi_t(k)$),
	 gives 
	\begin{equation}
	\label{eqn:sc-with-boundary-term}
	\begin{aligned}
	Q_k (t) & = Q_k(0) + \lambda_k t + 2 \displaystyle\sum_{j: Q_j (0) < Q_k (0)} \log |\lambda_k - \lambda_j| - 2 \displaystyle\sum_{j: Q_j (t) < Q_k (t)} \log |\lambda_k - \lambda_j| \\
	& \qquad + 	q_{N_1} (t) - q_{N_1} (0) + \mathfrak{C}(t) + O \big( R^{-1} \mathfrak{M} \log (|\lambda_k|+2) \big) .
	\end{aligned}
	\end{equation}
	
	\noindent 	Next we will show that in fact,
	\begin{equation}
		\label{eqn:Btbd}
		q_{N_1}(t)-q_{N_1}(0)+\mathfrak{C}(t)=0.
	\end{equation}
	Indeed, differentiating the definition of \(\mathfrak{C}(t)\) and using
	Lemma \ref{uknt} gives
	\[
		\mathfrak{C}'(t)
		=
		-\sum_{j=1}^N \lambda_j u_j(N_1;t)^2.
	\]
	By the spectral decomposition of \(\L(t)\), the last sum is the
	\((N_1,N_1)\)-entry of \(\L(t)\), hence equals \(b_{N_1}(t)\). Since
	\(\partial_t q_{N_1}(t)=b_{N_1}(t)\) and \(\mathfrak C(0)=0\), we obtain
	\(\mathfrak C(t)=q_{N_1}(0)-q_{N_1}(t)\), which proves
	\eqref{eqn:Btbd}. Equations \eqref{eqn:sc-with-boundary-term} and \eqref{eqn:Btbd} imply
	\eqref{eqn:sc}.

	Set
	\[
	a:=T^{-1}Q_k(t),
	\qquad
	L_k:=\log(|\lambda_k|+2).
	\]
	By the fact that $|\varphi_0(k)| \le T (\log N)^{10}$, Lemmas \ref{l0eigenvalues} and \ref{lem:Estr}, and Proposition \ref{prop:outside-strip-confinement},
		after restricting to smaller overwhelmingly probable event, we 
		may assume that \(a\in[-(\log N)^{15} ,(\log N)^{15}]\). 
		
	First, note that on the overwhelmingly probable event of Lemma \ref{lem:sumlambdak01_half}, we have
	\begin{equation}\label{eqn:Qttox}
	\left| \displaystyle\sum_{\substack{j\ne k:\,T^{-1} Q_j (t) < a}} \log |\lambda_k - \lambda_j| - \displaystyle\sum_{\substack{j\ne k:\,x_{\tau}^{(T)}(\lambda_j; T^{-1}Q_j (0)) < a}} \log |\lambda_k - \lambda_j| \right|   
	\leq R^{-1} \mathfrak{M} L_k. 
	\end{equation}
	Combining \eqref{eqn:sc} and \eqref{eqn:Qttox} gives
	\begin{align}
	Q_k (t) & = Q_k(0) + \lambda_k t + 2 \displaystyle\sum_{j: Q_j (0) < Q_k (0)} \log |\lambda_k - \lambda_j| - 2 \displaystyle\sum_{\substack{j\ne k:\,x_{\tau}^{(T)}(\lambda_j; T^{-1}Q_j (0))  < a}} \log |\lambda_k - \lambda_j| \notag \\
		& \qquad  + O \big( R^{-1} \mathfrak{M} \log (|\lambda_k|+2) \big). \label{eqn:scat2}
		\end{align}
	Recall the integral $\mathcal I_{\tau,k}^{(T)}$ from Proposition
	\ref{prop:scat_conc}. Applying that proposition
	with this value of \(a\) to \eqref{eqn:scat2} provides, on an event of overwhelming probability,
	\begin{equation}\label{eq:qkt-moving-cutoff}
	Q_k(t)
	=
	q_{\varphi_0(k)}(0)+\lambda_k t
	+2T\mathcal I_{\tau,k}^{(T)}(a)
	+
	 O\Bigl( R^{-1}\mathfrak M L_k\Bigr).
	\end{equation}
	Since $a=T^{-1}Q_k(t)$, expanding the definition of
	$\mathcal I_{\tau,k}^{(T)}(a)$ in \eqref{eq:qkt-moving-cutoff} gives
	\eqref{eqn:scattering-continuum-cutoff}, proving the second part of the proposition.
	\end{proof} 
%
%
%
\subsection{Proofs of Theorems \ref{thm:local-charge-testfn-intro} and \ref{thm:main}}
\label{subsec:mainpfs}

Before proceeding to the proofs of main theorems, we state several preliminary lemmas. The first describes a coupling between two 
Toda lattices; one on $\mathbb{Z}$, which we denote $\L^{\mathbb{Z}}$, whose initial entries are 
sampled from
 distributions satisfying Assumption \ref{ass:model_inf},
and another on $\llbracket N_1, N_2 \rrbracket$ with initial entries satisfying
 Assumption \ref{ass:model}, which we simply denote by $\L$. We denote 
 their time evolutions by $\L^{\Z}(t) = [L_{i j}^{\Z}(t)]_{i,j \in \Z}$ and $\L(t) = [L_{ij}(t)]_{i,j \in \llbracket N_1, N_2 \rrbracket}$, respectively. We omit its 
 proof, as it follows from \cite[Lemma 4.5 and Proposition 2.5(3)]{Agg25a}.\footnote{While
 \cite[Proposition 2.5(3)]{Agg25a} was stated under thermal equilibrium, it
 is directly verified from its proof that it also holds under Assumption
 \ref{ass:model}.}

\begin{lem}[{\cite{Agg25a}}]\label{lem:N1N1Zcouple}
Let $\L^{\mathbb{Z}} = \L^{\Z}(0)$ and $\L = \L(0)$ be as described above. 
Then, there exists a coupling between $\L^{\mathbb{Z}}$ and $\L$ under which,
on an event of overwhelming probability, 
\begin{equation}\label{eqn:supLLbd}
\sup_{0 \leq t \leq T \log N} \sup_{|i|, |j| \leq T (\log N)^2} | L_{i j}^{\Z}(t) - L_{i j}(t)| \leq e^{- T }. 
\end{equation}
\end{lem}

Next, we state several lemmas which will be useful analyzing 
convergence of the linear statistics appearing
 in Theorem \ref{thm:local-charge-testfn-intro}, Corollary \ref{cor:local-charge-intro},
 and Theorem \ref{thm:main}.  The first lemma is a reduction which allows
to study the sums appearing in 
Corollary \ref{cor:local-charge-intro} by studying
certain
linear statistics.  It is similar to \cite[Proposition 2.10(1)]{Agg25a}, whose
proof uses localization properties of the Lax matrix. In our setting, the
required localization is provided by \eqref{eqn:uMR}. Lemma \ref{lem:cl} is
proved in Section \ref{subsec:local-charge-intro}.
	
	\begin{lem}
\label{lem:cl}
Adopt Assumption \ref{ass:model}. Fix
\(m\in\mathbb Z_{\ge0}\) and $t \in [0, T \log N]$. There is a constant \(C>0\), depending
only on \(m\) and the constants in Assumption \ref{ass:model}, such
that, with overwhelming probability,
the following holds simultaneously for every interval
$J\subseteq
[-T(\log N)^{20},\,T(\log N)^{20}]$.
One has
\begin{equation}\label{eqn:cl}
\left|
\sum_{i:\,q_i(t)\in J}(\L(t)^m)_{ii}
-
\sum_{k:\,Q_k(t)\in J}\lambda_k^m
\right|
\le
T^{3/4}(\log N)^C.
\end{equation}
\end{lem}

The next lemmas are statements which together provide
quantitative error bounds for convergence of
 a certain class of linear statistics.

\begin{assumption}[Admissible mesoscopic test functions]
\label{ass:Fas}

Adopt Assumption \ref{ass:model}.
Fix \(m\in\mathbb Z_{\geq 0}\). Let \(I_T\subset\mathbb R\) be an interval satisfying
$|I_T|\leq T\log N$,
and set
$K_T:=T^{-1}I_T$.
Assume that
$K_T\subseteq
\bigl[-(\log N)^{19},(\log N)^{19}\bigr]$.

Let
$F:\mathbb R^2\longrightarrow\mathbb R$ 
be bounded and locally Lipschitz, with
$\supp F\subseteq \mathbb R\times K_T$.
We denote the argument of $F$ by $(\lambda,x)$ and therefore interpret
$\partial_\lambda F$ and $\partial_xF$ as the almost everywhere
(equivalently, weak) derivatives of $F$ in their first and second variables,
respectively. Define
\begin{align*}
\mathcal A(F)
&:=\sup_{\lambda,x\in\mathbb R}|F(\lambda,x)|,\\
\mathcal V(F)
&:=\sup_{\lambda\in\mathbb R}\int_{\mathbb R}
|\partial_xF(\lambda,x)|\,dx,\\
\mathcal D(F)
&:=\sup_{\lambda\in\mathbb R}\sup_{x\neq x'}
\frac{|F(\lambda,x)-F(\lambda,x')|}{|x-x'|} ,\\
\mathcal B(F)
&:=
\sup_{x\in\mathbb R}
\sup_{\lambda\neq\lambda'}
\frac{|F(\lambda,x)-F(\lambda',x)|}
{|\lambda-\lambda'|}.
\end{align*}
Assume that $\mathcal A(F)+\mathcal V(F)+\mathcal D(F)+\mathcal B(F)<\infty$.

\end{assumption}
%
%

In each of the following lemmas, we fix \(\tau\in[0,\log N]\) and set
\(t=T\tau\). Recall that $\bar q_i=T^{-1}q_i(0)$. Lemmas \ref{lem:Qtox} and \ref{lem:sumtoint} are proved in Section \ref{subsec:Qtox},
and Lemma \ref{lem:xTtox} is proved 
in Section \ref{subsec:xTtox}.

\begin{lem}\label{lem:Qtox}
Adopt Assumption \ref{ass:Fas} for some \(F : \mathbb R^2 \rightarrow \mathbb R \). 
There exists
a constant $C> 0$ such that with overwhelming probability,
\begin{equation}\label{eqn:FLQx}
\left| \sum_{i=N_1}^{N_2} F(\Lambda_i, x_\tau^{(T)}(\Lambda_i; \bar q_i)) -  \sum_{k=1}^{N} F(\lambda_k, Q_k(t) T^{-1})\right| \leq 
 \mathcal{D}(F) (\log N)^C   T^{3/4}.
\end{equation}
\end{lem}


\begin{lem}\label{lem:sumtoint}
Adopt Assumption \ref{ass:Fas} for some \(F : \mathbb R^2 \rightarrow \mathbb R \). There exists
a constant $C> 0$ such that with overwhelming probability,
\begin{multline}\label{eqn:Fconv}
\left| \sum_{i=N_1}^{N_2} F(\Lambda_i, x_\tau^{(T)}(\Lambda_i; \bar q_i)) - T \int_{-\infty}^{\infty}  \int_{-\infty}^{\infty}  F(\lambda, q) m_{\tau}^{(T)}(\lambda, q) d \lambda d q \right| \\
 \leq 
(\log N)^C
	\left[
	(\mathcal{A} + \mathcal{V})T^{3/5}
	+(\mathcal{B} + \mathcal{D}) T^{2/5}
	\right].
\end{multline}
\end{lem}

\begin{lem}\label{lem:xTtox}
Adopt Assumption \ref{ass:Fas} for some \(F : \mathbb R^2 \rightarrow \mathbb R \),
 and assume in addition that
the interval $K_T$ satisfies
$K_T \subset [-(\log N)^2, (\log N)^2]$. There exist constants
$c>0$ such that, with overwhelming probability,
\begin{equation}\label{eqn:FxTx}
\left| \sum_{i=N_1}^{N_2} F(\Lambda_i, x_\tau^{(T)}(\Lambda_i; \bar q_i)) -  \sum_{i=N_1}^{N_2} F(\Lambda_i, x_\tau(\Lambda_i; \bar q_i)) \right| \leq 
\mathcal D e^{-c(\log N)^2} .
\end{equation}
In addition, 
\begin{equation}\label{eqn:FmTm}
T \left|\int_{-\infty}^{\infty}  \int_{-\infty}^{\infty}  F(\lambda, q) m_{\tau}^{(T)}(\lambda, q) d \lambda d q -  \int_{-\infty}^{\infty}  \int_{-\infty}^{\infty}  F(\lambda, q) m_{\tau}(\lambda, q) d \lambda d q \right| \leq 
(\mathcal{A} +  \mathcal{D}) e^{-c(\log N)^2}  .
\end{equation}
\end{lem}

\begin{proof}[Proof of Theorem \ref{thm:local-charge-testfn-intro}]
Set $t:=T\tau$. Let us first replace the matrix $\L^{\Z}$, under which the
charges $\mathfrak k_i^{[n]}$ are defined in Theorem
\ref{thm:local-charge-testfn-intro}, by $\L$ (recall the notation defined at the
beginning of this subsection). In what follows, we let
$\bm q(t)=(q_i(t))_{i=N_1}^{N_2}$ and
$\bm q^{\Z}(t)=(q_i^{\Z}(t))_{i\in\Z}$ denote the Toda lattice position
variables associated with $\L(t)$ and $\L^{\Z}(t)$, respectively. It quickly
follows from Lemma \ref{lem:Estr} and \eqref{eqn:alph_bd}, with Lemma
\ref{l0eigenvalues} and \eqref{eqn:supLLbd} used to bound the derivatives of
$q_i$ and $q_i^{\Z}$, that, with overwhelming probability, any index $i$ for
which either $g(q_i(t)/T)\ne0$ or $g(q_i^{\Z}(t)/T)\ne0$ must satisfy
$|i|\le T(\log N)^2/2$. For such $i$, the combination of
\eqref{eqn:supLLbd}, Lemma \ref{l0eigenvalues} (at $A=\log N$), and the
tridiagonality of $\L(t)$ and $\L^{\Z}(t)$ yields
\[
\left|\mathfrak{k}_i^{[n]}(t)-(\L(t)^n)_{ii}\right|
\leq C_n(\log N)^n e^{-T},
\qquad
\left|g(q_i^{\Z}(t)/T)-g(q_i(t)/T)\right|
\leq \|g'\|_\infty T^{-1}e^{-cT}.
\]
Summing over $|i|\leq \frac12T(\log N)^2$ shows that replacing $\L^{\Z}$ and
$q^{\Z}$ by $\L$ and $q$ changes the normalized sum by $O(e^{-cT})$.
It suffices to prove
\[
\frac1T\sum_{i=N_1}^{N_2}g(q_i(t)/T)(\L(t)^n)_{ii}
\longrightarrow
\iint_{\mathbb R^2}g(x)\lambda^n m_\tau(\lambda,x)\,d\lambda\,dx.
\]

Choose \(K=[a,b]\) with \(\supp g\subset(a,b)\).
For \(y\in[a,b]\), set
\[
D_t(y):=
\sum_{i:\,q_i(t)\in[aT,yT]}(\L(t)^n)_{ii}
-
\sum_{k:\,Q_k(t)\in[aT,yT]}\lambda_k^n.
\]
Lemma \ref{lem:cl} then gives, with overwhelming probability,
\begin{equation}\label{eq:local-charge-cumulative}
\sup_{y\in[a,b]}|D_t(y)|
\le T^{3/4}(\log N)^C.
\end{equation}
Define the signed measure
\[
\nu_t:=
\sum_{i=N_1}^{N_2}(\L(t)^n)_{ii}\delta_{q_i(t)/T}
-
\sum_{k=1}^N\lambda_k^n\delta_{Q_k(t)/T}.
\]
Thus \(D_t(y)=\nu_t([a,y])\).  Since \(g\) vanishes in neighborhoods
of \(a\) and \(b\), integration by parts and
\eqref{eq:local-charge-cumulative} give
\begin{equation}
\label{eq:local-charge-particle-quasiparticle}
\begin{aligned}
\left|\sum_{i=N_1}^{N_2} g(q_i(t)/T)(\L(t)^n)_{ii}
-\sum_{k=1}^N g(Q_k(t)/T)\lambda_k^n\right|
&=\left|\int_a^b g(y)\,\nu_t(dy)\right|\\
&=\left|\int_a^b g'(y)D_t(y)\,dy\right|\\
&\le T^{3/4}(\log N)^C\int_{\mathbb R}|g'(x)|\,dx.
\end{aligned}
\end{equation}
On the event from Lemma \ref{l0eigenvalues} with \(A=\log N\), fix a smooth
function $\chi:\mathbb R\to[0,1]$ supported on $[-2,2]$, such that
$\chi=1$ on $[-1,1]$, and set
\(H_N(\lambda,x)=\chi(\lambda/\log N)\lambda^ng(x)\).  Then Assumption \ref{ass:Fas}
is satisfied with $F(\lambda, x) = H_N(\lambda,x)$ and
\(\sum_{\mathcal X\in\{\mathcal A,\mathcal V,\mathcal D,\mathcal B\}}
\mathcal X(H_N)\le C(\log N)^n\) and
\(H_N(\lambda_k,x)=\lambda_k^ng(x)\) for each $\lambda_k \in \eig \L(0)$ 
since we restrict to $\mathsf{BND}_{\L(0)}(\log N)$. Displays
\eqref{eq:local-charge-particle-quasiparticle}, \eqref{eqn:FLQx},
\eqref{eqn:Fconv}, and \eqref{eqn:FmTm}
from Lemmas \ref{lem:Qtox},
\ref{lem:sumtoint}, and \ref{lem:xTtox} with
\((F,\tau,I_T)=(H_N,\tau,[aT,bT])\), therefore give, with overwhelming
probability, after division by \(T\),
\[
\begin{aligned}
 &\left|\frac{1}{T}\sum_{i=N_1}^{N_2} g(q_i(t)/T)(\L(t)^n)_{ii}
-\int_{-\infty}^{\infty}  \int_{-\infty}^{\infty} \chi(\lambda/\log N)\lambda^ng(q)m_{\tau}(\lambda, q) d \lambda d q  \right|\\
&\qquad \qquad \leq C(\log N)^C(T^{-1/4}+T^{-2/5}+T^{-3/5}).
\end{aligned}
\]
The Gaussian bound \eqref{eqn:avrbd} removes \(\chi\) from the limiting
integral at cost \(O(e^{-c(\log N)^2})\). This completes the proof.
\end{proof}

\begin{proof}[Proof of Corollary \ref{cor:local-charge-intro}]
Set $t:=T\tau$. Following the argument in the beginning of the proof of
Theorem \ref{thm:local-charge-testfn-intro} (we omit further details), it
suffices to prove
\begin{equation}\label{eqn:corts}
\frac1{|J_T|}\sum_{i:\,q_i(t)\in J_T}(\L(t)^n)_{ii}
\longrightarrow
\int_{\mathbb R}\lambda^n m_\tau(\lambda,x)\,d\lambda.
\end{equation}

Put $K_T:=T^{-1}J_T$ and
$\rho_T:=T^{-1/4}(\log N)^{C_0}$, where $C_0$ is a sufficiently large
constant, to be fixed later. Since
$T^{a-1}\le |K_T|\le T^{a'-1}$,
$a>\frac56$,
we have $\rho_T=o(|K_T|)$.
For $r\ge0$, let
$K_T^{(r)}:=\{y\in\mathbb R:\operatorname{dist}(y,K_T)\le r\}$.
Moreover, \(K_T\) converges to the fixed point \(x\), and hence
$K_T^{(2\rho_T)}\subset[-(\log N)^2,(\log N)^2]$ for all sufficiently
large \(T\).

For \(k\in\llbracket1,N\rrbracket\), write
\[
Y_k
:=
x_\tau^{(T)}
\left(\lambda_k;\frac{Q_k(0)}{T}\right).
\]
We restrict to the events from Theorem \ref{ttn} and Lemma
\ref{l0eigenvalues}, with \(A=\log N\). By \eqref{eqn:Qx}, for $C_0$
sufficiently large,
\begin{equation}
\label{eq:cor-local-QY}
\left|\frac{Q_k(t)}{T}-Y_k\right|
\le\frac{\rho_T}{4}
\end{equation}
for every $k$ for which either $Q_k(t)/T\in K_T$ or
$Y_k\in K_T^{(\rho_T)}$. The required hypothesis
$\varphi_t(k)\in\llbracket N_1+T(\log N)^{10},
N_2-T(\log N)^{10}\rrbracket$ in Definition \ref{mt} follows, for
sufficiently large $T$ and after restricting to a smaller overwhelmingly
probable event, from Lemma \ref{lem:Estr} and Proposition
\ref{prop:outside-strip-confinement}.

Choose a smooth function
$\psi_T:\mathbb R\longrightarrow[0,1]$
such that
\[
\psi_T\equiv1\quad\text{on }K_T^{(\rho_T/4)},
\qquad
\supp\psi_T\subset K_T^{(\rho_T)},
\qquad
\|\psi_T'\|_\infty\le10\rho_T^{-1}.
\]
Such a function is obtained by smoothing the indicator of \(K_T\) on
scale \(\rho_T\). Also fix a smooth function
$\chi:\mathbb R\to[0,1]$ supported on $[-2,2]$, such that $\chi=1$ on
$[-1,1]$, and define
\[
F_T(\lambda,y)
:=
\chi(\lambda/\log N)\lambda^n\psi_T(y).
\]

On the event $\mathsf{BND}_{\L(t)}(\log N)$ (of Lemma \ref{l0eigenvalues}), \(\chi(\lambda_k/\log N)=1\)
for every \(k\). Thus \eqref{eq:cor-local-QY} gives
\begin{equation}
\label{eq:cor-local-discrete-smoothing}
\left|
\sum_{k:\,Q_k(t)\in J_T}\lambda_k^n
-
\sum_{k=1}^N F_T(\lambda_k,Y_k)
\right|
\le
2 (\log N)^n \sum_{k=1}^N
\mathbf 1_{K_T^{(5\rho_T/4)}\setminus K_T}(Q_k(t)T^{-1})
\leq T^{3/4}(\log N)^C.
\end{equation}
Indeed, the last inequality follows from \eqref{eqn:qbd-lower} and the fact
that the $(q_i(t))_{i=N_1}^{N_2}$ constitute a permutation of the
$(Q_k(t))_{k=1}^N$.

The function \(F_T\) satisfies Assumption \ref{ass:Fas} with
$I_T=TK_T^{(\rho_T)}$, and
$\mathcal A+\mathcal V+\mathcal B\le C(\log N)^C$ and
$\mathcal D\le C\rho_T^{-1}(\log N)^C$.
Lemmas \ref{lem:sumtoint} and \ref{lem:xTtox} give, with overwhelming
probability,
\begin{equation}
\label{eq:cor-local-sum-to-integral}
\begin{aligned}
\sum_{k=1}^N F_T(\lambda_k,Y_k)
=
T\iint_{\mathbb R^2}
F_T(\lambda,y)m_\tau(\lambda,y)\,d\lambda\,dy
+
O\left(
(\log N)^C
\left[
T^{3/5}+\rho_T^{-1}T^{2/5}
\right]
\right).
\end{aligned}
\end{equation}
Since $\rho_T^{-1}T^{2/5}=T^{13/20}(\log N)^{-C_0}$,
the error in \eqref{eq:cor-local-sum-to-integral} is
$O\left(T^{3/4}(\log N)^C\right)$.

In the integral on the right side of
\eqref{eq:cor-local-sum-to-integral}, we may by \eqref{eqn:avrbd} replace
$F_T(\lambda,y)=\chi(\lambda/\log N)\lambda^n\psi_T(y)$ with
$\lambda^n\psi_T(y)$, up to an error of $Ce^{-c(\log N)^2}$, and also
replace $\psi_T$ with $\mathbf1_{K_T}$, up to an error of
$C\rho_T(\log N)^C$.
Combining these estimates with
\eqref{eq:cor-local-discrete-smoothing} gives
\begin{equation}
\label{eq:cor-local-quasiparticle-sum}
\sum_{k:\,Q_k(t)\in J_T}\lambda_k^n
=
T\int_{K_T}\int_{\mathbb R}
\lambda^n m_\tau(\lambda,y)\,d\lambda\,dy
+
O\left(T^{3/4}(\log N)^C\right).
\end{equation}

By Lemma \ref{lem:cl}, with overwhelming probability the left side in the
display \eqref{eqn:corts} differs from the left side of
\eqref{eq:cor-local-quasiparticle-sum} by at most
$T^{3/4}(\log N)^C$. Since $|J_T|=T|K_T|$ and $a>5/6>3/4$, we have
$T^{3/4}(\log N)^C/|J_T|\longrightarrow0$. If
$f_\tau(y):=\int_{\mathbb R}\lambda^n m_\tau(\lambda,y)\,d\lambda$, then
$\frac{T}{|J_T|}\int_{K_T}f_\tau(y)\,dy
=\frac1{|K_T|}\int_{K_T}f_\tau(y)\,dy\longrightarrow f_\tau(x)$ because
$\sup_{y\in K_T}|y-x|\longrightarrow0$ and $f_\tau$ is continuous at $x$
by assumption. This proves the corollary.
\end{proof}

%

\begin{proof}[Proof of Theorem \ref{thm:main}]
Let \(K\) contain the \(x\)-projection of \(\supp f\).  Apply Lemmas
\ref{lem:Qtox}, \ref{lem:sumtoint}, and \ref{lem:xTtox} with
\((F,\tau,I_T)=(f,\tau,TK)\). Thus, by \eqref{eqn:FLQx},
\eqref{eqn:Fconv}, and \eqref{eqn:FmTm},
\[
\left|\frac1T\sum_{i=1}^N
f\bigl(\lambda_i,Q_i(T\tau)T^{-1}\bigr)
-\iint_{\mathbb R^2} f(\lambda,x)m_\tau(\lambda,x)\,d\lambda\,dx\right|
\le C(\log N)^C\bigl(T^{-1/4}+T^{-2/5}+T^{-3/5}\bigr),
\]
with overwhelming probability, and the right side tends to zero.
\end{proof}

\section{Proof of Proposition \ref{prop:qbdtfix}}
\label{sec:qbd}

Throughout this section, we adopt Assumption \ref{ass:model} and use the
notation of Definitions \ref{def:flandL}, \ref{def:quasi_intro}, and
\ref{mt}, as well as the characteristic flow
\(x_\tau^{(T)}\) from Section \ref{subsec:tc}. As specified in Section
\ref{notation}, for an index interval
$I=\llbracket a,b\rrbracket\subseteq\llbracket N_1,N_2\rrbracket$,
$\L^I(t)=\L^{[a,b]}(t)$ denotes the corresponding principal block of
$\L(t)$. In addition, for any such index interval and $s\geq0$, set
\begin{equation}\label{eqn:index-neighborhood}
I^{(s)} \coloneqq \{j \in \llbracket N_1, N_2 \rrbracket:
\operatorname{dist}(j,I)\leq s\}.
\end{equation}
The constants \(c_1,C_1>0\) are those from Definition \ref{mt}; they will
be chosen in Section \ref{subsec:qbd}. 

We first prove the improved spacing
bounds in the bulk, then extend them to the edge using the coupling
from Section \ref{subsec:coupl}, and then finally prove Proposition
\ref{prop:qbdtfix}.

\subsection{Improving \eqref{eqn:qbd} for $i, j$ in the bulk}
\label{subsec:qbd}

We first state two technical preliminary lemmas that facilitate the proof of
Proposition \ref{prop:qbdtfix}. The first is a spacing estimate 
for the reordered collection of particles. Its proof is given in
Section \ref{subsec:qbd0} below.

\begin{lem}
\label{lem:ordered-cloud-spacing}
Fix \(t=T\tau\in[0,T \log N]\).
Let
\[
y_{N_1}(t) \leq y_{N_1+1}(t)\leq \cdots \leq y_{N_2}(t)
\]
denote the increasing rearrangement of the multiset
$\{Q_k(t): k\in\llbracket 1,N\rrbracket\}
=
\{q_i(t): i\in\llbracket N_1,N_2\rrbracket\}$.
There are constants \(c_*,C_*\), independent of
\(c_1,C_1\), and \(\mathfrak A\), such that, with overwhelming probability, the following holds.
Let $i_0 =  - \lfloor T (\log N)^{16} \rfloor$, $i_1 = \lfloor   T (\log N)^{16} \rfloor$, 
and set
$M_*:=\frac{\mathfrak M}{2R}$.
If \(\mathfrak B(t)\) holds, then for every
\(u,v\in\llbracket N_1,N_2\rrbracket\) with \(u>v\),
as long as $y_u(t), y_v(t) \in [q_{i_0}(t) , q_{i_1}(t)]$, one has
\begin{equation}
\label{eqn:ordered-cloud-spacing}
\begin{aligned}
y_u(t)-y_v(t) &\geq c_*\mathfrak A(u-v), \qquad u-v\ge \frac{1}{4} M_*\\
y_u(t)-y_v(t) &\le
C_*\mathfrak A \max(u-v, M_*).
\end{aligned}
\end{equation}
\end{lem}

The next is an estimate bounding the local 
energy on intervals of mesoscopic size. Its proof is given in Section
\ref{subsec:local_energy_upper}.

\begin{lem}
\label{lem:local-lax-trace}
Fix \(t\in[0,T\log N]\). There is a constant \(C>0\), independent of \(\mathfrak A\), such that,
with overwhelming probability, on the event \(\mathfrak B(t)\), the
following holds. Set $M_*:=(2R)^{-1}\mathfrak M$. For every
$1\le L\le 30 R^{-1}\mathfrak M$ and every index interval \(I\) satisfying
$|I|=L$ and
$I^{(5M_*)}\subseteq
\llbracket N_1+\frac32T(\log N)^{10},
N_2-\frac32T(\log N)^{10}\rrbracket$ (where
$I^{(5M_*)}$ is defined by \eqref{eqn:index-neighborhood}),
one has
\[
\operatorname{Tr}\bigl((\L^{I}(t))^2\bigr)
\le
C\max\{L,M_*\}.
\]
\end{lem}

 \begin{lem}\label{lem:qbdlow}
Fix \(t\in[0,T\log N]\).
Suppose the constant \(c_1>0\) from \eqref{eqn:qbd-lower} is small enough
(independently of \(\mathfrak A\)).
On an event of overwhelming probability, if \(\mathfrak B(t)\) holds, then for all
$
i,j\in
\llbracket -T(\log N)^{15},\,T(\log N)^{15}\rrbracket$
with $i-j\ge \frac{\mathfrak M}{2R}$,
\begin{equation}\label{eqn:qlb}
q_i(t)-q_j(t)\ge 2c_1\mathfrak A(i-j).
\end{equation}
\end{lem}

\begin{proof}
 Set
$
M_*:=\frac{1}{2R}\mathfrak{M}$.
Assume $\mathfrak B(t)$ holds.  Thus Lemmas
\ref{lem:ordered-cloud-spacing} and \ref{lem:local-lax-trace} are
applicable at the fixed time $t$.
We first record the two inputs used below.

First, by Lemma \ref{lem:ordered-cloud-spacing}, 
with overwhelming probability the ordered particle positions
$y_{N_1}(t)\leq y_{N_1+1}(t)\leq \cdots \leq y_{N_2}(t)$
satisfy
\begin{equation}
\label{eqn:qbd-ordered-spacing}
y_u(t)-y_v(t)\ge c_*\mathfrak{A}(u-v)
\end{equation}
whenever $y_v, y_u \in [ q_{- \lfloor T (\log N)^{16} \rfloor}(t), q_{ \lfloor T (\log N)^{16} \rfloor}(t) ]$ and 
$u-v\ge \frac{1}{4} M_*$,
for a fixed constant $c_*>0$. We restrict to the event
on which this holds, to
$\cap_{r \geq 0} \mathsf{BND}_{\L(r)}(\log N)$, which ensures that
$|q_\ell'(s)|\leq\log N$ for every $s\in[0,t]$, and to the event in
\eqref{eq:stretch-good}. Note that with $i,j$ as in the statement,
$q_{\ell}(t) \in [q_{-\lfloor   T (\log N)^{16} \rfloor}(t), q_{\lfloor   T (\log N)^{16} \rfloor}(t)]$ 
for each $\ell \in \llbracket j, i \rrbracket$ by these two restrictions.

Second, it suffices to establish the lemma when
$L:=i-j\in[M_*,2M_*]$, so we assume this below.
Lemma \ref{lem:local-lax-trace} then implies that there exists
a constant $C_2>0$, independent of $\mathfrak A$, such that for every index
interval
$I\subseteq \llbracket N_1 + T(\log N)^{11},N_2-T (\log N)^{11}\rrbracket$
with $|I|\le L$, we have
\begin{equation}
\label{eqn:qbd-local-energy}
\sum_{j\in I} a_j(t)^2\le C_2 L.
\end{equation}

We now choose the constant $c_1>0$ in Definition \ref{mt} so small, and $\mathfrak{A}$ so large, that
\begin{equation}
\label{eqn:qbd-c1choice}
2c_1\le \frac{c_*}{8},
\qquad
\frac{c_*\mathfrak{A}}{4}\ge 1,
\qquad
e^{c_*\mathfrak{A}/4}>C_2.
\end{equation}
We claim that this choice yields the desired lower bound in \eqref{eqn:qlb}.

Fix $i>j$ as above and, seeking a contradiction, assume that
\begin{equation}
\label{eqn:qbd-contradiction-assumption}
q_i(t)-q_j(t)<2c_1\mathfrak{A}L.
\end{equation}
Let $\pi_t$ denote the permutation defined by
$q_\ell(t)=y_{\pi_t(\ell)}(t)$ for
$\ell\in \llbracket N_1,N_2\rrbracket$, and write
$p_s:=\pi_t(j+s)$ for $0\le s\le L$.

We claim that
\begin{equation}
\label{eqn:qbd-small-rank-gap}
p_L-p_0<\frac{L}{4}.
\end{equation}
If
$p_L-p_0\ge \frac{1}{4} M_*$,
then \eqref{eqn:qbd-ordered-spacing} gives (using the fact that $q_i(t), q_j(t) \in [q_{-\lfloor   T (\log N)^{16} \rfloor}(t), q_{\lfloor   T (\log N)^{16} \rfloor}(t)]$ 
on the event we are working on)
\[
q_i(t)-q_j(t)
=
y_{p_L}(t)-y_{p_0}(t)
\ge
c_*\mathfrak{A}(p_L-p_0).
\]
Together with \eqref{eqn:qbd-contradiction-assumption}, this yields
$p_L-p_0<\frac{2c_1}{c_*}L\le \frac{L}{4}$,
by \eqref{eqn:qbd-c1choice}. If instead
$p_L-p_0< \frac{1}{4} M_*$,
then, since $L\ge M_*$,
$p_L-p_0<\frac{L}{4}$.

Since the integers $p_0,p_1,\ldots,p_L$ are distinct, we have
$\max_{0\le s\le L} p_s-\min_{0\le s\le L} p_s\ge L$.
Combining this with \eqref{eqn:qbd-small-rank-gap} yields
\[
\bigl(p_0-\min_{0\le s\le L}p_s\bigr)
+
\bigl(\max_{0\le s\le L}p_s-p_L\bigr)
\ge
L-(p_L-p_0)
>
\frac{3L}{4}.
\]
Hence at least one of the two quantities on the left is at least $3L/8$.
Let us assume that
$p_0-\min_{0\le s\le L}p_s\ge 3L/8$, as the alternative case is addressed very similarly.
Choose $s_*\in \llbracket 0,L\rrbracket$ so that
$p_{s_*}=\min_{0\le s\le L}p_s$.
Since $L\ge M_*$,
$\frac{3L}{8}\ge \frac{1}{4} M_*$,
so \eqref{eqn:qbd-ordered-spacing} gives
\[
q_j(t)-q_{j+s_*}(t)
=
y_{p_0}(t)-y_{p_{s_*}}(t)
\ge
\frac{3c_*}{8}\mathfrak{A}L.
\]
Thus, in either case we have found an interval
$J\subseteq \llbracket j,i-1\rrbracket$, so in particular $|J| \leq L$, such that
\begin{equation}
\label{eqn:qbd-PVJ}
\sum_{m\in J} \bigl(q_m(t)-q_{m+1}(t)\bigr) \ge \frac{3c_*}{8}\mathfrak{A}L.
\end{equation}
Jensen's inequality gives
\[
\sum_{m\in J} a_m(t)^2
=
\sum_{m\in J} e^{q_m(t)-q_{m+1}(t)}
\ge
|J| \exp\left( \frac{1}{|J|}\sum_{m\in J} \bigl(q_m(t)-q_{m+1}(t)\bigr)\right) \geq |J| \exp\left( \frac{1}{|J|} \frac{3c_*}{8}\mathfrak{A}L \right) .
\]
Since
$ \frac{c_*\mathfrak{A}}{4}\ge 1$
by \eqref{eqn:qbd-c1choice}, the function $t\mapsto te^{\frac{3c_*}{8}\mathfrak{A}L /t}$ is decreasing on
$[1, L ]$. Therefore
\[
\sum_{m = j}^{i-1} a_m(t)^2
\ge
L \exp\left(  \frac{3c_*}{8}\mathfrak{A} \right) .
\]
By \eqref{eqn:qbd-c1choice}, this is strictly larger than $C_2 L$, contradicting
\eqref{eqn:qbd-local-energy}. This contradiction proves \eqref{eqn:qlb}.
\end{proof}

\begin{lem}
\label{lem:qbd-short-gap-upper}
Fix \(t\in[0,T \log N]\), and set
$M_*:=\frac{\mathfrak M}{2R}$.
Suppose the constant \(C_{1}>0\) from \eqref{eqn:qbd} is large enough (independently of $\mathfrak{A}$).
On an event of overwhelming probability, if $\mathfrak B(t)$ holds, then for all
$i,j\in
\llbracket -T(\log N)^{15},\,T(\log N)^{15}\rrbracket$
with $i > j$,
one has
\begin{equation}\label{eqn:qiqjup}
q_i(t)-q_j(t)\le \frac{1}{2}C_{1}\mathfrak A \max(\frac{\mathfrak M}{2R}, i-j).
\end{equation}
\end{lem}

\begin{proof}
If \(q_i(t)\le q_j(t)\), there is nothing to prove, so assume
\(q_i(t)>q_j(t)\). Set $L:=i-j$, let $\pi_t$ be the permutation from the
proof of Lemma \ref{lem:qbdlow}, and set $p_s:=\pi_t(j+s)$ for
$0\le s\le L$. Then $p_L>p_0$, and by
\eqref{eqn:ordered-cloud-spacing},
\[
q_i(t)-q_j(t)
=
y_{p_L}(t)-y_{p_0}(t)
\le
C_*\mathfrak A
\max\{p_L-p_0,M_*\}.
\]

Let
$B^-:=\{\ell<j:\ p_0<\pi_t(\ell)<p_L\},
B^+:=\{\ell>i:\ p_0<\pi_t(\ell)<p_L\}$.
If \(\ell\in B^-\), then \(q_\ell(t)>q_j(t)\). Hence the lower bound
\eqref{eqn:qbd-lower} implies that \(j-\ell\le 2 M_*\) (if this 
was not the case, then certainly we would have $q_\ell(t)<q_j(t)$); therefore
\(|B^-|\le 2 M_*\). Similarly, \(|B^+|\le 2 M_*\). Since every $\pi_t(\ell)$ strictly between \(p_0\) and \(p_L\) has
 either a label $\ell$ in \(\{j+1,\ldots,i-1\}\) or a label $\ell$
in \(B^-\cup B^+\),
\[
p_L-p_0
\le
L+4M_*.
\]
Combining this with the first display in the proof and choosing \(C_1\)
sufficiently large in terms of \(C_*\) gives
$q_i(t)-q_j(t)
\le
\frac{C_1}{2}\mathfrak A\max\{L,M_*\}$,
as desired.
\end{proof}

\subsection{Improving \eqref{eqn:qbd} when one of $i, j$ is near the edge}
\label{subsec:qbde}
\label{subsec:qbd1N}

The preceding lemmas establish the improved spacing bounds in a central
region. The next lemma extends them to the tails by coupling the
finite system to thermal equilibrium.

\begin{lem}
\label{lem:qbdlow-edge}
Fix \(t\in[0,T\log N]\), and suppose the constant \(c_1>0\) from
\eqref{eqn:qbd-lower} and the constant \(C_1>0\) from \eqref{eqn:qbd} are
sufficiently small and large, respectively, both independently of
$\mathfrak A$. On an event of overwhelming probability,
if \(\mathfrak B(t)\) holds, then the following holds. Let
$i,j\in
\llbracket N_1+T(\log N)^{4},\,N_2-T(\log N)^{4}\rrbracket$
with \(i>j\), and assume that at least one of \(i,j\) lies outside
$\llbracket -T(\log N)^{14},\,T(\log N)^{14}\rrbracket $.
Then
\begin{equation}\label{eqn:stretchtimet}
q_i(t)-q_j(t)
\le
\frac{C_1}{2}\mathfrak A
\max\left\{
\frac{\mathfrak M}{2R},\,i-j
\right\}.
\end{equation}
If in addition
$i-j\ge \frac{\mathfrak M}{2R}$,
then
\begin{equation}\label{eqn:stretchtimet-low}
q_i(t)-q_j(t)
\ge
2c_1\mathfrak A(i-j).
\end{equation}
\end{lem}

\begin{proof}
Assume $\mathfrak B(t)$ holds.
Recall from Section \ref{subsec:coupl} the coupled thermal equilibrium
evolution $\L^{\rm th}(s)$ and its Flaschka variables
$a_k^{\rm th}(s),b_k^{\rm th}(s)$. Let $q_k^{\rm th}(s)$ denote the
associated particle positions, as in \eqref{eqn:flv}, set
$\alpha_{\rm th}:=-2\mathbb E\log a_k^{\rm th}(0)$. By Lemma
\ref{lem:edge-thermal-coupling}, with overwhelming probability,
\begin{equation}\label{eqn:cplest}
\sup_{0\le s\le T\log N}
\sup_{\substack{k\in\llbracket N_1,N_2\rrbracket\\
|k|\ge2T(\log T)^5}}
\left(|a_k(s)-a_k^{\rm th}(s)|+|b_k(s)-b_k^{\rm th}(s)|\right)
\le e^{-T(\log N)^2}.
\end{equation}
Applying \cite[Lemma 3.11]{Agg25} with the parameter $T$ there equal to
$S:=T\log N$ here, with overwhelming probability,
\begin{equation}\label{eqn:ef311}
\left|q_u^{\rm th}(s)-q_v^{\rm th}(s)-\alpha_{\rm th}(u-v)\right|
\le |u-v|^{1/2}(\log N)^2,\qquad
u,v\in\llbracket N_1+T(\log N)^4,N_2-T(\log N)^4\rrbracket,
\end{equation}
simultaneously for $s\in[0,S]$.
The parameters chosen in Section \ref{subsec:coupl} ensure that
$c\mathfrak A\le\alpha_{\rm th}\le C\mathfrak A$.

If both \(i\) and \(j\) lie to the left of \(-T(\log N)^{14}\), or both lie
to the right of \(T(\log N)^{14}\), every index between them satisfies the
lower bound on $|k|$ in \eqref{eqn:cplest}. Applying \eqref{eqn:ef311} to
adjacent indices gives $a_k^{\rm th}(s)\ge e^{-C(\log N)^2}$ 
for $k \in \llbracket N_1+T(\log N)^4,N_2-T(\log N)^4-1\rrbracket$.
Consequently, \eqref{eqn:cplest} and \eqref{eqn:flv} give
$|(q_i(t)-q_j(t))-(q_i^{\rm th}(t)-q_j^{\rm th}(t))|=o(1)$, and hence
\[
q_i(t)-q_j(t)
=
\alpha_{\rm th}(i-j)
+O\bigl((i-j)^{1/2}(\log N)^2\bigr)+o(1).
\]
The display above implies
\eqref{eqn:stretchtimet} and, when \(i-j\ge\mathfrak M/(2R)\),
\eqref{eqn:stretchtimet-low}, after choosing \(c_1\) small and \(C_1\) large.

If both \(i\) and \(j\) lie in
\(\llbracket-T(\log N)^{15},T(\log N)^{15}\rrbracket\), the result follows
from Lemmas \ref{lem:qbdlow} and \ref{lem:qbd-short-gap-upper}. We may
therefore assume that $i>-T(\log N)^{14}$ and
$j<-T(\log N)^{15}$, as the alternative case is addressed very similarly.
Restrict further to the event in \eqref{eq:stretch-good} and to
$\cap_{s\ge0}\mathsf{BND}_{\L(s)}(\log N)$.
Then $i-j\ge \frac12T(\log N)^{15}$, and
\eqref{eq:stretch-good} and \eqref{eqn:alph_bd} give
\[
c\mathfrak A(i-j)
\le q_i(0)-q_j(0)
\le C\mathfrak A(i-j),
\]
while $\partial_tq_k=b_k$ and the boundedness event give
\(|q_k(t)-q_k(0)|\le T(\log N)^2\). The latter error is negligible on this
scale, proving both conclusions in this final case as well.
\end{proof}

Now we are in a position to prove Proposition \ref{prop:qbdtfix}.

\begin{proof}[Proof of Proposition \ref{prop:qbdtfix}]
On the intersection of the events in Lemmas \ref{lem:qbdlow},
\ref{lem:qbd-short-gap-upper}, and \ref{lem:qbdlow-edge}, apply Lemmas
\ref{lem:qbdlow} and \ref{lem:qbd-short-gap-upper} when both \(i\) and \(j\) lie
in \(\llbracket-T(\log N)^{14},T(\log N)^{14}\rrbracket\), and apply
Lemma \ref{lem:qbdlow-edge} otherwise.
\end{proof}

\subsection{Proof of Lemma \ref{lem:ordered-cloud-spacing}}
\label{subsec:qbd0}

\begin{proof}
The proof has two parts.  First we record upper and lower bounds for
the number of quasiparticles in a given interval, 
utilizing 
the property $\mathfrak{B}(t)$. Then we convert this estimate into
upper and lower bounds on spacing of $y_i$'s.

We claim that there are constants \(c_0,C_0>0\), independent of
$c_1,C_1$, and \(\mathfrak A\), such that, upon setting
$\theta:=(8C_0)^{-1}$, with overwhelming probability, for every interval
\(\mathcal I\subset [q_{i_0}(t) , q_{i_1}(t)]\) satisfying
$|\mathcal I|\ge \theta\mathfrak A M_*$,
one has
\begin{equation}
\label{eqn:Q-cloud-two-sided-count}
c_0\mathfrak A^{-1}|\mathcal I|
\le
\#\{k:Q_k(t)\in\mathcal I\}
\le
C_0\mathfrak A^{-1}|\mathcal I|.
\end{equation}
By a fine mesh argument and a union bound, it suffices to prove this for a single fixed interval $\mathcal I$,
which we now proceed to do.

Put
\[
\bar q_i \coloneqq q_i(0)/T, \qquad \Lambda_i:=\lambda_{\varphi_0^{-1}(i)}, \qquad
 Y_i(\tau):=x_\tau^{(T)}(\Lambda_i;\bar q_i),
\qquad 
I:=T^{-1}\mathcal I=[a,b].
\]
We restrict to the event in \eqref{eq:stretch-good}, the event
$\cap_{s\ge0}\mathsf{BND}_{\L(s)}(\log N)$ from Lemma
\ref{l0eigenvalues}, which in particular includes
$\mathsf{BND}_{\L(0)}(\log N)$, and the event in Proposition
\ref{prop:outside-strip-confinement}. By \eqref{eq:stretch-good} and the bound
$|q_\ell(t)-q_\ell(0)|\le T(\log N)^2$ (which holds on $\cap_{s\ge0}\mathsf{BND}_{\L(s)}(\log N)$), any
 index $\ell$ for which
$q_\ell(t)\in[q_{i_0}(t),q_{i_1}(t)]$ satisfies
$|\ell|\le2T(\log N)^{16}$. Hence $Q_k(t)\in\mathcal I$ implies
$|\varphi_t(k)|\le2T(\log N)^{16}$. Similarly, \eqref{eqn:vebd} gives
$|Y_{\varphi_0(k)}(\tau)-\bar q_{\varphi_0(k)}|
\le C\tau(1+|\lambda_k|)$, so if
$Y_{\varphi_0(k)}(\tau)\in I$ and $|\lambda_k|\le\log N$
(which holds on $\mathsf{BND}_{\L(0)}(\log N)$), then
\eqref{eq:stretch-good} gives
$|\varphi_0(k)|\le2T(\log N)^{16}$. When $|\varphi_0(k)|\le2T(\log N)^{16}$, Proposition
\ref{prop:outside-strip-confinement}, applied at
$m=\varphi_t(k)$ and $s'=0$, implies $\varphi_t(k) \in \llbracket N_1 + T (\log N)^{10}, N_2 - T (\log N)^{10} \rrbracket$,
as
required in Definition \ref{mt}. Thus, since $\mathfrak B(t)$ holds, property \ref{item:T3} of that
definition together with the discussion above gives
\begin{equation}\label{eqn:bsbd}
\left|T^{-1}Q_k(t)-Y_{\varphi_0(k)}(\tau)\right|
\le
T^{-1}\mathfrak M\log(|\lambda_k|+2)
\end{equation}
for each $k$ such that $Q_k(t) \in  \mathcal{I}$ or $Y_{\varphi_0(k)}(\tau) \in I$.

\smallskip
\emph{Step 1: Upper bound in \eqref{eqn:Q-cloud-two-sided-count}.}
For $n\ge0$, set $J_n=[n,n+1]\cup[-n-1,-n]$.
If \(Q_k(t)\in\mathcal I\) and \(\lambda_k\in J_n\), then
$Y_{\varphi_0(k)}(\tau)\in
[a-(n+2)\mathfrak M/T,b+(n+2)\mathfrak M/T]$ by
\eqref{eqn:bsbd}. Denote the latter interval by $I_n$. Choose smooth
functions $\chi_n,\psi_n:\mathbb R\to[0,1]$ such that $\chi_n=1$ on
$J_n$, $\operatorname{supp}\chi_n\subset
\{\lambda:\operatorname{dist}(\lambda,J_n)\le1/4\}$, and
$\|\chi_n'\|_\infty\le C$, while $\psi_n=1$ on $I_n$,
$\operatorname{supp}\psi_n\subset I_n+[-\mathfrak M/T,\mathfrak M/T]$, and
$\|\psi_n'\|_\infty\le CT/\mathfrak M$. Define
$h_n(\lambda,q):=\chi_n(\lambda)
\psi_n(x_\tau^{(T)}(\lambda;q))$.

The change of variables $x=x_\tau^{(T)}(\lambda;q)$
(which is valid since $q\mapsto x_\tau^{(T)}(\lambda;q)$ is invertible
by \eqref{eqn:Xtlin} and \eqref{eqn:avrbd2coarse}) gives the first bound
below, while \eqref{eq:flow-spec} gives the second:
\[
\|h_n\|_\infty+\sup_{\lambda\in\mathbb R}\int_{\mathbb R}
|\partial_qh_n(\lambda,q)|\,dq\le C,
\qquad
\|h_n\|_{\operatorname{Lip},\lambda}
\le C\left(1+\frac{\tau(n+2)T}{\mathfrak M}\right).
\]
For $n \leq \log N$, the support condition 
required for \eqref{eqn:hiTspace}, namely
$\supp h_n \subset \mathbb{R} \times [-(\log N)^{20}, (\log N)^{20}]$,
 follows from the assumption
 $\mathcal{I} \subset [q_{i_0}(t), q_{i_1}(t)]$,
\eqref{eq:stretch-good}, $|q_\ell(t)-q_\ell(0)|\le T(\log N)^2$, and
the first bound in \eqref{eqn:vebd}. Since $h_n=1$ on the set where 
$(\lambda, x_\tau^{(T)}(\lambda;q)) \in J_n \times I_n$,
\[
\#\{i:Y_i(\tau)\in I_n,\ \Lambda_i\in J_n\}
\le \sum_{i=N_1}^{N_2}h_n(\Lambda_i,\bar q_i).
\]
Applying \eqref{eqn:hiTspace} to the sum on the right, then using
\eqref{eqn:m0push} to push the integral forward and \eqref{eqn:avrbd} to
bound the resulting density $m_{\tau}^{(T)}$, gives
\begin{align*}
\#\{i:Y_i(\tau)\in I_n,\ \Lambda_i\in J_n\}
&\le C\mathfrak A^{-1}T
\left(|I|+\frac{(n+2)\mathfrak M}{T}\right)e^{-cn^2}\\
&\quad+(\log N)^C
\left(T^{3/5}+\frac{(n+2)T^{7/5}}{\mathfrak M}\right).
\end{align*}
On $\mathsf{BND}_{\L(0)}(\log N)$ the terms with
$n>\lceil\log N\rceil$ vanish. Summing over the remaining $n$ and using
$T|I|\ge\theta\mathfrak A\mathfrak M/(2R)$ gives, after choosing $C_0$
sufficiently large and then taking $\mathfrak A$ sufficiently large,
\[
\#\{k:Q_k(t)\in\mathcal I\}
\le
C_0\mathfrak A^{-1}T|I|.
\]
Since \(T|I|=|\mathcal I|\), this proves the upper bound in
\eqref{eqn:Q-cloud-two-sided-count}.

\smallskip
\emph{Step 2: Lower bound in \eqref{eqn:Q-cloud-two-sided-count}.}
Let $J_\star$ be the interval in Lemma \ref{lem:uniform-non-vacuum}, and
choose a compact interval $J_0$ such that $J_\star$ is 
contained in the interior of $J_0$. Write
$\ell:=|I|$. Recall that
$\ell\ge \theta\mathfrak A\frac{M_*}{T}
=
\frac{\theta\mathfrak A}{2R}\frac{\mathfrak M}{T}$.
Choose a smooth function $\chi_0$ with compact support contained in $J_0$
and satisfying
$0\le\chi_0\le1$ and
$\chi_0\equiv1$ on $J_\star$.
All constants $C>0$ below are allowed to depend on the fixed interval \(J_0\), 
and may change from line to line.

Set
$d_{\rm dyn}
:=
\frac{\mathfrak M}{T}
\max_{\lambda\in J_0}\log(|\lambda|+2)$.
Thus, for $\lambda_k\in J_0$, the right hand side of
\eqref{eqn:bsbd} is at most $d_{\rm dyn}$. Define
$I^\circ:=[a+2d_{\rm dyn},b-2d_{\rm dyn}]$.
Since $d_{\rm dyn}/\ell\le CR/(\theta\mathfrak A)$, we may, by first
choosing \(\mathfrak A\) sufficiently large, assume that
$|I^\circ|\ge\frac{9}{10}\ell$.

Set $\eta:=\ell/100$ and write $I^\circ=[u,v]$. Let \(\vartheta\) be a
smooth function satisfying $0\le\vartheta\le1$, with
$\vartheta(s)=0$ for $s\le0$ and $\vartheta(s)=1$ for $s\ge1$, and define
\[
\psi_{I^\circ,\eta}(x)
:=
\vartheta\left(\frac{x-u}{\eta}\right)
\vartheta\left(\frac{v-x}{\eta}\right).
\]
Then
\begin{equation}
\label{eq:smoothed-indicator-properties}
0\le\psi_{I^\circ,\eta}\le\mathbf 1_{I^\circ},
\qquad
\psi_{I^\circ,\eta}\equiv1
\quad\text{on }[u+\eta,v-\eta],
\end{equation}
and
\begin{equation}
\label{eq:smoothed-indicator-derivatives}
\|\psi_{I^\circ,\eta}'\|_\infty\le C\eta^{-1},
\qquad
\int_{\mathbb R}|\psi_{I^\circ,\eta}'(x)|\,dx\le C.
\end{equation}
Define
\begin{equation}
\label{eq:pulled-back-smoothed-indicator}
h_{I,\tau}(\lambda,q)
:=
\chi_0(\lambda)\,
\psi_{I^\circ,\eta}
\left(x_\tau^{(T)}(\lambda;q)\right).
\end{equation}
By the support properties of $\chi_0$ and $\psi_{I^\circ,\eta}$ and
\eqref{eqn:bsbd} (whose right hand side is upper bounded by $ d_{\rm dyn}$
by definition),
\begin{equation}\label{eqn:hupP}
\sum_{i=N_1}^{N_2}
h_{I,\tau}(\Lambda_i,\bar q_i)
\leq \left|\{i:\Lambda_i\in J_0,\ Y_i(\tau)\in I^\circ\}\right|
\leq \left| \{ k : Q_k(t) \in  \mathcal{I}  \}\right| .
\end{equation}

We verify the hypotheses of Assumption \ref{ass:admissible-test} in order
to apply Proposition \ref{prop:app-time-zero-dos}. First,
$0\le h_{I,\tau}\le1$, so $\|h_{I,\tau}\|_\infty\le1$. Differentiating in
the spatial coordinate gives
\[
\partial_qh_{I,\tau}(\lambda,q)
=
\chi_0(\lambda)
\psi_{I^\circ,\eta}'
\left(x_\tau^{(T)}(\lambda;q)\right)
\partial_qx_\tau^{(T)}(\lambda;q).
\]
For each fixed \(\lambda\), the map
$q\longmapsto x_\tau^{(T)}(\lambda;q)$
is strictly increasing and Lipschitz by \eqref{eqn:Xtlin} and
\eqref{eqn:avrbd2coarse}. Therefore the change of variables
$x=x_\tau^{(T)}(\lambda;q)$ and
\eqref{eq:smoothed-indicator-derivatives} give
\begin{equation}
\int_{\mathbb R}
|\partial_qh_{I,\tau}(\lambda,q)|\,dq
=
|\chi_0(\lambda)|
\int_{\mathbb R}
|\psi_{I^\circ,\eta}'(x)|\,dx
\le C.
\label{eq:pulled-test-qvariation}
\end{equation}
Consequently, $\sup_{\lambda\in\mathbb R}\int_{\mathbb R}
|\partial_qh_{I,\tau}(\lambda,q)|\,dq\le C$, so the constant \(A\) in
Assumption \ref{ass:admissible-test} may be taken to be $C$.

We next estimate the Lipschitz constant in \(\lambda\). On the 
compact interval \(J_0\), the characteristic flow estimate \eqref{eq:flow-spec} of Lemma \ref{lem:backward-flow-rectangle}
implies
$\sup_{\lambda\in J_0, q\in\mathbb R}
|\partial_\lambda x_\tau^{(T)}(\lambda;q)|
\le C \tau$. It follows that
\[
|\partial_\lambda h_{I,\tau}(\lambda,q)|
\le \|\chi_0'\|_\infty
+\|\psi_{I^\circ,\eta}'\|_\infty
|\partial_\lambda x_\tau^{(T)}(\lambda;q)|
\le C(1+\tau/\eta).
\]
Therefore, for some larger $C > 0$,
\begin{equation}
\label{eq:pulled-test-lambda-lip}
\|h_{I,\tau}\|_{\Lip,\lambda}
\le
C
\left(
1+\frac{\tau}{\ell}
\right).
\end{equation}

Finally, if \(h_{I,\tau}(\lambda,q)\ne0\), then \(\lambda\in J_0\) and
$x_\tau^{(T)}(\lambda;q)\in I^\circ\subset I$.
By the hypothesis that
\(T I\subset [q_{i_0}(t),q_{i_1}(t)]\), \eqref{eq:stretch-good}, and our
restriction to $\cap_{s\ge0}\mathsf{BND}_{\L(s)}(\log N)$, the
interval \(I\) lies inside
$[-C\mathfrak A(\log N)^{16},
 C\mathfrak A(\log N)^{16}]$.
Since \(J_0\) is fixed and \(\tau\le\log N\), \eqref{eqn:vebd} 
implies that if
$(\lambda,q)\in\operatorname{supp}h_{I,\tau}$, then
$|q|\le(\log N)^{20}$.
Thus $\operatorname{supp}h_{I,\tau}\subseteq
\mathbb R\times[-(\log N)^{20},(\log N)^{20}]$. In particular, 
 the
 support condition required for \eqref{eqn:hiTspace} holds,
and \(h_{I,\tau}\) satisfies Assumption
\ref{ass:admissible-test} with $A\le C$ and
$B\le C(1+\tau/|I|)$.

Thus, \eqref{eqn:hiTspace} from Proposition \ref{prop:app-time-zero-dos} implies that
\[
\begin{aligned}
\sum_{i=N_1}^{N_2}
h_{I,\tau}(\Lambda_i,\bar q_i)
\ge
T\iint_{\mathbb{R}^2}
h_{I,\tau}(\lambda,q)m_0^{(T)}(\lambda,q)
\,d\lambda\,dq
-
(\log N)^C
\left[
T^{3/5}
+
T^{2/5}
\left(
1+\frac{\tau}{|I|}
\right)
\right] .
\end{aligned}
\]
By \eqref{eqn:m0push} applied to $m_0^{(T)}$ and $m_\tau^{(T)}$,
Lemma \ref{lem:uniform-non-vacuum}, and the fact that $|I^\circ|\ge\frac{9}{10}\ell$,
\[
T \iint_{\mathbb{R}^2}
h_{I,\tau}(\lambda,q)m_0^{(T)}(\lambda,q)
\,d\lambda\,dq
=T \iint_{\mathbb{R}^2}
\chi_0(\lambda) \psi_{I^{\circ}, \eta}(x) m_{\tau}^{(T)}(\lambda,x)
\,d\lambda\,dx
\geq c\mathfrak{A}^{-1} |\mathcal{I}|
\geq \frac{c \theta}{R} \mathfrak M.
\]
whereas $C(\log N)^C (T^{3/5} + T^{2/5}\frac{T}{\mathfrak M}) < \frac{1}{100} \frac{c \theta}{R} \mathfrak M$ for $N$ large enough 
(here we use our assumption that $\delta < 1/200$ in Definition \ref{mt}) , so that 
\begin{equation}\label{eqn:hsumlwbd}
\sum_{i=N_1}^{N_2}
h_{I,\tau}(\Lambda_i,\bar q_i) \geq c\mathfrak{A}^{-1} |\mathcal{I}| 
\end{equation}
after decreasing the constant $c$. Now the lower bound in \eqref{eqn:Q-cloud-two-sided-count}  follows from \eqref{eqn:hsumlwbd} and \eqref{eqn:hupP}. 

\smallskip
\emph{Step 3: Lower bound in \eqref{eqn:ordered-cloud-spacing}.}
Set $D:=y_u(t)-y_v(t)$ and $L:=u-v$.

We first prove the lower bound. If
$D\ge \theta\mathfrak A M_*$,
then the interval \([y_v(t),y_u(t)]\) contains at least \(L+1\)
quasiparticles, so \eqref{eqn:Q-cloud-two-sided-count} gives
$L+1\le C_0\mathfrak A^{-1}D$. Hence
$D\ge C_0^{-1}\mathfrak A L$.

Suppose instead that \(D<\theta\mathfrak A M_*\). Enlarge
\([y_v(t),y_u(t)]\), within the larger bulk interval $[q_{i_0}(t) , q_{i_1}(t)]$, to an interval of length
\(\theta\mathfrak A M_*\). The larger interval still contains at least \(L+1\) points,
whereas the upper bound in \eqref{eqn:Q-cloud-two-sided-count} implies that it has at most \(C_0\theta M_*\) points.
Since \(L\ge M_*/4\) and $C_0\theta=1/8$, this is impossible.
Thus
$D\ge c_*\mathfrak A L$, with $c_* = C_0^{-1}$.

\smallskip
\emph{Step 4: Upper bound in \eqref{eqn:ordered-cloud-spacing}.}
If
$D\le 2\theta\mathfrak A M_*$,
then $D\le2\theta\mathfrak A\max\{L,M_*\}$, which gives the desired
upper bound after increasing $C_*$.
Otherwise choose an interval
$\mathcal I'\subset(y_v(t),y_u(t))$ of length \(D/2\). It contains at
most \(L-1\) quasiparticles and has length at least
\(\theta\mathfrak A M_*\). Therefore the
lower bound in \eqref{eqn:Q-cloud-two-sided-count} gives
$c_0\mathfrak A^{-1}D/2\le L-1$, and hence
$D\le 2c_0^{-1}\mathfrak A L$.
This proves the claim.

%
\end{proof}

\section{Localization estimates via eigenvalue matching}
	\label{sec:IStm2} 			

Throughout this section, we adopt Assumption \ref{ass:model} and the notation
of Definition \ref{def:quasi_intro}. Throughout this section, we also recall
$\delta\in(0,1/200)$,
$\mathfrak M\in[T^{3/4-\delta},T^{3/4}]$, $R=\mathfrak A^{1/10}$, and
$\mathfrak T=\mathfrak T_{\mathfrak M}$ from Definition \ref{mt}.
For $m\in\llbracket N_1+1,N_2-1\rrbracket$, set
\[
\H_-^m(t):=\L^{[N_1,m-1]}(t),
\qquad
\H_+^m(t):=\L^{[m+1,N_2]}(t).
\]
When $m=\varphi_t(k)$ in an argument where $t$ and $k \in \llbracket 1, N \rrbracket$ have been fixed, we
suppress the superscript $m$ from the above notation.

\subsection{Proof of Lemma \ref{lem:cutLbd}}
	\label{subsec:IStm2} 		

To motivate the construction below, fix $t \in [0, T \log N]$ and 
$m \in \llbracket N_1+2T(\log N)^{10},N_2-2T(\log N)^{10}\rrbracket$
(as in Lemma \ref{lem:AUpi} below). We will match the
eigenvalues of $\L(t)$ whose localization centers lie sufficiently far to the left of
$m$ with eigenvalues of $\L^{[N_1,m-1]}(t)$. Then to prove
 Lemma \ref{lem:cutLbd}, we will set $m = \varphi_t(k)$ 
 and argue that we may replace the
ordering of the centers by the ordering of the quasiparticles.

	\begin{lem}\label{lem:AUpi}
	There exists a constant $c > 0$ such that, for each fixed $t \in [0, T \log N]$, the following 
	holds with overwhelming probability. For 
	$m\in\llbracket N_1+2T(\log N)^{10},N_2-2T(\log N)^{10}\rrbracket$, let
		$\mathcal{A}_-^m \coloneqq  \llbracket N_1, m- 10 R^{-1} \mathfrak{M} \rrbracket $.
On the event $\mathfrak{B}(t)$, there exists an injective
map $\pi_-^m :  \mathcal{A}_-^m \rightarrow \eig \H_-^m(t) $ such that
\begin{equation}\label{eqn:materr}
	|\pi_-^m(i) - \Lambda_i(t)| \leq  e^{-c R^{-1} \mathfrak{M}}, \qquad \forall \text{ } i \in \mathcal{A}_-^m.
	\end{equation}
	\end{lem}
\begin{proof}
We restrict to the event of Proposition
\ref{prop:outside-strip-confinement} and the events
$\bigcap_{s\geq 0}\mathsf{BND}_{\L(s)}(\log N)$ and
$\mathsf{SEP}_{\L(0)}(e^{-(\log N)^2})$, which are overwhelmingly probable
by Lemmas \ref{l0eigenvalues} and \ref{eigenvalues0}. For each
$i\in\mathcal A_-^m$, we claim Assumption \ref{lmatrixl}
is satisfied with
$\mathcal I=\llbracket N_1,m-1\rrbracket$,
$\mathcal D=\llbracket m,N_2\rrbracket$, $\delta=0$, $\L=\L(t)$,
$\widetilde\L=\H_-^m(t)\oplus\bm 0$ (where $\bm 0$ acts on
$\llbracket m,N_2\rrbracket$), with eigenvalue $\lambda_j$ and localization
center $\varphi_j$ in that assumption given by $\Lambda_i(t)$ and
$i$ here, and with corresponding eigenvector given by
$\bm u_{\varphi_t^{-1}(i)}(t)$. In particular,
$\dist(i,\mathcal D)=m-i\ge10R^{-1}\mathfrak M>100$ for $N$ large enough. The required estimate
\eqref{eqn:lmatrix-boundary-decay} follows
from \eqref{eqn:uMR} when
$i\in\llbracket N_1+T(\log N)^{10},N_2-T(\log N)^{10}\rrbracket$, and
from \eqref{eqn:umphi} with $s=s'=t$, otherwise, since in the latter case, the
assumptions on $i$ and $m$ imply that $i\leq N_1+T(\log N)^{10}$ and
$m-i\geq T(\log N)^3$.

Lemma \ref{lem:full-to-truncation}, together with the definition of the
localization center bijection and the assumptions on $\mathfrak M$, therefore
associates to each $i\in\mathcal A_-^m$ an eigenvalue of $\H_-^m(t)$ satisfying
\eqref{eqn:materr}, after decreasing $c>0$ if necessary. These eigenvalues
are distinct by our restriction to
$\mathsf{SEP}_{\L(0)}(e^{-(\log N)^2})$, so this association defines the
required injection $\pi_-^m$.
\end{proof}

We next state three estimates used in the proof of Lemma \ref{lem:cutLbd};
in the proof, after using the previous lemma to match up most of the eigenvalues
of $\H_-$ to those of $\L(t)$ with
centers $\varphi' < \varphi$ which are not too close to $\varphi$, 
we must control the remaining errors contributing to the left hand side of \eqref{eqn:itoQ1}, 
and this is the purpose of the lemmas below.
Lemma \ref{lem:lsumbd}
controls sums of the form $\sum \log |\lambda_k - \lambda_j|$ and
$\sum \log |1 - \lambda_j/\lambda_k|$ over indices $j$ with
 $\varphi_t(j)$ ranging
 over short intervals of localization
centers. Lemma \ref{lem:lqebd} controls part of the contribution from 
the unmatched part of the spectrum of the
truncated matrix $\H_-$, and Lemma \ref{lem:logdivlambdak}
controls the remaining piece contributed by that unmatched part of the
spectrum. Lemmas \ref{lem:lsumbd}, \ref{lem:lqebd}, and
\ref{lem:logdivlambdak} are proved in Sections \ref{subsec:min},
\ref{subsec:unmatchedmu2pf}, and \ref{subsec:logsum_quadsum}, respectively.

\begin{lem}\label{lem:lsumbd}
Adopt Assumption \ref{ass:model} and fix $t\in[0,T\log N]$.
Let $r\ge\mathfrak A^{-1/10}$ and consider an index interval
$W\subseteq\llbracket N_1+T(\log N)^7,N_2-T(\log N)^7\rrbracket$ with
$|W|\le r\mathfrak M$. With overwhelming
probability, whenever $\mathfrak B(t)$ holds, the following two statements
hold.
\begin{enumerate}
\item For each $k\in\llbracket1,N\rrbracket$,
\begin{equation}\label{eqn:sumlambdak01_index}
\sum_{\substack{j\in\llbracket1,N\rrbracket:\,\varphi_t(j)\in W\\j\ne k}}
\bigl|\log|\lambda_k-\lambda_j|\bigr|
\le Cr\mathfrak M\log(|\lambda_k|+2).
\end{equation}
\item If $|\lambda_k|>10$,
\begin{equation}\label{eqn:sumellk}
\sum_{\substack{j\in\llbracket1,N\rrbracket:\,\varphi_t(j)\in W\\j\ne k}}
\left|\log|1-\frac{\lambda_j}{\lambda_k}|\right|
\le Cr\mathfrak M.
\end{equation}
\end{enumerate}
\end{lem}

\begin{lem}\label{lem:lqebd}
Adopt Assumption \ref{ass:model} and fix $t\in[0,T\log N]$. There is a
constant $C>0$, independent of $\mathfrak A$, such that the following holds
with overwhelming probability whenever $\mathfrak B(t)$ holds. Let
$m\in\llbracket N_1+2T(\log N)^{10},N_2-2T(\log N)^{10}\rrbracket$, let
$\mathcal A_-^m$ and $\pi_-^m$ be as in Lemma \ref{lem:AUpi}, and
define
$\mathcal U_-^m:=\eig\H_-^m(t)\setminus\pi_-^m(\mathcal A_-^m)$. Then
\begin{equation}\label{eqn:mu2U}
\sum_{\mu\in\mathcal U_-^m}\mu^2\le CR^{-1}\mathfrak M.
\end{equation}
\end{lem}

\begin{lem}\label{lem:logdivlambdak}
Adopt Assumption \ref{ass:model}, fix $t\in[0,T\log N]$.
For $k \in \llbracket 1, N \rrbracket$, let $\mathcal A_-, \pi_-, \mathcal U_-$ be 
defined as $\mathcal A_-^m, \pi_-^m, \mathcal U_-^m$ in 
Lemma \ref{lem:lqebd} with $m = \varphi_t(k)$.
There is a constant $C>0$, independent of $\mathfrak A$, such that, with
overwhelming probability, the following two statements hold whenever
$\mathfrak B(t)$ holds and $k\in\llbracket1,N\rrbracket$ satisfies
$\varphi=\varphi_t(k)\in
\llbracket N_1+2T(\log N)^{10},N_2-2T(\log N)^{10}\rrbracket$.
\begin{enumerate}
\item We have
\begin{equation}\label{eqn:less1}
\sum_{\substack{\mu\in\mathcal U_-\\|\mu-\lambda_k|<1}}
\left|\log|\lambda_k-\mu|\right|
\le CR^{-1}\mathfrak M\log(|\lambda_k|+2).
\end{equation}
\item If $|\lambda_k|>10$, then
\begin{equation}\label{enq:divlk1}
\sum_{\mu\in\mathcal U_-}
\left|\log|1-\frac{\mu}{\lambda_k}|\right|
\le CR^{-1}\mathfrak M.
\end{equation}
\end{enumerate}
\end{lem}

We now prove Lemma \ref{lem:cutLbd}, conditional on the lemmas above.

\begin{proof}[Proof of Lemma \ref{lem:cutLbd}]
Set $\varphi=\varphi_t(k)$. We restrict to the events that Lemmas
\ref{lem:Estr}, \ref{lem:AUpi}, \ref{lem:lsumbd}, \ref{lem:lqebd}, and
\ref{lem:logdivlambdak} hold, as well as
$\mathsf{SEP}_{\L(0)}(e^{-(\log N)^2})$, which holds with overwhelming
probability by Lemma \ref{eigenvalues0}. We use the notation of Lemma
\ref{lem:lqebd} for $\mathcal A_-^m$, $\pi_-^m$, and $\mathcal U_-^m$ with
$m=\varphi$, and suppress the superscript $m$. By
\eqref{eqn:materr} and eigenvalue separation,
\begin{equation}\label{eqn:match1}
\left|
\sum_{\mu\in\pi_-(\mathcal A_-)}\log|\lambda_k-\mu|
-
\sum_{j\in\mathcal A_-}\log|\lambda_k-\Lambda_j(t)|
\right|
\le e^{-(\log N)^{100}}.
\end{equation}

It remains to bound the unmatched contribution
\begin{equation}\label{eqn:Upm_bd}
\sum_{\mu\in\mathcal U_-}\left|\log|\lambda_k-\mu|\right|
\le CR^{-1}\mathfrak M\log(|\lambda_k|+2).
\end{equation}
Indeed, given \eqref{eqn:Upm_bd}, the bound
\eqref{eqn:sumlambdak01_index} from Lemma \ref{lem:lsumbd}, applied with
$r=21R^{-1}$, gives
\begin{equation}\label{eqn:windwbd}
\sum_{\substack{|j-\varphi|\le10R^{-1}\mathfrak M\\j\ne\varphi}}
\left|\log|\lambda_k-\Lambda_j(t)|\right|
\le CR^{-1}\mathfrak M\log(|\lambda_k|+2).
\end{equation}
Moreover, by \eqref{eqn:qbd-lower}, restriction to the boundedness event
$\cap_{s\ge0}\mathsf{BND}_{\L(s)}(\log N)$, and
\eqref{eq:stretch-good},
\[
\left|
\sum_{j:Q_j(t)<Q_k(t)}\log|\lambda_k-\lambda_j|
-
\sum_{\varphi_t(j)\in\mathcal A_-}\log|\lambda_k-\lambda_j|
\right|
\]
is bounded by the left-hand side of \eqref{eqn:windwbd}.
Thus \eqref{eqn:windwbd}, together with \eqref{eqn:match1} and \eqref{eqn:Upm_bd}, yields the
lemma.

To deduce \eqref{eqn:Upm_bd}, we use the bounds \eqref{eqn:mu2U} from Lemma
\ref{lem:lqebd} and \eqref{eqn:less1} from Lemma
\ref{lem:logdivlambdak}. Specifically, if
$|\mu-\lambda_k|>1$, then
\[
\left|\log|\lambda_k-\mu|\right|
\le C\bigl(\log(|\lambda_k|+2)+\mu^2\bigr),
\]
and $|\mathcal U_-|\le CR^{-1}\mathfrak M$ by construction. Thus
\eqref{eqn:mu2U} bounds the part of \eqref{eqn:Upm_bd} with
$|\mu-\lambda_k|>1$, while \eqref{eqn:less1} bounds the complementary
part.
\end{proof}

\subsection{Proof of Lemma \ref{lem:log_sum_dif}}
\label{subsec:ideas4}

\begin{proof}
Set $\varphi=\varphi_t(k)$. We first assume that
$\varphi\in\llbracket N_1+20T(\log N)^{10},
N_2-20T(\log N)^{10}\rrbracket$, and then briefly outline the more direct
argument for the complementary case at the end of the proof. Since
$|i-\varphi|<T(\log N)^3$, in the first case both $i$ and $\varphi$ lie in
$\llbracket N_1+2T(\log N)^{10},N_2-2T(\log N)^{10}\rrbracket$ for
sufficiently large $N$. Throughout the treatment of this case, we work on
the overwhelmingly probable events from Lemmas \ref{lem:AUpi},
\ref{lem:lsumbd}, \ref{lem:lqebd}, and \ref{lem:logdivlambdak}, and use the
notation $\mathcal A_-^m$, $\pi_-^m$, and $\mathcal U_-^m$ from Lemma
\ref{lem:lqebd} with $m\in\{i,\varphi\}$.

If $|\lambda_k|\le10$, further restriction to $\mathsf{SEP}_{\L(0)}(e^{-(\log N)^2})$ and \eqref{eqn:materr} give
\begin{align*}
&\sum_{\mu\in\eig\L^{[N_1,i-1]}(t)}\log|\lambda_k-\mu|
-\sum_{\mu\in\eig\L^{[N_1,\varphi-1]}(t)}\log|\lambda_k-\mu|\\
&\quad=
\sum_{j\in\mathcal A_-^i}\log|\lambda_k-\Lambda_j(t)|
-\sum_{j\in\mathcal A_-^\varphi}\log|\lambda_k-\Lambda_j(t)|\\
&\qquad\quad+
\sum_{\mu\in\mathcal U_-^i}\log|\lambda_k-\mu|
-\sum_{\mu\in\mathcal U_-^\varphi}\log|\lambda_k-\mu|+o(1)\\
&\quad\le
\sum_{i-10R^{-1}\mathfrak M\le\varphi_t(j)<\varphi}
\left|\log|\lambda_k-\lambda_j|\right|+CR^{-1}\mathfrak M \le C|\varphi-i|+CR^{-1}\mathfrak M
\le2C|\varphi-i|.
\end{align*}
For the first inequality, Lemma \ref{lem:lqebd}, applied with $m=i$, bounds
the positive part of the sum over $\mathcal U_-^i$; when
$|\lambda_k-\mu|>1$ we use
$\log|\lambda_k-\mu|\le C(1+\mu^2)$, while the remaining logarithms are
nonpositive. The bound \eqref{eqn:less1} from Lemma \ref{lem:logdivlambdak} bounds the
negative part of the sum over $\mathcal U_-^\varphi$. The last inequality
follows from the bound \eqref{eqn:sumlambdak01_index} from Lemma
\ref{lem:lsumbd} and
$R^{-1}\mathfrak M\le2|\varphi-i|$.

Suppose next that $|\lambda_k|>10$. Using
\(
\log|\lambda_k-x|=\log|\lambda_k|
+\log|1-x/\lambda_k|
\), we similarly obtain
\begin{align*}
&\sum_{\mu\in\eig\L^{[N_1,i-1]}(t)}\log|\lambda_k-\mu|
-\sum_{\mu\in\eig\L^{[N_1,\varphi-1]}(t)}\log|\lambda_k-\mu|\\
&\quad=-|\varphi-i|\log|\lambda_k|\\
&\qquad\quad+
\sum_{j\in\mathcal A_-^i}\log\left|1-\frac{\Lambda_j(t)}{\lambda_k}\right|
-\sum_{j\in\mathcal A_-^\varphi}\log\left|1-\frac{\Lambda_j(t)}{\lambda_k}\right|\\
&\qquad\quad+
\sum_{\mu\in\mathcal U_-^i}\log\left|1-\frac{\mu}{\lambda_k}\right|
-\sum_{\mu\in\mathcal U_-^\varphi}\log\left|1-\frac{\mu}{\lambda_k}\right|+o(1)\\
&\quad\le-|\varphi-i|\log|\lambda_k|
+\sum_{i-10R^{-1}\mathfrak M\le\varphi_t(j)<\varphi}
\left|\log\big|1-\frac{\lambda_j}{\lambda_k}\big|\right|
+CR^{-1}\mathfrak M
\le C|\varphi-i|.
\end{align*}
Here Lemma \ref{lem:lqebd} and $|\mathcal U_-^i|\le CR^{-1}\mathfrak M$
give
\[
\sum_{\mu\in\mathcal U_-^i}
\max\left\{\log\left|1-\frac{\mu}{\lambda_k}\right|,0\right\}
\le C\left(|\mathcal U_-^i|+\sum_{\mu\in\mathcal U_-^i}\mu^2\right)
\le CR^{-1}\mathfrak M.
\]
This bounds the positive part of the sum over $\mathcal U_-^i$. The bound
\eqref{enq:divlk1} from Lemma \ref{lem:logdivlambdak} bounds the sum over
$\mathcal U_-^\varphi$, and the bound \eqref{eqn:sumellk} from Lemma
\ref{lem:lsumbd}
bounds the sum $\sum_{i-10R^{-1}\mathfrak M\le\varphi_t(j)<\varphi} \left|\log\big|1-\frac{\lambda_j}{\lambda_k}\big|\right|$.

It remains to consider the case when $\varphi \notin \llbracket N_1 + 20T
(\log N)^{10}, N_2 - 20T (\log N)^{10} \rrbracket$. In this setting, the
proof follows from a comparison to thermal equilibrium, in which setting the
lemma follows from known results in \cite{Agg25a,Agg25}; so, we only briefly
outline the argument.

Since $\varphi \in \llbracket N_1 + T (\log N)^9, N_2 - T (\log N)^9
\rrbracket$ and $|i - \varphi| \le T (\log N)^3$, we have that $i$ and
$\varphi$ are either both in $\llbracket N_1 + T (\log N)^8, N_1 + 20T
(\log N)^{10} \rrbracket$ or both in $\llbracket N_2 - 20T (\log N)^{10},
N_2 - T (\log N)^8 \rrbracket$. Let us assume the former, as the proof in
the latter case is entirely analogous.

Recall the matrix $\mathbf{L}^{\mathrm{th}} (t)$ from Section
\ref{subsec:coupl}. By Lemma \ref{lem:edge-thermal-coupling},
\cite[Proposition 4.4]{Agg25a}, and \cite[Corollaries 5.5 and 5.6]{Agg25a}, it
is quickly verified that it suffices to prove the lemma with $\mathbf{L}(t)$
replaced by $\mathbf{L}^{\mathrm{th}} (0)$. Specifically, denoting $\mathbf{M}
= \mathbf{L}^{\mathrm{th}} (0)$ and $\eig \mathbf{M} =
(\widetilde{\lambda}_1, \cdots , \widetilde{\lambda}_N)$, it suffices to show
that
\begin{flalign*}
\displaystyle\sum_{\mu \in \eig \mathbf{M}^{[N_1,i-1]}} \log
|\widetilde{\lambda}_k - \mu| - \displaystyle\sum_{\mu \in \eig
\mathbf{M}^{[N_1, \varphi-1]}} \log |\widetilde{\lambda}_k - \mu| \le C
|i-\varphi|.
\end{flalign*}

Again using \cite[Corollaries 5.5 and 5.6]{Agg25a}, now to compare $\eig
\mathbf{M}^{[N_1,i-1]}$ and $\eig \mathbf{M}^{[N_1, \varphi-1]}$ to $\eig
\mathbf{M}$ (and estimating $|\log |\widetilde{\lambda}_k - \mu|| \le
(\log N)^2$ by Lemma \ref{eigenvalues0}), we find that it suffices to show
with overwhelming probability that
\begin{flalign*}
-\displaystyle\sum_{j=i}^{\varphi-1} \log |\widetilde{\lambda}_k -
\widetilde{\Lambda}_j| \le C |i-\varphi|,
\end{flalign*}

\noindent where $\widetilde{\Lambda}_j =
\widetilde{\lambda}_{\tilde{\varphi}^{-1} (j)}$ for any $\zeta$-localization
center bijection $\widetilde{\varphi} : \llbracket 1, N \rrbracket
\rightarrow \llbracket N_1, N_2 \rrbracket$ of $\mathbf{M}$ (that is, the
$\widetilde{\Lambda}_j$ are the eigenvalues of $\mathbf{M}$, indexed by their
localization centers). The above estimate is a quick consequence of \eqref{eqn:hiTconv} (with \eqref{eqn:vrhob}); we omit further
details.
\end{proof}

\subsection{Proof of Lemma \ref{lem:logdivlambdak}}
\label{subsec:logsum_quadsum}

\begin{proof}
We work on the events from Lemmas
\ref{lem:AUpi}, \ref{lem:lsumbd}, and \ref{lem:lqebd}, together with
$\mathsf{SEP}_{\L(0)}(e^{-(\log N)^2})$. Recall
$\mathcal U_-:=\eig\H_-(t)\setminus\pi_-(\mathcal A_-)$. In order to apply
Lemma \ref{lem:suff}, identify the index set $\mathcal A_-$ with the set of
eigenvalues
$\widehat{\mathcal A}_-:=\{\Lambda_a(t):a\in\mathcal A_-\}$, and let
$\widehat\pi_-(\Lambda_a(t)):=\pi_-(a)$.
Set
$\mathcal A_+:=\llbracket\varphi+10R^{-1}\mathfrak M,N_2\rrbracket$.
By Lemma \ref{lem:AUpi} and symmetry, there exists an injection
$\pi_+:\mathcal A_+\to\eig\H_+(t)$ satisfying
$|\pi_+(a)-\Lambda_a(t)|\le e^{-cR^{-1}\mathfrak M}$ for each $a\in\mathcal A_+$.
Define
\[
\mathcal U_+:=\eig\H_+(t)\setminus\pi_+(\mathcal A_+),
\qquad
\widehat{\mathcal A}_+:=\{\Lambda_a(t):a\in\mathcal A_+\},
\qquad
\widehat\pi_+(\Lambda_a(t)):=\pi_+(a).
\]
Also set
\[
\mathcal B:=\eig\L(t)\setminus
\bigl(\widehat{\mathcal A}_-\cup\widehat{\mathcal A}_+
\cup\{\lambda_k\}\bigr).
\]
Thus Lemma \ref{lem:suff} applies with $\lambda=\lambda_k$, the spectral
sets $\widehat{\mathcal A}_\pm$, and the injections $\widehat\pi_\pm$.

By \eqref{eqn:materr}, its symmetric counterpart, and restriction to
$\mathsf{SEP}_{\L(0)}(e^{-(\log N)^2})$,
\begin{equation}\label{eqn:matchlogapp}
\sum_{\sigma\in\{-,+\}}\sum_{a\in\mathcal A_\sigma}
\left|\log|\lambda_k-\pi_\sigma(a)|
-\log|\lambda_k-\Lambda_a(t)|\right|
\le e^{-c R^{-1}\mathfrak M}=: \mathcal E.
\end{equation}
Indeed, each summand is at most
$Ce^{-cR^{-1}\mathfrak M+(\log N)^2}$, and there are at most $N$ summands. 

We first prove \eqref{eqn:less1}. Since
$\log|\lambda_k-\mu|<0$ when $|\lambda_k-\mu|<1$,
\eqref{eq:unmatched-singular-bound} and
\eqref{eqn:matchlogapp} imply
\begin{align*}
&\sum_{\sigma\in\{-,+\}}
\sum_{\substack{\mu\in\mathcal U_\sigma\\|\mu-\lambda_k|<1}}
\left|\log|\lambda_k-\mu|\right|\\
&\quad\le
\sum_{\theta\in\mathcal B}\left|\log|\lambda_k-\theta|\right|
+\sum_{\sigma\in\{-,+\}}
\sum_{\substack{\mu\in\mathcal U_\sigma\\|\mu-\lambda_k|\ge1}}
\log|\lambda_k-\mu|
+C|\log\zeta|+C e^{-c R^{-1}\mathfrak M}.
\end{align*}
The set $\mathcal B$ is contained in
$\{\Lambda_a(t):0<|a-\varphi|\le10R^{-1}\mathfrak M\}$. 
The bound \eqref{eqn:sumlambdak01_index} from Lemma \ref{lem:lsumbd}, applied to this interval
with $r=21R^{-1}$, bounds the first term on the right by
$CR^{-1}\mathfrak M\log(|\lambda_k|+2)$. By Lemma
\ref{lem:lqebd} and symmetry,
\begin{equation}\label{eqn:qeupum}
\sum_{\mu\in\mathcal U_-}\mu^2+
\sum_{\mu\in\mathcal U_+}\mu^2\le CR^{-1}\mathfrak M.
\end{equation}
Moreover, $|\mathcal U_-|+|\mathcal U_+|\le CR^{-1}\mathfrak M$, and, for $|\lambda_k-\mu|>1$,
\[
\log|\lambda_k-\mu|
\le C\bigl(\log(|\lambda_k|+2)+\mu^2\bigr).
\]
These estimates bound the remaining terms by
$CR^{-1}\mathfrak M\log(|\lambda_k|+2)$, since
$|\log\zeta|\leq (\log N)^2$, and thus prove the first
part of the lemma.

Now suppose that $|\lambda_k|>10$. Since
$|\mathcal U_-|+|\mathcal U_+|=|\mathcal B|$,
we may subtract $\log|\lambda_k|$ from each logarithm in
\eqref{eq:unmatched-singular-bound}. Moreover, if
$|\mu-\lambda_k|<1$, then $|1-\mu/\lambda_k|<1$, so the normalized
logarithm is negative. Splitting the resulting identity according to
whether $|\lambda_k-\mu|<1$ gives 
\begin{align*}
&\sum_{\sigma\in\{-,+\}}
\sum_{\mu\in\mathcal U_\sigma}
\left|\log\big|1-\frac{\mu}{\lambda_k}\big|\right|\\
&\quad\le
\sum_{\theta\in\mathcal B}
\left|\log\big|1-\frac{\theta}{\lambda_k}\big|\right|
+2\sum_{\sigma\in\{-,+\}}
\sum_{\substack{\mu\in\mathcal U_\sigma\\|\mu-\lambda_k|\ge1}}
\left|\log\big|1-\frac{\mu}{\lambda_k}\big|\right|
+C|\log\zeta|+C\mathcal E.
\end{align*}
The bound \eqref{eqn:sumellk} from Lemma \ref{lem:lsumbd} bounds the first term on the
right by $CR^{-1}\mathfrak M$, while
the fact that for $|\mu-\lambda_k|\ge1$,
$\left|\log \big|1-\frac{\mu}{\lambda_k}\big|\right|
\le C(1+\mu^2)$
uniformly for $|\lambda_k|>10$, and \eqref{eqn:qeupum} 
bound the second by $CR^{-1}\mathfrak M$. This proves
\eqref{enq:divlk1}.
\end{proof}
	
\subsection{Proof of Lemma \ref{lem:lqebd}}
\label{subsec:unmatchedmu2pf}

First, we record a deterministic lemma which converts exponential
localization of eigenvectors into an upper bound on the sum of $\lambda_i^2$
over eigenvectors whose centers lie in a spatial block. Its proof is given in
Section \ref{subsec:local_energy_upper}.

\begin{lem}
\label{lem:local-energy-from-localization}
There is a constant $c>0$ such that
the following holds for each $C_0 > 1$. Let $\mathbf J$ be a real symmetric tridiagonal matrix
with rows and columns indexed by $\llbracket N_1,N_2\rrbracket$, diagonal
entries $b_i$, and off diagonal entries $a_i$. Let
$\eig\mathbf J=(\lambda_1,\lambda_2,\ldots,\lambda_N)$ and, for each
$j\in\llbracket1,N\rrbracket$, let
$\mathbf u_j=(u_j(N_1),\ldots,u_j(N_2))$ denote a unit eigenvector of
$\mathbf J$ with eigenvalue $\lambda_j$; fix a bijection
$\llbracket1,N\rrbracket\ni j\mapsto\varphi_j\in\llbracket N_1,N_2\rrbracket$.

Let $I=\llbracket\ell,m\rrbracket$, and for $D\ge0$ define
$I^{(D)}:=\{i\in\llbracket N_1,N_2\rrbracket:\dist(i,I)\le D\}$. Assume
that
\begin{equation}\label{eqn:leakage}
\max_{i\in I}|u_j(i)|\le e^{-R^{-1}\mathfrak M}
\qquad\text{whenever }
\varphi_j\notin I^{(C_0R^{-1}\mathfrak M)},
\end{equation}
and that $\max_j|\lambda_j|\le\log N$. Then
\begin{equation}
\label{eqn:local-energy-rM}
\sum_{i=\ell}^{m}b_i^2
+2\sum_{i=\ell}^{m-1}a_i^2
\le
\sum_{j:\,\varphi_j\in I^{(C_0R^{-1}\mathfrak M)}}\lambda_j^2
+c^{-1}e^{-cR^{-1}\mathfrak M}.
\end{equation}
\end{lem}

		\begin{proof}[Proof of Lemma \ref{lem:lqebd}]
	We restrict to the event of Proposition
	\ref{prop:outside-strip-confinement} and to
	$\mathsf{BND}_{\L(t)}(\log N)$, which is overwhelmingly probable by
	Lemma \ref{l0eigenvalues}. Put
		$I_-:=\llbracket N_1,m-\lceil20R^{-1}\mathfrak M\rceil\rrbracket$.
	Note that
		$I_-^{(5R^{-1}\mathfrak M)}
		\subset
		\mathcal{A}_-^m$.
		
	By \eqref{eqn:uMR}, which holds on $\mathfrak B(t)$, and Proposition
	\ref{prop:outside-strip-confinement}, the bound \eqref{eqn:leakage}
	holds with $C_0=5$. Thus, \eqref{eqn:local-energy-rM} yields
	\[
		\Tr\bigl((\L^{I_-}(t))^2\bigr)
		\le
		\sum_{i\in \mathcal{A}_-^m }\Lambda_i(t)^2+o(1).
	\]
	Using \eqref{eqn:materr}, injectivity of $\pi_-^m$, the bound
		$|\Lambda_i(t)|+|\pi_-^m(i)|\le C\log N$, and
		the display above, we obtain
	\[
		\sum_{\mu\in \pi_-^m(\mathcal A_-^m)}\mu^2
		\ge
		\sum_{i\in \mathcal A_-^m}\Lambda_i(t)^2-o(1)
		\ge
		\Tr\bigl((\L^{I_-}(t))^2\bigr)-o(1),
	\]
	where $o(1)$ is a quantity whose absolute value is upper bounded by $1$. Consequently,
	\begin{equation}\label{eqn:Umainbd}
	\begin{aligned}
		\sum_{\mu\in \mathcal U_-^m}\mu^2
		&=
		\Tr((\H_-^m(t))^2)-\sum_{\mu\in\pi_-^m(\mathcal A_-^m)}\mu^2\\
		&\le
		\Tr\bigl((\L^{[N_1,m-1]}(t))^2\bigr)
		-
		\Tr\bigl((\L^{I_-}(t))^2\bigr)
		+o(1).
	\end{aligned}
	\end{equation}
		The trace difference is supported on an index interval of length $O(R^{-1}\mathfrak M)$
		near $m$, and by Lemma \ref{lem:local-lax-trace},
	\[
		\Tr\bigl((\L^{[N_1,m-1]}(t))^2\bigr)
		-
		\Tr\bigl((\L^{I_-}(t))^2\bigr) \leq \Tr\bigl((\L^{[m - \lfloor  30 R^{-1}\mathfrak M \rfloor,m]}(t))^2\bigr)
		\le
		C R^{-1}\mathfrak M.
	\]
	Insertion of the above into \eqref{eqn:Umainbd} proves \eqref{eqn:mu2U}.
	\end{proof}

\section{Eigenvalue estimates}
			\label{sec:min}	

The goal of this section is to provide proofs of Lemmas
\ref{lem:sumlambdak01_half}, \ref{lem:local-lax-trace},
\ref{lem:lsumbd}, and \ref{lem:local-energy-from-localization}. Throughout the section, we adopt Assumption \ref{ass:model} and the notation
of Definition \ref{def:quasi_intro}. We also recall $\delta\in(0,1/200)$,
$\mathfrak M\in[T^{3/4-\delta},T^{3/4}]$, $R=\mathfrak A^{1/10}$, and
$\mathfrak T=\mathfrak T_{\mathfrak M}$ from Definition \ref{mt}. For \(\tau \geq 0\), we use the notation
$m_{\tau}^{(T)}$ from Section \ref{subsec:tc}, and set
\[
Y_i(\tau):=x_\tau^{(T)}(\Lambda_i;\bar q_i).
\]				
	\subsection{Proofs of Lemmas \ref{lem:sumlambdak01_half} and \ref{lem:lsumbd}}
	\label{subsec:min} 
	
We first establish an analogue of Lemma \ref{lem:lsumbd}, where
the condition $\varphi_t(i) \in W$ (for $t \geq 0 $ and an index $i$ as in \eqref{eqn:sumlambdak01_index}) is replaced by \(Y_i(\tau)\in I\) (with $\tau = t/T$) for a
spatial interval $I \subset \mathbb{R}$.

\begin{lem}\label{lem:sumlambdak01}
Let $\mathfrak M$ be as in Definition \ref{mt}, but assume more generally in
this lemma that
$0<\delta\le \frac{1}{16}$ and
$\mathfrak M\in[T^{3/4-\delta},T^{3/4}]$.
Fix \(\tau \in [0, \log N]\), and set
$
J_n:=[n,n+1]\cup[-n-1,-n]$.
There are constants \(c,C>0\) such that, with overwhelming probability,
the following holds.

For every \(r\in[\A^{-1/10},\infty)\), every
interval \(I\) satisfying $I\subset \mathbb{R},
|I|\le r \A \frac{\mathfrak M}{T}$,
and every \(k\in\llb1,N\rrb\),
\begin{equation}\label{eqn:sumlambdak01_shell}
\sum_{\substack{i\in\llb N_1,N_2\rrb:\,
Y_i(\tau)\in I,\ \Lambda_i\in J_n\\ i\ne\varphi_0(k)}}
\bigl|\log|\lambda_k-\Lambda_i|\bigr|
\le
C r \mathfrak M e^{-cn^2}\log(|\lambda_k|+2)+(\log N)^C .
\end{equation}
\end{lem}

	Before proving Lemma \ref{lem:sumlambdak01}, we state an intermediate lemma.
	The following lemma states that the $Y_{i}(\tau)$ with nearby eigenvalues 
	$\Lambda_i$ (see the beginning of the section for this notation) do not accumulate too densely in a given interval.
	The proof is postponed to Section \ref{subsec:weighted_rectangle}.

		\begin{lem}\label{lem:weighted_rectangle}
				Fix \(\tau\in[0,\log N]\). Fix intervals
				\(J\subseteq[-\log N,\log N]\) and \(I\subseteq\mathbb R\) with
				\(T^{-5/16}\le |I|\le N^2\), and put
				$\lambda_J:=\inf_{\lambda\in J}|\lambda|$.
				Then there are constants $c,C>0$ such that, with overwhelming
				probability,
				\begin{equation}\label{eq:rect-thin}
				\#\{i\in\llb N_1,N_2\rrb:\ Y_i(\tau)\in I,\ \Lambda_i\in J\}
				\le
				C\Big(\mathfrak A^{-1}T|I||J|e^{-c\lambda_J^2}+(\log N)^C\Big).
				\end{equation}
				\end{lem}

				\begin{proof}[Proof of Lemma \ref{lem:sumlambdak01}]
	It suffices to prove the estimate \eqref{eqn:sumlambdak01_shell} when
	$I\subseteq[-N^2/2,N^2/2]$, since Lemmas \ref{l0eigenvalues},
	\ref{lem:Estr}, and \ref{lem:vebd} imply that
	$Y_i(\tau)\in[-N^2/2,N^2/2]$ with overwhelming probability, for all
	$i\in\llbracket N_1,N_2\rrbracket$. Put
	$\eta:=\mathfrak M/T$. We may, and will, assume $|I| \geq \eta$.
	Since \(\mathfrak M\in[T^{3/4-\delta},T^{3/4}]\) and
	\(\delta\le 1/16\), we have
	$T^{-5/16}\le \eta\le T^{-1/4}$.
	Let \(\mathcal I_\eta\) be the family of all intervals with
	endpoints in \(\eta\mathbb Z\), which are contained in
	$
	\left[-N^2, N^2\right]$.
	The number of such intervals satisfies \(|\mathcal I_\eta| \leq N^6\).  Set \(\ell_0=\lceil \frac12\log T\rceil\). For
	\(0\le \ell\le\ell_0\), let \(\mathcal J_\ell\) be the family of intervals of
	length \(4e^{-\ell}\), with left endpoints in \(e^{-\ell-2}\mathbb Z\), that
	intersect \([-\log N-1,\log N+1]\).  The total number of these
	intervals is at most \(CT^{1/2}\log N\).

	By Lemmas \ref{l0eigenvalues} and \ref{eigenvalues0}, the event
	$\mathsf E_0:=
	\Bigl\{\max_{1\le j\le N}|\lambda_j|\le \log N\Bigr\}
	\cap
	\mathsf{SEP}_{\L}(e^{-(\log N)^2})$
	is overwhelmingly probable; moreover, so is the intersection $\mathsf{E}_1$ of events 
	\eqref{eq:rect-thin} (from Lemma \ref{lem:weighted_rectangle}) over all rectangles \(I_0\times J\) with
	\(I_0\in\mathcal I_\eta\), and \(J\) of the form $J = J_n$ or $J = \tilde J \cap J_n$ for \(0\le n\le\lceil\log N\rceil\), 
	where
	\(\tilde J \in\mathcal J_\ell\) for \(0\le\ell\le\ell_0\). We work on these events in what follows.

	Fix \(n\ge0\), \(k\), and an interval \(I\) as in the statement (recalling that
	we now assume $|I| \geq \eta$).
	If \(n>\lceil\log N\rceil\), the left side of
	\eqref{eqn:sumlambdak01_shell} is zero on \(\mathsf E_0\).  Assume then that
	\(n\le\lceil\log N\rceil\).  Note that \(I\) is covered by 
	an element \(I_0 \in \mathcal I_\eta\) of size at most $3 |I| $.
	On $\mathsf{E}_1$, the bound \eqref{eq:rect-thin}
	applied to $(I_0, J)$
	 for any $J = J_n$ or 
	$J = \tilde J \cap J_n$ of the 
	form described in the previous paragraph gives, using
	$\mathfrak A^{-1}T|I|\le r\mathfrak M$,
	\begin{equation}\label{eq:sumlambdak01-shell-count-new}
	\#\{i:Y_i(\tau)\in I,\ \Lambda_i\in  J\}
	\le
	Cr\mathfrak M e^{-cn^2}|J|
	+C(\log N)^C .
	\end{equation}

	We split the left side of \eqref{eqn:sumlambdak01_shell} into the following chunks:
		\[
	\widetilde{J}_1 \coloneqq \{\lambda: |\lambda-\lambda_k|>1 \} ,\qquad
	\widetilde{J}_2 \coloneqq \{\lambda: T^{-1/2}<|\lambda-\lambda_k|\le 1\} ,\qquad
	\widetilde{J}_3 \coloneqq \{\lambda: 0<|\lambda-\lambda_k|\le T^{-1/2}\}.
	\]
	For $\Lambda_i \in \widetilde{J}_1$,
	\[
	|\log|\lambda_k-\Lambda_i||\le C  \log(|\lambda_k|+2)+C(n+1) \leq 2 C (n+1) \log(|\lambda_k|+2),
	\]
	and \((n+1)e^{-cn^2}\le C e^{-c'n^2}\) for some $c'>0$.  Thus
	\eqref{eq:sumlambdak01-shell-count-new} bounds the contribution to the left
	hand side of \eqref{eqn:sumlambdak01_shell} from
	$\Lambda_i\in\widetilde{J}_1$ by
	\begin{align*}
	 \sum_{\substack{i\in\llb N_1,N_2\rrb:\,
	Y_i(\tau)\in I,\ \Lambda_i\in \widetilde{J}_1\cap J_n\\
	i\ne\varphi_0(k)}}
	\bigl|\log|\Lambda_i-\lambda_k|\bigr| 
	&\le
	2C(n+1)\log(|\lambda_k|+2)
	\#\{i:Y_i(\tau)\in I,\ \Lambda_i\in J_n\}\\
	&\le
	2 C^2 (n+1)\log(|\lambda_k|+2)
	\bigl(r\mathfrak M e^{-cn^2}+(\log N)^C\bigr)\\
	&\le
	C^3 r\mathfrak M\log(|\lambda_k|+2)
	e^{-c' n^2}
	+(\log N)^{C+3}.
	\end{align*}
 The right hand side of the last bound matches the
	right side of \eqref{eqn:sumlambdak01_shell}, after changing \(c,C\) there.

	For the range $\widetilde J_2$, define
	\[
	J_{k,\ell}^{+}:=[\lambda_k+e^{-(\ell+1)},\lambda_k+e^{-\ell}],
	\qquad
	J_{k,\ell}^{-}:=[\lambda_k-e^{-\ell},\lambda_k-e^{-(\ell+1)}],
	\]
	for \(0\le\ell<\ell_0\).  Choose
	\(\widetilde J_{k,\ell}^{\sigma}\in\mathcal J_\ell\), \(\sigma\in\{+,-\}\),
	containing \(J_{k,\ell}^{\sigma}\), but not containing $\lambda_k$, for each $\ell \leq \ell_0$;
	recall that $|\widetilde J_{k,\ell}^{\sigma}| = 4 e^{-\ell}$ by definition.
	  Applying \eqref{eq:sumlambdak01-shell-count-new} with
	\(J = \widetilde J_{k,\ell}^{\sigma}\cap J_n\) gives
	\[
	\#\{i:Y_i(\tau)\in I,\ \Lambda_i\in J_n\cap J_{k,\ell}^{\sigma}\}
	\le
	Cr\mathfrak M e^{-cn^2}e^{-\ell}
	+(\log N)^C .
	\]
	Multiplying by \(\ell+1\), since $|\log|\lambda_k-\Lambda_i|| \leq (\ell + 1)$ for $\Lambda_i \in \widetilde J_{k,\ell}^{\sigma}$, and summing over
	\(\sigma\in\{+,-\}\), \(0\le\ell<\ell_0\), yields the claimed bound for the
	contribution to the left side of \eqref{eqn:sumlambdak01_shell} from
	$\Lambda_i \in \widetilde{J}_2$.

	To bound the contribution to the left side of
	\eqref{eqn:sumlambdak01_shell} from $\Lambda_i\in\widetilde J_3$, 
	 choose
	\(\widehat J_k\in\mathcal J_{\ell_0}\) containing $\widetilde J_3$.
	Applying \eqref{eq:sumlambdak01-shell-count-new} with
	\(J=\widehat J_k\cap J_n\) gives
	\[
	\#\{i:Y_i(\tau)\in I,\ \Lambda_i\in J_n\cap \widetilde J_3\}
	\le
	Cr\mathfrak M e^{-cn^2}T^{-1/2}
	+(\log N)^C .
	\]
	On \(\mathsf E_0\),
	\(|\log|\lambda_k-\Lambda_i|| \leq (\log N)^2\).  Since \(T^{-1/2}(\log N)^2\le 1\) for large \(T\), this proves
	the contribution to the left side of \eqref{eqn:sumlambdak01_shell} from
	$\Lambda_i \in \widetilde{J}_3$ is again bounded by its right hand side.
	Adding the bounds for the three
	ranges proves \eqref{eqn:sumlambdak01_shell}, after increasing constants.
	\end{proof}

	\begin{proof}[Proof of Lemma \ref{lem:lsumbd}]
	Set $\tau=t/T$ and $J_n:=[n,n+1]\cup[-n-1,-n]$, and assume
	$\mathfrak B(t)$ holds.
	We prove the first part in detail. By partitioning $W$ if necessary, it
	suffices to consider $r\le\mathfrak A^{10}$. We may also assume
	$W\subseteq\llbracket N_1+T(\log N)^{15},
	N_2-T(\log N)^{15}\rrbracket$. Indeed, the portions outside this interval
	are handled by the coupling to thermal equilibrium from Lemma
	\ref{lem:edge-thermal-coupling}, through an argument
	similar to the one outlined
	at the end of the proof of Lemma \ref{lem:log_sum_dif}. The bound
	\eqref{eqn:qbd} gives real numbers $x_-<x_+$ such that
	\[
	\{q_a(t):a\in W\}\subseteq[x_-,x_+],
	\qquad
	x_+-x_-\le C\mathfrak A\bigl(|W|+\mathfrak M/R\bigr)
	\le Cr\mathfrak A\mathfrak M.
	\]
	If $\varphi_t(j)\in W$ and $\lambda_j\in J_n$, then
	$Q_j(t)=q_{\varphi_t(j)}(t)\in[x_-,x_+]$. Hence \eqref{eqn:Qx} gives
	\[
	Y_{\varphi_0(j)}(\tau)\in I_n:=
	\left[T^{-1}x_--C(n+1)\frac{\mathfrak M}{T},
	T^{-1}x_++C(n+1)\frac{\mathfrak M}{T}\right].
	\]
	The interval $I_n$ satisfies
	$|I_n|\le C(r\mathfrak A+n+1)\mathfrak M/T$. Thus, applying Lemma
	\ref{lem:sumlambdak01}, with the $r$ there equal to
	$\widetilde r_n=C(r+\mathfrak A^{-1}(n+1))$ here, and summing over
	$0\le n\le\lceil\log N\rceil$ gives
	\[
	\sum_{\substack{j:\,\varphi_t(j)\in W\\j\ne k}}
	\left|\log|\lambda_k-\lambda_j|\right|
	\le Cr\mathfrak M\log(|\lambda_k|+2)+(\log N)^C.
	\]
	The second term is absorbed into the first, since
	$r\ge\mathfrak A^{-1/10}$, $\log(|\lambda_k|+2)\ge\log2$, and
	$\mathfrak M\ge T^{11/16}$. This
	proves the first part.

	We now outline the similar argument for the second part. Fix
	$a:=|\lambda_k|>10$ and set
	\[
	\mathcal S_n:=
	\{j:\varphi_t(j)\in W,\ j\ne k,\ \lambda_j\in J_n\}.
	\]
	Lemma \ref{lem:weighted_rectangle} and Lemma \ref{lem:sumlambdak01},
	with $I$ and $r$ there given by 
	$I_n$ and $\widetilde r_n$ as in the
	first part of the proof, give
	\begin{equation}\label{eqn:2ests}
	\begin{aligned}
	|\mathcal S_n|
	&\le C\widetilde r_n\mathfrak M e^{-cn^2}+C(\log N)^C,\\
	\sum_{j\in\mathcal S_n}\left|\log|\lambda_k-\lambda_j|\right|
	&\le C\widetilde r_n\mathfrak M e^{-cn^2}\log(|\lambda_k|+2)+C(\log N)^C.
	\end{aligned}
	\end{equation}
	If $n+1\le a/2$, then
	\[
	\left|\log\big|1-\frac{\lambda_j}{\lambda_k}\big|\right|
	\le C\frac{n+1}{a},
	\qquad j\in\mathcal S_n,
	\]
	so the first estimate in \eqref{eqn:2ests} implies (upon summing over $n \leq a/2-1$) that the part of the sum
	in \eqref{enq:divlk1} over $|\lambda_j| \le a/2$
	is at most $Cr\mathfrak M+C(\log N)^C$. For the sum over $j \in \mathcal{S}_n$ with $n> a/2-1$, we use
	\[
	\left|\log\big|1-\frac{\lambda_j}{\lambda_k}\big|\right|
	\le\left|\log|\lambda_k-\lambda_j|\right|+\log a.
	\]
	The two estimates in \eqref{eqn:2ests} (using the Gaussian factor $e^{-c a^2}$, which remains after summing over $n$,
	to cancel out the growth of the factor $\log(|\lambda_k|+2)$), again give a bound of $Cr\mathfrak M+C(\log N)^C$. Combining these estimates proves the second part of Lemma \ref{lem:lsumbd}.
	\end{proof}

\begin{proof}[Proof of Lemma \ref{lem:sumlambdak01_half}]
Assume \(\mathfrak B(t)\) holds, fix \(k\) and \(a\) as in the statement, and set
$i_j:=\varphi_0(j),
L_k:=\log(|\lambda_k|+2),
J_n:=[n,n+1]\cup[-n-1,-n]$.
On an event of overwhelming probability, the following holds. Either $j$
is such that
$\varphi_t(j)\in\llbracket N_1+T(\log N)^{10},
N_2-T(\log N)^{10}\rrbracket$, in which case
\[
\left|T^{-1}Q_j(t)-Y_{i_j}(\tau)\right|
\le
\frac{\mathfrak M}{T}\log(|\lambda_j|+2),
\]
by \eqref{eqn:Qx}; otherwise, if
$\varphi_t(j) \notin \llbracket N_1 + T (\log N)^{10}, N_2 - T (\log N)^{10} \rrbracket$, then 
both $Y_{i_j}(\tau) - a$ and $Q_j(t)- a$ have the same sign.
Indeed, if $\varphi_t(j) < N_1 + T (\log N)^{10}$, then
on the event from
Proposition \ref{prop:outside-strip-confinement}, $\varphi_0(j) < N_1 + T (\log N)^{11}$, 
and $Y_{i_j}(\tau) < a$ follows
from the fact that $a > -(\log N)^{100}$ together with the fact that 
$Y_{i_j}(\tau) < -(\log N)^{100}$ by Lemmas \ref{lem:Estr} and \ref{lem:vebd},
after further restriction to an overwhelmingly probable event. 
In this same case $\varphi_t(j) < N_1 + T (\log N)^{10}$, the bound $Q_j(t)/T < a$ follows
from Lemmas \ref{lem:Estr} and \ref{l0eigenvalues} with $A = \log N$. Thus,  
if $\varphi_t(j) < N_1 + T (\log N)^{10}$, $\max(T^{-1}Q_j(t), Y_{i_j}(\tau)) < a$. 
Finally,
in the case $\varphi_t(j)>N_2-T(\log N)^{10}$ a very similar argument 
shows that $\min(T^{-1}Q_j(t), Y_{i_j}(\tau)) > a$.

Hence, because the two indicators below use the same cutoff $a$,
$\mathbbm 1_{\{T^{-1}Q_j(t)<a\}}
\neq
\mathbbm 1_{\{Y_{i_j}(\tau)<a\}}$
only if
\begin{equation}\label{eq:bootstrap-cutoff-strip}
\left|Y_{i_j}(\tau)- T^{-1}Q_j(t)\right|
\le
C\frac{\mathfrak M}{T}\log(|\lambda_j|+2).
\end{equation}
For \(n\ge0\), let us bound the contribution to the left side of
\eqref{eqn:sumlambdak01_half} from the indices satisfying $\mathbbm 1_{\{T^{-1}Q_j(t)<a\}}
\neq
\mathbbm 1_{\{Y_{i_j}(\tau)<a\}}$ and
\(\lambda_j\in J_n\); for such indices,
$\log(|\lambda_j|+2)\le C(n+2)$.
Thus Lemma \ref{lem:sumlambdak01}, with the $r$ there equal to
$r_n
:=
C\max\left\{
\mathfrak A^{-1/10},
\mathfrak A^{-1}(n+2)
\right\}$ here,
gives
\[
\sum_{\substack{j\ne k:\,\lambda_j\in J_n\\
\mathbbm 1_{\{T^{-1}Q_j(t)<a\}}
\ne
\mathbbm 1_{\{Y_{i_j}(\tau)<a\}}}}
\left|\log|\lambda_k-\lambda_j|\right|
\le
C\left(
\mathfrak A^{-1/10}
+
\mathfrak A^{-1}(n+2)
\right)
\mathfrak M e^{-cn^2}L_k
+
(\log N)^C.
\]
On the event of Lemma \ref{l0eigenvalues}, only
\(0\le n\le \lceil\log N\rceil\) contribute. Summing over such \(n\), and using
\(R=\mathfrak A^{1/10}\), yields
\[
\sum_{\substack{j\ne k:\\
\mathbbm 1_{\{T^{-1}Q_j(t)<a\}}
\ne
\mathbbm 1_{\{Y_{i_j}(\tau)<a\}}}}
\left|\log|\lambda_k-\lambda_j|\right|
\le
C R^{-1}\mathfrak M L_k
+
(\log N)^C.
\]
The final term is absorbed into the first, since
$L_k\ge\log 2$ and $R^{-1}\mathfrak M$ dominates every fixed power of
$\log N$.
This proves \eqref{eqn:sumlambdak01_half}.
\end{proof}



\subsection{Proof of Lemma \ref{lem:weighted_rectangle}}
\label{subsec:weighted_rectangle}

First we state the following lemma, proved in Appendix \ref{app:frozen-dos}. In this statement,
and in the proof of Lemma \ref{lem:weighted_rectangle}, which follows it, we 
denote $\L = \L(0)$, and for an
index interval $B \subset \llbracket N_1, N_2 \rrbracket$, recall from
Section \ref{notation}
that we denote by $\L^{B}$ the block of $\L$ indexed
by indices in $B$.
	
			\begin{lem}\label{lem:block_wegner_input}
			There are constants \(c,C>0\) such that the following holds. Let
			\(B\subseteq\llb N_1,N_2\rrb\) be an interval, let
			\(J\subseteq\mathbb R\) be an interval, and let
			\[
			X_B(J):=\#\bigl(\eig(\L^{B})\cap J\bigr),
			\qquad
			\lambda_J:=\inf_{\lambda\in J}|\lambda|.
			\]
			Then, 
			\[
			\E X_B(J)\le C |B|\,|J|\,e^{-c\lambda_J^2}.
			\]
			\end{lem}
			
			\begin{proof}[Proof of Lemma \ref{lem:weighted_rectangle}]
 Let
$\ell:=|I|,
h:=
\min\left\{
1,\frac{\ell}{C_0(\log N)^3},|J|
\right\}$,
where \(C_0\) is sufficiently large. Let $A\ge1$ be an integer and partition
\(J\) into disjoint intervals $(J_\alpha)_{\alpha=1}^A$ such that
$h/2\le|J_\alpha|\le h$. For every
\(\alpha\in\llbracket1,A\rrbracket\), define
\[
\mathcal P_{0,\alpha}
:=
\{(q,\lambda):
\lambda\in J_\alpha,\
x_\tau^{(T)}(\lambda;q)\in I\}.
\]
	By Lemma \ref{lem:uniform-non-vacuum}, the solution
	\(m_\tau^{(T)}\) satisfies the lower bound hypothesis 
	\eqref{eqn:mintlow} of Lemma
	\ref{lem:backward-flow-rectangle}, uniformly in \(T\). 
	Thus, Lemma
	\ref{lem:backward-flow-rectangle} gives an interval \(I_\alpha^+\) such that
	$\mathcal P_{0,\alpha}\subseteq I_\alpha^+\times J_\alpha$ and
\begin{equation}\label{eq:weighted-Ialpha-size}
|I_\alpha^+|
\le
C\left(
\ell+(\log N)^2|J_\alpha|
\right)
\le C\ell.
\end{equation}
Indeed, the first inequality follows by applying the lemma with
$L=\log N$ (which we may 
by the assumption $J \subseteq [-\log N, \log N]$) and inserting \(\tau\le\log N\) into its bound,
and the last inequality follows from the definition of \(h\).

Define the random index set
\[
S_{\alpha;0}
:=
\left\{
i\in\llbracket N_1,N_2\rrbracket:
q_i(0)/T
\in
I_\alpha^+
\right\}.
\]
If $S_{\alpha;0}\ne\varnothing$, set
$S_\alpha:=\llbracket\min S_{\alpha;0},\max S_{\alpha;0}\rrbracket$; if
$S_{\alpha;0}=\varnothing$, set $S_\alpha:=\varnothing$.
If \(Y_i(\tau)\in I\) and \(\Lambda_i\in J_\alpha\), then
\(\bar q_i\in I_\alpha^+\), and hence
\(i\in S_{\alpha;0}\subseteq S_\alpha\).
Lemma \ref{lem:Estr} implies, with overwhelming probability, that
\begin{equation}\label{eq:weighted-Salpha-size}
|S_\alpha|
\le
C\mathfrak A^{-1}T\ell.
\end{equation}

Choose
$D:=(\log N)^{5}$. Let $\mathfrak n\ge1$ be an integer and
$(B_m)_{m=1}^{\mathfrak n}$ be a family of consecutive disjoint intervals
such that
$\bigcup_{m=1}^{\mathfrak n}B_m=\llbracket N_1,N_2\rrbracket$ and
$D/2\le|B_m|\le D$. Let
$B_m^+:=\{i\in\llbracket N_1,N_2\rrbracket:\dist(i,B_m)\le D/4\}$ denote
the enlargement of $B_m$ by $D/4$ on both sides.

The proof proceeds by subdividing the Lax matrix $\L$ into the principal
blocks $\L^{B_m^+}$, approximating the relevant eigenvalues of $\L$ by
eigenvalues of these blocks, and using independence between the blocks and a
Bernstein estimate to bound the left side of \eqref{eq:rect-thin}. For each \(m\in\llbracket1,\mathfrak n\rrbracket\), set
\[
K_m
:=
\bigcup_{\alpha:\,B_m\cap S_\alpha\ne\varnothing}
J_\alpha
\]
and, recalling for $s \geq 0$ and $K \subset \mathbb{R}$ we set $K^{(s)} = \{t \in \mathbb{R} : \dist(t, K) \leq s\}$, let 
\begin{equation}\label{eqn:Xm}
X_m
:=
\#\left(
\operatorname{eig}(\L^{B_m^+})
\cap K_m^{(\varepsilon_D)}
\right),
\qquad
\varepsilon_D:=e^{-c_0 D},
\end{equation}
where $c_0>0$ is a sufficiently small constant.

We claim that, with overwhelming probability,
\begin{equation}\label{eq:weighted-target-to-block}
\#\{i\in B_m:Y_i(\tau)\in I,\ \Lambda_i\in J\}
\le X_m.
\end{equation}
Indeed, work on the events furnished by Lemmas \ref{bijectionl},
\ref{l0eigenvalues} with $A=\log N$, and \ref{eigenvalues0} with
$\delta=e^{-(\log N)^2}$. Suppose that \(i\in B_m\), \(Y_i(\tau)\in I\), and
\(\Lambda_i\in J_\alpha\). Then \(i\in S_\alpha\), and hence
\(J_\alpha\subseteq K_m\). On these events, Lemma
	\ref{lem:full-to-truncation}, applied with the $\widetilde\L$ there equal
	to $\L^{B_m^+}$ here, viewed as an $N\times N$ matrix by padding with
	zeros, and with the $\mathcal D$ there equal to
	$\llbracket N_1,N_2\rrbracket\setminus B_m^+$ here (using Lemma
\ref{l0eigenvalues} to verify \eqref{ll} and Lemma \ref{bijectionl} to verify
\eqref{eqn:lmatrix-boundary-decay}), gives an eigenvalue
\(\mu_i\in\operatorname{eig}(\L^{B_m^+})\) satisfying
\[
|\mu_i-\Lambda_i|\le C e^{-cD}.
\]
Thus \(\mu_i\in K_m^{(\varepsilon_D)}\), if we take \(c_0\) small enough.
Since $Ce^{-cD}<e^{-(\log N)^2}$ and the eigenvalues of $\L(0)$ are
separated by at least $e^{-(\log N)^2}$ by our restriction to the event of
Lemma \ref{eigenvalues0}, we have that the $(\mu_i)$ are mutually distinct
with overwhelming probability. Hence, the map from $\Lambda_i$ to $\mu_i$ is
an injection; together with the above, this confirms
\eqref{eq:weighted-target-to-block}.

Due to \eqref{eq:weighted-target-to-block}, we would like to bound
$\sum_{m=1}^{\mathfrak n}X_m$. However, since the $X_m$ depend on the
$S_\alpha$ through $K_m$, and the latter are random, the $(X_m)$ are not
quite independent. To remedy this, we define a variant $X_m^{(a)}$ of $X_m$
obtained from a deterministic variant $S_\alpha^{(a)}$ of $S_\alpha$.
Specifically, for any
$a\in[-N^2,N^2]\cap N^{-1}\mathbb Z$, define, recalling
$F_{\mathrm{str}}^{(T)}$ from \eqref{eqn:Fstr} and Section \ref{subsec:tc},
\begin{equation}\label{eq:weighted-Salpha-a}
S_\alpha^{(a)}
:=
\left\{i\in\llbracket N_1,N_2\rrbracket:
F_{\mathrm{str}}^{(T)}(i/T)+a\in(I_\alpha^+)^{(\ell)}\right\}.
\end{equation}

\noindent Following the above definitions of $K_m$ and $X_m$, set, recalling
$\varepsilon_D$ from \eqref{eqn:Xm},
\begin{equation}\label{eq:weighted-KX-a}
K_m^{(a)}
:=
\bigcup_{\alpha:\,B_m\cap S_\alpha^{(a)}\ne\varnothing}J_\alpha,
\qquad
X_m^{(a)}
:=
\#\left(\operatorname{eig}(\L^{B_m^+})
\cap(K_m^{(a)})^{(\varepsilon_D)}\right).
\end{equation}

We now claim that the following two statements hold. First, there exists some
$a\in[-N^2,N^2]\cap N^{-1}\mathbb Z$ such that $X_m\le X_m^{(a)}$ for all
$m\in\llbracket1,\mathfrak n\rrbracket$. Second, for any
$a\in[-N^2,N^2]\cap N^{-1}\mathbb Z$, we have with overwhelming probability
that
\begin{equation}\label{eqn:Xmalph}
\sum_{m=1}^{\mathfrak n}X_m^{(a)}
\le
C\mathfrak A^{-1}T\ell|J|e^{-c'\lambda_J^2}+(\log N)^C.
\end{equation}
These statements, together with \eqref{eq:weighted-target-to-block} summed
over $m\in\llbracket1,\mathfrak n\rrbracket$, imply the lemma.

To show the first statement, it suffices to show that there exists some
$a\in[-N^2,N^2]\cap N^{-1}\mathbb Z$ such that
$S_\alpha\subseteq S_\alpha^{(a)}$ for each
$\alpha\in\llbracket1,A\rrbracket$. If every $S_\alpha$ is empty, this is
immediate; otherwise, to that end, after relabeling the intervals, fix $i\in S_{1;0}$. There
exists $a\in[-N^2,N^2]\cap N^{-1}\mathbb Z$ such that
$|\bar q_i-F_{\mathrm{str}}^{(T)}(i/T)-a|\le1/N$, and we fix such a value of
$a$. Now fix $j\in S_\alpha$, for any
$\alpha\in\llbracket1,A\rrbracket$. Choose
$j_-,j_+\in S_{\alpha;0}$ with $j_-\le j\le j_+$. Note that
$\bar q_i\in I_1^+$ and $\bar q_{j_-},\bar q_{j_+}\in I_\alpha^+$. Applying Lemma
\ref{lem:backward-flow-rectangle}, this time with $I, J$ there given by 
 $I, J$ here, we obtain an interval $I_0$ which 
satisfies
$\cup_{\alpha=1}^A  I_\alpha^+  \subset I_0$, and which satisfies (using $\tau\le\log N$ and
$J\subseteq[-\log N,\log N]$) 
 $|I_0| \leq C(1+\ell)(\log N)^3$. Thus, 
\[
\max_{\sigma\in\{-,+\}}|\bar q_i-\bar q_{j_\sigma}|
\le C(1+\ell)(\log N)^3.
\]
Combining this with \eqref{eq:stretch-good} and \eqref{eqn:alph_bd} gives
$\max_{\sigma\in\{-,+\}}|j_\sigma-i|\le
T(1+\ell)(\log N)^4$, and hence the same bound for $|j-i|$. Therefore,
\eqref{eq:stretch-good} and $\ell\ge T^{-5/16}$ give, for all sufficiently
large $T$,
\[
\left|\bar q_i-\bar q_j-
\bigl(F_{\mathrm{str}}^{(T)}(i/T)-F_{\mathrm{str}}^{(T)}(j/T)\bigr)\right|
\le\frac{\ell}{20};
\]
indeed, the left side is at most
$C\mathfrak A T^{-1/2}(1+\ell)^{1/2}(\log N)^4=o(\ell)$ under these bounds.
It follows that
$|\bar q_j-F_{\mathrm{str}}^{(T)}(j/T)-a|\le\ell/10$. Moreover, the
monotonicity of $F_{\mathrm{str}}^{(T)}$, the inequalities
$j_-\le j\le j_+$, and the same estimate \eqref{eq:stretch-good} show that
$\dist(\bar q_j,I_\alpha^+)\le\ell/10$. Therefore
$\dist(F_{\mathrm{str}}^{(T)}(j/T)+a,I_\alpha^+)\le\ell/5<\ell$, and hence
$j\in S_\alpha^{(a)}$.

We next prove the second statement. For
$\alpha\in\llbracket1,A\rrbracket$, set
$\lambda_\alpha:=\inf_{\lambda\in J_\alpha}|\lambda|$. Since
$(K_m^{(a)})^{(\varepsilon_D)}$ is contained in the union of the intervals
$(J_\alpha)^{(\varepsilon_D)}$ over the indices $\alpha$ contributing to
$K_m^{(a)}$, Lemma \ref{lem:block_wegner_input} gives, after changing
$c,C>0$, for every
$m\in\llbracket1,\mathfrak n\rrbracket$,
\[
\mathbb E X_m^{(a)}
\le
CD
\sum_{\substack{\alpha\in\llbracket1,A\rrbracket:\,
S_{\alpha}^{(a)} \cap B_m \neq \varnothing}}
(|J_\alpha|+2\varepsilon_D)
e^{-c\lambda_\alpha^2}.
\]
For each \(\alpha\in\llbracket1,A\rrbracket\), the definition
\eqref{eq:weighted-Salpha-a}, the bound
\eqref{eq:weighted-Ialpha-size}, \eqref{eqn:alph_bd}, and \eqref{eqn:Fstr}
give
$|S_\alpha^{(a)}|\le C\mathfrak A^{-1}T\ell$. Hence the interval
\(S_\alpha^{(a)}\) intersects at most
$C\mathfrak A^{-1}T\ell/D$ blocks $B_m$. Consequently,
\[
\sum_{m=1}^{\mathfrak n}\mathbb E X_m^{(a)}
\le CD\sum_{\alpha=1}^{A}\frac{\mathfrak A^{-1}T\ell}{D}
(|J_\alpha|+e^{-c (\log N)^5})e^{-c\lambda_\alpha^2}
\le C\mathfrak A^{-1}T\ell|J|e^{-c'\lambda_J^2}=: \nu.
\]

For $\sigma\in\{0,1\}$, let
$\mathcal M_\sigma:=\{m\in\llbracket1,\mathfrak n\rrbracket:
m\equiv\sigma\pmod 2\}$. If $m,m'\in\mathcal M_\sigma$ are distinct, then
$B_m^+$ and $B_{m'}^+$ are disjoint. Thus the matrices
\(\L^{B_m^+}\), $m\in\mathcal M_\sigma$, depend on disjoint sets of
entries of $\L(0)$, which are mutually independent, so the variables
$(X_m^{(a)})_{m\in\mathcal M_\sigma}$ are independent. Moreover,
$0\le X_m^{(a)}\le CD$
and
$\operatorname{Var}(X_m^{(a)})
\le
CD\,\mathbb E X_m^{(a)}$.
Bernstein's inequality, applied separately for each $\sigma\in\{0,1\}$,
therefore gives
\[
\sum_{m=1}^{\mathfrak n} X_m^{(a)}
\le
 C\nu + (\log N)^{10}
\]
with overwhelming probability. Since there are only polynomially many
admissible values of $a$, a union bound makes this estimate simultaneous in
$a$. This proves \eqref{eqn:Xmalph}, and thus completes the proof.
\end{proof}

\subsection{Proofs of Lemmas \ref{lem:local-lax-trace} and \ref{lem:local-energy-from-localization}}
\label{subsec:local_energy_upper}

\begin{proof}[Proof of Lemma \ref{lem:local-energy-from-localization}]
Let $\mathbf P_I$ be the orthogonal projection onto
$\operatorname{span}\{\mathbf e_i:i\in I\}$, where for each
$i\in\llbracket N_1,N_2\rrbracket$ we have let $\mathbf e_i\in\mathbb R^N$
denote the coordinate vector, whose entry indexed by $i$ is equal to $1$ and
whose other entries are equal to $0$. It suffices to prove the deterministic bound
\begin{equation}
\label{eqn:local-energy-main-bound}
\sum_{i=\ell}^{m}b_i^2+2\sum_{i=\ell}^{m-1}a_i^2
\le
\sum_{j:\,\varphi_j\in I^{(D)}}\lambda_j^2
+\left(\sup_{j:\,\varphi_j\notin I^{(D)}}
\|\mathbf P_I\mathbf u_j\|_2^2\right)
\sum_{j:\,\varphi_j\notin I^{(D)}}\lambda_j^2.
\end{equation}
Indeed, take $D=C_0R^{-1}\mathfrak M$. The hypothesis
\eqref{eqn:leakage} gives
\[
\sup_{j:\,\varphi_j\notin I^{(D)}}
\|\mathbf P_I\mathbf u_j\|_2^2
\le |I|e^{-2R^{-1}\mathfrak M}.
\]
Since $|I|\le N$, $\max_j|\lambda_j|\le\log N$, and
$R^{-1}\mathfrak M>(\log N)^2$, the second term on the right side of
\eqref{eqn:local-energy-main-bound} is at most
$c^{-1}e^{-cR^{-1}\mathfrak M}$. Hence
\eqref{eqn:local-energy-main-bound} implies
\eqref{eqn:local-energy-rM}.

Now we show \eqref{eqn:local-energy-main-bound}. Since $\mathbf J$ is symmetric and tridiagonal,
\begin{equation}
\label{eqn:local-trace-identity}
\sum_{j=\ell}^{m} b_j^2
+2\sum_{j=\ell}^{m-1} a_j^2
=\|\mathbf P_I\mathbf J\mathbf P_I\|_{\mathrm{HS}}^2
\le\|\mathbf P_I\mathbf J\|_{\mathrm{HS}}^2
=\Tr(\mathbf P_I\mathbf J^2\mathbf P_I).
\end{equation}
Using the spectral decomposition of $\mathbf{J}$, we have
$\mathbf{J}^2=\sum_{j=1}^N\lambda_j^2\mathbf u_j\mathbf u_j^{\mathsf T}$.
Therefore
\[
\Tr(\mathbf P_I\mathbf J^2\mathbf P_I)
=
\sum_{j=1}^N\lambda_j^2
\Tr\bigl((\mathbf P_I\mathbf u_j)(\mathbf P_I\mathbf u_j)^{\mathsf T}\bigr)
=
\sum_{j=1}^N\lambda_j^2\|\mathbf P_I\mathbf u_j\|_2^2.
\]
Combining this with \eqref{eqn:local-trace-identity} gives
\begin{equation}
\label{eqn:local-energy-spectral-projection}
\sum_{i=\ell}^{m}b_i^2
+2\sum_{i=\ell}^{m-1}a_i^2
\le
\sum_{j=1}^N\lambda_j^2\|\mathbf P_I\mathbf u_j\|_2^2.
\end{equation}

We now split the last sum according to whether $\varphi_j\in I^{(D)}$.
For $\varphi_j\in I^{(D)}$, we use only the trivial bound
$\|\mathbf P_I\mathbf u_j\|_2^2\le\|\mathbf u_j\|_2^2=1$. Thus
\begin{equation}
\label{eqn:near-centers-bound}
\sum_{j:\,\varphi_j\in I^{(D)}}
\lambda_j^2\|\mathbf P_I\mathbf u_j\|_2^2
\le
\sum_{j:\,\varphi_j\in I^{(D)}}\lambda_j^2.
\end{equation}

For the complementary indices we use
\begin{equation}
\label{eqn:far-centers-bound}
\sum_{j:\,\varphi_j\notin I^{(D)}}
\lambda_j^2\|\mathbf P_I\mathbf u_j\|_2^2
\le
\left(\sup_{j:\,\varphi_j\notin I^{(D)}}
\|\mathbf P_I\mathbf u_j\|_2^2\right)
\sum_{j:\,\varphi_j\notin I^{(D)}}\lambda_j^2.
\end{equation}
Combining \eqref{eqn:local-energy-spectral-projection},
\eqref{eqn:near-centers-bound}, and \eqref{eqn:far-centers-bound} gives
\eqref{eqn:local-energy-main-bound}.
\end{proof}

To prove Lemma \ref{lem:local-lax-trace}, we require the following estimate on
the sums of squares of eigenvalues, whose localization centers are in a given
interval.

\begin{lem}
\label{lem:quadratic_lambda_sum}
Fix \(s=T\tau\in[0,T\log N]\). There is a constant \(C>0\) such
that, with overwhelming probability,
the following holds on the event \(\mathfrak B(s)\). For
 every index interval
$W\subseteq
\llbracket
N_1+T(\log N)^{7},
N_2-T(\log N)^{7}
\rrbracket$
with \(|W|\ge R^{-1}\mathfrak M\), one has
\begin{equation}\label{eqn:weighted_rectangle_index}
\sum_{j:\,\varphi_s(j)\in W}\lambda_j^2
\le
C|W|.
\end{equation}
\end{lem}

\begin{proof}
By a union bound
(and by summing the bounds obtained for smaller intervals), we may assume
that $W$ is fixed and $|W| = r \mathfrak M$ for some
$r\in[\mathfrak A^{-1/10},\mathfrak A^{10}]$.
We may also assume $W \subset \llbracket
N_1+T(\log N)^{15},
N_2-T(\log N)^{15}
\rrbracket$, as otherwise the lemma follows from Lemma
\ref{lem:edge-thermal-coupling}, through an argument
	similar to the one outlined
	at the end of the proof of Lemma \ref{lem:log_sum_dif}.

For \(n\ge0\), set $J_n:=[n,n+1]\cup[-n-1,-n]$.
Apply Lemma \ref{lem:weighted_rectangle} to the intervals
\(I,J_n\), where $I$ ranges over intervals such that $I \subset [-N^2/2, N^2/2]$,
$|I| \geq T^{-5/16}$, and both of its endpoints lie in $N^{-1} \mathbb{Z}$;
by a union bound, this gives an event
of overwhelming probability on which
\begin{equation}\label{eq:quadratic-shell-count}
\#\{i:Y_i(\tau)\in I,\ \Lambda_i\in J_n\}
\le
C\mathfrak A^{-1}T|I|e^{-cn^2}+C(\log N)^C
\end{equation}
simultaneously for each such interval $I$. In addition,
we restrict to the event of Lemma \ref{l0eigenvalues} with $A = \log N$.

By \eqref{eqn:qbd}, there are
\(x_-<x_+\) such that
\[
\{q_a(s):a\in W\}\subseteq[x_-,x_+],
\qquad
x_+-x_-\le Cr\mathfrak A\mathfrak M.
\]
For \(0\le n\le\lceil\log N\rceil\), let $I_n$ be the interval with left
endpoint $T^{-1}x_--C(n+2)\mathfrak M/T$ and right endpoint
$T^{-1}x_++C(n+2)\mathfrak M/T$.
Let \(\widetilde I_n\subset[-N^2/2, N^2/2]\) be an enlargement of \(I_n\) whose
endpoints lie in \(N^{-1}\mathbb Z\). It may be chosen so that
$|\widetilde I_n|\le C(r\mathfrak A+n+2)\mathfrak M/T$.

By \eqref{eqn:Qx},
which holds on $\mathfrak B(s)$, $\varphi_s(j) \in W, \lambda_j \in J_n$ implies
$Y_{\varphi_0(j)}(\tau)\in I_n\subseteq\widetilde I_n$. Hence, by
\eqref{eq:quadratic-shell-count},
\[
\sum_{j:\,\varphi_s(j)\in W}\lambda_j^2
\le C\sum_{n=0}^{\lceil\log N\rceil}(n+1)^2
\left[\mathfrak A^{-1}T|\widetilde I_n|e^{-cn^2}+(\log N)^C\right]
\le Cr\mathfrak M+(\log N)^C.
\]
Since $|W|=r\mathfrak M$ and $|W|\ge R^{-1}\mathfrak M$, the
polylogarithmic term is absorbed into $C|W|$. This proves
\eqref{eqn:weighted_rectangle_index}.
\end{proof}

\begin{proof}[Proof of Lemma \ref{lem:local-lax-trace}]
Fix $L$ and $I$ as in the statement, set
$\widehat L:=\max\{L,M_*\}$.
We work on the event from Lemma \ref{lem:quadratic_lambda_sum}, together with
the overwhelming probability events defined throughout the argument below.
For eigenvectors $\mathbf u_j(t)$ whose
	 localization centers $\varphi_t(j)$ lie in $\llbracket N_1 + T (\log N)^{10}, N_2 - T (\log N)^{10} \rrbracket \setminus I^{(5M_*)}$, \eqref{eqn:uMR}, which holds on $\mathfrak B(t)$, gives
\begin{equation}\label{eq:local-lax-leakage}
\max_{r\in I}|u_j(r;t)|
\le e^{\mathfrak M-5RM_*}
\le e^{-R^{-1}\mathfrak M},
\end{equation}
where we used $M_*=(2R)^{-1}\mathfrak M$. If, on the other hand, $\varphi_t(j) \notin \llbracket N_1 + T (\log N)^{10}, N_2 - T (\log N)^{10} \rrbracket$,
		then since $I^{(5M_*)}\subseteq
		\llbracket N_1+\frac32T(\log N)^{10},\,
		N_2-\frac32T(\log N)^{10}\rrbracket$, it must be
	the case that $\varphi_t(j)$ is at distance at least $T (\log N)^5$ from $I^{(5M_*)}$, and in this case
	Proposition \ref{prop:outside-strip-confinement} (upon further restriction 
	to an overwhelmingly probable event) again verifies
	\eqref{eq:local-lax-leakage}.
Lemma \ref{lem:local-energy-from-localization}, applied with $C_0=5/2$ and
with the $\mathbf J$ there equal to the Toda Lax matrix $\L(t)$ here, together with
	further restriction to the event $\cap_{t \geq 0} \mathsf{BND}_{\L(t)}(\log N)$ (from Lemma \ref{l0eigenvalues}), then gives
\[
\operatorname{Tr}\bigl((\L^{I}(t))^2\bigr)
\le
\sum_{\substack{j\in\llbracket1,N\rrbracket:\,
\varphi_t(j)\in I^{(5M_*)}}}\lambda_j^2
+ C e^{-c M_*} .
\]
Since
$R^{-1}\mathfrak M\le |I^{(5M_*)}|\le C\widehat L$, Lemma
\ref{lem:quadratic_lambda_sum} gives
\[
\sum_{\substack{j\in\llbracket1,N\rrbracket:\,
\varphi_t(j)\in I^{(5M_*)}}}\lambda_j^2
\le
C\widehat L.
\]
The combination of the preceding two bounds proves the claim.
\end{proof}

	\section{Proofs of remaining auxiliary results}
	\label{sec:fp}
	
	The goal of this section is to prove Proposition \ref{prop:cont} and Lemmas \ref{lem:cl},
	\ref{lem:Qtox}, and \ref{lem:sumtoint}. Throughout, we adopt Assumption
	\ref{ass:model} and the notation of Definition \ref{def:quasi_intro}.
	
	\subsection{Proof of Proposition \ref{prop:cont}}
	\label{subsec:proof-prop-cont}
	
	\begin{proof}[Proof of Proposition \ref{prop:cont}]
Let $
\mathcal G_N
:=
[0,T\log N]\cap N^{-1}\mathbb Z$,
and define the index intervals
$\mathcal B_{10}
:=
\llbracket
N_1+T(\log N)^{10},\,
N_2-T(\log N)^{10}
\rrbracket$, $
\mathcal B_9
:=
\llbracket
N_1+T(\log N)^9,\,
N_2-T(\log N)^9
\rrbracket$.
We restrict to $\bigcap_{r\ge0}\mathsf{BND}_{\L(r)}(\log N)$, which is overwhelmingly
probable by Lemma \ref{l0eigenvalues}, and to the event on which Proposition
\ref{prop:outside-strip-confinement} holds. For each
\(s\in\mathcal G_N\), let us also
restrict to the
events of Propositions \ref{vni_out} and
\ref{prop:qbdtfix}, with the time \(t\) there given by
$s$ here. Since $
|\mathcal G_N|
\le
C N T\log N
\le N^3$,
a union bound gives a single event of overwhelming probability on
which their conclusions hold at every \(s\in\mathcal G_N\) satisfying
\(\mathfrak{B}(s)\).

First, we provide a preliminary bound on localization centers. Simultaneously for every \(k\) and every
\(r\in[0,T\log N]\),
\begin{equation}
\label{eq:prop-cont-global-center}
|\varphi_r(k)-\varphi_0(k)|
\le
T(\log N)^3.
\end{equation}
Indeed, otherwise \eqref{eqn:umphi}, evaluated with $s'=0$ at
\(m=\varphi_r(k)\), would give 
$
|u_k(\varphi_r(k);r)|
\le
\exp\left(
-4RT(\log N)^3
\right)
<
\zeta$
for sufficiently large \(N\), which contradicts the definition of
\(\varphi_r(k)\) as a \(\zeta\)-localization center.

In particular, if $\varphi_{t'}(k)\in\mathcal B_{10}$
and \(|t-t'|\le2N^{-1}\), then
\[
|\varphi_t(k)-\varphi_{t'}(k)|
\le 2T(\log N)^3,
\]
and thus we obtain, for all sufficiently large \(N\),
\begin{equation}
\label{eq:prop-cont-weak-bulk}
\varphi_t(k)
\in
\llbracket
N_1+\tfrac12T(\log N)^{10},\,
N_2-\tfrac12T(\log N)^{10}
\rrbracket.
\end{equation}
This proves \eqref{eqn:lc_weak}. 

\medskip
\noindent
\emph{Proof of Item 1}: The Toda equations give $
\partial_rq_i(r)=b_i(r)$.
On the event 
$\bigcap_{r\ge0}\mathsf{BND}_{\L(r)}(\log N)$, 
$|b_i(r)|\le\log N$ for every \(i\) and \(r\ge0\). Consequently, if
\(|t-t'|\le2N^{-1}\), then
\[
|q_i(t)-q_i(t')|
\le
\int_{t\wedge t'}^{t\vee t'}|b_i(r)|\,dr
\le
\frac{2\log N}{N}.
\]
This proves \eqref{eqn:qbd_cont}.

\medskip
\noindent
\emph{Proof of Item 2}: Fix \(t\in[0,T \log N] \cap N^{-1} \mathbb{Z}\) satisfying property $\mathfrak{B}(t)$ and
\(t'\in[0,T\log N]\) with
$|t-t'|\le2N^{-1}$ and
$\varphi_{t'}(k)\in\mathcal B_{10}$.
Put $
p_t:=\varphi_t(k)$.
By \eqref{eq:prop-cont-weak-bulk} and the hypothesis on
\(\varphi_{t'}(k)\), we have \(p_t\in\mathcal B_9\). Proposition
\ref{vni_out}, applied at the mesh time \(t\), together with the trivial bound \(|u_k(\ell;t)|\le1\), implies
the global estimate
\begin{equation}
\label{eq:prop-cont-mesh-envelope}
|u_k(\ell;t)|
\le
\exp\left\{
\frac{\mathfrak M}{2}
-
R|\ell-p_t|
\right\}
\qquad
\text{for each }\ell\in\llbracket N_1,N_2\rrbracket.
\end{equation}

For each $r \in \mathbb{R}_{\geq0}$, let \(\mathbf{V}(r,t)\) denote the orthogonal matrix
with the property that for each $k \in \llbracket 1, N \rrbracket$,
$
\mathbf u_k(r)=\mathbf{V}(r,t) \mathbf u_k(t)$.
Let $\mathbf W_{p_t}$ denote the $N\times N$ diagonal matrix, with rows and
columns indexed by $\llbracket N_1,N_2\rrbracket$, whose $(j,j)$ entry is
equal to $e^{R|j-p_t|}$ for each
$j\in\llbracket N_1,N_2\rrbracket$. Thus,
\eqref{eq:prop-cont-mesh-envelope} is equivalent to
$\|\mathbf W_{p_t}\mathbf u_k(t)\|_\infty\le e^{\mathfrak M/2}$.
By \cite[Section 2]{Mos75}, the skew-symmetric, tridiagonal matrix
$\mathbf P(r)=[P_{i,j}(r)]$ with
$P_{i,i+1}(r)=a_i(r)/2$ satisfies
$\partial_r\mathbf L(r)=[\mathbf L(r),\mathbf P(r)]$. Hence, for the
matrix $\mathbf V(r,t)$ defined above,
$\partial_r\mathbf V(r,t)=-\mathbf P(r)\mathbf V(r,t)$. On
$\bigcap_{r\ge0}\mathsf{BND}_{\L(r)}(\log N)$,
\[
\left\|
\mathbf{W}_{p_t} \mathbf{P}(r) \mathbf{W}_{p_t}^{-1}
\right\|_{\ell^\infty\to\ell^\infty}
\le
C e^R\log N.
\]
Grönwall's inequality therefore gives
\begin{equation}
\label{eq:prop-cont-weighted-propagator}
\left\|
\mathbf{W}_{p_t} \mathbf{V}(r,t) \mathbf{W}_{p_t}^{-1}
\right\|_{\ell^\infty\to\ell^\infty}
\le
\exp\left\{
C e^R(\log N)|r-t|
\right\}.
\end{equation}
For \(r=t'\), we have \(|r-t|\le2N^{-1}\), so the right hand side is at
most \(2\) for sufficiently large \(N\). Combining
\eqref{eq:prop-cont-mesh-envelope} and
\eqref{eq:prop-cont-weighted-propagator}, we obtain
\begin{equation}
\label{eq:prop-cont-propagated-envelope}
|u_k(j;r)|
\le
2\exp\left\{
\frac{\mathfrak M}{2}
-
R|j-p_t|
\right\},
\qquad
r\in\{t,t'\}.
\end{equation}

Let  $p_{t'}:=\varphi_{t'}(k)$.
Since \(p_{t'}\) is a \(\zeta\)-localization center, $
|u_k(p_{t'};t')|\ge\zeta$.
Applying \eqref{eq:prop-cont-propagated-envelope} at \(j=p_{t'}\) gives, since
$\log\zeta^{-1}
\leq
100(\log N)^{3/2}\leq\mathfrak M/4$,
\begin{equation}
\label{eq:prop-cont-center-movement}
|p_{t'}-p_t|
\le
\frac{\mathfrak M}{R}.
\end{equation}

Taking \(r=t'\) in \eqref{eq:prop-cont-propagated-envelope}, and using
\eqref{eq:prop-cont-center-movement}, yields
\[
|u_k(j;t')|\le2\exp\left\{\frac{\mathfrak M}{2}
+R|p_{t'}-p_t|-R|j-p_{t'}|\right\}
\le2\exp\left\{\frac{3\mathfrak M}{2}
-R|j-\varphi_{t'}(k)|\right\}.
\]
After choosing \(\mathfrak A\) sufficiently large,
$2e^{3\mathfrak M/2}
\le
e^{\mathfrak A\mathfrak M}$, and
 \eqref{eqn:uMRct} follows.

\medskip
\noindent
\emph{Proof of Item 3}: We next prove \eqref{eqn:Qxct}. We again let $ t \in \mathcal{G}_N$
and assume $\mathfrak{B}(t)$ holds.
 We first record a consequence of Proposition
\ref{prop:qbdtfix}. If \(a,b\in\mathcal B_9\) and
$|a-b|\le\frac{2\mathfrak M}{R}$,
then, at the mesh time \(t\),
\begin{equation}
\label{eq:prop-cont-absolute-q-spacing}
|q_a(t)-q_b(t)|
\le
C\mathfrak A\frac{\mathfrak M}{R}.
\end{equation}

By \eqref{eq:prop-cont-center-movement}, 
$|p_t-p_{t'}|
\le
\frac{2\mathfrak M}{R}$.
Thus \eqref{eq:prop-cont-absolute-q-spacing}, applied with
\(a=p_t\) and \(b=p_{t'}\), gives
$|q_{p_t}(t)-q_{p_{t'}}(t)|
\le
C\mathfrak A\frac{\mathfrak M}{R}$.
Moreover, using \(\partial_rq_i(r)=b_i(r)\) and our restriction to
$\mathsf{BND}_{\L(0)}(\log N)$,
\[
|q_{p_{t'}}(t')-q_{p_{t'}}(t)|
\le
\frac{C\log N}{N}.
\]
Since
$Q_k(t)=q_{p_t}(t)$ and
$Q_k(t')=q_{p_{t'}}(t')$,
we conclude that
\begin{equation}
\label{eq:prop-cont-Q-raw}
|Q_k(t)-Q_k(t')|
\le
C\mathfrak A\frac{\mathfrak M}{R}
+\frac{C\log N}{N}.
\end{equation}
Recalling that \(R=\mathfrak A^{1/10}\),
and using \(\log(|\lambda_k|+2)\ge\log2\),
\eqref{eq:prop-cont-Q-raw} implies
\eqref{eqn:Qxct}.

\medskip
\noindent

To prove \eqref{eqn:Qxct2}, observe that
\[
\begin{aligned}
T\left|x_\tau^{(T)}(\lambda_k;\bar q_{\varphi_0(k)})
-x_{\tau'}^{(T)}(\lambda_k;\bar q_{\varphi_0(k)})\right|
&\le CT(1+\log N)|\tau-\tau'|
=C(1+\log N)|t-t'|
\le\frac{C\log N}{N},
\end{aligned}
\]
where the first inequality follows from \eqref{eqn:xtode},
\eqref{eqn:vebd}, and our restriction to
$\mathsf{BND}_{\L(0)}(\log N)$; the second follows from the definitions of
$t$ and $t'$; and the third follows from the fact that
$|t-t'|\le2N^{-1}$. This proves
\eqref{eqn:Qxct2} and completes the proof.
\end{proof}

	\subsection{Proof of Lemma \ref{lem:cl}}
	\label{subsec:local-charge-intro}

	\begin{proof}[Proof of Lemma \ref{lem:cl}]
Let us restrict to the overwhelming probability event 
given by the intersection of the events of Theorem
\ref{ttn}, Lemma \ref{l0eigenvalues} with $A = \log N$, and
Lemma \ref{lem:Estr}, and Proposition
\ref{prop:outside-strip-confinement}.  Put
\[
S:=\{i:q_i(t)\in J\}.
\]
If \(S=\varnothing\), then also
\(\{k:Q_k(t)\in J\}=\varnothing\), and there is nothing to prove.
Otherwise, let \(i_0=\min S\), \(i_1=\max S\), and
\(I=\llbracket i_0,i_1\rrbracket\).  The assumption on
\(J\), Lemma \ref{lem:Estr}, and the fact that for each $i \in \llbracket N_1, N_2 \rrbracket$,
\(|q_i(t)-q_i(0)|\le T(\log N)^2\) (which
holds by our restriction to $\cap_{r \geq 0}\mathsf{BND}_{\L(r)}(\log N)$) ensure that \(I\) lies in the
index range where \eqref{eqn:qbd-lower} applies.

If \(h\in I\setminus S\) and \(q_h(t)<\inf J\), then
\(h-i_0\le\mathfrak M/R\), since otherwise
\eqref{eqn:qbd-lower} would give \(q_h(t)>q_{i_0}(t)\).
Similarly, if \(q_h(t)>\sup J\), then
\(i_1-h\le\mathfrak M/R\).  Hence
$|I\setminus S|\le 2\left\lceil\frac{\mathfrak M}{R}\right\rceil$.
Since \(|L_{i,j}(t)|\le \log N\) by the restriction to $\mathsf{BND}_{\L(t)}(\log N)$,
\begin{equation}\label{eq:cl-hull-error}
\left|
\sum_{i\in S}(\L(t)^m)_{ii}
-
\sum_{i\in I}(\L(t)^m)_{ii}
\right|
\le C\mathfrak M(C\log N)^m.
\end{equation}

Let $\mathbf P_I$ be the orthogonal projection onto
$\operatorname{span}\{\mathbf e_i:i\in I\}$, where for each
$i\in\llbracket N_1,N_2\rrbracket$ we let $\mathbf e_i\in\mathbb R^N$
denote the coordinate vector whose entry indexed by $i$ is equal to $1$ and
whose other entries are equal to $0$. By the spectral
decomposition,
\[
\sum_{i\in I}(\L(t)^m)_{ii}
=
\sum_{k=1}^N\lambda_k^m\| \mathbf P_I \mathbf u_k(t)\|_2^2.
\]
Consequently, using restriction to $\mathsf{BND}_{\L(t)}(\log N)$,
\[
\begin{aligned}
&\left|
\sum_{i\in I}(\L(t)^m)_{ii}
-
\sum_{\varphi_t(k)\in I}\lambda_k^m
\right| \le
(C\log N)^m
\left[
\sum_{\varphi_t(k)\in I}\|\mathbf P_{I^{\complement}} \mathbf u_k(t)\|_2^2
+
\sum_{\varphi_t(k)\notin I}\|\mathbf P_I \mathbf u_k(t)\|_2^2
\right].
\end{aligned}
\]
Set \(D=\lceil2\mathfrak M/R\rceil\).
There are at most $5D$ indices $k$ for which
$\varphi_t(k)\in
\llbracket i_0-D,i_0+D\rrbracket\cup
\llbracket i_1-D,i_1+D\rrbracket$.
For the remaining $k$,
\eqref{eqn:uMR} (and
\eqref{eqn:umphi}, to handle localization centers near $N_1$ or $N_2$)
and summation of the geometric tails give a total
contribution \(O(1)\).  Thus the bracketed expression is at most
\(C(D+1)\le C\mathfrak M\), and hence
\begin{equation}\label{eq:cl-spectral-error}
\left|
\sum_{i\in I}(\L(t)^m)_{ii}
-
\sum_{i\in I}\Lambda_i(t)^m
\right|
\le C\mathfrak M(C\log N)^m.
\end{equation}
The bound on \(\left| I\setminus S \right|\) gives
\[
\left|
\sum_{i\in I}\Lambda_i(t)^m
-
\sum_{i\in S}\Lambda_i(t)^m
\right|
\le C\mathfrak M(C\log N)^m.
\]
Recalling the definition of $S$, combining the last estimate with \eqref{eq:cl-hull-error} and \eqref{eq:cl-spectral-error} proves \eqref{eqn:cl}.
\end{proof}

\subsection{Proof of Lemma \ref{lem:Qtox} and  Lemma \ref{lem:sumtoint}}
\label{subsec:Qtox}

For Lemmas \ref{lem:Qtox} and \ref{lem:sumtoint}, we fix
\(\tau\in[0,\log N]\), an interval \(I_T\), and a function \(F\) as in
Assumption \ref{ass:Fas}, write \(t=T\tau\) and \(K_T=T^{-1}I_T\), and
abbreviate \(\mathcal A=\mathcal A(F)\), \(\mathcal V=\mathcal V(F)\),
\(\mathcal D=\mathcal D(F)\), and \(\mathcal B=\mathcal B(F)\).
	\label{subsec:sumtoint}
	\begin{proof}[Proof of Lemma \ref{lem:Qtox}]
	We restrict to the intersection of the events in Theorem \ref{ttn},
	Lemma \ref{lem:Estr}, Proposition \ref{prop:outside-strip-confinement},
	and $\cap_{s\geq0}\mathsf{BND}_{\L(s)}(\log N)$.  Set
	$Y_k:=x_\tau^{(T)}(\lambda_k;\bar q_{\varphi_0(k)})$.
	The first sum in
	\eqref{eqn:FLQx} equals $\sum_{k=1}^N F(\lambda_k,Y_k)$.

	Let
	\[
	\mathcal S:=\{k:Y_k\in K_T\text{ or }Q_k(t)/T\in K_T\}.
	\]
	The bound \eqref{eqn:vebd}, applied to the effective velocity $v_\sigma^{\mathrm{eff},(T)}$ defined from Lemma \ref{lem:int_ghd} 
	with $m^{(T)}$ in place of $m$, Lemma
	\ref{lem:Estr}, the restriction to $\cap_{s\geq0}\mathsf{BND}_{\L(s)}(\log N)$,
	and Proposition \ref{prop:outside-strip-confinement} 
	show that $\varphi_t(k) \in \llbracket N_1+T(\log N)^{10},N_2-T(\log N)^{10}\rrbracket$
	if $k \in \mathcal{S}$.  
	
	Consequently, \eqref{eqn:Qx} and $|\lambda_k|\leq\log N$ give
	\begin{equation}\label{eqn:TxiQi}
	|Y_k-Q_k(t)/T|\leq \rho,
	\qquad
	\rho:=C\mathfrak M T^{-1}\log N,
	\qquad k\in\mathcal S.
	\end{equation}
	For each $k\in\mathcal S$, at least one of $Y_k$ and $Q_k(t)/T$
	belongs to $K_T$ by the definition of $\mathcal S$. Thus
	\eqref{eqn:TxiQi} implies that $Q_k(t)/T$ belongs to the
	$\rho$-enlargement of $K_T$ for every $k\in\mathcal S$. By
	\eqref{eqn:qbd-lower},
	\begin{equation}\label{eqn:mcSbd}
	|\mathcal S|
	\leq C T(|K_T|+\rho)
	\leq C\left(T|K_T|+\mathfrak M\log N\right).
	\end{equation}

	Outside $\mathcal S$ both summands on the left side of
	\eqref{eqn:FLQx} vanish by the support assumption on
	$F$.  The definition of $\mathcal D(F)$ in Assumption \ref{ass:Fas}, \eqref{eqn:TxiQi}, and
	\eqref{eqn:mcSbd} therefore give
	\[
	\left|\sum_{k=1}^N
	\bigl(F(\lambda_k,Y_k)-F(\lambda_k,Q_k(t)/T)\bigr)\right|
	\leq \mathcal D(F)\rho|\mathcal S|
	\leq \mathcal D(F)T^{3/4}(\log N)^C,
	\]
	where we used $|K_T|\leq2(\log N)^{19}$ and
	$\mathfrak M\leq T^{3/4}$.
	\end{proof}



	\begin{proof}[Proof of Lemma \ref{lem:sumtoint}]
	Throughout this proof, abbreviate $\mathcal X(F)=\mathcal X$ for each
	$\mathcal X\in\{\mathcal A,\mathcal V,\mathcal B,\mathcal D\}$. Set
	$L=\log N$ and
	\[
	\pi_L(\lambda):=\max(-L,\min(\lambda,L)),
	\qquad
	h_{\tau,L}(\lambda,q)
	:=
	F\bigl(\pi_L(\lambda),
	x_\tau^{(T)}(\pi_L(\lambda);q)\bigr).
	\]
	Also set
	$h_\tau(\lambda,q):=F(\lambda,x_\tau^{(T)}(\lambda;q))$.
	We restrict to the overwhelmingly probable event in Lemma \ref{l0eigenvalues} with
	$A=\log N$. On this event,
	\begin{equation}\label{eqn:hhL}
	h_{\tau,L}(\lambda,q)=h_\tau(\lambda,q)
	\qquad
	\text{for every $\lambda\in\eig\L(0)$ and $q\in\mathbb R$}.
	\end{equation}
	We will apply Proposition \ref{prop:app-time-zero-dos} to
	$h_{\tau,L}$, so we check the conditions of Assumption
	\ref{ass:admissible-test}.

	First,
	\[
	\partial_qh_{\tau,L}(\lambda,q) =
	\partial_xF\bigl(\pi_L(\lambda),
	x_\tau^{(T)}(\pi_L(\lambda);q)\bigr)
	\partial_qx_\tau^{(T)}(\pi_L(\lambda);q).
	\]
	By \eqref{eqn:Xtlin} and \eqref{eqn:avrbd2coarse}, the map
	$q\mapsto x_\tau^{(T)}(\pi_L(\lambda);q)$ is increasing. Therefore the
	substitution
	\[
	x=x_\tau^{(T)}(\pi_L(\lambda);q),
	\qquad
	dx=\partial_qx_\tau^{(T)}(\pi_L(\lambda);q)dq
	\]
	gives the exact identity
	\[
	\int_{\mathbb R}|\partial_qh_{\tau,L}(\lambda,q)|dq
	=
	\int_{\mathbb R}|\partial_xF(\pi_L(\lambda),x)|dx.
	\]
	In particular,
	\begin{equation}\label{eq:sumtoint-q-variation}
	\sup_{\lambda\in\mathbb R}
	\int_{\mathbb R}|\partial_qh_{\tau,L}(\lambda,q)|\,dq
	\le \mathcal V.
	\end{equation}

	Next, for almost every $\lambda\in\mathbb R$,
	\begin{equation}\label{eqn:cr}
	\begin{aligned}
	\partial_\lambda h_{\tau,L}(\lambda,q)
	={}&\pi_L'(\lambda)\Bigl[
	\partial_\lambda F\bigl(\pi_L(\lambda),
	x_\tau^{(T)}(\pi_L(\lambda);q)\bigr)\\
	&\quad+
	\partial_xF\bigl(\pi_L(\lambda),
	x_\tau^{(T)}(\pi_L(\lambda);q)\bigr)
	\partial_\lambda x_\tau^{(T)}(\pi_L(\lambda);q)
	\Bigr].
	\end{aligned}
	\end{equation}
	By \eqref{eq:flow-spec}, applied to $m^{(T)}$, the almost everywhere
	derivative $\partial_\lambda x_\tau^{(T)}(\lambda;q)$ satisfies
	\begin{equation}\label{eq:sumtoint-characteristic-lambda}
	\sup_{\substack{|\lambda|\le L\\q\in\mathbb R}}
	|\partial_\lambda x_\tau^{(T)}(\lambda;q)|
	\le
	C\tau(1+L).
	\end{equation}
	Since $|\pi_L'|\le1$ almost everywhere,
	\eqref{eq:sumtoint-characteristic-lambda} and
	\eqref{eqn:cr} give
	\begin{equation}\label{eq:sumtoint-lambda-lipschitz}
	\|h_{\tau,L}\|_{\Lip,\lambda}
	\le
	\mathcal B+
	C\mathcal D\,\tau(1+L)
	\le
	(\mathcal D+\mathcal B)(\log N)^C.
	\end{equation}

	Moreover, \(\|h_{\tau,L}\|_\infty\le \mathcal A\). The characteristic flow equation
	\eqref{eqn:xtode}, together with \eqref{eqn:vebd}, gives
	\begin{equation}\label{eqn:xtTqdif}
	|x_\tau^{(T)}(\lambda;q)-q|
	\le C\tau(1+L),
	\qquad |\lambda|\le L.
	\end{equation}
	As a consequence, for any $\lambda\in[-\log N,\log N]$, the set of
	$q$ for which $x_\tau^{(T)}(\lambda;q)$ belongs to
	$T^{-1}I_T\subset
	[-(\log N)^{19},(\log N)^{19}]$ is contained in
	\([-(\log N)^{20},(\log N)^{20}]\) for all sufficiently large
		\(N\), and the latter interval is contained in $[N_1/T,N_2/T]$.
	Hence \(h_{\tau,L}\) satisfies Assumption
		\ref{ass:admissible-test} with
		\[
		A_h\le C(\mathcal A+\mathcal V),
		\qquad
		B_h\le C(\mathcal B+\mathcal D)(\log N)^C.
		\]

		Proposition \ref{prop:app-time-zero-dos}, specifically
		\eqref{eqn:hiTspace} (whose additional condition on $I$ is verified by
		the preceding support bound), now gives, with overwhelming probability,
	\begin{align}
	&\left|
	\sum_{i=N_1}^{N_2}h_{\tau,L}(\Lambda_i,\bar q_i)
	-
	T\iint_{\mathbb R^2}
	h_{\tau,L}(\lambda,q)m_0^{(T)}(\lambda,q)\,d\lambda\,dq
	\right| \notag\\
	&\qquad\le
	C(\log N)^C
	\left[
	(\mathcal{A} + \mathcal{V})T^{3/5}
	+(\mathcal{B} + \mathcal{D}) T^{2/5}
	\right].
	\label{eq:sumtoint-dos}
	\end{align}

	By the equality of $h_{\tau,L}$ and $h_\tau$ on the eigenvalues \eqref{eqn:hhL}, the sum in \eqref{eq:sumtoint-dos} is the sum in
	\eqref{eqn:Fconv}. On the other hand, Lemma
	\ref{rho0estimate}, \eqref{eqn:alph_bd}, and \eqref{eqn:xtTqdif} imply
	\[
	T\iint_{\mathbb R^2}
		|h_{\tau,L}(\lambda,q)-h_\tau(\lambda,q)|
		m_0^{(T)}(\lambda,q)\,d\lambda\,dq
		\le
		C\mathcal A (\log N)^C e^{-cL^2}.
		\]
	Finally, \eqref{eqn:m0push} gives
	\[
	\iint_{\mathbb R^2}
	h_\tau(\lambda,q)m_0^{(T)}(\lambda,q)\,d\lambda\,dq
	=
	\iint_{\mathbb R^2}
	F(\lambda,x)m_\tau^{(T)}(\lambda,x)\,d\lambda\,dx.
	\]
		Combining the previous two estimates with \eqref{eq:sumtoint-dos} proves \eqref{eqn:Fconv}.
		\end{proof}

\subsection{Proof of Lemma \ref{lem:xTtox}}
\label{subsec:xTtox}

\begin{proof}[Proof of Lemma \ref{lem:xTtox}]
Set $L_0:=  (\log T)^5$. By Assumption
\ref{ass:model}, Definition \ref{rho0}, and \eqref{eqn:alph_bd},
\begin{equation}
\label{eq:xTtox-initial-agreement}
m_0^{(T)}(\lambda,q)=m_0(\lambda,q),
\qquad |q|\le L_0/2.
\end{equation}
Work on the event from Lemma \ref{l0eigenvalues} with $A=\log N$. If
either $x_\tau^{(T)}(\Lambda_i;\bar q_i)$ or
$x_\tau(\Lambda_i;\bar q_i)$ belongs to $K_T$, then \eqref{eqn:vebd}
and $K_T\subset[-(\log N)^2,(\log N)^2]$ give
$|\bar q_i|\le C(\log N)^2$. For sufficiently large $T$, the hypotheses
of Lemma \ref{lem:thermal-completion} therefore hold with
$q = \bar q_i$, $x_0=0$, and
$L=L_0/2$, and that lemma gives
\[
\left|x_\tau^{(T)}(\Lambda_i;\bar q_i)
-x_\tau(\Lambda_i;\bar q_i)\right|
\le Ce^{-c(\log N)^2}.
\]
This, together with Assumption \ref{ass:Fas}, proves
\eqref{eqn:FxTx}.

For the second assertion, \eqref{eqn:m0push}, applied once to
$(m_\tau^{(T)},x_\tau^{(T)})$ and once to $(m_\tau,x_\tau)$, gives
\[
\begin{aligned}
&T\left|\iint_{\mathbb R^2}F(\lambda,x)
\bigl(m_\tau^{(T)}(\lambda,x)-m_\tau(\lambda,x)\bigr)
\,d\lambda\,dx\right|\\
&\qquad=
T\left|\iint_{\mathbb R^2}
\Bigl[
F\bigl(\lambda,x_\tau^{(T)}(\lambda;q)\bigr)m_0^{(T)}(\lambda,q)
-F\bigl(\lambda,x_\tau(\lambda;q)\bigr)m_0(\lambda,q)
\Bigr]\,d\lambda\,dq\right|.
\end{aligned}
\]
In the region $|\lambda|\le\log N$ where
either point $x_\tau^{(T)}(\lambda; q)$ or $x_\tau(\lambda;q)$ belongs to $K_T$, 
the bound \eqref{eqn:vebd} again gives
$|q|\le C(\log N)^2$, so (because \eqref{eq:xTtox-initial-agreement} holds) the
application of Lemma \ref{lem:thermal-completion} from the previous paragraph applies. The Gaussian bound
\eqref{eqn:avrbd}, and its analogue for $m^{(T)}$, controls the integral over the
 region $|\lambda| > \log N$. Consequently,
\[
T\left|\iint_{\mathbb R^2}F(\lambda,q)
\left(m_\tau^{(T)}(\lambda,q)-m_\tau(\lambda,q)\right)
\,d\lambda\,dq\right|
\le C(\mathcal A+\mathcal D)e^{-c(\log N)^2},
\]
which, after decreasing $c$, is \eqref{eqn:FmTm}.
\end{proof}
	
\appendix

	\section{Proof of Lemma \ref{lem:Estr}}
\label{app:Estr}

\begin{proof}

Since \(a_i(0)=e^{(q_i(0)-q_{i+1}(0))/2}\), we have
\[
q_{i+1}(0)-q_i(0)=-2\log a_i(0).
\]
For \(i > j\), by \eqref{eqn:alpha-regularity} and \eqref{eqn:Fstr},
\[
\bar q_i-\bar q_j-(F(i/T)-F(j/T))
=
T^{-1}\sum_{m=j}^{i-1}
\bigl(-2\log a_m(0)-\alpha(m/T)\bigr) + O_{\alpha,i,j}(T^{-1} ) ,  
\]
where $|O_{\alpha}(T^{-1} )| \leq C \mathfrak{A} |i-j| T^{-2} $ with a constant $C>0$ independent of $\mathfrak{A}$, and the analogous formula holds for \(i<j\).
This term is negligible compared to the right hand side in the event \eqref{eq:stretch-good} 
because $N \leq 2 T^2$, so it will ultimately be absorbed into that error. Write
\[
\xi_i:=-2\log a_i(0) - \mathbb{E}[-2\log a_i(0)] = -2\log a_i(0)-\alpha(i/T)
\]
for the centered stretch increments.

The random variables \(\xi_i\) are independent, and by
\eqref{munu} and \eqref{eqn:alph_bd} they admit a uniform
sub exponential moment: there exists \(\varepsilon_0>0\)
such that
\[
\sup_{i\in\llb N_1,N_2\rrb}\E e^{\varepsilon_0 |\xi_i| / \mathfrak{A}}<\infty.
\]

Bernstein's inequality for sub-exponential variables, \cite[Theorem 2.8.1]{Ver18}, applied to the partial
sums above and combined with a union bound over \(i,j\in\llb N_1,N_2\rrb\),
implies that the event \eqref{eq:stretch-good}
holds with overwhelming probability.

\end{proof}

\section{Proof of Lemma \ref{lem:vebd}}
\label{app:vebd}

\begin{proof}[Proof of Lemma \ref{lem:vebd}]
The first bound follows as in the proof of \cite[Corollary 3.8]{Agg25}.

Now we prove the second bound, using the first assertion. Fix
\(t\ge0, x \in \mathbb{R}\), and suppress \(t,x\) from the notation:
\[
m(\mu):=m_t(\mu,x),
\qquad
v(\lambda):=\ve_t(\lambda;x),
\qquad
J(\lambda):=1+(\mathbf Tm)(\lambda).
\]
After changing the constants in the Gaussian exponent, the bound
\eqref{eqn:avrbd} and the first bound on \(v\) from this lemma imply for
some $C,c>0$
\begin{equation}\label{eqn:mv-gaussian}
0\le m(\mu)\le Ce^{-c\mu^2},
\qquad
|m(\mu)v(\mu)|\le Ce^{-c\mu^2}.
\end{equation}
Moreover, by \eqref{eqn:J-fixed-point-lower},
\begin{equation}\label{eqn:J-positive-ve-proof}
J(\lambda)\ge c>0
\qquad
\text{for every }\lambda\in\mathbb R.
\end{equation}

We first record an elementary continuity estimate that follows from
\eqref{eqn:Tdf}. If \(g\) satisfies $|g(\mu)|\le Ce^{-c\mu^2}$, then,
whenever \(0<\delta:=|a-b|\le\frac12\),
\begin{equation}\label{eqn:log-potential-modulus}
|(\mathbf Tg)(a)-(\mathbf Tg)(b)|
\le
C\delta\log\frac1\delta.
\end{equation}

Next, note that \eqref{eqn:effective-velocity-linear-identity}, in the
notation above, reads
\begin{equation}\label{eqn:ve-linear-form}
J(\lambda)v(\lambda)
=
\lambda+\bigl(\mathbf T(mv)\bigr)(\lambda).
\end{equation}
Set
$F(\lambda)
:=
\lambda+\bigl(\mathbf T(mv)\bigr)(\lambda)$.
By \eqref{eqn:mv-gaussian} and
\eqref{eqn:log-potential-modulus}, for
\(a,b\in[-L,L]\) and \(0<|a-b|\le\frac12\),
\begin{align}
|F(a)-F(b)|
&\le
C|a-b|
\left(1+\log\frac1{|a-b|}\right),
\label{eqn:F-log-lip}\\
|J(a)-J(b)|
&\le
C|a-b|
\left(1+\log\frac1{|a-b|}\right),
\label{eqn:J-log-lip}\\
|F(b)|
&\le C(1+L).
\label{eqn:F-growth}
\end{align}
Since \(v=F/J\) and $J\ge c$, these estimates imply
\begin{equation}\label{eqn:ve-log-lip}
|v(a)-v(b)|
\le
C(1+L)|a-b|
\left(1+\log\frac1{|a-b|}\right).
\end{equation}
Using
$r\left(1+\log\frac1r\right)\le Cr^{1/2}$, $0<r\le\frac12$,
and the linear growth bound when \(|a-b|>\frac12\), we obtain the
preliminary H\"older estimate
\begin{equation}\label{eqn:ve-holder-half}
|v(a)-v(b)|
\le
C(1+L)|a-b|^{1/2},
\qquad
a,b\in[-L,L].
\end{equation}

We now use the cancellation in \eqref{eqn:ve-linear-form} to upgrade
\eqref{eqn:ve-holder-half} to a Lipschitz estimate. Subtracting
\eqref{eqn:ve-linear-form} at \(a\) and \(b\) gives the identity
\begin{equation}\label{eqn:ve-difference-identity}
J(b)\bigl(v(a)-v(b)\bigr)
=a-b
+
2\int_{\mathbb R}
\left(
\log|a-\mu|-\log|b-\mu|
\right)
m(\mu)\bigl(v(\mu)-v(a)\bigr)\,d\mu.
\end{equation}
Let
$\delta:=|a-b|\le\frac12,
r:=\mu-a,
h:=b-a$.
We divide the integral in \eqref{eqn:ve-difference-identity} into
three regions.

On \(|r|\le2\delta\), the H\"older estimate
\eqref{eqn:ve-holder-half}, applied on
\([-L-1,L+1]\), gives
\[
|v(a+r)-v(a)|
\le C(1+L)|r|^{1/2}.
\]
Moreover, by the change of variables \(r=\delta s\),
\[
\begin{aligned}
&\int_{|r|\le2\delta}
\left|
\log|r|-\log|r-h|
\right||r|^{1/2}\,dr\\
&\qquad=
\delta^{3/2}
\int_{|s|\le2}
\left|
\log|s|-\log|s-h/\delta|
\right||s|^{1/2}\,ds
\le C\delta^{3/2}.
\end{aligned}
\]
The last integral is finite because the logarithmic singularities are
locally integrable. Hence the contribution from this region is at most
$C(1+L)\delta^{3/2}
\le
C(1+L)\delta$.

On \(2\delta<|r|\le1\), the mean value theorem gives
\[
\left|
\log|r|-\log|r-h|
\right|
\le
C\frac{\delta}{|r|}.
\]
Using \eqref{eqn:ve-holder-half} again, the contribution from this
region is bounded by
\[
C(1+L)\delta
\int_{2\delta}^1r^{-1/2}\,dr
\le
C(1+L)\delta.
\]

Finally, on \(|r|>1\), the same kernel estimate gives
$\left|
\log|r|-\log|r-h|
\right|
\le
C\frac{\delta}{|r|}$.
The first bound on \(v\) in \eqref{eqn:vebd} yields
\[
\begin{aligned}
|v(a+r)-v(a)|
&\le
C(1+|a+r|)+C(1+|a|)\\
&\le
2C\bigl(|r|+1+L\bigr).
\end{aligned}
\]
Since \(|r|>1\), $\frac{|v(a+r)-v(a)|}{|r|}
\le C(1+L)$.
Therefore, using \(\int_{\mathbb R}m(\mu)\,d\mu\le C\), the contribution
from the final region is also at most
$C(1+L)\delta$.

Inserting these three estimates into
\eqref{eqn:ve-difference-identity} and using \(J(b)\ge c\), we obtain
\[
|v(a)-v(b)|
\le
C(1+L)|a-b|
\]
whenever \(|a-b|\le\frac12\). When
\(|a-b|>\frac12\), the first estimate in \eqref{eqn:vebd} gives
\[
\begin{aligned}
|v(a)-v(b)|
&\le
C(1+|a|)+C(1+|b|)\\
&\le
4C(1+L)|a-b|.
\end{aligned}
\]
This proves the second estimate in \eqref{eqn:vebd}.
\end{proof}

\section{General tridiagonal matrix facts}
\label{app:time-zero-dos}

\subsection{Transfer matrix formulas for eigenvectors}
\label{app:transfer}
	In this section we let $\mathbf{M}$ denote a symmetric, tridiagonal matrix. Next we discuss more precise behavior for the eigenvectors of $\mathbf{M}$ through the method of transfer matrices, following \cite[Section 4]{APO}. For any integer $k \in \llbracket N_1, N_2-1 \rrbracket$ and real number $E \in \mathbb{R}$, define the \emph{transfer matrix} $\mathbf{S}_k (E) = \mathbf{S}_k (E; \mathbf{M}) \in \Mat_{2\times2}$ by setting
	\begin{flalign}
		\label{ks} 
		\mathbf{S}_k (E)= \left[ \begin{array}{cc} 0 & 1 \\ -M_{k,k+1}^{-1} M_{k-1,k} & M_{k,k+1}^{-1} (E - M_{k,k}) \end{array} \right].
	\end{flalign}
	
	\noindent Moreover, for any subset $\mathcal{K} = (k_1, k_2, \ldots , k_m) \subseteq \llbracket N_1, N_2 - 1 \rrbracket$ with $k_1 < k_2 < \cdots < k_m$, define $\mathbf{S}_{\mathcal{K}} (E) = \mathbf{S}_{\mathcal{K}} (E; \mathbf{M}) \in \Mat_{2 \times 2}$ by setting
	\begin{flalign}
		\label{k2s} 
		\mathbf{S}_{\mathcal{K}} (E) = \mathbf{S}_{k_m} (E) \cdot \mathbf{S}_{k_{m-1}} (E) \cdots \mathbf{S}_{k_1} (E).
	\end{flalign} 
	
	The following (standard) lemma indicates how the above transfer matrices can be used to evaluate eigenvector entries of $\mathbf{M}$.
	
	\begin{lem}[{\cite[Lemma 3.9]{Agg25a}}]
		
		\label{smatrixu} 
		
		Let $\mu \in \eig \mathbf{M}$, and let $\bm{u} = (u_{N_1}, u_{N_1+1}, \ldots , u_{N_2}) \in \mathbb{R}^N$ denote an eigenvector of $\mathbf{M}$ with eigenvalue $\mu$. For each $k \in \llbracket N_1, N_2 \rrbracket$, set $\bm{w}_k = (u_{k-1}, u_k) \in \mathbb{R}^2$, where $u_{N_1-1} = 0$. Then $\mathbf{S}_{\llbracket i,j \rrbracket} (\mu) \cdot \bm{w}_i = \bm{w}_{j+1}$ for any $i, j \in \llbracket N_1, N_2 - 1 \rrbracket$ with $i \le j$.
		
	\end{lem}

	The next lemma expresses entries of $\mathbf{S}_{\llbracket i, j \rrbracket} (E)$ through eigenvalues of truncations of $\mathbf{M}$. In what follows, for any integers $k, \ell \in \llbracket N_1, N_2 \rrbracket$ with $k \le \ell$, we let $\mathbf{M}^{[k, \ell]}$ denote the $(\ell-k+1) \times (\ell-k+1)$ matrix obtained by restricting $\mathbf{M}$ to rows and columns indexed by $i, j \in \llbracket k, \ell \rrbracket$; stated equivalently, for each $i, j \in \llbracket k, \ell \rrbracket$, the $(i, j)$ entry of $\mathbf{M}^{[k,\ell]}$ is equal to the $(i,j)$ entry of $\mathbf{M}$.
	
	\begin{lem}[{\cite[Equation (1.65)]{OINL}}]
		
		\label{sij} 
		
		Set $\eig \mathbf{M}^{[i,j]} = ( \mu_1^{[i,j]}, \mu_2^{[i, j]}, \ldots , \mu_{j-i+1}^{[i, j]} )$ for any $i, j \in \llbracket N_1, N_2 \rrbracket$ with $i \le j$. Then, for any $i, j \in \llbracket N_1, N_2 - 1 \rrbracket$, denoting $\ell = j - i + 1$ we have
		\begin{flalign*}
			\mathbf{S}_{\llbracket i, j \rrbracket} (E) = \left[ \begin{array}{cc} - M_{i-1,i}  \displaystyle\prod_{k=i}^{j-1} M_{k,k+1}^{-1}  \displaystyle\prod_{h=1}^{\ell-2} ( E - \mu_h^{[i+1,j-1]} ) & \displaystyle\prod_{k=i}^{j-1} M_{k,k+1}^{-1}  \displaystyle\prod_{h=1}^{\ell-1} ( E - \mu_h^{[i,j-1]} ) \\
				- M_{i-1,i}  \displaystyle\prod_{k=i}^{j} M_{k,k+1}^{-1}   \displaystyle\prod_{h=1}^{\ell-1} ( E - \mu_h^{[i+1,j]} ) & \displaystyle\prod_{k=i}^{j} M_{k,k+1}^{-1} \displaystyle\prod_{h=1}^{\ell} ( E - \mu_h^{[i, j]} )  \end{array} \right].
		\end{flalign*}
		
	\end{lem} 
	
	\subsection{Proof of Lemma \ref{lem:suff}}
\label{app:suffpf}

First, we require an intermediate lemma, which is an exact identity. We again let $\mathbf{M}$ denote a symmetric, tridiagonal matrix.
\begin{lem}[{\cite[Equation (4)]{DSRL}}]
		
		\label{n1n2u}
		
		Let $\mu \in \eig \mathbf{M}$, and let $\bm{u} = (u_{N_1}, \ldots , u_{N_2}) \in \mathbb{R}^N$ denote a unit eigenvector of $\mathbf{M}$ with eigenvalue $\mu$. Then,
		\begin{flalign*}
			\log |u_{N_1}| + \log |u_{N_2}| = \displaystyle\sum_{i=N_1}^{N_2-1} \log |M_{i,i+1}| - \displaystyle\sum_{\mu' \in \eig \mathbf{M} \setminus \{ \mu \}} \log |\mu - \mu'|.
		\end{flalign*}
		
	\end{lem}

	\begin{proof}[Proof of Lemma \ref{lem:suff}]
	Let $u$ be a unit eigenvector for $\lambda$ such that
	$|u(\varphi)|\ge\zeta$.
	Choose $(q_i)_{i=N_1}^{N_2}$ so that
	$L_{i,i+1}=\exp((q_i-q_{i+1})/2)$. Adding the two identities in Lemma \ref{lem:vkn1_formula} at $i=\varphi$ gives
	\begin{align*}
	2\log|u(N_1)|+2\log|u(N_2)|
	={}&4\log|u(\varphi)|+q_{N_1}-q_{N_2}\\
	&-2\sum_{\mu\in\eig \H_-}\log|\lambda-\mu|
	-2\sum_{\mu\in\eig \H_+}\log|\lambda-\mu|.
	\end{align*}
	On the other hand, Lemma \ref{n1n2u} and
	$2\sum_{i=N_1}^{N_2-1}\log L_{i,i+1}=q_{N_1}-q_{N_2}$ yield
	\[
	2\log|u(N_1)|+2\log|u(N_2)|
	=q_{N_1}-q_{N_2}
	-2\sum_{\theta\in\eig\L\setminus\{\lambda\}}\log|\lambda-\theta|.
	\]
	Comparing the last two displays, we obtain
	\begin{equation}\label{eq:principal-minor-log-identity}
	\sum_{\mu\in\eig \H_-}\log|\lambda-\mu|
	+
	\sum_{\mu\in\eig \H_+}\log|\lambda-\mu|
	-
	\sum_{\theta\in\eig\L\setminus\{\lambda\}}\log|\lambda-\theta|
	=2\log|u(\varphi)|.
	\end{equation}
	Rearranging \eqref{eq:principal-minor-log-identity} gives
	\begin{align*}
	-\sum_{\mu\in\mathcal U_-}\log|\lambda-\mu|
	-\sum_{\mu\in\mathcal U_+}\log|\lambda-\mu|
	&=-\sum_{\theta\in\mathcal B}\log|\lambda-\theta|
	-2\log|u(\varphi)|\\
	&\qquad\quad+
	\sum_{\theta\in\mathcal A_-}
	\bigl(\log|\lambda-\pi_-(\theta)|-\log|\lambda-\theta|\bigr)\\
	&\qquad\quad+
	\sum_{\theta\in\mathcal A_+}
	\bigl(\log|\lambda-\pi_+(\theta)|-\log|\lambda-\theta|\bigr).
	\end{align*}
	Now use \eqref{eq:matched-log-error} and
	$-2\log|u(\varphi)|\leq-2\log\zeta$ to obtain
	\eqref{eq:unmatched-singular-bound}.
	\end{proof}

\subsection{Proof of Lemma \ref{lem:full-to-truncation}}
\label{app:ftt}

\begin{proof}[Proof of Lemma \ref{lem:full-to-truncation}]
Recall that $\mathcal I=\llbracket m,\ell\rrbracket
=\llbracket N_1,N_2\rrbracket\setminus\mathcal D$. For each
$i\in\llbracket N_1,N_2\rrbracket$, let $\mathbf e_i\in\mathbb R^N$ be the
standard coordinate vector indexed by $i$, and let $\mathbf P$ denote the
orthogonal projection onto
$\operatorname{span}\{\mathbf e_i:i\in\mathcal I\}$. In addition, throughout
the proof, let $K=\log N$ and $p=\varphi_j$, and set
$\partial\mathcal I:=\{a\in\llbracket N_1,N_2\rrbracket:
\max\{\dist(a,\mathcal D),\dist(a,\mathcal I)\}\le1\}$, which is exactly
the set of indices covered by \eqref{eqn:lmatrix-boundary-decay}.

Let $\mathbf u=\mathbf u_j$ and let
$\mathbf v=\mathbf P\mathbf u=(v(i))_{i=N_1}^{N_2}$. We have
\begin{equation}\label{eqn:tildeLlambdav}
 (\widetilde \L - \lambda_j \mathbf{I}) \mathbf v =
\mathbf P (\widetilde \L - \L) \mathbf P\mathbf u
+ (\mathbf{I}-\mathbf P) \widetilde \L\mathbf P\mathbf u
- \mathbf P\L (\mathbf{I} - \mathbf P)\mathbf u.
\end{equation}
Note that 
\begin{equation}\label{eqn:tildeLlambdav1}
\|\mathbf P (\widetilde \L - \L)\mathbf P\mathbf u\|_2
\leq \sqrt{3}\delta\left(\sum_{i=N_1}^{N_2}
\sum_{|k-i|\leq1}|u_j(k)|^2\right)^{1/2}\leq3\delta,
\end{equation}
since $\|\mathbf u\|_2=1$ and by the second bound in \eqref{ll};
\begin{equation}\label{eqn:tildeLlambdav2}
\|(\mathbf{I}-\mathbf P)\widetilde \L\mathbf P\mathbf u\|_2
\leq C K\max_{a\in\partial\mathcal I}|u_j(a)|
\leq C K e^{-c_0 \dist(p, \mathcal{D})},
\end{equation}
by the first bound in \eqref{ll} and \eqref{eqn:lmatrix-boundary-decay}, as
$|a-p|\ge\dist(p,\mathcal D)-1$ for $a\in\partial\mathcal I$.
By similar reasoning,
\begin{equation}\label{eqn:tildeLlambda3}
\|\mathbf P\L(\mathbf{I}-\mathbf P)\mathbf u\|_2
\leq C K\max_{a\in\partial\mathcal I}|u_j(a)|
\leq C K e^{-c_0 \dist(p, \mathcal{D})}.
\end{equation}
Thus, combining \eqref{eqn:tildeLlambdav1}-\eqref{eqn:tildeLlambda3} with \eqref{eqn:tildeLlambdav}, we obtain
\begin{equation}\label{eqn:bd}
\|(\widetilde \L-\lambda_j\mathbf I)\mathbf v\|_2
\leq C(Ke^{-c_0\dist(p,\mathcal D)}+\delta).
\end{equation}
Let 
\begin{equation}
\mathbf w = \frac{\mathbf v}{\|\mathbf v\|_2}=(w(i))_{i=N_1}^{N_2}.
\end{equation}
Since $|v(p)|\geq\zeta$ and $\|\mathbf v\|_2\leq1$, we know
$|w(p)|\geq\zeta$. In addition,
\begin{equation}\label{eqn:tildLwbd}
 \|(\widetilde \L-\lambda_j\mathbf I)\mathbf w\|_2
 \leq C\zeta^{-1}(Ke^{-c_0\dist(p,\mathcal D)}+\delta)=:\alpha,
\end{equation}
since $\|\mathbf v\|_2\geq\zeta$.

Now, let $\boldsymbol\psi_a=(\psi_a(i))_{i=N_1}^{N_2}$,
$a\in\llbracket1,N\rrbracket$, denote an eigenbasis of $\widetilde \L$, such that
$\widetilde \L\boldsymbol\psi_a=\mu_a\boldsymbol\psi_a$. We have the
decomposition
\begin{equation}
\mathbf w=\sum_{a=1}^N c_a\boldsymbol\psi_a.
\end{equation}
Let $R = 2 \alpha / \zeta$, and let
$\mathbf Q = \mathbf{1}_{[\lambda_j-R,\lambda_j+R]}(\widetilde \L)$ denote the
spectral projection of $\widetilde \L$ onto
$[\lambda_j-R,\lambda_j+R]$. We have
\begin{equation}\label{eqn:IQz2}
\|(\mathbf I-\mathbf Q)\mathbf w\|_2^2
\leq \sum_{a : |\mu_a - \lambda_j| > R} |c_a|^2
\leq R^{-2} \sum_{a =1}^N |c_a|^2|\mu_a-\lambda_j|^2
\leq R^{-2} \alpha^2 = \zeta^2/4.
\end{equation}
In the third inequality, we used \eqref{eqn:tildLwbd} to see that
$\|(\widetilde \L-\lambda_j\mathbf I)\mathbf w\|_2^2
= \sum_{a =1}^N |c_a|^2 |\mu_a-\lambda_j|^2 \leq \alpha^2$.
In particular, \eqref{eqn:IQz2} implies
$|((\mathbf I-\mathbf Q)\mathbf w)(p)|\leq\zeta/2$; combined with
$|w(p)|\geq\zeta$, this gives
\begin{equation}
|(\mathbf Q\mathbf w)(p)| \geq \zeta / 2.
\end{equation}
Since $\mathbf Q\mathbf w=
\sum_{a:|\mu_a-\lambda_j|\leq R}c_a\boldsymbol\psi_a$,
Cauchy--Schwarz gives
\[
\frac{\zeta^2}{4}\leq |(\mathbf Q\mathbf w)(p)|^2
\leq \left(\sum_{a:|\mu_a-\lambda_j|\leq R}|c_a|^2\right)
\left(\sum_{a:|\mu_a-\lambda_j|\leq R}|\psi_a(p)|^2\right)
\leq \sum_{a:|\mu_a-\lambda_j|\leq R}|\psi_a(p)|^2.
\]
Thus, at least one $a$ satisfying $|\mu_a - \lambda_j| \leq R$ must satisfy $|\psi_a(p)|^2 \geq \zeta^2/(4 N)$. This completes the proof.
\end{proof}

\section{Proof of Lemma \ref{lem:thermal-tail-quasiparticle}}
\label{app:thermal_edge}

\begin{proof}[Proof of Lemma \ref{lem:thermal-tail-quasiparticle}]
Throughout this proof, we adopt the notation of Definition
\ref{def:quasi_intro} and write $i_k:=\varphi_0(k)$. By symmetry under reflection, 
it
suffices to consider
$i_k\in\llbracket N_1+\lfloor H_T/2\rfloor,-\lfloor H_T/2\rfloor\rrbracket$.
Let $\L^{\rm th}(s)$ be the coupled thermal Toda evolution introduced
immediately before Lemma \ref{lem:edge-thermal-coupling}, and denote
$\eig\L^{\rm th}(s)=(\eta_1,\eta_2,\ldots,\eta_N)$. For any $s\ge0$, let
$\varphi_s^{\rm th}:\llbracket1,N\rrbracket\to\llbracket N_1,N_2\rrbracket$
denote a $\zeta$-localization center bijection for $\L^{\rm th}(s)$, as in
\cite[Assumption 1.12]{Agg25}, and define
$Q_\ell^{\rm th}(s):=q_{\varphi_s^{\rm th}(\ell)}^{\rm th}(s)$ as in
\cite[Equation (1.15)]{Agg25}. We also set
$\alpha_{\rm th}:=-2\mathbb E\log a_0^{\rm th}$.

We work on the events from Lemmas \ref{l0eigenvalues}, with $A=\log N$,
\ref{eigenvalues0}, with $\delta=e^{-(\log N)^2}$, \ref{lem:Estr}, and
\ref{lem:edge-thermal-coupling}, and from Proposition
\ref{prop:outside-strip-confinement}. For the thermal evolution, we also
restrict to the corresponding overwhelmingly probable 
boundedness and eigenvalue-separation events (from \cite[Lemma 3.13]{Agg25})
\[
\bigcap_{s\ge0}\mathsf{BND}_{\L^{\rm th}(s)}(\log N),
\qquad
\mathsf{SEP}_{\L^{\rm th}}\bigl(e^{-(\log N)^2}\bigr),
\]
and on the particle-spacing events from \cite[Lemmas 3.11]{Agg25} on which
\begin{equation}\label{eqn:thmspace}
\begin{aligned}
\left|q_i^{\rm th}(0)-q_j^{\rm th}(0)-\alpha_{\rm th}(i-j)\right|
&\le \sqrt{|i-j|}(\log N)^2,
&&i,j\in\llbracket N_1,N_2\rrbracket,\\
\left|q_i^{\rm th}(s)-q_j^{\rm th}(s)-\alpha_{\rm th}(i-j)\right|
&\le \sqrt{|i-j|}(\log N)^2,
&&\begin{gathered}
0\le s\le T\log N,\\[-2mm]
i,j\in\llbracket N_1+T(\log N)^4,N_2-T(\log N)^4\rrbracket.
\end{gathered}
\end{aligned}
\end{equation}
We finally restrict to the event from
\cite[Lemma 3.22]{Agg25} on which 
any two $e^{-150 (\log N)^{3/2}}$-localization centers $\varphi, \varphi'$ for 
an eigenvalue $\lambda$ of $\L^{\rm th} (s)$ satisfy $|\varphi-\varphi'| \leq 2 (\log N)^3$.

\smallskip
\emph{Step 1: Confinement of localization centers.} Equation
\eqref{eqn:umphi}, with $s'=0$ and $m=\varphi_t(k)$, and
Definition \ref{def:loccent} imply that
$|\varphi_t(k)-i_k|\le T(\log N)^3$: Indeed, if not, then \eqref{eqn:umphi} implies
$|u_k(\varphi_t(k);t)|<\zeta$, contradicting that $\varphi_t(k)$ is a
$\zeta$-localization center. Since $T(\log N)^3=o(H_T)$, for all
sufficiently large $N$ we have
$\varphi_t(k)\in\llbracket N_1+\lfloor H_T/3\rfloor,
-\lfloor H_T/3\rfloor\rrbracket$.

\smallskip
\emph{Step 2: Match \(\lambda_k\) to one thermal eigenvalue.} Set $m_-:=-2T(\log T)^5$ and
$\mathcal D_-:=\llbracket m_-,N_2\rrbracket$. For $s\in\{0,t\}$, set
$p_s:=\varphi_s(k)$, so that $p_0=i_k$. The assumption on $i_k$ gives
$\dist(p_0,\mathcal D_-)\ge cH_T$, and Step 1 gives the same bound for
$p_t$.

Apply Lemma \ref{lem:full-to-truncation} with $\L=\L(s)$,
$\widetilde\L=\L^{\rm th}(s)$, and $\mathcal D=\mathcal D_-$. The second 
bound in \eqref{ll}, with $\delta = e^{-T(\log N)^2}$, is provided by
 Lemma \ref{lem:edge-thermal-coupling}; the first bound in \eqref{ll}
follow from the events
of Lemma \ref{l0eigenvalues} and
\cite[Lemma 3.13]{Agg25} which we have restricted to. Finally, if $|r-m_-|\le1$, then
$|r-p_s|\ge cH_T>T(\log N)^4$, so Proposition
\ref{prop:outside-strip-confinement} gives the localization estimate
required in Assumption \ref{lmatrixl}.

The lemma therefore gives, at each $s\in\{0,t\}$, an eigenvalue
$\eta=\eta_\ell\in\eig\L^{\rm th}(s)$ within $e^{-(\log N)^3}$ of
$\lambda_k$. It also shows that $p_s$ is an
$N^{-1}\zeta$-localization center for this eigenvalue. Since the spectrum of
$\L^{\rm th}(s)$ is independent of $s$ and since we restricted
to $\mathsf{SEP}_{\L^{\rm th}}\bigl(e^{-(\log N)^2}\bigr)$, the $\eta \in \eig\L^{\rm th}(s)$ such that $|\lambda_k-\eta|\le e^{-(\log N)^3}$ is
unique; at both $s=0$ and $s=t$, $p_s$
is an $N^{-1}\zeta$-localization center for $\eta$ (with respect to $\L^{\rm th}(s)$).

\smallskip
\emph{Step 3: Compare $p_s$ to $\varphi_s^{\rm th}(\ell)$.} Recall that $p_s$ is an
$N^{-1}\zeta$-localization center for $\eta$ at $s\in\{0,t\}$, while
$\varphi_s^{\rm th}(\ell)$ is a $\zeta$-localization center for the same
eigenvalue. For sufficiently large $N$,
$N^{-1}\zeta\ge e^{-150(\log N)^{3/2}}$, and Step 1 shows that $p_s$
is at
distance at least $\lfloor H_T/3\rfloor$, and hence at least $T(\log N)^5$, from the endpoints
$N_1, N_2$. Thus, \cite[Lemma 3.22]{Agg25}, applied at $s=0$ and
$s=t$, gives
\begin{equation}
\label{eq:thermal-center-index-match}
\left|
p_s-\varphi_s^{\rm th}(\ell)
\right|
\le
(\log N)^3,
\qquad s\in\{0,t\}.
\end{equation}

Applying the thermal particle-spacing event \eqref{eqn:thmspace}
 with $i=p_s$ and
$j=\varphi_s^{\rm th}(\ell)$, and then using
\eqref{eq:thermal-center-index-match}, gives
\begin{equation}
\label{eq:thermal-center-position-match}
\left|
q_{p_s}^{\rm th}(s)
-
Q_{\ell}^{\rm th}(s)
\right|
\le
(\log N)^5,
\qquad
s\in\{0,t\}.
\end{equation}

\smallskip
\emph{Step 4: Compare displacements of $Q_k$ and $Q_\ell^{\rm th}$.} By the definition of the coupling, $a_r(0)=a_r^{\rm th}(0)$ for every
$r\in\llbracket N_1,-\lfloor T(\log T)^5\rfloor-1\rrbracket$.
Consequently, $q_a(0)-q_a^{\rm th}(0)$ is independent of
$a\in\llbracket N_1,-\lfloor T(\log T)^5\rfloor\rrbracket$; denote this
constant by $c_-$. Both $p_t$ and $i_k$ belong to this interval for all
sufficiently large $N$. Since $\partial_s q_a(s)=b_a(s)$ and
$\partial_s q_a^{\rm th}(s)=b_a^{\rm th}(s)$, Lemma
\ref{lem:edge-thermal-coupling} gives, for the indices used below,
\begin{equation}
\label{eq:fixed-index-position-coupling}
\left|
q_a(t)-q_a^{\rm th}(t)-c_-
\right|
\le
t e^{-T(\log N)^2}.
\end{equation}
Apply this at the index \(a=p_t\), and use
\eqref{eq:thermal-center-position-match} at \(s=t\). At time zero,
use $p_0=i_k$, for which $q_{i_k}(0)-q_{i_k}^{\rm th}(0)=c_-$, and use
\eqref{eq:thermal-center-position-match} at $s=0$. Since
$Q_k(t)=q_{p_t}(t)$ and $Q_k(0)=q_{i_k}(0)$, we obtain
\begin{equation}
\label{eq:actual-thermal-displacement-match}
\left|
\bigl(Q_k(t)-Q_k(0)\bigr)
-
\bigl(
Q_{\ell}^{\rm th}(t)
-
Q_{\ell}^{\rm th}(0)
\bigr)
\right|
\le
(\log N)^C.
\end{equation}

\smallskip
\emph{Step 5: Apply the thermal quasiparticle trajectory theorem.} Let
$v_{\rm eff}^{\rm th}$ denote the thermal effective velocity from
\cite[Definition 1.8]{Agg25}. Apply \cite[Theorem 1.13]{Agg25} with $T$ there given by
 $T\log N$ here. Its scale
assumption $T\log N\le N(\log N)^{-6}$ holds for all sufficiently large
$T$ because $N=\lfloor2T^2\rfloor$. By
\eqref{eq:thermal-center-index-match}, the distance from
$\varphi_0^{\rm th}(\ell)$ to $\{N_1,N_2\}$ is at least $\lfloor H_T/3\rfloor$, which
is larger than the required buffer $T(\log N)^6$. Thus
\begin{equation}
\label{eq:thermal-linear-trajectory}
\left|
Q_{\ell}^{\rm th}(t)
-
Q_{\ell}^{\rm th}(0)
-
t v_{\rm eff}^{\rm th}(\eta)
\right|
\le
(T\log N)^{1/2}(\log N)^{35}
\le
 T^{1/2}(\log N)^C.
\end{equation}

By Step 2, $|\lambda_k-\eta|\le e^{-(\log N)^3}$. On the spectral
boundedness events we have restricted to,
$|\lambda_k|+|\eta|\le C\log N$, so the derivative estimate in
\cite[Corollary 3.8]{Agg25} gives
$t|v_{\rm eff}^{\rm th}(\lambda_k)-v_{\rm eff}^{\rm th}(\eta)|\le1$.
Combining this with
\eqref{eq:actual-thermal-displacement-match} and
\eqref{eq:thermal-linear-trajectory}, we conclude that
\begin{equation}
\label{eq:actual-thermal-linear-trajectory}
\left|
Q_k(t)-Q_k(0)
-
t v_{\rm eff}^{\rm th}(\lambda_k)
\right|
\le
C T^{1/2}(\log N)^C.
\end{equation}

\smallskip
\emph{Step 6: Identify the finite-\(T\) characteristic in the thermal
region.} Let $\widetilde\varrho_{\rm th}$ be the thermal density of states
$\varrho$ from \cite[Definition 1.2]{Agg25}, with parameters $\beta=1$
and $\theta=\mathfrak A^{-1}$; equivalently, it is the limiting
empirical eigenvalue density of the thermal Lax matrix $\L^{\rm th}(0)$ 
used above. Set
$m_{\rm th}(\lambda,x):=\alpha_{\rm th}^{-1}
\widetilde\varrho_{\rm th}(\lambda)$.
This profile is independent of $x$. Formula \eqref{eqn:vedf} for its
effective velocity agrees with \cite[Definition 1.8]{Agg25}, and hence its
effective velocity is $v_{\rm eff}^{\rm th}$. Its characteristic flow is therefore
\begin{equation}
\label{eq:homogeneous-thermal-characteristic}
x_\tau^{\rm th}(\lambda;q)
=
q+\tau v_{\rm eff}^{\rm th}(\lambda).
\end{equation}

The bound $\widetilde\alpha^{(T)}(y)\ge\mathfrak K^{-1}\mathfrak A$ from
\eqref{eqn:alph_bd} and the inequality
$i_k/T\le-\frac12(\log N)^{10}$ show that
$F_{\rm str}^{(T)}(i_k/T)$ is at distance at least
$c\mathfrak A(\log N)^{10}$ from
$F_{\rm str}^{(T)}\bigl(-(\log T)^5\bigr)$. Thus, by Assumption
\ref{ass:model} and the event of Lemma
\ref{lem:Estr}, $\bar q_{i_k}$ is contained in
a neighborhood of radius $c\mathfrak A(\log N)^{10}$ throughout which
$m_0^{(T)}$ and $m_{\rm th}$ agree.

We apply Lemma \ref{lem:thermal-completion} to these two profiles with
$x_0=q=\bar q_{i_k}$, $L=c_0\mathfrak A(\log N)^{10}$ for a sufficiently
small constant $c_0>0$, $s=\tau$, and $\lambda=\lambda_k$. Indeed,
$\tau\le\log N<(\log T)^2$, while $|\lambda_k|\le\log N$ and hence
$\tau(1+|\lambda_k|)\le C(\log N)^2\ll L$. We obtain
\begin{equation}
\label{eq:finiteT-thermal-characteristic}
\left|
x_\tau^{(T)}
\left(\lambda_k;\bar q_{i_k}\right)
-
x_\tau^{\rm th}
\left(\lambda_k;\bar q_{i_k}\right)
\right|
\le
Ce^{-c(\log N)^2}.
\end{equation}

Finally, $Q_k(0)=q_{i_k}(0)=T\bar q_{i_k}$, while
\eqref{eq:homogeneous-thermal-characteristic} gives
$T x_\tau^{\rm th}(\lambda_k;\bar q_{i_k})
=Q_k(0)+t v_{\rm eff}^{\rm th}(\lambda_k)$. Combining this identity with
\eqref{eq:actual-thermal-linear-trajectory} and
\eqref{eq:finiteT-thermal-characteristic} proves
\eqref{eq:thermal-tail-quasiparticle}.
\end{proof}

	\section{Time zero density of states and linear statistics}
\label{app:time-zero-dos}

The purpose of this appendix is to prove certain facts about the
time zero Lax matrix $\L(0)$ and its density of states, 
and give quantitative error bounds for the 
convergence of linear statistics. In
particular, we provide a proof of Lemma \ref{lem:Estr} in 
Section \ref{app:Estr};
Proposition \ref{prop:outside-strip-confinement} in Section \ref{app:lemLocR};  
Lemmas \ref{lem:frozen-dos-intro}, \ref{rho0estimate}, and
\ref{lem:block_wegner_input} in Section \ref{app:frozen-dos};
Lemma \ref{lem:uniform-non-vacuum} in Section \ref{app:vac};
Proposition \ref{prop:app-time-zero-dos} in Section \ref{app:time-zero-dos-proof};
and Proposition \ref{prop:scat_conc} in Section \ref{app:scat_conc}.

Throughout the appendix,
unless otherwise specified, we adopt
Assumption \ref{ass:model}, and we adopt the notation 
of Definition \ref{def:quasi_intro}.
However, until Section \ref{app:scat_conc}, we suppress the superscript
$(T)$ on objects constructed from the $T$-dependent entry laws; for example,
we will simply write $m_0$ rather than $m_0^{(T)}$, 
and write $p_{\dg}(b;y)$ and $ \mu_{\sd}(a,y)$ rather
 than $p_{\dg}^{(T)}(b;y)$ and $\mu_{\sd}^{(T)}(a,y)$.
In fact, we may drop the assumption that 
$p_{\dg}(b;y)$ and $\mu_{\sd}(a;y)$ are given by 
the thermal equilibrium densities for $|y| > (\log T)^{5}$
(as specified in Assumption \ref{ass:model}), 
and all of the results of this appendix remain valid.

\subsection{Lax density of states with homogeneous entries}
\label{app:frozen-dos}

For \(y\in\mathbb R\) and \(K\ge1\), let \(\mathbf{S}_y^{(K)}\) denote the
\(K\times K\) symmetric tridiagonal matrix whose diagonal entries 
are sampled from the density
\(\mu_{\dg}(\cdot,y)\), and whose off diagonal are sampled from
\(\mu_{\sd}(\cdot,y)\), all independently.  Set
\[
\nu_y^{(K)}:=\frac1K\sum_{\lambda\in\eig \mathbf{S}_y^{(K)}}\delta_\lambda .
\]
First we record a useful continuity result for linear statistics as a function 
of the parameter $y$.

\begin{lem}
\label{lem:app-hilbert-schmidt-spectral-stability}
Let \(\mathbf{A},\mathbf{B}\) be \(K\times K\) real symmetric matrices. For each
 \(\mathbf{M}\in\{\mathbf{A},\mathbf{B}\}\), let
\(\eig \mathbf{M}=(\lambda_1(\mathbf{M}),\lambda_2(\mathbf{M}),\ldots,
\lambda_K(\mathbf{M}))\). Then
\[
\frac1K\sum_{\ell=1}^K
|\lambda_\ell(\mathbf{A})-\lambda_\ell(\mathbf{B})|^2
\le
\frac1K\Tr(\mathbf{A}-\mathbf{B})^2 .
\]
Consequently, for every Lipschitz function \(h:\mathbb R\to\mathbb R\),
\[
\left|
\frac1K\Tr h(\mathbf{A})-\frac1K\Tr h(\mathbf{B})
\right|
\le
\|h\|_{\mathrm{Lip}}
\left(\frac1K\Tr(\mathbf{A}-\mathbf{B})^2\right)^{1/2}.
\]
\end{lem}

\begin{proof}
The first bound is the Hoffman-Wielandt inequality; see, for example,
\cite[Lemma 2.1.19]{AGZ10}.  The second bound
follows from Cauchy Schwarz:
\[
\frac1K\sum_{\ell=1}^K |h(\lambda_\ell(\mathbf{A}))-h(\lambda_\ell(\mathbf{B}))|
\le \|h\|_{\Lip}
\left(\frac1K\sum_{\ell=1}^K|\lambda_\ell(\mathbf{A})-\lambda_\ell(\mathbf{B})|^2\right)^{1/2}.
\]\qedhere
\end{proof}

The next lemma shows that linear statistics of $\mathbf{S}_y^{(K)}$
converge in probability, and in $L^1$, to a deterministic limit.
\begin{lem}
\label{lem:app-frozen-dos}
For each \(y\in\mathbb R\), there is a probability measure
\(\widetilde\nu_y\) such that, for every bounded continuous \(h\),
\[
\frac1K\Tr h(\mathbf{S}_y^{(K)})
\longrightarrow
\int_{\mathbb R} h(\lambda)\,\widetilde\nu_y(d\lambda)
\]
in probability and in \(L^1\) as \(K\to\infty\).  
\end{lem}

\begin{proof}
Truncate the diagonal entries to \([-M,M]\) and the off diagonal entries to
\([M^{-1},M]\), obtaining \(\mathbf{S}_y^{(K,M)}\).  For fixed \(M\), this is an
ergodic Jacobi matrix (that is, satisfies
\cite[Hypotheses H.5.1 and H.5.6]{OINL}) with bounded coefficients bounded away from zero, so
\cite[Lemma 5.12]{OINL} provides a probability measure \(\widetilde\nu_y^{(M)}\)
satisfying
\[
\frac1K\Tr h(\mathbf{S}_y^{(K,M)})
\to \int_{\mathbb{R}} h\,d\widetilde\nu_y^{(M)}
\]
in probability and \(L^1\), for every bounded continuous \(h\).  The
bounds \eqref{munu} give
\begin{equation}\label{eqn:Mtoinf}
\lim_{M\to\infty}\sup_{K,y} 
\frac1K\E\left[ \Tr(\mathbf{S}_y^{(K)}-\mathbf{S}_y^{(K,M)})^2 \right] =0.
\end{equation}
Taking the limit as $K \rightarrow \infty$ the second bound in 
 Lemma \ref{lem:app-hilbert-schmidt-spectral-stability} and using 
 the convergence above, 
we see that if \(h\) is bounded
and Lipschitz, then for each $y \in \mathbb{R}$
\[
\left|\int_{\mathbb R} h\,d\widetilde\nu_y^{(M)}
-\int_{\mathbb R} h\,d\widetilde\nu_y^{(M')}\right|
\le \|h\|_{\Lip}\,\varepsilon_{M,M'},
\qquad \varepsilon_{M,M'}\to0
\]
as \(M,M'\to\infty\). Thus, there exists a unique weak limit of the measures
$\widetilde\nu_y^{(M)}$ as $M$ tends to $\infty$; we denote this limit by
\(\widetilde\nu_y\).

We now show convergence of linear statistics of 
$\mathbf{S}_y^{(K)}$ (the untruncated random matrices) to
integrals against this limiting measure.  For a
bounded Lipschitz \(h\), write
\[
\left|\frac1K\Tr h(\mathbf{S}_y^{(K)})
-\int_{\mathbb R} h\,d\widetilde\nu_y\right|
\le A_{K,M}+B_{K,M}+C_M,
\]
where
\[
A_{K,M}:=\left|\frac1K\Tr h(\mathbf{S}_y^{(K)})
-\frac1K\Tr h(\mathbf{S}_y^{(K,M)})\right|,
\]
\[
B_{K,M}:=\left|\frac1K\Tr h(\mathbf{S}_y^{(K,M)})
-\int_{\mathbb R} h\,d\widetilde\nu_y^{(M)}\right|,
\]
and
\[
C_M:=\left|\int_{\mathbb R} h\,d\widetilde\nu_y^{(M)}
-\int_{\mathbb R} h\,d\widetilde\nu_y\right|.
\]
By \eqref{eqn:Mtoinf}, the first term is small in \(L^1\), uniformly in \(K\), once \(M\) is large.  For fixed \(M\), the second term tends
to zero in probability and \(L^1\) as \(K\to\infty\).  The third term tends to
zero by the definition of \(\widetilde\nu_y\).  Letting first \(K\to\infty\)
and then \(M\to\infty\), we obtain the claimed convergence for bounded
Lipschitz \(h\).

It remains only to pass from bounded Lipschitz test functions to bounded continuous
test functions. Using \eqref{munu}, this follows from a routine approximation argument.
\end{proof}

\begin{lem}
\label{lem:app-frozen-couple}
There is a constant $C> 0$ such that 
the following holds. Let $A > 1$.
For every integer $K\geq1$ and each $y,y'\in\mathbb R$ with $y\le y'$ satisfying (with $\mathcal{P}$
as in Assumption \ref{ass:entry-law-profiles})
$[y,y']
\cap
\bigl( \mathcal P \cup \{ - (\log T)^5, (\log T)^5 \} \bigr) =\varnothing$,
 there is a coupling between
$\mathbf{S}_y^{(K)}$ and $\mathbf{S}_{y'}^{(K)}$ such that, for every
bounded Lipschitz function \(h\) with
$\max(\|h\|_{\infty}, \| h \|_{\Lip}) \leq A$,
\begin{equation}
\label{eq:app-frozen-profile-coupling}
\frac1K \mathbb{E} \left|\Tr h(\mathbf{S}_y^{(K)})
- \Tr h(\mathbf{S}_{y'}^{(K)}) \right|
\leq C A |y-y'|.
\end{equation}
Moreover, let $B=\llbracket b,b+K-1\rrbracket\subseteq\llbracket N_1,N_2\rrbracket$
be an index interval satisfying
$[b/T,(b+K-1)/T]\cap\bigl(\mathcal P\cup\{-(\log T)^5,(\log T)^5\}\bigr)=\varnothing$.
There is a coupling between $\mathbf{S}_{b/T}^{(K)}$ and the corresponding $K\times K$ block
$\L^B$ of $\L(0)$ such that
\begin{equation}
\label{eq:app-frozen-block-coupling}
\frac1K\mathbb{E} \left|\Tr h(\mathbf{S}_{b/T}^{(K)}) - \Tr h(\L^B) \right|
\leq C A K/T .
\end{equation}
Finally, denoting $\widetilde \nu_y$ the family of measures from Lemma \ref{lem:app-frozen-dos}, we have 
\begin{equation}\label{eqn:nuydif}
\left| \int_{-\infty}^{\infty} h(\lambda)\,\widetilde\nu_y(d\lambda)-  \int_{-\infty}^{\infty} h(\lambda)\,\widetilde\nu_{y'}(d\lambda) \right| 
\leq C  A |y-y'|   .
\end{equation}
\end{lem}

\begin{proof}
Use the coupling afforded by the bounds
\eqref{eqn:W2bd} independently for every diagonal
and off diagonal entry.  Then
\[
\frac1K\E\Tr(\mathbf{S}_y^{(K)}-\mathbf{S}_{y'}^{(K)})^2\le C|y-y'|^2.
\]
Lemma \ref{lem:app-hilbert-schmidt-spectral-stability}
gives \eqref{eq:app-frozen-profile-coupling}. The argument for \eqref{eq:app-frozen-block-coupling}
is similar:
if \(B=\llbracket b,b+K-1\rrbracket\) and \(y=b/T\), then
\(|i/T-y|\le K/T\) for \(i\in B\), so the same coupling gives
\[
\frac1K\E\Tr(\L^B-\mathbf{S}_y^{(K)})^2\le C(K/T)^2.
\]
Lemma \ref{lem:app-hilbert-schmidt-spectral-stability} then gives
\eqref{eq:app-frozen-block-coupling}. Finally, let
\(K\to\infty\) in \eqref{eq:app-frozen-profile-coupling} and use Lemma
\ref{lem:app-frozen-dos} to obtain \eqref{eqn:nuydif}.
\end{proof}

\begin{prop}
\label{prop:app-frozen-wegner}
For each $y\in\mathbb R$ and $K\ge1$, the measure
$\widetilde\nu_{y,K}:=\E\,\frac1K\sum_{\lambda \in \eig \mathbf{S}_y^{(K)}}
\delta_{\lambda}$ is absolutely continuous with respect to Lebesgue measure.
Denoting its density by $r_K(\lambda;y)$, there exists a constant $c>0$,
independent of $K$, such that
\begin{equation}
\label{eq:app-frozen-wegner-density}
0\le r_K(\lambda;y)\le c^{-1} e^{-c\lambda^2},
\qquad (\lambda,y)\in\mathbb R^2.
\end{equation}
Similarly,
if \(B\subset\llbracket N_1,N_2\rrbracket\) is an index interval with
\(|B|=K\), then
$\nu_{B,K}:=\E\,\frac1K\sum_{\lambda \in \eig \L^B(0)}
\delta_{\lambda}$
is absolutely continuous with respect to Lebesgue measure. Denoting its density
by \(r_B(\lambda)\), we have
\[
0\le r_B(\lambda)\le c^{-1} e^{-c\lambda^2},
\qquad \lambda\in\mathbb R.
\]
\end{prop}

\begin{proof}
We only show the first statement of the proposition, as the proof of the latter is entirely analogous.
It is enough to prove that, for every interval \(J\subset\R\),
\begin{equation}\label{eqn:desired}
\E\,\frac1K\Tr\1_J(\mathbf{S}_y^{(K)})
\le C\int_J e^{-c\lambda^2}\,d\lambda .
\end{equation}
Indeed, this estimate implies absolute continuity of the averaged spectral
measure and gives \eqref{eq:app-frozen-wegner-density}.

Let \(\mathbf A\) be a real symmetric matrix, let \(\mathbf u\) be a unit vector,
let \(\H_s=\mathbf A+s|\mathbf u\rangle\langle\mathbf u|\), and let
\(w : \mathbb{R} \rightarrow \mathbb{R}_{\geq 0}\) be bounded. We claim that,
for every interval \(J\),
\begin{equation}
\label{eqn:min_type}
\int_\R w(s)\,\langle \mathbf u,\1_J(\H_s)\mathbf u\rangle\,ds
\le \|w\|_\infty |J|.
\end{equation}
We show this for \(w\equiv1\), since 
the general case is similar. Let $z=E+\mathrm{i}\eta$, where
$E,\eta\in\mathbb R$ and $\eta > 0$. First, we record the resolvent identity
\[
\langle \mathbf u,(\H_s-z)^{-1}\mathbf u\rangle
=\frac{m(z)}{1+s\,m(z)},\qquad
m(z):=\langle \mathbf u,(\mathbf A-z)^{-1}\mathbf u\rangle.
\]
Write \(m(E+\mathrm{i}\eta)=a+\mathrm{i}b\), where \(b>0\).  Then
\[
\Im\frac{m(E+\mathrm{i}\eta)}{1+s\,m(E+\mathrm{i}\eta)}
=\frac{b}{(1+sa)^2+(sb)^2},
\]
and hence
\[
\int_\R \Im\frac{m(E+\mathrm{i}\eta)}{1+s\,m(E+\mathrm{i}\eta)}\,ds=\pi .
\]
Therefore,
\[
\int_\R\int_J
\Im\langle \mathbf u,(\H_s-E-\mathrm{i}\eta)^{-1}\mathbf u\rangle\,dE\,ds
=\pi |J|.
\]
Moreover,
\[
\langle \mathbf u,\1_J(\H_s)\mathbf u\rangle
\le \liminf_{\eta\downarrow0}\frac1\pi
\int_{J_\varepsilon}
\Im\langle \mathbf u,(\H_s-E-\mathrm{i}\eta)^{-1}\mathbf u\rangle\,dE
\]
for any \(\varepsilon>0\), where \(J_\varepsilon\) is the \(\varepsilon\)-neighborhood
of \(J\).  Integrating in \(s\) and then letting
\(\varepsilon\downarrow0\), gives the estimate for \(w\equiv1\).  The general
bounded \(w\) case is similar.

We now prove \eqref{eqn:desired} for intervals
\(J=[E,E+\delta]\), where \(0<\delta\le1\). We assume below that $|E| \ge 10$, as the proof when $|E| < 10$ 
is entirely analogous. 
We further
assume that $E\ge10$, as the case when $E\le-10$ is very similar. Put
\[
\mathbf P:=\1_J(\mathbf{S}_y^{(K)}),\qquad X:=\Tr \mathbf P .
\]
Write the diagonal entries of \(\mathbf{S}_y^{(K)}\) as \((\mathbf{S}_y^{(K)})_{i i} =\xi_i\), $1 \leq i \leq K$ and the
off diagonal entries as \((\mathbf{S}_y^{(K)})_{i, i+1} =\eta_i>0\), \(1\le i\le K-1\).  Since every
eigenvalue counted by \(\mathbf P\) is at least \(E\),
\begin{equation}\label{eqn:EXTrsyP}
E  X\le \Tr(\mathbf{S}_y^{(K)}\mathbf P).
\end{equation}
Set \(p_i:=P_{ii}\), \(\eta_0=\eta_K=0\).  As \(\mathbf P\) is a projection,
\[
|P_{i,i+1}|\le \sqrt{p_i p_{i+1}}\le \frac{p_i+p_{i+1}}2,
\]
and therefore, with $\xi_i^+ \coloneqq \max(\xi_i, 0)$
\begin{equation}\label{eqn:TrsyP}
\Tr(\mathbf{S}_y^{(K)}\mathbf P)
\le
\sum_{i=1}^K \bigl(\xi_i^+ +\eta_{i-1}+\eta_i\bigr)p_i .
\end{equation}
Let \(t=E/6\).  For all nonnegative \(x_1,x_2,x_3\),
\[
x_1+x_2+x_3
\le
\frac E2+\sum_{\ell=1}^3 x_\ell \1_{\{x_\ell\ge t\}}.
\]
Applying this with \(x_1=\xi_i^+\), \(x_2=\eta_{i-1}\), \(x_3=\eta_i\) in \eqref{eqn:TrsyP}, and 
substituting the resulting bound back into \eqref{eqn:EXTrsyP}, 
we get
\[
X\le \frac2E\left[
\sum_{i=1}^K \xi_i^+\1_{\{\xi_i\ge t\}}p_i
+\sum_{i=1}^{K-1}\eta_i\1_{\{\eta_i\ge t\}}(p_i+p_{i+1})
\right].
\]

We estimate the expectation of the right side using \eqref{eqn:min_type} in the
diagonal variables.  For the diagonal term, condition on all entries except
\(\xi_i\).  The weight
\[
w_i(s):=s\1_{\{s\ge t\}}p_{\dg}(s;y)
\]
satisfies
\begin{equation}\label{eqn:winfbd}
\|w_i\|_\infty\le C t e^{-c t^2}
\end{equation}
by the Gaussian bound on $p_{\dg}(s;y)$ in \eqref{munu}.  Hence by \eqref{eqn:min_type} and \eqref{eqn:winfbd},
\[
\E\bigl[\xi_i^+\1_{\{\xi_i\ge t\}}p_i\bigr]
\le C\delta\,t e^{-c t^2}.
\]
For the off diagonal term, let $\mathbb E_i$ denote the expectation with respect
to $\xi_i$, conditioning on all other entries of $\mathbf S_y^{(K)}$.
The bound \eqref{eqn:min_type} with weight
\(w = p_{\dg}(\cdot;y)\) gives $\mathbb E_i[p_i]\le C\delta$.
Thus, using the second bound for
\(\mu_{\sd}\) in \eqref{munu},
\[
\E\bigl[\eta_i\1_{\{\eta_i\ge t\}}p_i\bigr]
\le C\delta\,\E[\eta_i\1_{\{\eta_i\ge t\}}]
\le C\delta e^{-c t^2}.
\]
The same estimate holds with \(p_{i+1}\), using $\mathbb E_{i+1}$ instead.
Summing over \(i\) yields
\[
\E X\le C K\delta e^{-cE^2}.
\]

Hence
for every interval \(J\) of length at most one,
\[
\E\,\frac1K\Tr\1_J(\mathbf{S}_y^{(K)})
\le C |J| \exp\{-c\,\dist(J,0)^2\}.
\]
Partitioning an arbitrary interval into subintervals of length at most one
and adjusting \(c,C\), this becomes
\[
\E\,\frac1K\Tr\1_J(\mathbf{S}_y^{(K)})
\le C\int_J e^{-c\lambda^2}\,d\lambda ,
\]
as claimed.
\end{proof}

\begin{proof}[Proof of Lemma \ref{lem:block_wegner_input}]
By Proposition \ref{prop:app-frozen-wegner}, 
\[
\mathbb E X_B(J)
\le C|B|\int_J e^{-c\lambda^2}\,d\lambda
\le C|B||J|e^{-c\lambda_J^2},
\]
as desired.
\end{proof}

Next, we obtain a similar bound for the density of $\widetilde\nu_y$, which 
will allow us to prove Lemma \ref{lem:frozen-dos-intro}.
\begin{lem}
\label{lem:app-frozen-dos-regularity}
For each $y\in\mathbb R$, the measure $\widetilde\nu_y$ from Lemma
\ref{lem:app-frozen-dos} is absolutely continuous with respect to Lebesgue
measure. Denoting its density by $\widetilde\varrho_0(\lambda;y)$, there
exists a constant $c>0$ such that
\[
0\le \widetilde\varrho_0(\lambda;y)\le c^{-1}e^{-c\lambda^2}
\]
for all $(\lambda,y)\in\mathbb R^2$.
\end{lem}

\begin{proof}
Fix \(y\).  
Taking $K \rightarrow \infty$ in 
Proposition \ref{prop:app-frozen-wegner} 
implies that \(\widetilde\nu_y\) is absolutely continuous with respect to Lebesgue
measure, and its Radon Nikodym derivative
\(\widetilde\varrho_0(\lambda;y)\) satisfies
\[
0\le \widetilde\varrho_0(\lambda;y)\le C e^{-c\lambda^2}
\]
for almost every \(\lambda\) with respect to Lebesgue measure. Modifying the density on a null set for each
\(y\), we may take the displayed bound to hold for all \(\lambda\).
\end{proof}

\begin{proof}[Proof of Lemma \ref{lem:frozen-dos-intro}]
This follows from Lemmas \ref{lem:app-frozen-dos} and \ref{lem:app-frozen-dos-regularity}. 
\end{proof}

\begin{proof}[Proof of Lemma \ref{rho0estimate}]
By Definition \ref{rho0} and Lemma \ref{lem:app-frozen-dos-regularity}, there 
are constants \(c,C>0\)
such that $0\le \varrho_0(\lambda;q)\le C e^{-c\lambda^2}$ uniformly in
\((\lambda,q)\in\mathbb R^2\). Moreover, Definition \ref{rho0} and
\eqref{eqn:alph_bd} give $0\le m_0(\lambda,q)
=\alpha(q)^{-1}\varrho_0(\lambda,q)
\le C\mathfrak A^{-1}e^{-c\lambda^2}$.
This completes the proof.
\end{proof}

The proof of the following lemma is entirely analogous to that of
\cite[Lemma 5.7]{Agg25} (while the latter is only stated under thermal
equilibrium, it is quickly verified using Lemma \ref{gijexponential} that its
proof also holds under Assumption \ref{ass:model}); therefore, it is omitted.
\begin{lem}
\label{lem:app-frozen-moving-tests}
Fix \(A,B > 1\). Let \(h : \mathbb{R} \rightarrow \mathbb{R}\) be a bounded Lipschitz
function satisfying $\|h\|_\infty\le A$ and $\|h\|_{\Lip}\le B$.
Then, for some \(C>0\), uniformly for \(y\in\mathbb R\),
\[
\left| \mathbb{E}\left[ \Tr h(\mathbf{S}_y^{(K)}) \right]
-
 K \int_{-\infty}^{\infty} h(\lambda)\widetilde\varrho_0(\lambda,y)\,d\lambda \right|
 \leq  A(\log N)^C+C B K^{-500}
\]
\end{lem}

\subsection{Proof of Lemma \ref{lem:uniform-non-vacuum}}
\label{app:vac}

\begin{proof}
Throughout this section, we increase
\(\mathfrak A_0=\mathfrak A_0(\mathfrak K,\epsilon)\), if necessary, and
assume that \(\mathfrak A\geq \mathfrak A_0\). All constants below are
uniform in \(T\), \(y\), \(x\), and \(\tau\). Here, as at the beginning of
this appendix, \(T\) is the parameter indexing the entry laws, whose
superscript \((T)\) is suppressed.

We first prove a uniform lower bound for the density of states $\widetilde\varrho_0(\lambda;y)$.
By \eqref{eqn:mudgc}, there exists \(b_\star\in[-\mathfrak K,\mathfrak K]\)
such that
$p_{\dg}(b,y)\geq \mathfrak K^{-1}$,
$b\in I_\star:=[b_\star-\mathfrak K^{-1},b_\star+\mathfrak K^{-1}]$, $y \in \mathbb{R}$.
There is consequently \(c>0\) such that
\begin{equation}
\label{eq:finiteT-diagonal-lower}
p_{\dg}(b,y)\geq c,
\qquad
b\in I_\star,\quad y\in\mathbb R,\quad T\geq2.
\end{equation}

We next record that the off diagonal laws $\mu_{\sd}(\cdot, y)$ place a uniformly positive
amount of mass on small values. Namely, there are
\(c_0,p_0>0\) such that
\begin{equation}
\label{eq:small-hopping-probability}
\inf_{T\geq2}\inf_{y\in\mathbb R}
\mu_{\sd}
\bigl((0,e^{-c_0\mathfrak A}],y\bigr)
\geq p_0.
\end{equation}
To see this, let \(a\) be distributed according to
\(\mu_{\sd}(\cdot,y)\), and set
\[
\widetilde\alpha
=
-2\int_0^\infty\log a'\,\mu_{\sd}(da',y),
\qquad
Z:=2(-\log a)_+,
\qquad
Y:=\frac{Z}{\widetilde\alpha}.
\]
Since
$\widetilde\alpha
=
\mathbb E Z-\mathbb E\bigl[2(\log a)_+\bigr]$,
we have \(\mathbb E Y\geq1\). Moreover, by \eqref{munu},
\[
\mathbb E e^{\epsilon Y/2}
=
\mathbb P(a>1)
+
\int_0^1 a^{-\epsilon/\widetilde\alpha}\,
\mu_{\sd}(da,y)
\leq 1+\mathfrak K.
\]
Thus, we may choose \(M=M(\mathfrak K,\epsilon)\) sufficiently large that
$\mathbb E\left[Y\mathbf 1_{\{Y>M\}}\right]\leq\frac14$. 
It follows that
\[
\frac34
\leq
\mathbb E\left[Y\mathbf 1_{\{Y\leq M\}}\right]
\leq
\frac12+M\mathbb P(Y\geq\tfrac12),
\]
and therefore
$\mathbb P(Y\geq\tfrac12)\geq\frac1{4M}$.
On the event \(Y\geq\frac12\),
\[
a\leq e^{-\widetilde\alpha/4}
\leq e^{-\mathfrak A/(4\mathfrak K)}
\]
by \eqref{eqn:alph_bd}. 

Fix \(y \in \mathbb{R}\), and let \(\L^{(y)} =[L_{i, j}^{(y)}]_{i, j \in \mathbb{Z}}\) be the
doubly infinite tridiagonal symmetric matrix with diagonal and
off diagonal entries $L_{i i}^{(y)}  = b_i, L_{i,i+1}^{(y)} = a_i$ sampled from $\mu_{\dg}(\cdot, y)$
and $\mu_{\sd}(\cdot, y)$, respectively, with all entries
on and above the diagonal
mutually independent. Let
\[
G_{00}(z)
:=
\left\langle\delta_0,
(\L^{(y)}-z)^{-1}\delta_0
\right\rangle,
\qquad z\in\mathbb C_+.
\]
Removing the \(0\)th row and column separates the matrix into left and right
principal restrictions: \(\L_-^{(y)} = [L_{i, j}^{(y)}]_{i, j \in \mathbb{Z}_{<0}}\) is the restriction of \(\L^{(y)}\)
to \(\mathbb Z_{<0}\), and \(\L_+^{(y)}= [L_{i, j}^{(y)}]_{i, j \in \mathbb{Z}_{>0}}\) is its restriction to
\(\mathbb Z_{>0}\). Write
\[
g_-(z)
=
\left\langle\delta_{-1},
(\L_-^{(y)}-z)^{-1}\delta_{-1}
\right\rangle,
\qquad
g_+(z)
=
\left\langle\delta_{1},
(\L_+^{(y)}-z)^{-1}\delta_{1}
\right\rangle.
\]
The Schur complement formula gives
\begin{equation}
\label{eq:frozen-schur}
G_{00}(z)
=
\frac1{
b_0-z-a_{-1}^2g_-(z)-a_0^2g_+(z)
}.
\end{equation}

We will use the elementary uniform fractional moment estimate
\begin{equation}
\label{eq:cavity-fractional-moment}
\sup_{y,z\in\mathbb C_+}
\mathbb E|g_\pm(z)|^{1/2}
\leq C.
\end{equation}
Indeed, condition on every entry of $\L_{\pm}^{(y)}$
except the diagonal entry \(b_{\pm}\) at site $\pm1$. The same Schur complement
identity writes
\[
g_\pm(z)=\frac1{b-\zeta}
\]
for some \(\zeta\in\mathbb C_+\) independent of \(b\). Since the
diagonal densities \(p_{\dg}\) are uniformly bounded by \eqref{munu},
\[
\int_{\mathbb{R}} |g_\pm(z)|^{1/2}p_{\dg}(b,y)\,db
\leq
C\int_{|b-\Re\zeta|\leq1}
|b-\Re\zeta|^{-1/2}\,db
+
\int_{|b-\Re\zeta|>1}p_{\dg}(b,y)\,db
\leq C,
\]
which proves \eqref{eq:cavity-fractional-moment}.

Set
\[
\delta_{\mathfrak A}:=e^{-c_0\mathfrak A},
\qquad
r:=\frac{\mathfrak K^{-1}}8,
\qquad
\Sigma(z):=a_{-1}^2g_-(z)+a_0^2g_+(z).
\]
The variables \(a_{-1}\) and \(a_0\) are independent of the two
Green functions $g_{\pm}(z)$. Therefore, by
\eqref{eq:small-hopping-probability},
\eqref{eq:cavity-fractional-moment}, and Markov's inequality,
\[
\mathbb P\bigl(|\Sigma(z)|\leq r\bigr)\geq
\mathbb P(a_{-1},a_0\leq\delta_{\mathfrak A})-
\sum_{\sigma\in\{-,+\}}\mathbb P\left(a_\sigma\leq\delta_{\mathfrak A},\ a_\sigma^2|g_\sigma(z)|>r/2\right)
\geq p_0^2-C\delta_{\mathfrak A}.
\]
After increasing \(\mathfrak A\), we consequently have
\begin{equation}
\label{eq:self-energy-small}
\inf_{y,z\in\mathbb C_+}
\mathbb P\bigl(|\Sigma(z)|\leq r\bigr)
\geq \frac{p_0^2}{2}
=:p_\Sigma>0.
\end{equation}

Define the smaller interval
$J_\star
:=
\left[b_\star-\frac{\mathfrak K^{-1}}4,
b_\star+\frac{\mathfrak K^{-1}}4\right]$.
Recall, from the beginning of the proof, that $b_\star\in[-\mathfrak K,\mathfrak K]$.
Fix \(\lambda\in J_\star\) and \(0<\eta\leq r\), and put
\(z=\lambda+\mathrm{i}\eta\). Write
\[
z+\Sigma(z)=u+\mathrm{i}v.
\]
Since $g_{\pm}(z)$ have nonnegative imaginary part, \(v>0\).
On the event \(\{|\Sigma(z)|\leq r\}\), we have
\begin{equation}
\label{eq:self-energy-real-imaginary-bounds}
|u-\lambda|\leq r,
\qquad
0<v\leq\eta+|\Sigma(z)|\leq2r=\frac{\mathfrak K^{-1}}4.
\end{equation}
In particular, the distance from \(u\) to either endpoint of \(I_\star\)
is at least \(5\mathfrak K^{-1}/8\), so
\(|b-u|/v\geq5/2\) for either endpoint \(b\) of \(I_\star\), by
\eqref{eq:self-energy-real-imaginary-bounds}. Computing the conditional
expectation of
\eqref{eq:frozen-schur}, conditional on all variables except \(b_0\), and using
\eqref{eq:finiteT-diagonal-lower}, we obtain on this event
\[
\begin{aligned}
\mathbb E [\Im G_{00}(\lambda+\mathrm{i}\eta) | \{b_i\}_{i \neq 0} , \{a_i\}]
&=
\int_{\mathbb R}
p_{\dg}(b,y)
\frac{v}{(b-u)^2+v^2}\,db \\
&\geq
c\int_{I_\star}
\frac{v}{(b-u)^2+v^2}\,db
\geq
2c\arctan\left(\frac52\right).
\end{aligned}
\]
Averaging over the remaining variables and using
\eqref{eq:self-energy-small} gives
\begin{equation}
\label{eq:poisson-dos-lower}
\frac1\pi\mathbb E\Im G_{00}(\lambda+\mathrm{i}\eta)
\geq c>0,
\qquad
\lambda\in J_\star,\quad 0<\eta\leq r.
\end{equation}

Let \(\widetilde\nu_y\) denote the probability measure
from Lemma \ref{lem:app-frozen-dos}.
By \cite[Lemma 5.12]{OINL}, the left hand side of
\eqref{eq:poisson-dos-lower} is its Poisson
transform:
\[
\frac1\pi\mathbb E\Im G_{00}(\lambda+\mathrm{i}\eta)
=
\int_{\mathbb R}
\frac1\pi\frac{\eta}{(\lambda-E)^2+\eta^2}\,
\widetilde\nu_y(dE).
\]
By Lemma \ref{lem:app-frozen-dos-regularity},
\(\widetilde\nu_y\) has a density
\(\widetilde\varrho_0(\lambda,y)\). Thus, letting \(\eta\downarrow0\) in
\eqref{eq:poisson-dos-lower} yields
\begin{equation}
\label{eq:frozen-dos-uniform-lower}
\widetilde\varrho_0(\lambda,y)\geq c,
\qquad
T\geq2,\quad y\in\mathbb R,\quad\lambda\in J_\star.
\end{equation}
By \eqref{eqn:alph_bd}, 
$\widetilde\alpha(y)\leq C\mathfrak A$
uniformly in \(T\) and \(y\).
Hence, by \eqref{eq:frozen-dos-uniform-lower} and the definition of \(m_0\),
\begin{equation}
\label{eq:m0T-uniform-lower}
m_0(\lambda,x)
\geq c\mathfrak A^{-1},
\qquad
x\in\mathbb R,\quad\lambda\in J_\star.
\end{equation}

We now prove that the lower bound in \eqref{eq:m0T-uniform-lower}
continues to hold for \(m_\tau\), uniformly for every \(\tau\geq0\).
By Lemma \ref{lem:Jlow} and \eqref{eqn:Jtbd} from Lemma
\ref{lem:fixed_point_derivative0}, since \(J_\star\) is compact,
\[
0<J_0(\lambda,x)\leq C,
\qquad x\in\mathbb R,\quad\lambda\in J_\star.
\]
We claim that
\begin{equation}
\label{eq:varpi0-uniform-lower}
\varpi_\tau(\lambda;X)
\geq c\mathfrak A^{-1},
\qquad
\tau \geq 0, \quad X\in\mathbb R,\quad\lambda\in J_\star.
\end{equation}
Indeed, the relation
\eqref{eqn:vpdef} together with
\eqref{eq:m0T-uniform-lower} gives \eqref{eq:varpi0-uniform-lower} at $\tau = 0$;
by \eqref{eqn:wtX},
the same lower bound holds for $\varpi_\tau(\lambda;X)$ for every \(\tau\geq0\).

Finally, recall that \(m_\tau\) is
characterized by \eqref{eqn:mt}. By Lemma \ref{lem:Jlow},
$J_\tau(\lambda,x)\geq c>0$,
hence using \eqref{eq:varpi0-uniform-lower}, we obtain
\[
m_\tau(\lambda,x)
=
\varpi_\tau
\bigl(\lambda;X_\tau(\lambda,x)\bigr)
J_\tau(\lambda,x) 
\geq
c\mathfrak A^{-1}
\]
for every \(\tau\geq0\), \(x\in\mathbb R\), and
\(\lambda\in J_\star\). Decreasing $c$ so that $c\leq\mathfrak K^{-1}/4$ and replacing
$J_\star$ by $[\Lambda_\star-c,\Lambda_\star+c]$, where
$\Lambda_\star:=b_\star$, proves the lemma.
\end{proof}

\subsection{Concentration estimates for Lax matrix functionals}
\label{app:inhomogeneous-dos}

In this section we establish three results (Lemmas
\ref{lem:time-zero-site-concentration},
\ref{lem:app-time-zero-block-replacement}, and
\ref{lem:replace-initial-position}) providing bounds for functionals of
\(\L(0)\); they will later be used to show Proposition
\ref{prop:app-time-zero-dos} in Section
\ref{app:time-zero-dos-proof} below. We begin with the following lemma,
which is the counterpart in our setting of \cite[Lemma 4.6]{Agg25}; its
proof is also similar, so we only outline it.

\begin{lem}
\label{lem:time-zero-site-concentration}
Let \(h\) satisfy Assumption \ref{ass:admissible-test}. 
Then, with overwhelming probability,
\begin{equation}
\label{eq:time-zero-site-concentration}
\left|
\sum_{i=N_1}^{N_2}h(\Lambda_i,i/T)
-
\mathbb E\sum_{i=N_1}^{N_2}h(\Lambda_i,i/T)
\right|
\le
 T^{1/2}(\log N)^C A
+
T^{-90}B.
\end{equation}
\end{lem}

\begin{proof}[Proof of Lemma \ref{lem:time-zero-site-concentration} (Outline)]
Set $Z_h:=\sum_{i=N_1}^{N_2}h(\Lambda_i,i/T)$.
Only indices $i$ with $i/T\in I$ can contribute to $Z_h$, and
\begin{equation}
\label{eq:number-relevant-site-labels}
\#\{i\in\llbracket N_1,N_2\rrbracket:i/T\in I\}
\le
C T(\log N)^{20},
\end{equation}
which implies
\begin{equation}
\label{eq:Zh-deterministic-bound}
|Z_h|
\le
C T(\log N)^{20}A.
\end{equation}
We regard \(Z_h=Z_h((\xi_k)_k)\) as a function of the mutually
independent pairs
\(\xi_k:=(a_k(0),b_k(0))\),
\(k\in\llbracket N_1,N_2\rrbracket\), with the convention
\(a_{N_2}(0)=0\).

Fix
$L:= \lceil (\log N)^{5} \rceil$. Let
$\mathcal I_{\rm supp}
:=
\{i\in\llbracket N_1,N_2\rrbracket:i/T\in I\}$,
and let
$\mathcal I_{\rm supp}^+
:=
\{k\in\llbracket N_1,N_2\rrbracket:
\dist(k,\mathcal I_{\rm supp})\le 10L\}$.
Partition \(\mathcal I_{\rm supp}^+\) into $m$ consecutive nonempty index intervals
$D_1,\ldots,D_m$
of length between $L/2$ and \(L\), and let
$D_0:=\llbracket N_1,N_2\rrbracket\setminus\mathcal I_{\rm supp}^+$.
Denoting \(\xi(D_i):=\{\xi_k:k\in D_i\}\) for each
\(i\in\llbracket0,m\rrbracket\), observe that
\((\xi(D_i))_{0\leq i\leq m}\) decomposes the independent random
variables \((a_i,b_i)\) into a disjoint union.
By \eqref{eq:number-relevant-site-labels},
\begin{equation}
\label{eq:number-influence-groups}
m+1\le  T(\log N)^C.
\end{equation}

In what follows, we recall the influence
\(\operatorname{Infl}_{\xi(D_i)}(Z_h;p)\) from
\cite[Definition 5.1]{Agg25}. We claim that, for each
\(a\in\{0,1,\ldots,m\}\),
\begin{equation}
\label{eq:site-cloud-influence}
\operatorname{Infl}_{\xi(D_a)}(Z_h;p)
\le
(\log N)^C A+T^{-100}B,
\end{equation}
with $p:=e^{-c_0(\log N)^2}$ for some constant $c_0>0$;
\eqref{eq:site-cloud-influence} is the counterpart of
\cite[Lemma 5.3]{Agg25}. Given \eqref{eq:site-cloud-influence}, the proof
of the lemma is very similar to that of \cite[Lemma 4.6]{Agg25} given
\cite[Lemma 5.3]{Agg25}, as a consequence of the concentration bound
\cite[Lemma 5.2]{Agg25}; we therefore omit it, so it remains to verify
\eqref{eq:site-cloud-influence}.

Define a new matrix  \(\widetilde\L(0)\) as follows; starting with $\L(0)$, for every \(k\in D_a\), replace \(a_k(0)\) and \(b_k(0)\) by independent
copies \(\widetilde a_k(0)\) and \(\widetilde b_k(0)\), independent also
of all the original entries, and leave every entry with index outside
\(D_a\) unchanged. Write
\(\eig\widetilde\L(0)=(\widetilde\lambda_1,\widetilde\lambda_2,\ldots,
\widetilde\lambda_N)\), let
\(\widetilde\varphi_0:\llbracket 1,N\rrbracket\to\llbracket N_1,N_2\rrbracket\)
denote a \(\zeta\)-localization center bijection for \(\widetilde\L(0)\), and set
\(\widetilde\Lambda_i:=\widetilde\lambda_{\widetilde\varphi_0^{-1}(i)}\)
for each \(i\in\llbracket N_1,N_2\rrbracket\). Finally, set
\[
\widetilde Z_h
=
\sum_{i=N_1}^{N_2}
h(\widetilde\Lambda_i,i/T).
\]
By \cite[Definition 5.1]{Agg25}, it suffices to show
\begin{equation}\label{eqn:tildZhZh}
|\widetilde Z_h-Z_h| < (\log N)^C A+T^{-100}B
\end{equation}
with probability at least $1-p$.

The following proof of \eqref{eqn:tildZhZh} is similar to that of
\cite[Lemma 5.3]{Agg25}; first, fix
\(a\in\llbracket0,m\rrbracket\) and restrict to the events of
Lemma \ref{bijectionl}, Lemma \ref{l0eigenvalues} with \(A=\log N\),
and Lemma \ref{eigenvalues0} with
\(\delta=e^{-(\log N)^2}\) for both \(\L(0)\) and
\(\widetilde\L(0)\); their intersection has probability
at least $1-p$ for $c_0$ small enough. Applying Lemma \ref{lem:full-to-truncation} 
with $\L = \L(0)$, $\widetilde\L = \widetilde\L(0)$, and
$\mathcal{D}=D_a$, using Lemma \ref{bijectionl} to verify its hypotheses,
gives a bijection
\(\psi:\llbracket N_1,N_2\rrbracket\longrightarrow
\llbracket N_1,N_2\rrbracket\) such that, whenever
\(\dist(i,D_a)\geq3L\),
\begin{equation}
\label{eq:resampling-matching}
|\Lambda_i-\widetilde\Lambda_{\psi(i)}|
\le
e^{(\log N)^2-cL},
\qquad
|\psi(i)-i|\le(\log N)^2.
\end{equation}
For \(a=0\), every \(i\in\mathcal I_{\rm supp}\) is at distance at least
\(10L\) from \(D_0\). For \(a\ge1\), at most \(CL\) indices are within
distance \(3L\) of \(D_a\).

Since \(\psi\) is a bijection,
$\widetilde Z_h
=
\sum_{i=N_1}^{N_2}
h(\widetilde\Lambda_{\psi(i)},\psi(i)/T)$; moreover,
the indices within distance \(3L\) of \(D_a\) contribute at most
$CL A$
to $|Z_h-\widetilde Z_h|$. For every remaining index $i \in \mathcal I_{\rm supp}$, \eqref{eq:resampling-matching} gives
\begin{equation}\label{eq:resampling-spectral-error}
|h(\Lambda_i,i/T)
-
h(\widetilde\Lambda_{\psi(i)},\psi(i)/T)|
\le
e^{-c L} B
+
|h(\widetilde\Lambda_{\psi(i)},i/T)
-
h(\widetilde\Lambda_{\psi(i)},\psi(i)/T)|.
\end{equation}

It remains to control the second term on the right side of \eqref{eq:resampling-spectral-error}. For each
$\ell \in \llbracket 1, m \rrbracket$, let 
$s_{\ell}$ be such that
\(D_\ell=\llbracket s_\ell,s_\ell+|D_{\ell}|-1\rrbracket\), and define
\[
D_\ell^+
:=
\llbracket
s_\ell-2L,
s_\ell+3L-1
\rrbracket, \qquad
C_\ell^+
:=
\llbracket
\frac{s_\ell-2L}{T},
\frac{s_\ell+3L}{T}
\rrbracket
\qquad
g_\ell(\lambda)
:=
\int_{C_\ell^+}
|\partial_qh(\lambda,q)|\,dq.
\]
Because
\(|D_\ell|\leq L\), the interval \(D_\ell^+\) enlarges \(D_\ell\) by at
least \(2L\) on each side. If $\dist(i, D_a) \geq 3 L$, then
\(|\psi(i)-i|\le(\log N)^2\), so if \(i\in D_\ell\), then
the interval with endpoints \(i/T\) and \(\psi(i)/T\) is contained in
\(C_\ell^+\). Hence
\begin{equation}
\label{eq:spatial-shift-gell}
|h(\lambda,i/T)-h(\lambda,\psi(i)/T)|
\le g_\ell(\lambda),
\qquad i\in D_\ell.
\end{equation}
Moreover, by Assumption \ref{ass:admissible-test},
\begin{equation}
\label{eq:gell-uniform-bound}
0\le g_\ell(\lambda)\le A.
\end{equation}

On the event we have restricted to, Lemma
\ref{lem:full-to-truncation}, applied with \(\mathcal{I} =D_\ell^+\),
matches each eigenvalue
\(\widetilde\Lambda_{\psi(i)}\), with \(i\in D_\ell\),
injectively to an eigenvalue
$\mu_{\ell,i}\in\eig\widetilde\L^{D_\ell^{+}}(0)$
such that
$|\widetilde\Lambda_{\psi(i)}-\mu_{\ell,i}|
\le e^{-cL}$.
Therefore, for each $i \in D_{\ell}$ with $\dist(i, D_a) \geq 3 L$,
\[
|h(\widetilde\Lambda_{\psi(i)},i/T)
-
h(\widetilde\Lambda_{\psi(i)},\psi(i)/T)| \le
|h(\mu_{\ell,i},i/T)
-
h(\mu_{\ell,i},\psi(i)/T)|
+
2e^{-cL}B
\le
g_\ell(\mu_{\ell,i})
+
2 e^{-c L}B,
\]
and thus,
\begin{equation}
\label{eq:spatial-shift-block-traces}
\sum_{i : \dist(i, D_a)\geq3 L }
|h(\widetilde\Lambda_{\psi(i)},i/T)
-
h(\widetilde\Lambda_{\psi(i)},\psi(i)/T)|
\le
\sum_{\ell=1}^m X_\ell
+
e^{-c L}B,
\end{equation}
where
\[
X_\ell
:=
\Tr g_\ell(\widetilde\L^{D_\ell^{+}}(0)).
\]

We now prove that
\begin{equation}
\label{eq:block-trace-sum-bound}
\sum_{\ell=1}^m X_\ell
\le
C(\log N)^C A
\end{equation}
with probability at least \(1-p\), so long as $c_0$ is chosen small enough, 
since this is sufficient to prove \eqref{eqn:tildZhZh}, and thus the lemma,
by \eqref{eq:spatial-shift-block-traces},
\eqref{eq:resampling-spectral-error}, and the content above
\eqref{eq:resampling-spectral-error}.

First, Proposition \ref{prop:app-frozen-wegner} gives
\[
\mathbb EX_\ell
\le
C|D_\ell^+|
\int_{\mathbb R}
g_\ell(\lambda)e^{-c\lambda^2}\,d\lambda.
\]
Since \(cL\le |D_\ell|\le |D_\ell^+|\le CL\) for each
\(\ell\in\llbracket1,m\rrbracket\) and since, given
\(\ell\in\llbracket1,m\rrbracket\), there are at most \(C\) indices
\(\ell'\) for which \(C_\ell^+\cap C_{\ell'}^+\) is nonempty, we have
\[
\sum_{\ell=1}^m\mathbb EX_\ell
\le
CL\sum_{\ell=1}^m
\iint_{\mathbb R\times C_\ell^+}
|\partial_qh(\lambda,q)|
e^{-c\lambda^2}\,dq\,d\lambda
\leq
10 CL\iint_{\mathbb R^2}
\left|\partial_q h(\lambda,q)\right|
e^{-c\lambda^2}
\,dq\,d\lambda.
\]
By Assumption \ref{ass:admissible-test}, the double integral in the
right hand side above is upper bounded by $CA$,
and consequently
\begin{equation}
\label{eq:block-trace-total-mean}
\sum_{\ell=1}^m\mathbb E[X_\ell]
\le
CL A.
\end{equation}

Let \(\mathcal C_1,\mathcal C_2,\ldots,\mathcal C_{10}\) denote pairwise
disjoint subsets of \(\llbracket1,m\rrbracket\) whose union is
\(\llbracket1,m\rrbracket\), such that
\(D_\ell^+\cap D_{\ell'}^+\) is empty for any distinct
\(\ell,\ell'\in\mathcal C_i\) with \(i\in\llbracket1,10\rrbracket\).
As at the end of the proof of Lemma \ref{lem:weighted_rectangle}, the
random variables \((X_\ell)_{\ell\in\mathcal C_i}\) are mutually independent
for any \(i\in\llbracket1,10\rrbracket\). Moreover, by
\eqref{eq:gell-uniform-bound}, \(0\le X_\ell\le CLA\), so for any
\(\mathcal C\in\{\mathcal C_1,\mathcal C_2,\ldots,\mathcal C_{10}\}\),
\[
\sum_{\ell \in \mathcal{C}} \Var(X_\ell)
\le
\sum_{\ell \in \mathcal{C}} \mathbb EX_\ell^2
\le
CL A\sum_{\ell \in \mathcal{C}} \mathbb EX_\ell
\le
C L^2A^2.
\]
Bernstein's inequality, applied to each \(\mathcal C\), gives with
overwhelming probability
\[
\sum_{\ell\in\mathcal C}X_\ell
\le
\sum_{\ell\in\mathcal C}\mathbb EX_\ell+LA(\log N)^2
\le
C LA(\log N)^2.
\]
Summing over all ten such \(\mathcal C\) then gives
\eqref{eq:block-trace-sum-bound} and thus the lemma.
\end{proof}

The following lemma allows us to approximately decompose 
(the expected value of) a linear statistic of $\L(0)$ 
into (the expected value of) a sum of linear statistics of independent blocks of $\L(0)$.
Recall the set of discontinuity points $\mathcal{P}$ from Assumption~\ref{ass:entry-law-profiles}.
For this, let
\[
K=\lfloor T^{2/5}\rfloor,
\qquad
\widetilde{\mathcal B}_j:=\llbracket jK,(j+1)K-1\rrbracket,
\qquad j\in\mathbb Z.
\]
Let \(j_1<j_2<\cdots<j_M\) be precisely the integers \(j\) for which
\(\widetilde{\mathcal B}_j\subseteq\llbracket N_1,N_2\rrbracket\) and
\[
\left[\frac{jK}{T},\frac{(j+1)K}{T}\right]
\cap\bigl(\mathcal P\cup\{-(\log T)^5,(\log T)^5\}\bigr)
=\varnothing.
\]
Set
\[
\mathcal B_r:=\widetilde{\mathcal B}_{j_r},
\qquad
y_r:=\frac{\min \mathcal B_r}{T},
\qquad r\in\llbracket1,M\rrbracket.
\]

\begin{lem}
\label{lem:app-time-zero-block-replacement}
Let $h : \mathbb{R}^2 \rightarrow \mathbb{R}$ satisfy 
Assumption \ref{ass:admissible-test}, and adopt the notation 
above for $\mathcal B_r, y_r, 1 \leq r \leq M$.
Then
\begin{equation}
\label{eq:app-time-zero-block-replacement}
\left|
\mathbb{E}\left[ \sum_{i=N_1}^{N_2}h(\Lambda_i,i/T) \right]
-
\mathbb E \left[ \sum_{r=1}^{M}\Tr h(\L^{\mathcal B_r}(0),y_r) \right]
\right| 
\le
(\log N)^C
\left(
(K + \frac{T}{K})A
+
B
\right).
\end{equation}
\end{lem}

\begin{proof}
Set
$D=(\log N)^{5}$. The indices $k \in \llbracket N_1,N_2\rrbracket \setminus \left( \cup_{i=1}^M \mathcal B_i \right)$
either satisfy \(k \in \widetilde{\mathcal B}_j\) for some $j$ 
satisfying $\left[\frac{jK}{T},\frac{(j+1)K}{T}\right] \cap \{N_1, N_2 \} \neq \varnothing$,
or
\(k \in \widetilde{\mathcal B}_j\) for some $j$
such that $\left[\frac{jK}{T},\frac{(j+1)K}{T}\right]$ intersects
\(\mathcal P\cup\{-(\log T)^5,(\log T)^5\}\). Hence, there are at most
\(CK\) such indices, and their contribution to the first sum on the left
side of \eqref{eq:app-time-zero-block-replacement} is at most
\(CK\|h\|_\infty\).

By Assumption \ref{ass:admissible-test}, only indices
\(r\in\llbracket1,M\rrbracket\) for which \(\dist(y_r,I)\le K/T\) give a
nonzero contribution to the second sum on the left side of
\eqref{eq:app-time-zero-block-replacement}. The set
\(\mathcal J_h\subseteq\llbracket1,M\rrbracket\) of such indices satisfies
\(N_h:=|\mathcal J_h|\le CTK^{-1}(\log N)^{20}\).

Let us restrict to the overwhelmingly probable events
from Lemmas \ref{bijectionl},
\ref{l0eigenvalues} (with $A = \log N$), and \ref{eigenvalues0} (with
$\delta = e^{-(\log N)^2}$). For every \(r\in\mathcal J_h\), Lemma
\ref{lem:full-to-truncation}, applied with $\L = \L(0)$,
$\widetilde \L = \L^{\mathcal B_r}(0)$, and
$\mathcal{D} = \mathcal B_r^{\complement}$, injectively matches each \(\Lambda_i\),
with \(i\in\mathcal B_r\) and \(\dist(i,\mathcal B_r^{\complement})\ge D\), to an
eigenvalue \(\mu_{r,i}\in\eig\L^{\mathcal B_r}(0)\) such that
\[
|\Lambda_i-\mu_{r,i}|\le e^{-cD}.
\]
Consequently, at most \(CD\) eigenvalues \(\Lambda_i\) of \(\L(0)\) with
\(i\in\mathcal B_r\) are not matched to those of
\(\L^{\mathcal B_r}(0)\), and at most \(CD\) eigenvalues of
\(\L^{\mathcal B_r}(0)\) are not matched to those of \(\L(0)\). We refer
to such eigenvalues as unmatched.

The unmatched terms from both sums in
\eqref{eq:app-time-zero-block-replacement} contribute at most
\begin{equation}\label{eqn:hinfbd}
C D N_h\|h\|_\infty
\le
(\log N)^C\frac{T}{K}\|h\|_\infty.
\end{equation}
For the matched terms,
\[
|h(\Lambda_i,i/T)-h(\mu_{r,i},y_r)|
\leq
B|\Lambda_i-\mu_{r,i}|
+
\int_{y_r}^{i/T}
|\partial_qh(\mu_{r,i},q)|\,dq.
\]
Summing the first term over all \(r\in\mathcal J_h\) and all matched
\(i\in\mathcal B_r\) gives
\begin{equation}\label{eqn:Tlgcmatch}
 T (\log N)^C \,e^{-cD}B
\le
B.
\end{equation}

For the error term with the spatial derivative, 
\[
\sum_{r\in\mathcal J_h}\sum_{\substack{i\in \mathcal B_r\\
\dist(i,\mathcal B_r^{\complement})\ge D}}
\int_{y_r}^{i/T}
|\partial_qh(\mu_{r,i},q)|\,dq
\le
 \sum_{r\in\mathcal J_h}\int_{y_r}^{y_r+K/T}
\Tr\left|\partial_qh(\L^{\mathcal B_r}(0),q)\right|\,dq.
\]
Taking expectations and applying Proposition
\ref{prop:app-frozen-wegner} to the nonnegative function
\(\lambda\mapsto|\partial_qh(\lambda,q)|\), we obtain
\[
\mathbb E\Tr
\left|\partial_qh(\L^{\mathcal B_r}(0),q)\right|
\le
CK\int_{\mathbb R}
|\partial_qh(\lambda,q)|e^{-c\lambda^2}\,d\lambda.
\]
Therefore
\[
\begin{aligned}
\mathbb E\sum_{r\in\mathcal J_h} \int_{y_r}^{y_r+K/T}
\Tr\left|\partial_qh(\L^{\mathcal B_r}(0),q)\right|\,dq &\le
CK\sum_{r\in\mathcal J_h}\int_{y_r}^{y_r+K/T}
\int_{\mathbb R}
|\partial_qh(\lambda,q)|e^{-c\lambda^2}\,d\lambda\,dq\\
&\le CK \iint_{\mathbb R^2}
|\partial_qh(\lambda,q)|e^{-c\lambda^2}\,d\lambda\,dq.
\end{aligned}
\]
By Assumption \ref{ass:admissible-test}, the last double integral is at
most \(CA\).

Combining the matched and unmatched contributions from \eqref{eqn:Tlgcmatch} and \eqref{eqn:hinfbd}, and the contribution above
proves
\[
\left|
\mathbb E\sum_{i=N_1}^{N_2} h(\Lambda_i,i/T)
-
\mathbb E\sum_{r=1}^{M}\Tr h(\L^{\mathcal B_r}(0),y_r)
\right|
\le
(\log N)^C
\left(
KA
+
\frac{T}{K}A
+
B
\right),
\]
as desired.
\end{proof}

\begin{lem}
\label{lem:replace-initial-position}
Let \(h: \mathbb{R}^2 \rightarrow \mathbb{R} \) satisfy Assumption \ref{ass:admissible-test} and 
$\supp h \subset \mathbb{R} \times [-(\log N)^{20}, (\log N)^{20}]$.
Then, with overwhelming probability we have
\begin{equation}
\label{eq:replace-initial-position}
\left| \sum_{i=N_1}^{N_2}
\bigl(h(\Lambda_i,\bar q_i)-h(\Lambda_i,F_{\rm str}(i/T))\bigr)
\right|
\le
C \mathfrak A^{1/2}T^{1/2}(\log N)^C
A +B.
\end{equation}
\end{lem}

\begin{proof}
Set
\[
D:=(\log N)^{20},\qquad
F_i:=F_{\rm str}(i/T),\qquad
\delta:=T^{-1/2}(\log N)^{14}.
\]
Let \(\mathsf E_{\rm str}\) denote the event in
\eqref{eq:stretch-good}. By \eqref{eqn:alph_bd} and
\eqref{eqn:Fstr},
\begin{equation}
\label{eq:replacement-F-bilipschitz}
\mathfrak K^{-1}\mathfrak A(v-u)
\le
F_{\rm str}(v)-F_{\rm str}(u)
\le
\mathfrak K\mathfrak A(v-u),
\qquad u,v\in\mathbb R,\quad u\le v.
\end{equation}
Moreover, since \(q_0(0)=F_{\rm str}(0)=0\), on
\(\mathsf E_{\rm str}\),
\begin{equation}
\label{eq:replacement-stretch-at-zero}
|\bar q_i-F_i|
\le
\mathfrak A\frac{\sqrt{|i|}}{T}(\log N)^2,
\qquad i\in\llbracket N_1,N_2\rrbracket.
\end{equation}
If either \(h(\Lambda_i,\bar q_i)\) or \(h(\Lambda_i,F_i)\) is
nonzero, then either \(|\bar q_i|\le D\) or \(|F_i|\le D\).
Combining \eqref{eq:replacement-F-bilipschitz} and
\eqref{eq:replacement-stretch-at-zero} gives
\[
\mathfrak K^{-1}\mathfrak A|i|
\le
DT+\mathfrak A\sqrt{|i|}(\log N)^2 \leq 2 D T.
\]
Consequently, for every such index,
\begin{equation}\label{eqn:qbfonestr}
|\bar q_i-F_i|
\le
C\mathfrak A^{1/2}T^{-1/2}(\log N)^{12}
\le \delta
\end{equation}
for all sufficiently large \(N\). In particular, \(|F_i|\le D+\delta\).

Next, set
\[
J:=\left\lceil\frac{2D+4\delta}{\delta}\right\rceil,
\qquad
z_j:=-D-2\delta+j\delta,
\qquad
C_j:=[z_j,z_j+\delta),
\quad j\in\llbracket0,J-1\rrbracket.
\]
Thus the intervals \(C_j\) cover every possible value of \(F_i\) 
for which one of the
two summands in \eqref{eq:replace-initial-position} can be nonzero.
Define the sets
\[
S_j:=\{i\in\llbracket N_1,N_2\rrbracket:F_i\in C_j\}, \quad j\in\llbracket0,J-1\rrbracket.
\]
Discard the empty intervals and write
\(\mathcal J:=\{j\in\llbracket0,J-1\rrbracket:S_j\ne\varnothing\}\).
The lower bound in \eqref{eq:replacement-F-bilipschitz} shows that each
\(S_j\) is an index interval, that is, a finite set of consecutive integers, and
\begin{equation}
\label{eq:Sj-size}
|S_j|\le C(T\delta+1).
\end{equation}

Set \(C_j^+:=[z_j-2\delta,z_j+3\delta]\) and
\(g_j(\lambda):=\int_{C_j^+}|\partial_qh(\lambda,q)|\,dq\).
If \(i\in S_j\), then on \(\mathsf E_{\rm str}\), by 
\eqref{eqn:qbfonestr} both \(F_i\) and
\(\bar q_i\) belong to \(C_j^+\). Hence
\begin{equation}
\label{eq:replacement-cell-reduction}
|h(\lambda,\bar q_i)-h(\lambda,F_i)|
\le g_j(\lambda).
\end{equation}

Next, set \(D_0=(\log N)^5\), and let
\[
\widehat S_j
:=
\{k\in\llbracket N_1,N_2\rrbracket:\operatorname{dist}(k,S_j)\le D_0\}.
\]
Let us restrict further to the
event from Lemma \ref{bijectionl} for
\(\zeta=e^{-100(\log N)^{3/2}}\), the event
\(\mathsf{BND}_{\L(0)}(\log N)\) from Lemma
\ref{l0eigenvalues}, and the event
\(\mathsf{SEP}_{\L(0)}(e^{-(\log N)^2})\) from Lemma
\ref{eigenvalues0}. This intersection is overwhelmingly probable. Apply Lemma
\ref{lem:full-to-truncation} with
\(\mathcal I=\widehat S_j\),
\(\mathcal D=\llbracket N_1,N_2\rrbracket\setminus\widehat S_j\),
\(\L=\L(0)\), and
\(\widetilde\L=\L^{\widehat S_j}(0)\), where the last matrix is viewed
as an \(N\times N\) matrix by padding it with zeros outside
\(\widehat S_j\). The parameter \(\delta\) from Assumption
\ref{lmatrixl} is equal to \(0\) here. Since \(i\) is a
\(\zeta\)-localization center of
\(\Lambda_i\) and
\(\operatorname{dist}(i,\mathcal D)\ge D_0\) for \(i\in S_j\), the lemma
gives an eigenvalue
\(\mu_i^{(j)}\in\operatorname{eig}\L^{\widehat S_j}(0)\) such that
\(|\Lambda_i-\mu_i^{(j)}|\le e^{-cD_0}\).
Since
\(e^{-cD_0}<e^{-(\log N)^2}\), the eigenvalue-separation event shows
that the map \(S_j \ni i\mapsto\mu_i^{(j)}\) is injective. Moreover,
\eqref{eq:replacement-cell-reduction} and the definition of the 
Lipschitz constant $\|h\|_{\Lip,\lambda}$
 give, for every \(i\in S_j\),
\[
|h(\Lambda_i,\bar q_i)-h(\Lambda_i,F_i)|
\le
g_j(\mu_i^{(j)})
+2 e^{-cD_0} \|h\|_{\Lip,\lambda}.
\]
Consequently,
\begin{equation}
\label{eq:full-to-block-weighted}
\begin{aligned}
\sum_{i\in S_j}|h(\Lambda_i,\bar q_i)-h(\Lambda_i,F_i)|
&\le
\sum_{i\in S_j}g_j(\mu_i^{(j)})
+2|S_j|  e^{-cD_0} \|h\|_{\Lip,\lambda}\\
&\le
\sum_{\mu\in\operatorname{eig}\L^{\widehat S_j}(0)}g_j(\mu)
+e^{-cD_0}\|h\|_{\Lip,\lambda}.
\end{aligned}
\end{equation}

Set \(X_j:=\sum_{\mu\in\operatorname{eig}\L^{\widehat S_j}(0)}
g_j(\mu)\) for \(j\in\mathcal J\). It remains to upper bound
\(\sum_{j\in\mathcal J}X_j\). We proceed similarly to the proof of the
upper bound in \eqref{eqn:Xmalph}, by applying Bernstein's inequality
after dividing the indices in \(\mathcal J\) into parity classes. To that
end, Proposition
\ref{prop:app-frozen-wegner} gives
\begin{equation}
\label{eq:block-wegner-gj}
\mathbb E
\sum_{\mu\in\operatorname{eig}\L^{\widehat S_j}(0)}
g_j(\mu)
\le
C|\widehat S_j|
\int_{\mathbb R}g_j(\lambda)e^{-c\lambda^2}\,d\lambda
.
\end{equation}

By Assumption \ref{ass:admissible-test} and the fact that each
\(q\in\mathbb R\) belongs to at most six of the intervals \(C_j^+\),
\[
\begin{aligned}
\sum_{j \in \mathcal{J}}
\int_{\mathbb R}g_j(\lambda)e^{-c\lambda^2}\,d\lambda
=
\sum_{j \in \mathcal{J}}
\iint_{\mathbb R\times C_j^+}
|\partial_qh(\lambda,q)|e^{-c\lambda^2}
\,dq\,d\lambda
&\le 6\iint_{\mathbb R^2}
|\partial_qh(\lambda,q)|e^{-c\lambda^2}\,dq\,d\lambda \\
&\le CA.
\end{aligned}
\]
Since \(|\widehat S_j|\le C(T\delta+D_0)\), summing
\eqref{eq:block-wegner-gj} over \(j\in\mathcal J\) therefore gives
\[
\sum_{j\in\mathcal J}\mathbb E[X_j]
\le C(T\delta+D_0)A 
\le C T^{1/2}(\log N)^C A.
\]

Since \(0\le X_j\le |\widehat S_j|\|g_j\|_\infty
\le CT^{1/2}(\log N)^C A\), the preceding bound also
implies
\begin{equation}\label{eqn:sig2bd}
\sum_{j\in\mathcal J}\mathbb E[X_j^2]
\le  T(\log N)^C A^2.
\end{equation}

For \(\sigma\in\{0,1\}\), define the parity class
\(\mathcal J_\sigma:=\{j\in\mathcal J:j\equiv\sigma\pmod 2\}\).
If \(j,j'\in\mathcal J_\sigma\) are distinct, then
\eqref{eq:replacement-F-bilipschitz} gives
\(\operatorname{dist}(S_j,S_{j'})\ge c \mathfrak{A}^{-1}T\delta\). Since
\(D_0=(\log N)^5=o(T\delta)\), the enlarged blocks
\(\widehat S_j\), \(j\in\mathcal J_\sigma\), are pairwise disjoint for
all sufficiently large \(N\). The matrices
\(\L^{\widehat S_j}(0)\), \(j\in\mathcal J_\sigma\), therefore depend
on disjoint sets of the mutually independent entries of \(\L(0)\), so
the variables \((X_j)_{j\in\mathcal J_\sigma}\) are independent.
Bernstein's inequality, applied separately for
\(\sigma\in\{0,1\}\) using \eqref{eqn:sig2bd}, gives, with
overwhelming probability,
\begin{equation}\label{eqn:XjJbd}
\sum_{j \in \mathcal{J}}X_{j }
\le CT^{1/2}(\log N)^C A.
\end{equation}
Finally,
\(|\mathcal J|\le C D/\delta
\le C T^{1/2}(\log N)^6\),
and since \(D_0=(\log N)^5\), the sum over \(j\) of the last term in
\eqref{eq:full-to-block-weighted} is at most \(B\) for all sufficiently
large \(N\). Thus, summing \eqref{eq:full-to-block-weighted} over \(j \in \mathcal{J}\) 
and using \eqref{eqn:XjJbd} proves
\eqref{eq:replace-initial-position}.
\end{proof}

%

\subsection{Proof of Proposition \ref{prop:app-time-zero-dos}}
\label{app:time-zero-dos-proof}

\begin{proof}
Recall the blocks \(\mathcal B_1,\ldots,\mathcal B_M\) and the points
\(y_r\) defined above Lemma
\ref{lem:app-time-zero-block-replacement}. Set
\(K:=\lfloor T^{2/5}\rfloor\), and write
\[
Z(h):=\sum_{i=N_1}^{N_2}h(\Lambda_i,i/T), \qquad \mathfrak B(h):=
\sum_{r=1}^{M}\Tr h(\L^{\mathcal B_r}(0),y_r),
\]
and
\[
\mathcal R(h):=
K\sum_{r=1}^{M}
\int_{\mathbb R}
h(\lambda,y_r)\widetilde\varrho_0(\lambda,y_r)\,d\lambda, \qquad  \mathcal I(h):=
T\iint_{\mathbb R^2}
h(\lambda,y)\widetilde\varrho_0(\lambda,y)\,d\lambda\,dy.
\]
We first prove \eqref{eqn:hiTconv}. We compare the above four quantities in turn.

\medskip
\noindent
\textit{Step 1: Bounding \(\left|Z(h)-\mathbb E[\mathfrak B(h)]\right|\).}
By \eqref{eq:time-zero-site-concentration} from Lemma
\ref{lem:time-zero-site-concentration} and the estimate from Lemma
\ref{lem:app-time-zero-block-replacement}, with overwhelming probability,
\begin{equation}
\label{eq:dos-proof-concentration-block}
\begin{aligned}
\left|Z(h)-\mathbb E[\mathfrak B(h)]\right|
&\le |Z(h)-\mathbb E[Z(h)]|
+\left|\mathbb E[Z(h)]-\mathbb E[\mathfrak B(h)]\right|\\
&\le C(\log N)^C\left(AT^{3/5}+B\right).
\end{aligned}
\end{equation}

\medskip
\noindent
\textit{Step 2: Bounding \(\left|
\mathbb E[\mathfrak B(h)]-\mathcal R(h)
\right|\).}
For each \(r\in\llbracket1,M\rrbracket\), define
\(h_r(\lambda):=h(\lambda,y_r)\). By Assumption
\ref{ass:admissible-test}, \(\|h_r\|_\infty\le A\) and
\(\|h_r\|_{\Lip}\le B\).
By the definition of the blocks \(\mathcal B_r\),
$[y_r,y_r+K/T]\cap
\bigl(\mathcal P\cup\{-(\log T)^5,(\log T)^5\}\bigr)=\varnothing$. By
\eqref{eq:app-frozen-block-coupling} from
Lemma \ref{lem:app-frozen-couple}, applied with its parameter \(A\) set to
\(A+B\), 
\begin{equation}
\label{eq:dos-proof-frozen-coupling}
\left|
\mathbb E\Tr h_r(\L^{\mathcal B_r}(0))
-
\mathbb E\Tr h_r(\mathbf{S}_{y_r}^{(K)})
\right|
\le
C(A+B)\frac{K^2}{T}.
\end{equation}
Moreover, Lemma \ref{lem:app-frozen-moving-tests} gives
\begin{equation}
\label{eq:dos-proof-frozen-moving-test}
\left|
\mathbb E\Tr h_r(\mathbf{S}_{y_r}^{(K)})
-
K\int_{\mathbb R}
h(\lambda,y_r)\widetilde\varrho_0(\lambda,y_r)\,d\lambda
\right|
\le
C\left(A(\log N)^C + B T^{-90}\right).
\end{equation}

We wish to sum the errors in the previous two displays in order to bound
\(\left|\mathbb E[\mathfrak B(h)]-\mathcal R(h)\right|\).
Only indices \(r\in\llbracket1,M\rrbracket\) such that \(y_r\in I\)
contribute. By the assumption on \(I\), the number \(N_h\) of such
indices satisfies \(N_h\le CT(\log N)^C/K\).
Summing the preceding two estimates
\eqref{eq:dos-proof-frozen-coupling} and \eqref{eq:dos-proof-frozen-moving-test}
over these indices therefore
gives
\begin{align}
\left|
\mathbb E[\mathfrak B(h)]-\mathcal R(h)
\right|
&\le
 \frac{T (\log N)^C}{K}
\left[
(A+B)\frac{K^2}{T}
+
A
\right]  = (\log N)^C
\left[
(A+B)K+A\frac TK
\right].
\label{eq:dos-proof-freezing}
\end{align}

\medskip
\noindent
\textit{Step 3: Bounding \(|\mathcal R(h)-\mathcal I(h)|\).}
Set
$f(y):=
\int_{\mathbb R}
h(\lambda,y)\widetilde\varrho_0(\lambda,y)\,d\lambda$.
For \(y\in[y_r,y_r+K/T]\), write
\[
f(y)-f(y_r)=E_{r,1}(y)+E_{r,2}(y),
\]
where
\[
E_{r,1}(y):=
\int_{\mathbb R}
\bigl(h(\lambda,y)-h(\lambda,y_r)\bigr)
\widetilde\varrho_0(\lambda,y)\,d\lambda
\]
and
\[
E_{r,2}(y):=
\int_{\mathbb R}h(\lambda,y_r)
\left[
\widetilde\varrho_0(\lambda,y)
-
\widetilde\varrho_0(\lambda,y_r)
\right]d\lambda.
\]

Recall by the construction of $y_1,\dots, y_M$ that 
\[
\left[N_1/T,N_2/T\right]\setminus
\bigcup_{r=1}^M[y_r,y_r+K/T]
\]
has length at most \(CK/T\). Since
\(\operatorname{supp}f\subseteq[N_1/T,N_2/T]\) and
\(\|f\|_\infty\le A\) by Assumption \ref{ass:admissible-test}, the
contribution to \(\mathcal I(h)\) from the set
$(\lambda,y)\in\mathbb R\times
\left(\left[N_1/T,N_2/T\right]\setminus
\bigcup_{r=1}^M[y_r,y_r+K/T]\right)$
is at most \(CAK\). It remains to estimate 
\(|\mathcal R(h)-T \iint_{\mathcal{S}} h(\lambda, y) \widetilde \varrho_0(\lambda, y) d\lambda d y |  \), with
\(\mathcal{S} = \mathbb R\times\bigcup_{r=1}^M[y_r,y_r+K/T]\).

By the Gaussian upper bound for
\(\widetilde\varrho_0\) from Lemma \ref{lem:app-frozen-dos-regularity},
\[
|E_{r,1}(y)|
\le
C\int_{y_r}^{y}
\int_{\mathbb R}
|\partial_qh(\lambda,s)|e^{-c\lambda^2}
\,d\lambda\,ds.
\]
Hence, by Assumption \ref{ass:admissible-test},
\begin{align}
T\sum_{r=1}^{M} \int_{y_r}^{y_r+K/T}|E_{r,1}(y)|\,dy
&\le
CK
\iint_{\mathbb R^2}
|\partial_qh(\lambda,s)|e^{-c\lambda^2}
\,d\lambda\,ds \le CKA.
\label{eq:dos-proof-riemann-h}
\end{align}

For the second term, apply \eqref{eqn:nuydif} from Lemma
\ref{lem:app-frozen-couple} to the function
\(\lambda\mapsto h(\lambda,y_r)\).
Since the supremum and Lipschitz
norms of $h(\cdot, y_r)$ are bounded by \(A\) and \(B\), respectively,
\[
|E_{r,2}(y)|
\le
C(A+B)|y-y_r|.
\]
There are at most \(C T(\log N)^C/K\) indices $r$ for which 
$E_{r,2}(y) \neq 0$, and therefore
\begin{equation}
T\sum_{r=1}^{M}\int_{y_r}^{y_r+K/T}|E_{r,2}(y)|\,dy
\le
 T\frac{T(\log N)^C}{K}
(A+B)\left(\frac KT\right)^2
\le
(\log N)^C K(A+B).
\label{eq:dos-proof-riemann-dos}
\end{equation}
Combining the estimate for the omitted cells with
\eqref{eq:dos-proof-riemann-h} and
\eqref{eq:dos-proof-riemann-dos}, we obtain
\begin{equation}
\label{eq:dos-proof-riemann}
|\mathcal R(h)-\mathcal I(h)|
\le
(\log N)^C K(A+B).
\end{equation}

Combining \eqref{eq:dos-proof-concentration-block},
\eqref{eq:dos-proof-freezing}, and \eqref{eq:dos-proof-riemann}, and
recalling that \(K=\lfloor T^{2/5}\rfloor\), gives, with overwhelming
probability,
\begin{equation}
\label{eq:dos-proof-site-general}
|Z(h)-\mathcal I(h)|
\le
(\log N)^C
\left(
A T^{3/5}+B T^{2/5}
\right),
\end{equation}
and thus proves \eqref{eqn:hiTconv}.

\medskip
We now prove \eqref{eqn:hiTspace}. Define
\(\widehat h(\lambda,y):=h(\lambda,F_{\rm str}(y))\). Then
\(\widehat h(\Lambda_i,i/T)=h(\Lambda_i,F_{\rm str}(i/T))\), and it is
quickly verified, using the fact that
\(F_{\rm str}'\ge\mathfrak K^{-1}\mathfrak A>1\) by
\eqref{eqn:alph_bd}, that \(\widehat h\) satisfies Assumption
\ref{ass:admissible-test} with the same parameters \(A\) and \(B\).
By \eqref{eq:replace-initial-position} from Lemma
\ref{lem:replace-initial-position}, with overwhelming
probability,
\begin{equation}
\label{eq:dos-proof-position-replacement}
\left|
\sum_{i=N_1}^{N_2}h(\Lambda_i,\bar q_i)
-
\sum_{i=N_1}^{N_2}\widehat h(\Lambda_i,i/T)
\right|
\le
C\mathfrak A^{1/2}T^{1/2}(\log N)^C A+B.
\end{equation}
Applying \eqref{eq:dos-proof-site-general} with \(h\) replaced by
\(\widehat h\) gives
\begin{equation}
\label{eq:dos-proof-composed-site}
\left|
\sum_{i=N_1}^{N_2}\widehat h(\Lambda_i,i/T)
-
T\iint_{\mathbb R^2}
\widehat h(\lambda,y)
\widetilde\varrho_0(\lambda,y)\,d\lambda\,dy
\right|
\le
C(\log N)^C
\left(
A T^{3/5}+B T^{2/5}
\right).
\end{equation}

Finally, by Definition \ref{rho0} and \eqref{eqn:Fstr}, changing
variables \(q=F_{\rm str}(y)\) gives
\[
\iint_{\mathbb R^2}
\widehat h(\lambda,y)
\widetilde\varrho_0(\lambda,y)\,d\lambda\,dy=
\iint_{\mathbb R^2}
h(\lambda,F_{\rm str}(y))
\widetilde\varrho_0(\lambda,y)\,d\lambda\,dy =
\iint_{\mathbb R^2}
h(\lambda,q)m_0(\lambda,q)\,d\lambda\,dq.
\]
Since \(\mathfrak A\) is fixed,
\(\mathfrak A^{1/2}T^{1/2}\le T^{3/5}\) for all sufficiently large
\(T\). Combining the preceding identity with
\eqref{eq:dos-proof-position-replacement} and
\eqref{eq:dos-proof-composed-site} proves
\eqref{eqn:hiTspace}.
\end{proof}

\subsection{Proof of Proposition \ref{prop:scat_conc}}
\label{app:scat_conc}

Recall that, at the beginning of Appendix \ref{app:time-zero-dos}, we
suppressed the superscript \((T)\) on objects constructed from the
\(T\)-dependent entry laws. In this subsection, we restore this
superscript on all hydrodynamic limit quantities defined under Assumption
\ref{ass:model}, as in Section \ref{subsec:tc}.

\begin{proof}
Put $\rho:=T^{-3/10}$.
We restrict to the overwhelmingly probable events supplied by
Lemma \ref{l0eigenvalues} with $A=\log N$; Lemma \ref{eigenvalues0} with
$\delta=e^{-(\log N)^2}$; Lemma \ref{lem:sumlambdak01}, applied at
times $0$ and $\tau$ with $\delta=1/20$, $\mathfrak M=T^{7/10}$,
and for each
 $J_n$, $0\le n\le\lceil\log N\rceil$; and to the event \eqref{eq:stretch-good}.
After relabelling \(j=\varphi_0^{-1}(i)\), the difference of summations 
in the left hand side of the bound in
the statement is
\[
S_{k,a}
:=
\sum_{\substack{i\ne i_k}}
\log|\lambda_k-\Lambda_i|
\left(
\mathbbm 1_{\{\bar q_i<\bar q_{i_k}\}}
-
\mathbbm 1_{\{x_\tau^{(T)}(\Lambda_i;\bar q_i)<a\}}
\right).
\]

Let \(\vartheta:\mathbb R\to[0,1]\) be a smooth nonincreasing function
satisfying $\vartheta(r)=1$ for $r\le-1$ and
 $\vartheta(r)=0$ for $r \geq 1$, and set
\[
\vartheta_\rho(r):=\vartheta(r/\rho).
\]
Let also \(\chi : \mathbb{R} \rightarrow \mathbb{R}\) be a smooth function satisfying
$\chi=1$ on $[-1,1]$ and $\supp\chi\subset[-2,2]$,
and define
\[
\chi_N(\lambda)
:=
\chi\left(\frac{\lambda}{2\log N}\right).
\]
For \(|\xi|\le\log N\), put
\[
\ell_{\xi,\rho}(\lambda)
:=
\chi_N(\lambda)
\log\bigl(|\lambda-\xi|\vee\rho\bigr).
\]
Then
\begin{equation}
\label{eq:scat-regularized-log-bounds}
\|\ell_{\xi,\rho}\|_\infty\le C\log N,
\qquad
\|\ell_{\xi,\rho}\|_{\Lip}\le C\rho^{-1}(\log N)^C.
\end{equation}

For any \(u,a,\lambda,q\in\mathbb R\) and
\(\xi\in[-\log N,\log N]\), define the regularized test function (which
might be viewed as an approximation to the summand in \(S_{k,a}\), when
\(u=\bar q_{i_k}\), \(\xi=\lambda_k\), \(\lambda=\Lambda_i\), and
\(q=\bar q_i\))
\begin{equation}
\label{eq:scat-regularized-test}
h_{\xi,u,a}^{(T)}(\lambda,q)
:=
\ell_{\xi,\rho}(\lambda)
\left[
\vartheta_\rho(q-u)
-
\vartheta_\rho
\left(x_\tau^{(T)}(\lambda;q)-a\right)
\right].
\end{equation}
First, we claim that on an
event of overwhelming probability, for all
\(\xi\in[-\log N,\log N]\) and
\(u,a\in[-(\log N)^{15},(\log N)^{15}]\),
\begin{equation}
\left|
\sum_{i=N_1}^{N_2}
h_{\xi,u,a}^{(T)}(\Lambda_i,\bar q_i)
-
T\iint_{\mathbb R^2}
h_{\xi,u,a}^{(T)}(\lambda,q)m_0^{(T)}(\lambda,q)
\,d\lambda\,dq
\right| 
\le
C T^{7/10}(\log N)^C.
\label{eq:scat-smooth-dos}
\end{equation}

We first verify the assumptions of Proposition \ref{prop:app-time-zero-dos} for 
arbitrary fixed $u, a \in  [-(\log N)^{15}, (\log N)^{15}]$. The 
map \(q\mapsto x_\tau^{(T)}(\lambda;q)\) is increasing and invertible by \eqref{eqn:Xtlin} and Lemma \ref{lem:Jlow}. Hence
\[
\int_{\mathbb R}\left|\partial_qh_{\xi,u,a}^{(T)}(\lambda,q)\right|\,dq
\le 2|\ell_{\xi,\rho}(\lambda)|
\|\vartheta_\rho'\|_{L^1(\mathbb R)}\le C\log N.
\]
For the second term in \eqref{eq:scat-regularized-test} we used the change of variables
\(x=x_\tau^{(T)}(\lambda;q)\). Also,
\[
\|h_{\xi,u,a}^{(T)}\|_\infty\le C\log N,
\]
and, since \eqref{eq:flow-spec} gives, uniformly for \(|\lambda|\le4\log N\),
$\left|
\partial_\lambda x_\tau^{(T)}(\lambda;q)
\right|
\le(\log N)^3$,
we have (using \eqref{eq:scat-regularized-log-bounds})
\[
\left\|
h_{\xi,u,a}^{(T)}
\right\|_{\Lip,\lambda}
\le
C\rho^{-1}(\log N)^C
=
C T^{3/10}(\log N)^C.
\]
In addition, \eqref{eqn:vebd} shows that the support of
\(h_{\xi,u,a}^{(T)}\) lies in 
$\mathbb{R} \times [-(\log N)^{20},(\log N)^{20}]$.

We now restrict further to the intersection of events from Proposition
\ref{prop:app-time-zero-dos}. More precisely, take deterministic
\(T^{-100}\)-meshes 
\[
\xi\in[-\log N,\log N] \cap T^{-100}\mathbb{Z} ,
\qquad
u,a\in[-(\log N)^{15},(\log N)^{15}] \cap T^{-100}\mathbb{Z} .
\]
The resulting
family has cardinality at most $T^{301}$, and the estimates from the preceding paragraph 
show that Assumption 
\ref{ass:admissible-test} is satisfied by 
$h_{\xi,u,a}^{(T)}$ for each $\xi, u, a$ in the family, 
with the same parameters $A = (\log N)^2$ and $B = T^{3/10} (\log N)^C$.
	Thus, by \eqref{eqn:hiTspace}  (whose additional condition on $I$ follows from
	the preceding support bound), there
is a constant $C>0$ and
an overwhelmingly probable event on which
for each \(a,u,\xi\) in the meshes specified above,
\begin{equation}
\left|
\sum_{i=N_1}^{N_2}
h_{\xi,u,a}^{(T)}(\Lambda_i,\bar q_i)
-
T\iint_{\mathbb R^2}
h_{\xi,u,a}^{(T)}(\lambda,q)m_0^{(T)}(\lambda,q)
\,d\lambda\,dq
\right| 
\le
(\log N)^C
\left(
T^{3/5}
+
T^{7/10}
\right).  \label{eq:scat-smooth-dos2}
\end{equation}
Moreover, by the Lipschitz dependence of
\eqref{eq:scat-regularized-test} on \((\xi,u,a)\) (together with the bound
on $\varrho_0^{(T)}$ in Lemma \ref{rho0estimate}) we
conclude that \eqref{eq:scat-smooth-dos2}
holds not only on the mesh but for arbitrary values of
\(\xi\in[-\log N,\log N]\) and
\(u,a\in[-(\log N)^{15},(\log N)^{15}]\),
up to an error upper bounded by $2 (\log N)^C T^{7/10}$. Consequently, the
claimed bound \eqref{eq:scat-smooth-dos} holds with overwhelming
probability. By our restriction to \eqref{eq:stretch-good} and the event of 
Lemma \ref{l0eigenvalues}, $|\lambda_k| \leq \log N$ and 
$ \bar q_{i_k} \in [-(\log N)^{15},(\log N)^{15}]$, so we may specialize this bound to $\xi=\lambda_k$ and
 $u = \bar q_{i_k}$.

We next bound $\left|S_{k,a}
-
\sum_{\substack{i\ne i_k}}
h_{\lambda_k,\bar q_{i_k},a}^{(T)}(\Lambda_i,\bar q_i)
\right|$. Note that
\begin{equation}\label{eqn:1u0}
\left|
\mathbbm 1_{\{\bar q_i<\bar q_{i_k}\}}
-
\vartheta_\rho(\bar q_i-\bar q_{i_k})
\right|
\le
\mathbbm 1_{\{|\bar q_i-\bar q_{i_k}|\le\rho\}},
\end{equation}
and
\begin{equation}\label{eqn:1ut}
\left|
\mathbbm 1_{\{x_\tau^{(T)}(\Lambda_i;\bar q_i)<a\}}
-
\vartheta_\rho\bigl(x_\tau^{(T)}(\Lambda_i;\bar q_i)-a\bigr)
\right|
\le
\mathbbm 1_{\{|x_\tau^{(T)}(\Lambda_i;\bar q_i)-a|\le\rho\}}.
\end{equation}
Here \(\rho=\mathfrak M/T\), so the intervals in both applications
of Lemma
\ref{lem:sumlambdak01} 
below have length \(2\rho\leq2\mathfrak A\mathfrak M/T\), as required.
Two applications of \eqref{eqn:sumlambdak01_shell} with \(r=2\) and summed over
$0 < n  < \lceil \log N \rceil$, therefore give, uniformly for
$a\in[-(\log N)^{15},(\log N)^{15}]$,
\begin{align}
&\sum_{\substack{i\ne i_k\\|\bar q_i-\bar q_{i_k}|\le\rho}}
\bigl|\log|\lambda_k-\Lambda_i|\bigr|
+
\sum_{\substack{i\ne i_k\\|x_\tau^{(T)}(\Lambda_i;\bar q_i)-a|\le\rho}}
\bigl|\log|\lambda_k-\Lambda_i|\bigr| \notag\\
&\qquad\le
C\left(T\rho+(\log N)^C\right)\log(|\lambda_k|+2)
\le
T^{7/10}(\log N)^C\log(|\lambda_k|+2).
\label{eq:scat-spatial-boundary}
\end{align}
This implies 
$$
\left| 
\sum_{\substack{i\ne i_k}} h_{\lambda_k,\bar q_{i_k},a}^{(T)}(\Lambda_i,\bar q_i) - 
\sum_{\substack{i\ne i_k}} \ell_{\lambda_k,\rho}(\Lambda_i)
\left(
\mathbbm 1_{\{\bar q_i<\bar q_{i_k}\}}
-
\mathbbm 1_{\{x_\tau^{(T)}(\Lambda_i;\bar q_i)<a\}}
\right) \right| \leq T^{7/10}(\log N)^C\log(|\lambda_k|+2).$$

It remains to bound the error coming from replacing $\ell_{\lambda_k,\rho}(\Lambda_i)$ with 
$\log|\lambda_k - \Lambda_i|$ in the second summation above. Whenever
$|\Lambda_i-\lambda_k|\ge\rho$,
and \(|\Lambda_i|\le\log N\), one has
$\ell_{\lambda_k,\rho}(\Lambda_i)
=
\log|\lambda_k-\Lambda_i|$, so it suffices
to bound the error coming from indices $i$ with $|\Lambda_i-\lambda_k| < \rho$.
By \eqref{eqn:xtode} and the first bound in \eqref{eqn:vebd}, if
\(\mathbbm 1_{\{\bar q_i<\bar q_{i_k}\}}-
\mathbbm 1_{\{x_\tau^{(T)}(\Lambda_i;\bar q_i)<a\}}\ne0\), then
\(\bar q_i\in K_0:=[-(\log N)^{20},(\log N)^{20}]\). For
$m\in\mathbb Z$, set
\[
J_m^\rho
:=
[(m-2)\rho,(m+3)\rho]\cap[-\log N,\log N],
\]
There are at most \(C\rho^{-1}\log N\) indices \(m\in\mathbb Z\) for
which
\([m\rho,(m+1)\rho]\cap[-\log N,\log N]\) is nonempty; we only consider
the associated \(J_m^\rho\) in the below. We now further restrict to the
intersection of the events in \eqref{eq:rect-thin} from Lemma
\ref{lem:weighted_rectangle}, applied at time $0$ with $I=K_0$ and
$J=J_m^\rho$ for every such $m$; this remains an overwhelmingly
probable event. For each $k$, choose $m$ such that
$\lambda_k\in[m\rho,(m+1)\rho]$. On the event from Lemma
\ref{l0eigenvalues}, every $\Lambda_i$ satisfying
$|\Lambda_i-\lambda_k|\le2\rho$ belongs to $J_m^\rho$. Hence
\eqref{eq:rect-thin} gives
\begin{equation}\label{eqn:numK0xi}
\#\left\{
i:\bar q_i\in K_0,\
0<|\Lambda_i-\lambda_k|\le2\rho
\right\}
\le
C T|K_0|\rho+C(\log N)^C.
\end{equation}
By eigenvalue separation,
\[
\bigl|\log|\lambda_k-\Lambda_i|\bigr|
\le(\log N)^2
\qquad (i\ne i_k).
\]
It follows from \eqref{eqn:1u0}, \eqref{eqn:1ut}, \eqref{eq:scat-spatial-boundary}, \eqref{eqn:numK0xi} and the bound above,
and adding the omitted \(i=i_k\) term (as in \eqref{eq:scat-smooth-dos}),
that
\begin{equation}
\label{eq:scat-log-regularization-discrete}
\left|
S_{k,a}
-
\sum_{i=N_1}^{N_2}
h_{\lambda_k,\bar q_{i_k},a}^{(T)}(\Lambda_i,\bar q_i)
\right|
\le
 T^{7/10}(\log N)^C\log(|\lambda_k|+2).
\end{equation}

Next, we identify the integral in \eqref{eq:scat-smooth-dos2}
with the integral $\mathcal{I}_{\tau, k}^{(T)}(a)$ in the statement of
the proposition. By Remark \ref{rem:pushrem},
$N_\tau^{(T)}(\lambda,x_\tau^{(T)}(\lambda;q))=N_0^{(T)}(\lambda,q)$. Writing
$\Phi_\tau^{(T)}(a,\lambda):=(x_\tau^{(T)}(\lambda;\cdot))^{-1}(a)$,
we have
$N_\tau^{(T)}(\lambda,a)=N_0^{(T)}(\lambda,\Phi_\tau^{(T)}(a,\lambda))$. Since
$q\mapsto x_\tau^{(T)}(\lambda;q)$ is increasing, moreover,
$\mathbbm 1_{\{x_\tau^{(T)}(\lambda;q)<a\}}
=\mathbbm 1_{\{q<\Phi_\tau^{(T)}(a,\lambda)\}}$.
It follows that
\begin{align}
N_0^{(T)}(\lambda,\bar q_{i_k})-N_\tau^{(T)}(\lambda,a)
&=
\int_{\mathbb R}
\left[
\mathbbm 1_{\{q<\bar q_{i_k}\}}
-
\mathbbm 1_{\{x_\tau^{(T)}(\lambda;q)<a\}}
\right]
m_0^{(T)}(\lambda,q)\,dq.
\label{eq:scat-primitive-halfspace}
\end{align}
Consequently,
\begin{align}
\mathcal I_{\tau,k}^{(T)}(a)
&=
\iint_{\mathbb R^2}
\log|\lambda_k-\lambda|
\left[
\mathbbm 1_{\{q<\bar q_{i_k}\}}
-
\mathbbm 1_{\{x_\tau^{(T)}(\lambda;q)<a\}}
\right]
m_0^{(T)}(\lambda,q)\,dq\,d\lambda.
\label{eq:scat-I-halfspace}
\end{align}

Using \eqref{eqn:1u0}, \eqref{eqn:1ut}, the bound on
$m_0^{(T)}$ from Lemma \ref{rho0estimate},
\eqref{eq:backflow-space}, and the estimate
\[
\int_{\mathbb R}
|\log|\lambda_k-\lambda||e^{-c\lambda^2}\,d\lambda
\le
C\log(|\lambda_k|+2),
\]
we obtain
\begin{align}
\left|
T\iint_{\mathbb R^2}
h_{\lambda_k,\bar q_{i_k},a}^{(T)}(\lambda,q)m_0^{(T)}(\lambda,q)\,d\lambda\,dq
-
T\mathcal I_{\tau,k}^{(T)}(a)
\right| 
&\le
C T\rho(\log N)^C\log(|\lambda_k|+2) \notag\\
&=
C T^{7/10}(\log N)^C\log(|\lambda_k|+2).
\label{eq:scat-continuum-regularization}
\end{align}

Combining
\eqref{eq:scat-smooth-dos},
\eqref{eq:scat-log-regularization-discrete}, and
\eqref{eq:scat-continuum-regularization}, we conclude that
\[
\left|
S_{k,a}
-
T\mathcal I_{\tau,k}^{(T)}(a)
\right|
\le
C T^{7/10}(\log N)^C\log(|\lambda_k|+2).
\]
The right hand side above is at most
\(T^{3/4-1/100}\log(|\lambda_k|+2)\), which proves the proposition.
\end{proof}

	\section{Proofs of Lemma \ref{lem:time-0-locR} and Proposition \ref{prop:outside-strip-confinement}}
\label{app:lemLocR}

\subsection{Proof of Lemma \ref{lem:time-0-locR}}
\label{app:time0loc}

	\begin{proof}
We prove the case \(i<\varphi_0(j) =: \varphi\); the other case follows by
 symmetry. Let
$d:=\varphi-i$, and $\H_m:=\L^{[N_1,m-1]}(0)$.
If 
$d\le \frac{T^{2/3}}{20R}$,
then the desired estimate follows from \(|u_j(i)|\le1\). We
therefore assume throughout that
$d>\frac{T^{2/3}}{20R}$.

Lemmas \ref{smatrixu} and \ref{sij} give the exact identity
\begin{equation}
\label{eq:time-zero-transfer-short}
\begin{aligned}
2\log|u_j(i)|
&=
q_i(0)-q_\varphi(0)+2\log|u_j(\varphi)|\\
&\quad
+2\sum_{\mu\in\eig \H_i}\log|\lambda_j-\mu|
-2\sum_{\nu\in\eig \H_\varphi}\log|\lambda_j-\nu|.
\end{aligned}
\end{equation}
We claim, first that with overwhelming probability,
\begin{equation}
\label{eq:time-zero-log-reduction-short}
\sum_{\mu\in\eig \H_i}\log|\lambda_j-\mu|
-
\sum_{\nu\in\eig \H_\varphi}\log|\lambda_j-\nu|
\le
-\sum_{a=i}^{\varphi-1}
\log|\lambda_j-\Lambda_a|
+
C(\log N)^C,
\end{equation}
and second that 
with overwhelming probability,
\begin{equation}\label{eq:time-zero-singular-block}
-\sum_{a=i}^{\varphi-1}
\log|\lambda_j-\Lambda_a|
\le Cd.
\end{equation}

The proof follows from the two previous claims and \eqref{eq:time-zero-transfer-short} as follows. 
Lemma \ref{lem:Estr} and \eqref{eqn:alph_bd} give with overwhelming probability
\[
q_\varphi(0)-q_i(0)
\ge \mathfrak A \left( d-\sqrt d(\log N)^2 \right).
\]
Using \(2\log|u_j(\varphi)|\le0\) in
\eqref{eq:time-zero-transfer-short}, and then applying
\eqref{eq:time-zero-log-reduction-short} and
\eqref{eq:time-zero-singular-block}, we obtain
\[
2\log|u_j(i)|
\le
-\mathfrak A d+C d+\mathfrak{A}\sqrt d(\log N)^C.
\]
Finally, choose \(\mathfrak A\) sufficiently large that \(C\le\mathfrak A/4\)
and \(\mathfrak A/4\ge20R\).
Then,
\[
2\log|u_j(i)|
\le
-\frac12\mathfrak A d
\le
-40Rd.
\]
Dividing by \(2\) gives the desired result.

Next, we prove \eqref{eq:time-zero-log-reduction-short}. 
Let $\mathcal A_m\coloneqq \{\Lambda_a:N_1\le a\le m-D\}$, 
with $D = \lfloor (\log N)^5 \rfloor$.
On the events of Lemmas
\ref{bijectionl}, \ref{l0eigenvalues}, and \ref{eigenvalues0},
with parameters $A=\log N$ and $\delta=e^{-(\log N)^2}$
in the latter two, Lemma \ref{lem:full-to-truncation}
(whose hypotheses \eqref{ll} and \eqref{eqn:lmatrix-boundary-decay} follow
from Lemmas \ref{l0eigenvalues} and \ref{bijectionl}),
applied with $\L = \L(0)$ and $\widetilde \L$ given by the (zero-padded) minor $\H_m$, supplies a map
$\pi_m:\mathcal A_m\to\eig\H_m$ satisfying
$|\pi_m(\Lambda_a)-\Lambda_a|\le e^{-cD}$, for $m \in \{\varphi, i\}$.
Lemma \ref{eigenvalues0} makes $\pi_m$ injective. For each
$m \in \{\varphi, i\}$, define
$\mathcal U_m:=\eig\H_m\setminus\pi_m(\mathcal A_m)$. For every
$a\in\llbracket N_1,m-D\rrbracket$,
$\left| \log|\lambda_j-\Lambda_a| -
\log|\lambda_j-\pi_m(\Lambda_a)|  \right| \leq e^{-c (\log N)^2}$ by
our restriction to $\mathsf{SEP}_{\L}(e^{-(\log N)^2})$ and since
$|\pi_m(\Lambda_a)-\Lambda_a|\le e^{-cD}$. To show \eqref{eq:time-zero-log-reduction-short},
we will track the (positive part of the) error from replacing its left hand side with 
$$\sum_{N_1 \leq a  < i} \log|\lambda_j - \Lambda_a| - \sum_{N_1 \leq a  < \varphi} \log|\lambda_j - \Lambda_a| = -\sum_{a=i}^{\varphi-1} \log|\lambda_j - \Lambda_a|.$$
Our restriction to $\mathsf{SEP}_{\L}(e^{-(\log N)^2})$ 
bounds the total contribution 
from the $O(D)$ indices in the two sums on the left hand side above
within distance $D$ of either $i$ or $\varphi$
by $(\log N)^C$. Finally, $|\mathcal U_i|\le D$, and our restriction to
$\mathsf{BND}_{\L(0)}(\log N)$ gives
$\sum_{\mu\in\mathcal U_i}\log|\lambda_j-\mu|\le CD\log\log N$.
It follows that
\[
\sum_{\mu\in\eig \H_i}\log|\lambda_j-\mu|
-
\sum_{\nu\in\eig \H_\varphi}\log|\lambda_j-\nu|
\leq
-\sum_{a=i}^{\varphi-1}
\log|\lambda_j-\Lambda_a|
- \sum_{\mu \in \mathcal U_{\varphi}} \log|\lambda_j-\mu| 
+(\log N)^C.
\]

This gives \eqref{eq:time-zero-log-reduction-short}, 
except that it remains only to control the 
contribution from $\mathcal U_{\varphi}$.
Put \(\H_\varphi^+:=\L^{[\varphi+1,N_2]}(0)\). Define $\mathcal A_+\coloneqq \{\Lambda_a:a\geq\varphi+D\}$.
Lemma \ref{lem:full-to-truncation} again supplies an injection
\(\pi_+:\mathcal A_+\to\eig \H_\varphi^+\); set $\mathcal{U}_+ = \eig \H_{\varphi}^+ \setminus \pi_+(\mathcal A_+)$. 
Apply Lemma \ref{lem:suff} with
$\mathcal A_-:=\mathcal A_\varphi$, $\pi_-:=\pi_\varphi$, and
$\mathcal U_-:=\mathcal U_\varphi$, and let $\mathcal B$ be the resulting
set from its statement. Then $|\mathcal U_\pm|+|\mathcal B|\leq CD$, and
Lemma \ref{eigenvalues0} gives
\[
-\sum_{\theta\in\mathcal B}\log|\lambda_j-\theta|\leq C D(\log N)^2.
\]
The preceding application of Lemma \ref{lem:full-to-truncation}, its
analog for $\pi_+$, and Lemma \ref{eigenvalues0} show that the 
error $\mathcal E$ in \eqref{eq:matched-log-error} is at most
$e^{-cD}$. Moreover, Lemma \ref{l0eigenvalues}, applied with
$A=\log N$, gives
$\sum_{\mu\in\mathcal U_+}\log|\lambda_j-\mu|\leq CD\log\log N$.
Since $|u_j(\varphi)|\geq\zeta\geq
e^{-200(\log N)^{3/2}}$, Lemma \ref{lem:suff} therefore yields
\[
-\sum_{\mu\in\mathcal U_-}\log|\lambda_j-\mu|
\leq C(\log N)^C.
\]
This shows \eqref{eq:time-zero-log-reduction-short}.

To obtain \eqref{eq:time-zero-singular-block}, we use 
\eqref{eqn:hiTconv} from Proposition \ref{prop:app-time-zero-dos}, and we use
a smoothing argument together with Minami estimates to deal with the singularity of the logarithm. 
We omit details on this point, as they are similar to the proof of Proposition \ref{prop:scat_conc}.
\end{proof}

\subsection{Proof of Proposition \ref{prop:outside-strip-confinement}}
\label{app:outside-strip}

\begin{proof}[Proof of Proposition \ref{prop:outside-strip-confinement}]
Following \cite[Section 2]{Mos75}, let
\(\mathbf P(s)=[P_{ij}(s)]_{i,j\in\llbracket N_1,N_2\rrbracket}\) be the
skew-symmetric tridiagonal matrix defined by
\[
P_{i,i+1}(s)=\frac{a_i(s)}2,
\qquad
P_{i+1,i}(s)=-\frac{a_i(s)}2,
\qquad i\in\llbracket N_1,N_2-1\rrbracket,
\]
with all other entries equal to zero. Then the Toda Lax equation is
$\partial_s\mathbf L(s)=[\mathbf L(s),\mathbf P(s)]$.
Let \(\mathbf V(s)=[V_{mn}(s)]_{m,n\in\llbracket N_1,N_2\rrbracket}\)
solve
\[
\partial_s \mathbf V(s)=-\mathbf P(s)\mathbf V(s),
\qquad
\mathbf V(0)=\mathbf I.
\]
Then
$\mathbf L(s)=\mathbf V(s)\mathbf L(0)\mathbf V(s)^{\mathsf T}$.
Recalling the notation \(\mathbf u_k(s)=(u_k(n;s))_{n=N_1}^{N_2}\) 
for eigenvectors from
Definition \ref{def:quasi_intro}, we note that
$\mathbf u_k(s)=\mathbf V(s)\mathbf u_k(0)$.
Thus,
\begin{equation}
\label{eq:quantitative-propagation-convolution}
u_k(m;s)
=
\sum_{n=N_1}^{N_2}V_{mn}(s)u_k(n;0).
\end{equation}

The bound in
\cite[Equation (5.5)]{Agg25a} gives, on an event of overwhelming probability,
\[
|V_{mn}(s)|
\le
\sum_{\ell\ge |m-n|}(2es\log N/\ell)^\ell.
\]
If $D\ge \frac{1}{4}T(\log N)^{3}$ and
$|m-n|\ge\frac D2$,
then, since \(s\le T\log N\),
$\frac{2es\log N}{\ell}
\le
\frac{C}{\log N}$.
For sufficiently large
\(N\), \(\frac{C}{\log N} \leq e^{-32R-1}\). Consequently,
\begin{equation}
\label{eq:V-fast-tail}
|V_{mn}(s)|
\le
2e^{-32R|m-n|}
\qquad
\text{if }|m-n|\ge\frac D2.
\end{equation}

Set $D\coloneqq |m-\varphi_0(k)|$, and suppose $D \geq T(\log N)^3/2$. Split the sum in
\eqref{eq:quantitative-propagation-convolution} into the terms for which
$|m-n|\ge D/2$ and those for which $|m-n|<D/2$.
For the first case, using \(|u_k(n;0)|\le1\) and
\eqref{eq:V-fast-tail},
\[
\sum_{|m-n|\ge D/2}
|V_{mn}(s)u_k(n;0)|
\le
C\sum_{r\ge D/2}e^{-32Rr}
\le
C e^{-10RD}.
\]
For the second case, if \(|m-n|<D/2\), then
$|n-\varphi_0(k)|\ge D-|m-n|>D/2$.
In particular, after restricting to the overwhelmingly probable event of 
Lemma \ref{lem:time-0-locR}, and
using \(|V_{mn}(s)|\le1\) and
$T^{2/3}\le R|n-\varphi_0(k)|$ for sufficiently large $N$,
\[
\sum_{|m-n|<D/2}
|V_{mn}(s)u_k(n;0)|
\le
\sum_{|n-\varphi_0(k)|\ge D/2}
e^{-19R|n-\varphi_0(k)|}
\le C e^{-9RD}.
\]

Combining the preceding estimates, we obtain
$|u_k(m;s)|\le Ce^{-8RD}$ for sufficiently large \(N\). This implies the $s' = s$
case of the proposition. Moreover, the previous display implies that 
$$|\varphi_s(k) - \varphi_0(k)| \leq \frac{1}{2} T (\log N)^3.$$
Therefore, if $d \coloneqq |\varphi_s(k) -m| > T (\log N)^3$, 
then $|\varphi_0(k) -m| \geq T (\log N)^3/2$,  and thus we have
$
|u_k(m;s)|
\le
e^{-4Rd}$,
as claimed.
\end{proof}

			\bibliography{bibliotek.bib}
	\bibliographystyle{abbrv}

\end{document}